\documentclass[11pt,
]{article}
\usepackage{amssymb,amstext,amsmath,amsthm,latexsym,mathrsfs,mathbbol}

\usepackage{makeidx}      

\makeindex

\lccode`\-=`\-
\usepackage[T2A]{fontenc}
\usepackage[cp1251]{inputenc}
\usepackage[english,russian]{babel}

\newcommand{\nek}{\newcommand}

\DeclareMathAlphabet{\cur}{U}{eur}{m}{n}

\nek{\skr}{\mathscr}

\nek{\bfit}{\bfseries\itshape}
\nek{\bfsl}{\bfseries\slshape}
\nek{\sfbs}{\mdseries\sffamily\itshape}        

\nek{\vyk} [1] {}

\nek{\imar}[1]{\marginpar[
\flushright\footnotesize\sl
$\mtho\longrightarrow$\\ \vspace{-1ex}{#1}$\mtho\dashv$\vspace*{1ex}]
{
\flushleft\footnotesize\sl
$\mtho\longleftarrow$\\ \vspace{-1ex}{#1}$\mtho\dashv$\vspace*{1ex}}}

\nek{\parf}{\section}
\nek{\punk}{\subsection}
\renewcommand{\thesubsection}{\thesection\asbuk{subsection}}

\vyk{
\chardef\eurl ='025
\chardef\eurM ='026
\nek{\fla} {{\cur{\eurl}}} 
\nek{\fmu} {{\cur{\eurM}}} 
\chardef\eurF = '011
\nek{\fFi} {{\cur{\eurF}}}
}

\newcounter{enuf}
\nek{\enufi}{\addtocounter{enuf}{1}}
\nek{\fenu}{
\def\theenumi{$\mtho(\fnsymbol{enuf})$}
\def\labelenumi{\theenumi}
\enufi\itsep
}
\newcounter{enuc}
\nek{\enuci}{\addtocounter{enuc}{1}}
\nek{\cenu}{
\def\theenumi{$\mtho\arabic{enuc}^\circ$}
\def\labelenumi{$\mtho\arabic{enuc}^\circ.$}
\enuci\itsep
}

\newcounter{enuF}
\nek{\enuFi}{\addtocounter{enuF}{1}}
\nek{\Fenu}{
\def\theenumi{\psur(\arabic{enuF})\psur}
\def\labelenumi{\theenumi}
\nek{\ifla}{\addtocounter{enuF}1\itla}
\enuFi\itsep
}

\nek{\itsep}{\itemsep=0.4ex plus 0.15ex minus 0.15ex}
\nek{\tenu}[1]{
\def\theenumi{#1}
\def\labelenumi{\theenumi}\itsep
}
\nek{\tenui}[1]{
\def\theenumii{#1}
\itsep
}

\theoremstyle{plain}

\newtheorem{theore}             {Theorem} 
\newtheorem{corollar}  [theore]{Corollary}
\newtheorem{propo}     [theore]{Proposition}
\newtheorem{lemm}      [theore]{Lemma}
\newtheorem{cla}       [theore]{Claim} 
\newtheorem{clt}       {Claim} [theore]

\theoremstyle{definition}
\newtheorem{defn}      [theore]{Definition}
\newtheorem{prim}      [theore]{Example}
\newtheorem{rem}       [theore]{Remark}

\newtheorem*{prF}{{\bf Proof}}               
\newtheorem*{prC}{{\bf Proof of the claim}}  
\newtheorem*{con}{{\bfit Construction\/}{\bf.}}    
\newtheorem*{aq} {{\bf Acknowledgements\/}}          

\newtheorem{que}  {Question}
\nek{\thsp}{\hspace{0.1ex plus \mathsurround}}
\nek{\bpro}{\begin{propo}}
\nek{\epro}{\end{propo}}
\nek{\bcl} {\begin{cla}}
\nek{\ecl} {\end{cla}}
\nek{\bct} {\begin{clt}}
\nek{\ect} {\end{clt}}
\nek{\bco} {\begin{con}}
\nek{\eco} {\qeDD{Construction}\end{con}}
\nek{\bcor}{\begin{corollar}}
\nek{\ecor}{\end{corollar}}
\nek{\baq} {\begin{aq}}
\nek{\eaq} {\end{aq}}
\nek{\bex} {\begin{prim}}
\nek{\eex} {\qeD\end{prim}}
\nek{\bdf} {\begin{defn}} 
\nek{\eDf} {\end{defn}}
\nek{\edf} {\qeD\end{defn}}
\nek{\edF} {\end{defn}}
\nek{\ble} {\begin{lemm}}
\nek{\ele} {\end{lemm}}
\nek{\bte} {\begin{theore}}
\nek{\ete} {\end{theore}}
\nek{\bre} {\begin{rem}} 
\nek{\ere} {\qeD\end{rem}} 
\nek{\bqe} {\begin{que}} 
\nek{\eqe} {\qeD\end{que}} 
\nek{\bpf} {\begin{prF}} 
\nek{\epf} {\qed\end{prF}} 
\nek{\ePf} {\end{prF}} 
\nek{\bpc} {\begin{prC}} 
\nek{\epc} {\qeDD{Claim}\end{prC}} 
\nek{\qeD} {\hfill$\mtho\Box$}
\nek{\qeDD} [1] 
{\hfill\hbox{$\mtho\Box$~({\sl #1\/}\hspace{0.1ex})}}
\nek{\epF} [1] {\qeDD{#1}\end{prF}} 

\nek{\bde}{\begin{description}}
\nek{\ede}{\end{description}}
\nek{\ben}{\begin{enumerate}}
\nek{\een}{\end{enumerate}}
\nek{\bit}{\begin{itemize}}
\nek{\eit}{\end{itemize}}
\nek{\bay}{\begin{array}}
\nek{\eay}{\end{array}}
\nek{\bmp}{\begin{minipage}}
\nek{\emp}{\end{minipage}}
\nek{\fF} {{\bf F}}
\nek{\fG} {{\bf G}}
\nek{\Gd} {\fG_\da}
\nek{\Fs} {\fF_\sg}
\nek{\fS} {{\bf S}}
\nek{\fT} {{\bf T}}
\nek{\fH} {{\bf H}}
\nek{\fa} {{\bf a}}
\nek{\fb} {{\bf b}}
\nek{\ZFC}{{\bf ZFC}}
\nek{\ZC} {{\bf ZC}}
\nek{\zhc}{{\bf ZFHC}}
\nek{\pli}{{\bf I}}
\nek{\plj}{{\phantom{I}\bf I}}
\nek{\pld}{{\bf II}}
\nek{\gai} [2] {\fH^\cZ_{#1#2}}
\nek{\gad} [2] {\fG^{#1}_{#2}} 
\nek{\iai} [1] {\fH_{#1}}
\nek{\iad} [1] {\fG_{#1}}
\nek{\hi} {\fH}
\nek{\hd} {\fG}
\nek{\etc} {{\sl etc}}
\nek{\iesp}{\hspace{0.3ex}}
\nek{\pw} {\hbox{a.\iesp e.}}
\nek{\ie} {\hbox{\sl i.\iesp e.}}
\nek{\eg} {\hbox{\sl e.\iesp g.}}
\nek{\ea} {\hbox{\sl et.\hspace{0.3ex}al.}}
\nek{\noo}
{\hbox{\sl w.\iesp l.\iesp o.\iesp g.}}
\nek{\pv} {\hbox{a.\iesp a.}}
\nek{\vrt} {\hbox{w.\iesp r.\iesp t.}}
\nek{\lsc} {\hbox{l.\iesp s.\iesp c.}}
\nek{\er} {{ER}}
\nek{\bm} {{BM}}
\nek{\dd}[1]{$\mtho\hspace{0.2ex}{#1}$-\hspace{0.0ex}}
\nek{\dw}{\dd\om}
\nek{\lis}[1] {\mathop{\tt lim\hspace{0.2ex}sup}_{#1}}
\nek{\len}[1] {\mathop{\tt lh}{#1}}
\nek{\Ord}  {{\tt{Ord}}}
\nek{\Exh}  {{\tt{Exh}}}
\nek{\Nul}  {{\tt{Null}}}
\nek{\Mod}  {\mathop{\tt{Mod}}}
\nek{\Aut}  {\mathop{\tt{Aut}}}
\nek{\card} {\mathop{\tt card}}
\nek{\lh}   {\mathop{\tt lh}}
\nek{\pr}   {\mathop{\tt pr}}
\nek{\sr}   {\mathop{\tt sr}}
\nek{\ran}  {\mathop{\tt ran}}
\nek{\dom}  {\mathop{\tt dom}}
\nek{\fld}  {\mathop{\tt field}}
\nek{\otp}  {\mathop{\tt otp}}
\nek{\Max}  {\mathop{\tt Max}}
\nek{\tsup} {\mathop{\tt sup}}
\nek{\tinf} {\mathop{\tt inf}}
\nek{\tmin} {\mathop{\tt min}}
\nek{\tmax} {\mathop{\tt max}}
\nek{\tlim} {\mathop{\tt lim}}
\nek{\tlis} {\mathop{\tt lim\hspace{0.3ex}sup}}
\nek{\tlii} {\mathop{\tt lim\hspace{0.3ex}inf}}
\nek{\Fin}  {{\tt Fin}} 
\nek{\bFin} {{\bf Fin}} 
\nek{\maxi}[1] {\Max^\xi_{#1}}
\nek{\HC} {{\rm HC}}
\nek{\ccc}{{\sc ccc}}
\nek{\al} {\alpha}
\nek{\ba} {\beta}
\nek{\ga} {\gamma}
\nek{\Da} {\Delta}
\nek{\da} {\delta}
\nek{\ka} {\kappa}
\nek{\la} {\lambda}
\nek{\La}{\Lambda}
\nek{\sg} {\sigma}
\nek{\Sg} {\Sigma}
\nek{\vpi}{\varphi}
\nek{\vpy}{\vpi_\iy}
\nek{\vt} {\vartheta}
\nek{\vT} {\Theta}
\nek{\ovt}{{\overline\vt}}
\nek{\ovi}{{\overline\vpi}}
\nek{\ops}{{\overline\psi}}
\nek{\ve} {\varepsilon}
\nek{\om} {\omega}
\nek{\Om} {\Omega}
\nek{\lom}{^{<\om}}
\nek{\za}{\zeta}
\nek{\tpi}{\tau_\vpi}
\nek{\omi} {\om_1}
\nek{\omm} [1] {\om^{\om^{#1}}}
\nek{\bse} {2\lom}
\nek{\nse} {\dN\lom}
\nek{\alo} {{\aleph_0}}
\nek{\sd}   {\mathbin{\Da}}
\nek{\bigd} {\mathbin{\hbox{\large$\mtho\Da$}}}
\nek{\sqe} {\fmu\hspace{0.05ex}}
\newcommand{\fs}[2]{{\bf\iSg}^{#1}_{#2}}
\newcommand{\fp}[2]{{\bf\iPi}^{#1}_{#2}}

\newcommand{\iSg}{{\mathchar"7106}}
\newcommand{\iPi}{{\mathchar"7105}}
\newcommand{\iDa}{{\mathchar"7101}}
\newcommand{\is}[2]{\iSg^{#1}_{#2}}
\newcommand{\ip}[2]{\iPi^{#1}_{#2}}
\newcommand{\id}[2]{\iDa^{#1}_{#2}}
\nek{\BBB}{\hspace{0.05ex}}
\nek{\dA}{{\BBB{\mathbb A}\BBB}}
\nek{\dC}{{\BBB{\mathbb C}\BBB}}
\nek{\dF}{{\BBB{\mathbb F}\BBB}}
\nek{\dN}{{\BBB{\mathbb N}\BBB}}
\nek{\dP}{{\BBB{\mathbb P}\BBB}}
\nek{\dQ}{{\BBB{\mathbb Q}\BBB}}
\nek{\dqp}{\dQ^+}
\nek{\dR}{{\BBB{\mathbb R}\BBB}}
\nek{\dS}{{\BBB{\mathbb S}\BBB}}
\nek{\dT}{{\BBB{\mathbb T}\BBB}}
\nek{\dV}{{\BBB{\mathbb V}\BBB}}
\nek{\dZ}{{\BBB{\mathbb Z}\BBB}}
\nek{\dX}{{\BBB{\mathbb X}\BBB}}
\nek{\dY}{{\BBB{\mathbb Y}\BBB}}
\nek{\dyn} {\dY^\dN}
\nek{\dvp} {\dV^+}
\nek{\dn}{2^\dN}
\nek{\dntn}{2^{\dN\ti\dN}}
\nek{\dnqn}{{(2^\dN)}{}^\dN}
\nek{\pnqn}{{\pn}{}^\dN}
\nek{\bn}{\dN^\dN}
\nek{\rn} {\dR^\dN}
\nek{\tm} [1] {2^{\mxi\xi}}
\nek{\nn}{{\dN\ti\dN}}
\nek{\ccs} {}
\nek{\cA}{{\ccs{\skr A}\ccs}}
\nek{\cD}{{\ccs{\skr D}\ccs}}
\nek{\cE}{{\ccs{\skr E}\ccs}}
\nek{\cF}{{\ccs{\skr F}\ccs}}
\nek{\cS}{{\ccs{\skr S}\ccs}}
\nek{\cP}{{\ccs{\skr P}\ccs}}
\nek{\cW}{{\ccs{\skr W}\ccs}}
\nek{\cX}{{\ccs{\skr X}\ccs}}
\nek{\cI} {{\skr I}} 
\nek{\cJ} {{\skr J}} 
\nek{\cO} {{\skr O}} 
\nek{\cZ} {{\skr Z}}
\nek{\zo} {\cZ_0}
\nek{\zw} {\cZ_{\hbox{\small\rm w}}}
\nek{\xn} {\cX^\dN}
\nek{\an} {\dA^\dN}
\nek{\cd} [1] {\cD_{#1}}
\nek{\pws}  [1] {\cP(#1)}
\nek{\cp}  [2] {\cP_{\hspace*{-0.4ex}\tt cnt}^{#1}(#2)}
\nek{\pwf} [1] {\cP_{\hspace*{-0.4ex}\tt fin}(#1)}
\nek{\pwc} [1] {\cP_{\hspace*{-0.4ex}\tt ctbl}(#1)}
\nek{\pnn}{\cP(\nn)}
\nek{\pn}{\cP(\dN)}
\nek{\ps}{\cP(\dS)}
\nek{\pz}{\cP(\dZ)}
\nek{\mm}{{\BBB{\mathfrak M}\BBB}}
\nek{\gE}{{\BBB{\mathfrak E}\BBB}}
\nek{\shi} {{\mathfrak s}}
\newcommand{\gc}{{\BBB{\mathfrak c}\BBB}}
\nek{\ilo}[1] {{[0,n_{#1})}}
\nek{\ii} [1] {{[n_{#1},\infty)}}
\nek{\ir} [2] {{[#1,#2)}}
\nek{\iry}[1] {{[#1,\iy)}}
\nek{\iq} [2] {\ir{\nu_{#1}}{\nu_{#2}}}
\nek{\iqy}[1] {\ir{\nu_{#1}}\iy}
\nek{\iqo}[1] {\ir0{\nu_{#1}}}
\nek{\iqn}[1] {\iq{#1}{#1+1}}
\nek{\opl} {\oplus}
\nek{\ap}  {\cdot}
\nek{\cj}  {\mathbin{\hspace{0.2ex}\&\hspace{0.2ex}}}
\nek{\dm}  {$$}
\nek{\sus} {\mathopen{\exists\hspace{0.35ex}}}
\nek{\kaz} {\mathopen{\forall\hspace{0.35ex}}}
\nek{\imp} {\Longrightarrow} 
\nek{\eqv} {\Longleftrightarrow} 
\nek{\ti}  {\times} 
\nek{\mo}  {\models} 
\nek{\sq}  {\subseteq}
\nek{\qs}  {\supseteq}
\nek{\su}  {\subset}
\nek{\sneq}{\subsetneqq}
\nek{\we}  {{\mathbin{\hspace*{0.2ex}^\wedge}}}
\nek{\obr} {^{-1}}
\nek{\dif} {\smallsetminus}
\nek{\res} {\mathbin{\restriction}}
\nek{\lef} {\preccurlyeq}
\nek{\gef} {\succcurlyeq}
\nek{\pu}  {\emptyset}
\nek{\iy}  {\infty}
\nek{\piy} {+\iy}
\nek{\nin} {\not\in}
\nek{\limp}{\,\imp\,}
\nek{\leqv}{\,\eqv\,}
\nek{\onto}{\stackrel{{\rm onto}}{\longrightarrow}}
\nek{\ang} [1] {\langle #1\rangle}
\nek{\stk} [2] {\ang{#1\hspace{0.3ex};\hspace{0.1ex}#2}}
\nek{\sis} [2] {\ans{#1}_{#2}}
\nek{\ans} [1] {\{\hspace{0.01ex}#1\hspace{0.01ex}\}}

\nek{\zz} {\linebreak[0]} 
\nek{\ens} [2] {\ans{{#1\hspace{0.5ex}{:}}\zz\hspace{0.5ex}#2}}

\nek{\suh} [1] {[\hspace{0.3pt}#1\hspace{0.3pt}]_{\sq}}
\nek{\itla} {\item\label}

\nek{\aeq} {\mathbin{\|}}
\nek{\laeq}{\,\aeq\,} 

\nek{\tz} {\mathbin{;}}

\nek{\seq}[2] {(#1)_{#2}}

\nek{\kb}[2]{#1^{(#2)}}

\nek{\ssty} {\textstyle}
\nek{\isum} [3] {{{\ssty\ugl\sum_{#1}\,#3}\mid{#2}\ugr}}

\nek{\ifi}  {{\cur{Fin}\hspace*{0.1ex}}}
\nek{\frt}  {{\cur{Fr}\hspace*{0.1ex}}}
\nek{\ibo}  {{\cur{Bou}\hspace*{0.1ex}}}
\nek{\fifi} {\ifi\ti\ifi}
\nek{\fio} {\ifi\ti0}
\nek{\ofi} {0\ti\ifi}

\nek{\xip} {{\xi+1}}
\nek{\etp} {{\eta+1}}

\nek{\vx} [1] {^{(#1)}}

\nek{\dop}[1] {{#1}^{\complement}}

\nek{\skl} {\hbox{\mtho\large$($}}
\nek{\skp} {\hbox{\mtho\large$)$}}

\nek{\ugl} {\hbox{\mtho\large$\langle$}}
\nek{\ugr} {\hbox{\mtho\large$\rangle$}}

\nek{\df} [1] {\dop {#1}}
\nek{\pl} [1] {{#1}^+}

\nek{\bbW} {\hbox{\mtho\boldmath $W$}}

\nek{\bo}{{\bf O}}

\nek{\Zo}  {{\cZ_0}}
\nek{\Eo}  {\rE_{\text{\sf0}}}
\nek{\Fo}  {\rF_0}
\nek{\Ec}  {\mathbin{\overline{\rE}}}
\nek{\Eco} {\mathbin{\overline{\Eo}}}
\nek{\rzo}  {\rZ_{\text{\sf0}}}
\nek{\neo} {\mathbin{\not{{\hspace{-0.4ex}\rE}}_0}}
\nek{\nfo} {\mathbin{\not{{\hspace{-0.4ex}\rF}}_0}}
\nek{\nrE} [1] {\mathbin{\not{{\hspace{-0.4ex}\rE}}_{#1}}}
\nek{\rD}  {\mathbin{\sf D}}
\nek{\rR}  {\mathbin{\sf R}}
\nek{\rT}  {\mathbin{\sf T}}
\nek{\rtd} {\rT_2}
\nek{\rZ}  {\mathbin{\sf Z}}
\nek{\rdi} {\mathbin{{\sf D}_I}}
\nek{\rda} {\mathbin{{\sf D}_A}}
\nek{\rF}  {\mathbin{\sf F}}
\nek{\rE}  {\mathbin{\sf E}}
\nek{\rP}  {\mathbin{\sf P}}
\nek{\ei}   {\rE^\iy}
\nek{\rfi}  {\rF^\iy}
\nek{\nD}  {\mathbin{{\not\hspace{-0.35ex}\sf D}}}
\nek{\nE}  {\mathbin{{\not\hspace{-0.35ex}\sf E}}}
\nek{\nR}  {\mathbin{{\not\hspace{-0.35ex}\sf R}}}
\nek{\nF}  {\mathbin{{\not\hspace{-0.25ex}\sf F}}}
\nek{\reo} {\rE_{\hspace{-1.0pt}\Zo}}
\nek{\rez} {\rE_{\hspace{-1.0pt}\cZ}}
\nek{\nrz} {\mathbin
{\not{{\hspace{-0.4ex}\rE}}_{\hspace{-1.0pt}\cZ}}}
\nek{\nre} {\mathbin
{\not{{\hspace{-0.4ex}\rE}}_{\hspace{-1.0pt}\Zo}}}
\nek{\reff}{\rE_{\fifi}}
\nek{\reb} {\le_{\rm B}}
\nek{\eqb} {\approx_{\rm B}}
\nek{\dzo} {\dd{\Zo}}
\nek{\dde} {\dd{\rE}}
\nek{\ddec}{\dd{\Ec}}
\nek{\dep} {\dd{\rE'}}
\nek{\ddf} {\dd{\rF}}
\nek{\ddv} {\dd{\vt}}
\nek{\ddz} {\dd\cZ} 
\nek{\der} {\er-}
\nek{\dee} {\dd{\rE,\rE'}}
\nek{\Def} {\dd{\rE,\rF}}
\nek{\ddd} [2] {\dd{#1,#2}}

\nek{\fdo}{\hbox{\raisebox{0.2ex}{\mtho\tiny$\bullet$}}}

\nek{\fdt}{\hbox{\raisebox{-0.25ex}{\LARGE\bf.}}}
\nek{\bdot}[1] {\raisebox{-0.07ex}{\mtho$\stackrel{\fdt}{#1}$}}
\nek{\doa} {{\bdot a}} 
\nek{\dox} {{\bdot x}}
\nek{\dog} {\raisebox{-0.28ex}{\mtho$\stackrel{\fdt}g$}}
\nek{\doxl}{\dox_{\tt left}}
\nek{\doxr}{\dox_{\tt right}}
\nek{\rkb} [1] {|#1|_{\rm CB}}
\nek{\rkt} [2] {{^{\om}\hspace*{-1.3pt}|#1|_{#2}}}
\nek{\rko} [2] {{|#1|_{#2}}}
\nek{\rkT} [1] {{^{\om}\hspace*{-1.3pt}|#1|}}
\nek{\rkO} [1] {{|#1|}}
\nek{\kt}  [1] {{^{#1}\hspace*{-0.8pt}T}}
\nek{\bsf} [1] {I_{#1}}
\nek{\fri} [1] {\mathbin{\displaystyle\rE^{\tt fr}_{#1}}}
\nek{\fre} [1] {\frt_{{#1}}}

\nek{\fps} [3] {{\prod_{#3}{#2}\,/\,{#1}}}
\nek{\fpt} [2] {{\prod{\sis{#2}{}}\,/\,{#1}}}
\nek{\fpd} [2] {{#1}\otimes{#2}}

\nek{\sto} {[s_0,t_0]}
\nek{\ostp}{[s',t']}
\nek{\ost} {[s,t]}
\nek{\osu} {[s,u]} 
\nek{\otu} {[t,u]}
\nek{\ouv} {[u,v]}

\nek{\zs} {{\tilde s}}
\nek{\zt} {{\tilde t}}
\nek{\zu} {{\tilde u}}
\nek{\zv} {{\tilde v}}
\nek{\zn} {{\tilde n}}
\nek{\zm} {{\tilde m}}

\nek{\ff}[2] {F_{#1}^{#2}}

\nek{\spa} {\dS}
\nek{\spn} {\spa^\dN}
\nek{\spp} {\dY}

\nek{\poq}{\underline}
\nek{\nad}{\overline}

\nek{\sm} {\hbox{\mtho{\large$\Sg$}}}
\nek{\smp} [2] {\sm_{#1}^{#2-1}}
\nek{\smy} [1] {\sm_{#1}^{\iy}}

\nek{\sui} [1] {\cS_{\ans{#1}}}
\nek{\sun} {{\sui{1/n}}}
\nek{\srn} {{\sui{r_n}}}
\nek{\ern} {\rE_{\ans{r_n}}}
\nek{\erpn} {\rE_{\ans{r'_n}}}
\nek{\eun} {\rE_{\ans{1/n}}}
\nek{\suo} {{\sui0}}
\nek{\eso} {\rE_{\hspace{-1.0pt}\suo}}
\nek{\nso} {\mathbin
{\not{{\hspace{-0.4ex}\rE}}_{\hspace{-1.0pt}\suo}}}

\nek{\gal} [3] {{\tt Gal}^{#1}_{#2}(#3)} 

\nek{\nbd} [1] {{\skr O}_1(#1)}

\nek{\aH} {H^\ast}
\nek{\aB} {B^\ast}

\nek{\ovl} [1] {\overline{#1}} 
\nek{\ovg} [1] {\ovl{g_{#1}}}
\nek{\ovp} [1] {\ovl{\ga_{#1}}}

\nek{\plo} {+1} 

\nek{\yk} [1] {k_{#1}}

\nek{\mtho}{\mathsurround=0mm}
\nek{\msur}{\hspace{-1\mathsurround}}
\nek{\psur}{\hspace{0.3\mathsurround}}
\nek{\dsur}{\hspace{-0.3\mathsurround}}
\nek{\hsur}{\hspace{-0.5\mathsurround}}
\nek{\noi}{\noindent}
\nek{\vom}{\vspace{1mm}}
\nek{\vtm}{\vspace{2mm}}

\nek{\uv}{{\bf V}}

\nek{\wA} {{\widehat A}}
\nek{\ha} {{\hat a}}
\nek{\hb} {{\hat b}}
\nek{\he} {{\hat\ve}}
\nek{\hT} {{\hat t}}
\nek{\hl} {{\hat l}}
\nek{\hs} {{\hat s}}
\nek{\hm} {{\hat m}}
\nek{\hn} {{\widehat n}}
\nek{\hk} {{\widehat k}}

\nek{\fo} {{\mathbf 0}}
\nek{\fr} {{\mathbf 1}}

\nek{\vnu} {{\vec \nu}} 
\nek{\dvn} {\dd\vnu} 
\nek{\wnu} {\cW_{\vnu}} 
\nek{\ewn} {\rE_{\vnu}}
\nek{\ret} {\rE_{T}}

\nek{\ci} [1] {I_{#1}}

\nek{\vex}{{\vec x}}
\nek{\vey}{{\vec y}}

\nek{\dns} [1] {d_{\vnu}^{#1}}
\nek{\den} {d_{\vnu}}
\nek{\din} {d_\vnu}

\nek{\poo}{=_{\tt df}}

\nek{\dpi}{d_\vpi}

\nek{\nol}[1] {{\bf 0}_{#1}}
\nek{\edi}[1] {{\bf 1}_{#1}}

\nek{\nrn} {_{\ans{r_n}}}
\nek{\vpr} {\vpi\nrn}
\nek{\dpr} {d\nrn}
\nek{\drn} {\dd{\ans{r_n}}}

\nek{\okr} [2] {{\skr O}_{#1}(#2)}

\nek{\nid}{\gE}

\nek{\bus}{\begin{equation}}   
\nek{\eus}{\end{equation}} 

\nek{\zid} [2] {\cI_{#1/#2}}
\nek{\zidef} {\zid\rE\rF}

\nek{\zfo} [2] {\dP_{#1/#2}}
\nek{\zfoef} {\zfo\rE\rF}

\nek{\peo} {\dP_{\Eo}}
\nek{\pep} {\dP'_{\Eo}}
\nek{\ieo} {\cI_{\Eo}}

\nek{\renek}{\renewcommand}

\renek{\thesubsection}
{\thesection.\alph{subsection}}
\renek{\parf}[1]{\section{{\protect\boldmath#1}}}
\renek{\punk}[1]{\subsection{{\protect\boldmath#1}}%
\setcounter{clt}1}

\theoremstyle{plain}
\newtheorem{cort}  [clt]{Corollary}  
\nek{\bcot}{\begin{cort}}
\nek{\ecot}{\end{cort}}

\newtheorem{prot}  [clt]{Proposition}  
\nek{\bprt}{\begin{prot}}
\nek{\eprt}{\end{prot}}

\newtheorem{lemt}  [clt]{Lemma}
\nek{\blt}{\begin{lemt}}
\nek{\elt}{\end{lemt}}

\theoremstyle{definition}

\newtheorem{remt}  [clt]{Remark}  
\nek{\bret}{\begin{remt}}
\nek{\eret}{\qeD\end{remt}}

\newtheorem{vop}       [theore]{Question}
\nek{\bqu} {\begin{vop}}
\nek{\equ} {\qeD\end{vop}}

\nek{\bup} {\begin{upr}}
\nek{\eup} {\qeD\end{upr}}

\newtheorem{uprt}       [clt]{Exercise}
\nek{\bupt} {\begin{uprt}}
\nek{\eupt} {\qeD\end{uprt}}

\nek{\btb}{\begin{tabular*}}
\nek{\etb}{\end{tabular*}}

\newtheorem{deft}  [clt]{Definition}
\nek{\bdt}{\begin{deft}\thsp\rm}
\nek{\edt}{\qeD\end{deft}}

\SetMathAlphabet{\cur}{bold}{U}{eur}{b}{n}

\renek{\Gd} {\fG_\fda}
\renek{\Fs} {\fF_\fsg}

\nek{\ubf}{\fontseries{b}\selectfont}

\mathchardef\alphA ="710B
\mathchardef\betA ="710C
\mathchardef\gammA ="710D
\mathchardef\deltA ="710E
\mathchardef\vartA ="7123
\mathchardef\kpA   ="7114
\mathchardef\mU    ="7116
\mathchardef\nU    ="7117
\mathchardef\rhO   ="711A
\mathchardef\sigmA ="711B

\nek{\fal} {{\cur{\alphA}}}
\nek{\fba} {{\cur{\betA}}}
\nek{\fsg} {{\cur{\sigmA}}}
\nek{\fda} {{\cur{\deltA}}}

\nek{\dds}{\dd\fsg}

\nek{\nmp} {\Longleftarrow}

\nek{\tsc}[1]{\hbox{\footnotesize\sc{#1}}}

\nek{\ddi} {\dd{\cI}}
\nek{\ddj} {\dd{\cJ}}

\nek{\ren} {\le_{\tsc{c}}}

\nek{\eui}[1] {{\text{\it{EU}}}_{\ans{#1}}}

\nek{\sqa} {\sq^\ast}
\nek{\sua} {\su^\ast}

\nek{\prf} [1] {\hbox{\S\hspace{0.3ex}\ref{#1}}}
\nek{\prff}[1] {\S\S~\ref{#1}}
\nek{\pff} [1] {\S~{#1}}

\nek{\mrn} {\mu\nrn}

\renek{\dop} [1] {\complement #1}
\renek{\df}  [1] {{#1}^\complement}
\nek{\doP}  [1] {{#1}^\complement}

\nek{\ima} {\hbox{\hspace{0.1ex}''}}

\nek{\aprb}{\approx_{\tsc b}}
\nek{\ismb}{\cong_{\tsc b}}

\nek{\incs} {<_{\tsc i}}
\nek{\inc} {\le_{\tsc i}}
\nek{\eqi} {\sim_{\tsc i}}

\nek{\reas} {<_{\tsc a}}
\nek{\rea} {\le_{\tsc a}}
\nek{\eqa} {\sim_{\tsc a}}

\nek{\reaas} {<_{\tsc{aa}}}
\nek{\reaa} {\le_{\tsc{aa}}}
\nek{\eqaa} {\sim_{\tsc{aa}}}

\nek{\rebs} {<_{\tsc b}}
\renek{\reb} {\le_{\tsc b}}
\renek{\eqb} {\sim_{\tsc b}}

\nek{\reBs}{<_{\tsc{bm}}}
\nek{\reB} {\le_{\tsc{bm}}}
\nek{\eqB} {\sim_{\tsc{bm}}}

\nek{\rds} {<_{\tsc{b},\Da}}
\nek{\rd} {\le_{\tsc{b},\Da}}
\nek{\eqd} {\approx_{\tsc{b},\Da}}

\nek{\rbas} {<_{\tsc{b,ba}}}
\nek{\rba} {\le_{\tsc{b,ba}}}
\nek{\eqba} {\approx_{\tsc{b,ba}}}

\nek{\rbasp} {<_{\tsc{b,ba}}^+}
\nek{\rbap} {\le_{\tsc{b,ba}}^+}
\nek{\eqbab} {\approx_{\tsc{b,ba}}^+}

\nek{\nab} [1] {\nabla(#1)}
\nek{\orb} {\le_{\tsc{rb}}}
\nek{\srb} {<_{\tsc{rb}}}
\nek{\erb} {\sim_{\tsc{rb}}}

\nek{\orbpp} {\le_{\tsc{rb}}^{++}}
\nek{\srbpp} {<_{\tsc{rb}}^{++}}
\nek{\erbpp} {\sim_{\tsc{rb}}^{++}}

\nek{\orbp} {\le_{\tsc{rb}}^{+}}
\nek{\srbp} {<_{\tsc{rb}}^{+}}
\nek{\erbp} {\sim_{\tsc{rb}}^{+}}

\nek{\ork} {\le_{\tsc{rk}}}
\nek{\srk} {<_{\tsc{rk}}}
\nek{\erk} {\sim_{\tsc{rk}}}

\nek{\obe} {\le_{\tsc{be}}}
\nek{\sbe} {<_{\tsc{be}}}
\nek{\ebe} {\sim_{\tsc{be}}}

\nek{\odl} {\le_{\sd}}
\nek{\sdl} {<_{\sd}}
\nek{\edl} {\sim_{\sd}}

\nek{\rei} {\rE_{\cI}}
\nek{\rej} {\rE_{\cJ}}

\nek{\eeb} {\sim_{\tsc{b}}}

\nek{\obep} {\le_{\tsc{be}}^+}
\nek{\sbep} {<_{\tsc{be}}^+}
\nek{\ebep} {\sim_{\tsc{be}}^+}

\nek{\seb} {<_{\tsc{b}}}

\nek{\supp} {\mathop{\tt supp}}

\nek{\atp} {\mathop{\tt at}^+}
\nek{\atm} {\mathop{\tt at}^-}

\nek{\hv} [2] {||#1||_{#2}}

\nek{\vmu} {{\vec \mu}}

\nek{\bel} [1] {\mathrel{\text{\boldmath\mtho$\ell$}}^{#1}}
\nek{\beL} [1] {\mathrel{\text{\rm\bf L}}^{#1}}
\nek{\bem}     {\mathrel{\text{\bfit m}}}
\nek{\fco}{\mathrel{\text{\rm\bf c}_{\text{\sf0}}}}
\nek{\fc} {\mathrel{\text{\rm\bf c}}}
\nek{\fvt}{\mathrel{\text{\boldmath\mtho$\vt$}}}

\nek{\beli} {\bel\iy}

\nek{\oin} {{[0,1]^\dN}}

\nek{\iz} [2] {\cI_{#1}^{#2}}
\nek{\jz} [1] {\iz{}{}(#1)}
\nek{\iw} [2] {\cW_{#1}^{#2}}
\nek{\jw} [1] {\iw{}{}(#1)}
\nek{\ib} [2] {\cB^{#1}_{#2}}
\nek{\jb} [1] {\ib{}{#1}}

\nek{\ovu} {{\overline u}}

\nek{\nr}[2]{\nor{#1}_{#2}} 

\nek{\TS}{\textstyle}
\nek{\DS}{\displaystyle}

\nek{\Ba}{B^\ast}

\nek{\resi} [1] {\mathop{\restriction_{#1}}}

\nek{\gP}{{\BBB{\mathfrak P}\BBB}}
\nek{\gF}{{\BBB{\mathfrak F}\BBB}}
\nek{\gJ}{{\BBB{\mathfrak J}\BBB}}
\nek{\dG}{{\BBB{\mathbb G}\BBB}}
\nek{\dH}{{\BBB{\mathbb H}\BBB}}

\nek{\ac} {\cdot} 

\nek{\curle}{\preccurlyeq}
\nek{\cle}{\curle}
\nek{\cge}{\succcurlyeq}
\nek{\cl} {\prec}
\nek{\resic} [1] {\resi{\cle #1}}
\nek{\rec} {\resic}
\nek{\recb} [1] {\rec{(#1)}}
\nek{\rsic} [1] {\resi{\cl #1}}
\nek{\rc} {\rsic}
\nek{\rcb} [1] {\rc{(#1)}}

\nek{\kai} {\forall^\iy\hspace{0.1ex}}
\nek{\exi} {\exists^\iy\hspace{0.1ex}}

\nek{\ovw}{\nad w}

\nek{\il}[2] {\ir{n_{#1}}{n_{#2}}}
\nek{\ia}[2] {\ir{a_{#1}}{a_{#2}}}
\renek{\ij}[2] {\ir{j_{#1}}{j_{#2}}}
\nek{\cB}{{\BBB{\skr B}\BBB}}
\nek{\cG}{{\BBB{\skr G}\BBB}}
\nek{\cL}{{\BBB{\skr L}\BBB}}
\nek{\cN}{{\BBB{\skr N}\BBB}}
\nek{\cU}{{\BBB{\skr U}\BBB}}

\nek{\pnd}{\pn^\dN}

\nek{\anp} [2] {\ang{#1}^{#2}}

\nek{\ta}{\tau}

\nek{\lev} {\mathop{\tt{lev}}}
\nek{\glu} {\mathop{\tt{dep}}}
\nek{\dia} {\mathop{\tt{diam\hspace{0.15ex}}}}

\nek{\Ei}{\rE_{\text{\sf 1}}}
\nek{\Ed}{\rE_{\text{\sf 2}}}
\nek{\Et}{\rE_{\text{\sf 3}}}
\nek{\Ey}{\rE_\iy}

\nek{\npi}{\nu_\vpi}
\nek{\nsi}{\nu_\psi}

\nek{\Ii}{\cI_1}
\nek{\Id}{\cI_2}
\nek{\It}{\cI_3}

\nek{\emb}  {\sqsubseteq_{\tsc{b}}}
\nek{\emn}  {\sqsubseteq_{\tsc{c}}}
\nek{\embi} {\sqsubseteq_{\tsc{b}}^{\rm i}}
\nek{\emni} {\sqsubseteq_{\tsc{c}}^{\rm i}}

\nek{\sio}{\cS_0}
\nek{\esn}{\rE_{\ans{1/n}}}
\nek{\nsn}{\mathbin{\not\hspace*{-0.3ex}\esn}}

\nek{\req}[2]{{\DS|_{#1}^{#2}}}
\nek{\rlq}[1]{\resi{\ge#1}}
\nek{\rmq}[1]{\resi{<#1}}

\nek{\Dij} {\dd{\cI,\cJ}}
\nek{\deo} {\dd{\Eo}}

\nek{\dc}[2] {{\bf W}^{#1}_{#2}}

\nek{\Za}{{\nad Q}}
\nek{\Ya}{{\widetilde Y}}

\nek{\Ga}{\Gamma}
\nek{\rG}  {\mathbin{\sf G}}

\nek{\inva}{{\tt{inv}}}

\nek{\tP}{{P^\ast}}
\nek{\app} {{\hspace{0.2ex}{\cdot}\hspace{0.2ex}}}

\nek{\lland}{\,\land\,}

\nek{\seqv} {\hspace{0.3ex}\Leftrightarrow\hspace{0.3ex}}
\nek{\simp} {\hspace{0.3ex}\Rightarrow\hspace{0.3ex}}

\nek{\eqn}{\equiv_n}

\nek{\pit}{\tilde\pi}
\nek{\tve}{\tilde\ve}

\nek{\dsu}[2] {{#1\oplus#2}}

\nek{\ef} {\dd{\rE,\rF}}

\nek{\isi} {\cong}
\nek{\isa} {\mathrel{\hspace{0.2ex}\isi^\ast}}

\nek{\ske} [3] {\equiv_{#1#2}^{#3}}
\nek{\skab}[1] {\ske AB{#1}}
\nek{\spab}[1] {\ske{A'}{B'}{#1}}

\nek{\Indent}{\hspace*{3ex}}
\nek{\fsur}{\hspace{0.5\mathsurround}}

\nek{\rdm}  {\rD_{\tt max}}

\nek{\co} {\dd{c_0}}
\nek{\lv} {{\normalfont\scshape lv}-}

\nek{\susi} {\exists^{\iy}\hspace*{0.2ex}}
\nek{\kazi} {\forall^{\iy}\hspace*{0.2ex}}

\nek{\igi}{$\hbox{\mtho\boldmath$\fal$}$}
\nek{\igd}{$\hbox{\mtho\boldmath$\fba$}$}

\nek{\cg}[1] {{\tt Choq}(#1)}
\nek{\cgs}[1]{{\tt Choq}^{\rm s}(#1)}

\nek{\dnn}{\dN^\dN}
\nek{\nnn}{(\dnn){\vphantom{\dN}}^\dN}

\nek{\isg} {{S_\iy}}

\nek{\uset}{{Universal sets}}

\nek{\ler}[2] {\mathbin{\sim^{#1}_{#2}}}
\nek{\aer}[2] {\mathbin{\rE_{#1}^{#2}}}
\nek{\ergx}{\aer\dG\dX}
\nek{\egx} {\ergx}
\nek{\lo} [3] {\cO(#3,#1,#2)}
\nek{\sym}{\ler} 
\nek{\ong} {1_\dG}
\nek{\rr} [2] {\rR^{#1}_{#2}}

\nek{\rav}[1] {\rD(#1)}

\nek{\toq}{\circle*{0.5}}
\nek{\tob}{\circle*{1.0}}

\nek{\stob}{\circle{3.0}}
\nek{\ktob}{\kras\circle{3.0}}

\nek{\cob}{\circle{1.5}}
\nek{\mtir}{\line(-1,0){2}}

\nek{\bon} [2] {\cO_{#1}(#2)}

\nek{\dnnn} {\dnnp\dN}
\nek{\dnnp} [1] {(\dnn){\vphantom{\dN}}^{#1}}

\nek{\prift}{\sf}
\nek{\Penu}{{\prift$\id11$ Enumeration}}
\nek{\Refl}{{\prift Reflection}}
\nek{\Uset}{{\prift Universal Sets}}
\nek{\Kres}{{\prift Kreisel Selection}}
\nek{\Cenu}{{\prift Countable-to-1 Enumeration}}
\nek{\Cuni}{{\prift Countable-to-1 Uniformization}}
\nek{\Cpro}{{\prift Countable-to-1 Projection}}
\nek{\Sepa}{{\prift Separation}}
\nek{\Redu}{{\prift Reduction}}

\nek{\dpf} {\mathord{{\dP^2\hspace*{-0.3ex}}\res\rF}}
\nek{\dpe} {\mathord{{\dP^2\hspace*{-0.3ex}}\res\rE}}

\nek{\ek}[2] {[#1]_{{#2}}}

\nek{\eke}[1] {\ek{#1}{\rE}}
\nek{\ekeo}[1] {\ek{#1}{\Eo}}
\nek{\ekec}[1] {\ek{#1}{\Eco}}
\nek{\ekco}[1] {\ek{#1}{\Ec}}
\nek{\ekfo}[1] {\ek{#1}{\Fo}}
\nek{\ekf}[1] {\ek{#1}{\rF}}
\nek{\ekg}[1] {\ek{#1}{G}}

\nek{\ur}{_{\tt right}}
\nek{\ul}{_{\tt left}}

\nek{\cont}{\hbox{\mtho\large$\gc$}}

\nek{\mem} {\dd\in}

\nek{\bL}{{\bf L}}
\nek{\cli} {{\sc cli}}

\nek{\di} [1] {{#1}^\#}
\nek{\drE} {\mathbin{\di\rE}}
\nek{\drF} {\mathbin{\di\rF}}

\nek{\kon} {\hbox{\mtho\large${\mathfrak c}$}}
\nek{\bk} [1] {{\cur B}_{#1}}

\nek{\wtau}{{\widehat\tau}}

\nek{\gra}[1]{\dd{#1}``grainy''}
\nek{\grap}{``grainy''}

\nek{\xE}[2] {\mathbin{\rR^{#2}_{\ge #1}}}

\nek{\moq} [1] {\Mod_{#1}}
\nek{\mox} {\moq} 
\nek{\loa} [1] {j_{#1}}
\nek{\ism} [1] {\cong_{#1}}
\nek{\izm} [2] {\cong_{#1}^{#2}}
\nek{\aut} [1] {\Aut_{#1}}

\nek{\hfn} {{\rm HF}(\dN)}
\nek{\tce} [1] {{\rm TC}_\ve(#1)}
\nek{\ihf} {\simeq_{\hfn}}

\nek{\symr}{\equiv}
\nek{\rrt} [3] {\symr_{#2#3}^{#1}}
\nek{\rrq} [5] {{#4}\symr_{#2#3}^{#1}{#5}}
\nek{\nrq} [5] {{#4}\not\symr_{#2#3}^{#1}{#5}}
\nek{\rrQ} [5] {{#4}\symr_{#2\,,\,#3}^{#1}{#5}}
\nek{\rro} [5] {\ang{#2,#4}\symr^{#1}\ang{#3,#5}}
\nek{\rrO} [1] {\symr^{#1}}

\nek{\lww} {\cL_{\omi\om}}

\nek{\forc}
{\mathrel{{|\hspace{-0.2ex}|\hspace{-1.1ex}-\hspace{-1.7ex}-}}}

\newlength{\dxii}
\nek{\fC} {{\bf C}}
\nek{\pg} {\fC_\dG}
\nek{\px} {\fC_\dX}

\nek{\gen} {gen.\ }
\nek{\hg}  {h.\hspace{0.4ex}gen.\ }
\nek{\hgp} {h.\hspace{0.4ex}gen.}

\nek{\rfy} {\rF^\iy}

\nek{\incl} [1] {\mathop{\text{\sc Int}}\overline{#1}}  

\nek{\PP} {pinned}
\nek{\PPP}{Pinned}

\nek{\bap}{{\bar p}}
\renek{\wtau} {{\widehat p}}

\nek{\zO} {,\linebreak[0]}
\nek{\zi} {,\linebreak[0]\,}
\nek{\zd} {,\linebreak[0]\:}
\nek{\zT} {,\linebreak[0]\;}
\nek{\yo} {,\linebreak[0]}
\nek{\yi} {,\linebreak[0]\,}
\nek{\yd} {,\linebreak[0]\:}
\nek{\yt} {,\linebreak[0]\;}

\nek{\prit} [1] {[{{\rm #1}}]}
\nek{\nor} [1] {\|#1\|}

\nek{\fap} {f.\hspace{0.1ex}a.\hspace{0.1ex}p.\hspace{0.1ex}m.}

\nek{\bpr} [1] {\bpf[{{\sl#1}\/}]}

\nek{\eqr} {equivalence relation}

\nek{\fx} {{\bf x}}

\nek{\rzd} {^{\text{\tt red}}}

\nek{\srez} [2] {{\bf S}_{#1}(#2)}

\nek{\qc}  [1] {\resi{> #1}}
\nek{\qec} [1] {\resi{\ge #1}}
\nek{\rme} [1] {\resi{\le#1}}

\nek{\alex} {<_{\text{\tt alex}}}
\nek{\lex} {<_{\text{\tt lex}}}
\nek{\act} {<_{\text{\tt act}}}

\nek{\bbo} {\mathbb 0}

\nek{\fz} {{\mathbf z}}

\nek{\dln} {2\lom}

\nek{\fK} {{\bf K}}

\nek{\Ks} {\fK_\fsg}

\nek{\dva}{{\ans{0,1}}}

\nek{\rtn}{\dR^\dN}
\nek{\ntn}{\dN^\dN}
\nek{\ztn}{\dZ^\dN}

\nek{\snos} [1] {\footnote{\ #1}}

\nek{\rH}  {\mathbin{\sf H}}

\nek{\fras}[2] {\text{\footnotesize$\DS\frac{#1}{#2}$}}
\nek{\fral}[2] {\text{\large$\frac{#1}{#2}$}}

\nek{\renu}{\tenu{{\rm(\roman{enumi})}}}

\nek{\ergg}{\aer\dG\dG}

\nek{\rit} [1] {{\it#1\/}}

\nek{\lap} [1] {<<#1>>}

\nek{\mto} {\mapsto}

\nek{\pgcr} {\addtocounter{page}1}

\renek{\pgcr}{}

\makeatletter
\if@twoside%
\else\fi%
\makeatother
         
\begin{document}

\title{{\sc Varia}\\  Ideals and Equivalence Relations}

\author{Vladimir Kanovei}

\date{Feb 2005, Caltech}
\maketitle
                                                                  
\begin{abstract}
A selection of basic results on Borel reducibility of ideals
and \er s, especially those with comparably short proofs.
This is an unfinished text as yet.
Some proofs have missing parts and loose ends. 
\snos
{{\tt kanovei@mccme.ru} and {\tt vkanovei@math.uni-wuppertal.de}
are my contact addresses.}
\end{abstract}

\newpage

{\small
\def\contentsname{\large\bf Contents}
\tableofcontents
}

\parf{Reducibility}
\label{red}

There are several reasonable ways to compare \er s, 
usually formalized in terms of existence of a 
{\it reduction\/}, \ie, a map of certain kind which 
allows to derive one of the \er s from the other one.  
{\it Borel reducibility\/} $\reb$ is the key one, yet 
there are several special types of $\reb,$ in particular, 
those induced by a low-level maps, useful in many cases.
Generally, the most of research on reducibility of Borel 
\er s or ideals is concentrated around the following 
notions of reducibility.

\punk{Borel reducibility}
\label{bored}

If $\rE$ and $\rF$ are \er s on Polish spaces 
resp.\ $\dX,\,\dY,$ then 
\bit
\item[$\ast$]\msur
$\rE\reb\rF$ ({\it Borel reducibility\/})
\index{reducibility!Borel}%
means that there is a Borel map 
\index{zzEleF@$\rE\reb\rF$}%
$\vt:\dX\to\dY$ (called {\it reduction\/}) 
\index{reduction@ reduction of \er s}%
such that ${x\rE y}\eqv {\vt(x)\rF\vt(y)}$ for all 
$x,\,y\in\dX$;

\item[$\ast$]\msur
$\rE\eqb\rF$ iff $\rE\reb\rF$ and $\rF\reb\rE$ 
({\it Borel bi-reducibility\/});  
\index{reducibility!Borel bi-reducibility}%
\index{zzEeqF@$\rE\eqb\rF$}%

\item[$\ast$]\msur
\index{zzElesF@$\rE\rebs\rF$}%
$\rE\rebs\rF$ iff $\rE\reb\rF$ but not $\rF\reb\rE$ 
\index{reducibility!Borel strict}%
({\it strict Borel reducibility\/});  

\item[$\ast$]\msur 
$\rE\emb\rF$ means that there is a Borel 
\index{zzEemF@$\rE\emb\rF$}%
\index{embedding@embedding of \er s}%
{\it embedding\/}, \ie, a $1-1$ reduction;

\item[$\ast$]\msur
$\rE\aprb\rF$ iff $\rE\emb\rF$ and $\rF\emb\rE$ 
(a rare form, \cite[\pff0]{sinf});
\index{zzEemeF@$\rE\aprb\rF$}%

\item[$\ast$]\msur 
$\rE\embi\rF$ means that there is a Borel 
{\it invariant\/} embedding, \ie, an embedding $\vt$ 
\index{embedding!invariant}%
\index{zzEemiF@$\rE\embi\rF$}%
such that $\ran\vt=\ens{\vt(x)}{x\in\dX}$ is an 
\index{set!invariant}%
\ddf{\it invariant\/} set 
(meaning that the \ddf{\it saturation\/} 
$\ekf{\ran\vt}=\ens{y'}{\sus x\:(y\rF\vt(x))}$ 
equals $\ran\vt$);

\item[$\ast$]\msur 
$\rE\ren\rF,\msur$ $\rE\emn\rF,\msur$ $\rE\emni\rF$ 
\index{zzErenF@$\rE\ren\rF$}%
\index{zzEemnF@$\rE\emn\rF$}%
\index{zzEemniF@$\rE\emni\rF$}%
mean that there 
is a \poq{continuous\/} resp.\ reduction, embedding, 
invariant embedding.
\eit 
Sometimes they write ${\dX/\rE}\reb{\dY/\rF}$ instead 
of $\rE\reb\rF$.

\bde
\item[{\bfsl Borel reducibility of ideals\/}:]
$\cI\reb\cJ$ iff $\rei\reb\rej.$ 
Thus it is required that there is a Borel map 
$\vt:\pws A\to\pws B$ such that ${x\sd y}\in \cI$ 
iff ${\vt(x)\sd \vt(y)}\in \cJ.$ 
(Here $\cI$ is an ideal on $A$ and $\cJ$ is an ideal 
on $B.$) 

Versions $\cI\ren\cJ,\msur$ $\cI\emb\cJ,\msur$ 
$\cI\emn\cJ$ have the corresponding meaning. 
\ede

\punk{``Algebraic'' Borel reducibility}
\label{abored}

This is a more special version of Borel reducibility 
of ideals, characterized by the property that the 
reduction must respect a chosen algebraic structure. 
We shall be especially interested in the Boolean 
algebra structure and a weaker \dd\sd group structure 
of sets of the form $\cP(A).$ 
Let $\cI,\,\cJ$ be ideals on resp.\ $A,\,B$.

\bde
\item[{\bfsl Borel BA reducibility\/}:]
$\cI\rba\cJ$ if there is a Borel \ddj approximate  
Boolean algebra homomorphism $\vt:\pws A\to\pws B$ 
with $x\in\cI\eqv \vt(x)\in\cJ$. 
 
A version: 
$\cI\rbap\cJ$ if there is a set $A\in\pl\cJ$ with 
$\cI\rba(\cJ\res A)$.
\ede
Here, $\vt:\pws A\to\pws B$ is an \ddj{\it approximate\/} 
Boolean algebra homomorphism if the sets 
$\skl\vt(x)\cup\vt(y)\skp\sd\vt(x\cup y)$ and 
$\vt(\dop x)\sd \dop(\vt(x))$ always belong to $\cJ$ 
whenever $x,\,y\sq A.$ 
Let further a \ddj approximate \dd\sd{\it homomorphism\/} 
be any map $\vt:\pws A\to\pws B$ such that 
$(\vt(x)\sd\vt(y))\sd\vt(x\sd y)$ 
always belongs to $\cJ.$
This leads to a weaker reducibility:  
\bde
\item[{\bfsl Borel \dd\sd reducibility\/}:]
$\cI\rds\cJ$ iff there is a Borel \ddj approximate  
\dd\sd homo\-mor\-phism $\vt:\pws A\to\pws B$ 
such that $x\in\cI\eqv \vt(x)\in\cJ$.
\ede

\punk{Borel, continuous, and Baire measurable reductions}
\label{b-c}

Many properties of Borel reductions 
hold for a bigger family of Baire measurable 
(\bm, for brevity) maps. 
Any reducibility definition in 
\prff{bored},~\ref{abored} 
admits a weaker \bm\ version, which claims that the reduction 
postulated to exist is only \bm, not necessarily Borel. 
Such a version will be denoted with a subscript BM instead of 
B, for instance, $\rE\reB\rF$ means that there is a \bm\ 
reduction, \ie, a \bm\ map 
$\vt:\dX=\dom\rE\to\dY=\dom\rF$ 
such that ${x\rE y}\eqv {\vt(x)\rF\vt(y)}$ for all 
$x,\,y\in\dX$.

On the other hand, a \poq{continuous} reducibility can 
sometimes be derived. 

\vyk{
\ble[{{\rm Louveau ?}}]
\label{l:bc}
{} \ 
If\/ $\cI\zd\cJ$ are Borel ideals on\/ $\dN,$ 
and\/ $\cI\reB\cJ,$ then\/ 
$\cI\ren{\dsu\cJ\cJ}$ 
{\rm(via a continuous reduction)}, 
that is~\footnote
{\ See the definition of $\dsu\cI\cJ$ in \prf{opeid}. 
${\dsu\cJ\cJ}\ren\cJ$ for many practically 
interesting ideals.}, 
there exist continuous maps\/ 
$\vt_0$ and\/ $\vt_1:\pn\to\pn$ such that, for any\/ 
$x,\,y\in\pn,$ 
$x\sd y\in\cI$ iff both\/ $\vt_0(x)\sd\vt_0(y)\in\cJ$ 
and\/ $\vt_1(x)\sd\vt_1(y)\in\cJ$. 
\ele
\bpf
Let $\vt:\pn\to\pn$ witness $\cI\reB\cJ.$  
As any \bm\ map, $\vt$ is continuous on a dense $\Gd$ set 
$D=\bigcap_i D_i\sq\pn,$ all $D_i$ being dense open and 
$D_{i+1}\sq D_i.$  
We can define a sequence 
$0=n_0<n_1<n_2<...$ and, for every $i,$ a set 
$u_i\sq \il{i}{i+1}$ such that 
$x\cap \il{i}{i+1}=u_i\imp x\in D_i.$~\footnote
{\ Sets like $u_i$ are called {\it stabilizers\/},
they are of much help in study of Borel ideals.}  
Let  
\dm
\textstyle
N_0=\bigcup_i \il{2i}{2i+1}\,,\quad 
N_1=\bigcup_i \il{2i+1}{2i+2}\,,\quad 
U_0=\bigcup_i u_{2i}\,,\quad U_1=\bigcup_i u_{2i+1}\,.
\dm
Now set 
$\vt_0(x)=\vt((x\cap N_0)\cup S_1)$ and 
$\vt_1(x)=\vt((x\cap N_1)\cup S_0)$ for $x\sq\dN$.
\epf
}

\ble[{{\rm Louveau ?}}]
\label{l:bc}
{} \ 
If\/ $\cI$ is a Borel ideal on a countable\/ $A,$
$\rE$ an \eqr\ on a Polish\/ $\dX,$ and\/ 
$\rei\reB\rE,$ then\/ 
$\rei\ren{\rE\ti\rE}$ 
{\rm(via a continuous reduction)}, 
that is,
there exist continuous maps\/ 
$\vt_0\zd\vt_1:\pws A\to\dX$ such that, for any\/ 
$x,\,y\in\pn,$ 
$x\sd y\in\cI$ iff both\/ $\vt_0(x)\rE\vt_0(y)$ 
and\/ $\vt_1(x)\rE\vt_1(y)$. 
\ele
\bpf
We \noo\ suppose that $A=\dN.$ 
Let $\vt:\pn\to\dX$ witness that $\rei\reB\rE.$  
Then $\vt$ is continuous on a dense $\Gd$ set 
$D=\bigcap_i D_i\sq\pn,$ all $D_i$ dense open and 
$D_{i+1}\sq D_i.$  
A sequence $0=n_0<n_1<n_2<\dots$ and, for any $i,$ a set 
$u_i\sq \il{i}{i+1}$ can be easily defined, by
induction on $i,$ so that 
$x\cap \il{i}{i+1}=u_i\imp x\in D_i.$\footnote
{\ Sets like $u_i$ are called {\it stabilizers\/}, 
\index{stabilizer}%
they are of much help in study of Borel ideals.}  
Let  
\dm
\textstyle
N_1=\bigcup_i \il{2i}{2i+1}\,,\;\; 
N_2=\bigcup_i \il{2i+1}{2i+2}\,,\;\; 
U_1=\bigcup_i u_{2i}\,,\;\; U_2=\bigcup_i u_{2i+1}\,.
\dm
Now set 
$\vt_1(x)=\vt((x\cap N_1)\cup U_2)$ and 
$\vt_2(x)=\vt((x\cap N_2)\cup U_1)$ for $x\sq\dN$.
\epf

The following question should perhaps be answered in the 
negative in general and be open for some particular cases.

\bqu
Suppose that $\rE\reb \rF$ are Borel \er s. 
Does there always exist a \poq{continuous} reduction ? 
\equ

\punk{Reducibility via maps between the underlying sets}
\label{reli}

This is an even more special kind of Borel reducibility. 
Let $\cI,\,\cJ$ be ideals on resp.\ $A,\,B,$ as above.

\bde
\item[{\bfsl Rudin--Keisler order\/}:]
$\cI\ork\cJ$ iff there exists a function $b:\dN\to\dN$ 
(a {\it Rudin--Keisler\/} reduction) 
such that $x\in\cI\eqv b\obr(x)\in\cJ$. 

\item[{\bfsl Rudin--Blass order\/}:]
$\cI\orb\cJ$ iff there is a \poq{finite-to-one} 
function $b:\dN\to\dN$ 
(a {\it Rudin--Blass\/} reduction) 
with the same property.

A version: $\cI\orbp\cJ$ allows $b$ to be defined on a 
proper subset of $\dN,$ in other words, we have pairwise 
disjoint finite non-empty sets $w_k=b\obr(\ans{k})$ 
such that $x\in\cI\leqv w_x=\bigcup_{k\in x}w_k\in\cJ.$

Another version: 
$\cI\orbpp\cJ$ requires that, in addition, 
the sets $w_k=b\obr(\ans{k})$ satisfy 
$\tmax w_k<\tmin w_{k+1}$.
\ede

There is a ``clone'' of the Rudin--Blass order which 
applies in a much more general situation. 
Suppose that $X=\prod_{k\in\dN}X_k$ and 
$Y=\prod_{k\in\dN}Y_k,$ 
$0=n_0<n_1<n_2<...,$ and 
$H_i:X_i\to\prod_{n_i\le k<n_{i+1}}Y_k$ for any $i.$ 
Then, we can define 
\dm
\Psi(x)=H_0(x_0)\cup H_1(x_1)\cup H_2(x_2)\cup...\in Y
\dm
for each $x=\sis{x_i}{i\in\dN}\in X.$ 
Maps $\Psi$ of this kind were called {\it additive\/} 
by Farah~\cite{f-co}. 
More generally, if, in addition, $0=m_0<m_1<m_2<...,$ and 
$H_i:\prod_{m_i\le j<m_{i+1}}X_j\to
\prod_{n_i\le k<n_{i+1}}Y_k$ for any $i,$ 
then we can define 
\dm
\Psi(x)=H_0(x\res\ir{m_0}{m_1})
\cup H_1(x\res\ir{m_1}{m_2})
\cup H_2(x\res\ir{m_2}{m_3})\cup...\in Y
\dm
for each $x\in X.$ 
Farah calls maps $\Psi$ of this kind 
{\it asymptotically additive\/}. 
All of them are Borel functions $X\to Y,$ 
provided all sets $X_j$ and $Y_k$ are finite. 

Suppose now that $\rE$ and $\rF$ are \er s on resp.\ 
$X=\prod_kX_k$ and $Y=\prod_kY_k$.

\bde
\item[{\bfsl Additive reducibility\/}:]
$\rE\rea\rF$ if there is an additive reduction $\rE$ 
to $\rF$.
$\rE\reaa\rF$ if there is an asymptotically 
additive reduction $\rE$ to $\rF$.
\ede

\ble[{{\rm Farah~\cite{f-co}}}] 
\label{arb} 
Suppose that\/ $\cI$ and\/ $\cJ$ are Borel ideals on\/ 
$\dN.$ 
Then\/ $\cI\orbpp\cJ$ iff\/ $\rei\rea\rej$.
\ele
(By definition $\rei$ and $\rej$ are \er s on $\pn,$ 
yet we can consider them as \er s on 
$\dn=\prod_{k\in\dN}\ans{0,1},$ 
as usual, which yields the intended 
meaning for $\rei\rea\rej.$)
\bpf
If $\cI\orbpp\cJ$ via a sequence of finite sets $w_i$ with 
$\tmax w_i<\tmin w_{i+1}$ then we put $n_0=0$ and 
$n_i=\tmin w_i$ for $k\ge1,$ 
so that $w_i\sq \il{i}{i+1},$ and, for any $i,$ 
put $H_i(0)=\il{i}{i+1}\ti\ans0$ and let $H_i(1)$ be 
the characteristic function of $w_i$ within $\il{i}{i+1}.$ 
Conversely, if $\rei\rea\rej$ via a sequence 
$0=n_0<n_1<n_2<...$ and a family of maps 
$H_i:\ans{0,1}\to2^{\il{i}{i+1}}$ then 
$\cI\orbpp\cJ$ via the sequence of sets  
$w_i=\ens{k\in\il{i}{i+1}}{H_i(0)(k)\ne H_i(1)(k)}$.
\epf

The following definition is taken from \cite{jkl}. 
Let $\cI,\,\cJ$ be ideals on $\dN$.

\bde
\item[{\bfsl Reducibility via inclusion\/}:] 
$\cI\inc\cJ$ if there is a map $b:\dN\to\dN$ such  
that $x\in\cI\imp b\obr(x)\in\cJ.$ 
(Note $\imp$ instead of $\eqv$!)
\ede

In particular if $\cI\sq\cJ$ then $\cI\inc\cJ$ via 
$b(k)=k.$ 
It follows that this order is not fully compatible with 
$\reb$ because $\sui{1/n}\sq\Zo$ while the summable 
ideal $\sui{1/n}$ and the density-0 ideal 
$\Zo$ are known to be \dd\reb incomparable.

\punk{Isomorphism}
\label{1isom}

Let $\cI,\,\cJ$ be ideals on resp.\ $A,\,B.$ 
{\it Isomorphism\/} $\cI\isi\cJ$ means that there is a 
\index{isomorphism!icj@$\cI\isi\cJ$}%
\index{zzicj@$\cI\isi\cJ$}%
bijection $\ba:A\onto B$ such that we have 
${x\in \cI}\leqv{\ba\ima x\in\cJ}$ for 
all $x\sq A$. 
 
Sometimes they use a weaker definition: let
\index{isomorphism!icja@$\cI\isa\cJ$}%
\index{zzicja@$\cI\isa\cJ$}%
$\cI\isa\cJ$ mean that there are sets $A'\in\df\cI$ 
and $B'\in\df\cJ$ such that 
${\cI\res A'}\isi{\cJ\res B'}.$ 
Yet this implies $\cI\isi\cJ$ in most usual cases, 
the only notable exception (among nontrivial ideals), 
is produced by the ideals $\cI=\ifi$ and 
$\cJ=\dsu\ifi{\pn}\isi\ens{x\sq\dN}{x\cap D\in\ifi},$ 
where $D$ is an infinite and coinfinite set~\footnote
{\ Kechris~\cite{rig} called ideals $\cJ$ of this kind 
{\it trivial variations of\/ $\ifi$}.}
: then $\cI\isa\cJ$ but not $\cI\isi\cJ$.

\punk{Remarks}

\imar{check this subsection once again}%
The following shows simple relationships between different 
reducibilities:   %
\dm
{\cI\orb\cJ}\simp{\cI\ork\cJ}\simp{\cI\obe\cJ}
\simp{\cI\obep\cJ}\simp{\cI\odl\cJ}\simp{\cI\reb\cJ}.
\dm
For instance if 
$b:\dN\to\dN$ witnesses $\cI\ork\cJ$ then 
$\vt_b(X)=b\obr(X)$ witnesses $\cI\obe\cJ.$ 
Note that any $\vt_b$ is an exact Boolean algebra 
homomorphism $\pn\to\pn;$ moreover, it is known that 
any BM Boolean algebra homomorphism $\pn\to\pn$ is 
$\vt_b$ for an appropriate $b:\dN\to\dN.$ 
{\it Approximate\/} homomorphisms are liftings of 
homomorphisms into quotients of $\pn,$ thus, any 
\ddj approximate $\vt:\pn\to\pn$ induces the map 
$\Theta(X)=\ens{\vt(X)\sd Y}{Y\in\cJ},$ which is a 
homomorphism $\pn\to\pn/\cJ.$ 
Farah~\cite{aq}, and Kanovei and Reeken \cite{kr} 
demonstrated that in some important cases 
(of ``nonpatological'' P-ideals and, generally, for all 
Fatou, or Fubini, ideals) 
we have ${\cI\ork\cJ}\eqv{\cI\obe\cJ}.$ 
On the other hand ${\cI\ork\cJ}\nmp{\cI\obe\cJ}$ fails 
for rather artificial P-ideals. 

The right-hand end is the most intrigueing: 
is there a pair of Borel ideals $\cI,\,\cJ$ such that  
$\cI\reb\cJ$ but not $\cI\odl\cJ\,?$ 
If we actually have the equivalence then the 
whole theory of Borel reducibility for Borel ideals 
can be greatly simplified 
because reduction maps which are \dd\sd homomorphisms 
are much easier to deal with.

\parf{Introduction to ideals}
\label{somi}

As many interesting \er s appear as $\rei$ for a 
Borel ideal $\cI,$ we take space to discuss a few 
basic items related to Borel ideals. 
We begin with several examples and notation, 
and then continue with 
some important types of ideals.

\bit
\item
$\ifi=\ens{x\sq\dN}{x\,\text{ is finite}},$ the 
ideal of all finite sets;

\item
$\Ii=\fio=
\ens{x\sq\dN^2}{\ens{k}{\seq xk\ne\pu}\in\ifi}$;
\imar{gde vvedeno $\seq xk$?}%

\item
$\Id=\sui{1/n}=
\ens{x\sq \dN}{\sum_{n\in x}\frac1{n+1}}<\piy,$ 
the {\it summable ideal\/};

\item
$\It=\ofi=\ens{x\sq\dN^2}{\kaz k\:(\seq xk\in\ifi)}$;

\item
$\cZ_0=\eui{1}=\ens{x\sq\dN}
{\tlim_{n\to\piy}\frac{\#(x\cap\ir0n)}n=0},$ 
the {\it density ideal\/}.
\eit

\punk{Notation}
\label{opeid}

\bit
\item
For any ideal $\cI$ on a set $A,$ we define 
$\pl\cI=\cP(A)\dif\cI$ (\ddi{\it positive\/} sets) and 
$\df\cI=\ens{X}{\dop X\in\cI}$ ({\it the dual filter\/}). 
Clearly $\pu\ne\df\cI\sq\pl\cI$.

\item
If $B\sq A,$ then we put  
$\cI\res B=\ens{x\cap B}{x\in\cI}$.

\item
If $\cI,\,\cJ$ are ideals on resp.\ $A,\,B,$ 
then $\dsu\cI\cJ$ (the {\it disjoint sum\/}) is the 
ideal of all sets 
$x\sq C=(\ans0\ti A)\cup(\ans1\ti B)$ with 
$\seq x0\in \cI$ and $\seq x1\in\cJ$ 
(where $\seq xi=\ens{c}{\ang{i,c}\in x},$ as usual).

If the sets $A,\,B$ are disjoint then $\dsu\cI\cJ$ 
can be equivalently defined as the 
ideal of all sets $x\sq A\cup B$ with 
$x\res A\in \cI$ and $x\res B\in\cJ$.

\item
The {\it Fubini product\/} $\fps{\cI}{\cJ_a}{a\in A}$ 
\index{Fubini product!of ideals}%
of ideals $\cJ_a$ on sets 
$B_a,$ over an ideal $\cI$ on a set $A$ 
is the ideal on the set 
$B=\ens{\ang{a,b}}{a\in A\land b\in B_a},$ which 
consists of all sets $y\sq B$ such that the set 
$\ens{a}{(y)_a\nin \cJ_a}$ belongs to $\cI,$ where 
$(y)_a=\ens{b}{\ang{a,b}\in y}$ (the cross-section). 

\item
In particular, the Fubini product $\fpd\cI\cJ$ of two
ideals $\cI,\cJ$ on sets resp.\ $A,B,$ is equal to
$\fps\cI{\cJ_a}{a\in A},$ 
where $\cJ_a=\cJ\zd\kaz a.$ 
Thus $\fpd\cI\cJ$ consists of all sets $y\sq A\ti B$ 
such that  
$\ens{a}{(y)_a\nin \cJ}\in \cI$. 
\eit

\punk{P-ideals and submeasures}
\label{pideals}

Many important Borel ideals belong to the class of 
P-ideals. 

\bdf
\label{api}
An ideal $\cI$ on $\dN$ is 
a {\it P-ideal\/} if for any sequence 
of sets $x_n\in\cI$ there is a set $x\in\cI$ such that 
$x_n\sqa x$ (\ie, $x_n\dif x\in\ifi$) for all $n$;
\edf

For instance, the ideals $\ifi,\,\Id,\,\It,\,\Zo$ 
(but not $\Ii$!) are P-ideals. 

This class admits several apparently different 
but equivalent characterizations, one of which is 
connected with submeasures. 
\bit
\item
A {\it submeasure\/} on a set $A$ is any map 
\index{submeasure}%
$\vpi:\cP(A)\to[0,\piy],$ satisfying $\vpi(\pu)=0,$ 
$\vpi(\ans a)<\piy$ for all $a,$
and $\vpi(x)\le \vpi(x\cup y)\le \vpi(x)+\vpi(y)$.

\item
A submeasure $\vpi$ on $\dN$ is 
{\it lover semicontinuous\/}, 
\index{submeasure!lover semicontinuous}%
\index{submeasure!lsc@\lsc}%
or \lsc\ for brevity, if we have 
$\vpi(x)=\tsup_n\vpi(x\cap\ir0n)$ for all $x\in\pn$. 
\eit
 
To be a {\it measure\/}, a submeasure $\vpi$ has to 
satisfy, in addition, that 
$\vpi(x\cup y)=\vpi(x)+\vpi(y)$ whenever $x,\,y$ 
are disjoint. 
Note that any \dds additive measure is \lsc, 
but if $\vpi$ is \lsc\ then $\vpy$ is not 
necessarily \lsc\ itself.

Suppose that $\vpi$ is a submeasure on $\dN.$ 
Define the {\it tailsubmeasure\/} 
\index{tailsubmeasure}%
\index{submeasure!tailsubmeasure}%
$\vpy(x)=\hv x\vpi=\tinf_n(\vpi(x\cap\iry n)).$ 
The following ideals are considered:
\dm
\bay{rcllll}
\Fin_\vpi &=& \ens{x\in\pn}{\vpi(x)<\piy} &&&;\\[1ex]
\index{zzfinf@$\Fin_\vpi$}%

\Nul_\vpi &=& \ens{x\in\pn}{\vpi(x)=0} &&&;\\[1ex]
\index{zznulf@$\Nul_\vpi$}%

\index{zzexhf@$\Exh_\vpi$}%
\Exh_\vpi &=& \ens{x\in\pn}{\vpy(x)=0}&=
& \Nul_{\vpy}&.\\[0.5ex]
\eay
\dm

\bex
\label{exh:e}
$\ifi=\Exh_\vpi=\Nul_\vpi,$ where $\vpi(x)=1$ for any 
$x\ne \pu.$ 
We also have $\ofi=\Exh_\psi,$ where 
$\psi(x)=\sum_k\,2^{-k}\,\vpi(\ens{l}{\ang{k,l}\in x})$ 
is \lsc.
\eex

It turns out (Solecki, see Theorem~\ref{sol} below) that 
analytic P-ideals are the same as ideals of the form 
$\Exh_\vpi,$ where $\vpi$ is a 
\lsc\ submeasure on $\dN.$ 
It follows that any analytic P-ideal is $\fp03$.

\punk{Polishable ideals} 
\label{punk:poli}

There is one more characterization of Borel P-ideals. 
Let $T$ be the ordinary Polish product topology on $\pn.$ 
Then $\pn$ is a Polish group in the sense of $T$ and the 
symmetric difference as the operation, and any ideal 
$\cI$ on $\dN$ is a subgroup of $\pn$. 

\bdf
\label{polI}
An ideal $\cI$ on $\dN$ is {\it polishable\/} if 
there is a Polish group topology $\tau$ on $\cI$ which 
produces the same Borel subsets of $\cI$ as $T\res\cI$.
\edf

The same Solecki's theorem (Theorem~\ref{sol}) proves 
that, for analytic ideals, to be a P-ideal is the same 
as to be polishable. 
It follows (see Example~\ref{exh:e}) that, for instance, 
$\ifi$ and $\It=\ofi$ are polishable, but $\Ii=\fio$ 
is not. 
The latter will be shown directly after the next lemma.

\ble
\label{sol:?}
Suppose that an ideal\/ $\cI\sq\pn$ is polishable. 
Then there is only one Polish group topology\/ 
$\tau$ on\/ $\cI.$ 
This topology\/ refines\/ $T\res\cI$ 
and is metrizable by a\/ \dd\sd invariant metric. 
If\/ $Z\in\cI$ then\/ $\tau\res\cP(Z)$ coincides 
with\/ $T\res\cP(Z).$ 
In addition, $\cI$ itself is\/ \dd TBorel.
\ele
\bpf
Let $\tau$ witness that $\cI$ is polishable. 
The identity map 
$f(x)=x{:}\,\stk\cI\tau\to\stk\pn T$ 
is a \dd\sd homomorphism and is Borel-measurable 
because all \dd{(T\res\cI)}open sets are \dd\tau Borel, 
hence, by the Pettis theorem (Kechris~\cite[??]{dst}), 
$f$ is continuous. 
It follows that all \dd{(T\res\cI)}open subsets of 
$\cI$ are \dd\tau open, and that $\cI$ is 
\dd TBorel in $\pn$ because $1-1$ continuous 
images of Borel sets are Borel. 

A similar ``identity map'' argument shows that 
$\tau$ is unique if exists. 

It is known (Kechris~\cite[]{dst}) 
that any Polish group topology admits a 
left-invariant compatible metric, which, in this case, 
is right-invariant as well since $\sd$ is an 
abelian operation.

Let $Z\in\pn.$ 
Then $\cP(Z)$ is \dd Tclosed, hence, 
\dd\tau closed by the above, subgroup of $\cI,$ and 
$\tau\res\cP(Z)$ is a Polish group topology on $\cP(Z).$ 
Yet $T\res\cP(Z)$ is another Polish group topology on 
$\cP(Z),$ with the same Borel sets. 
The same ``identity map'' argument proves  
that $T$ and $\tau$ coincide on $\cP(Z)$.
\epf

\bex
\label{fio}
$\cI_1=\fio$ is not polishable. 
Indeed we have $\fio=\bigcup_n W_n,$ where  
$W_n=\ens{x}{x\sq\ans{0,1,...,n}\ti\dN}.$  
Let, on the contrary, $\tau$ be a Polish group topology 
on $\cI_1.$ 
Then $\tau$ and the ordinary topology $T$ coincide on 
each set $W_n$ by the lemma, in particular, 
each $W_n$ remains \dd\tau nowhere dense 
in $W_{n+1},$ hence, in $\cI_1,$ a contradiction with 
the Baire category theorem for $\tau$.
\eex

\punk{Some $\Fs$ ideals}

Any sequence $\sis{r_n}{n\in\dN}$ of positive reals $r_n$ 
with $\sum r_n=\piy$ defines the ideal 
\dm
\sui{r_n}=\ens{X\sq\dN}{\sum_{n\in X}r_n<\piy}=
\ens{X}{\mrn(X)<\piy}\,,
\dm 
where $\mrn(X)=\sum_{n\in X}r_n.$
These ideals are called {\it summable ideals\/};
all of them are $\Fs.$ 
References \cite{mat72,maz91,aq}. 
Any summable ideal is easily a P-ideal:  
indeed, $\sui{r_n}=\Exh_\vpi,$ where 
$\vpi(X)=\sum_{n\in X}r_n$ is a \dds additive measure. 

Summable ideals are perhaps the easiest to study among 
all P-ideals. 
Further entries: 
1) Farah \cite[\pff1.12]{aq} on summable 
ideals under $\obe,$ 
2) Hjorth: \dd\reb structure 
of ideals \dd\reb reducible to summable ideals, in 
\cite{h-ban}.

\ble[{{\rm Folklore ?}}]
\label{sumi}
Suppose that\/ $r_n\ge 0,\msur$ $r_n\to0,$ and\/ 
$\sum_nr_n=\piy.$ 
Then any\/ summable ideal\/ $\cI$ satisfies\/ 
$\cI\orbpp\sui{r_n}.$
\ele
\bpf
Let $I=\sui{p_n},$ where $p_n\ge0$ 
(no other requirements !). 
Under the assumptions of the lemma we can associate a 
finite set $w_n\sq\dN$ to any $n$ so that 
$\tmax w_n<\tmin w_{n+1}$ and 
$|r_n-\sum_{j\in w_n}r_i|<2^{-n}.$ 
\epf

Farah \cite[\pff 1.10]{aq} defines a non-summable $\Fs$ 
P-ideal as follows. 
Let $I_k=\ir{2^k}{2^{k+1}}$ and 
$\psi_k(s)=k^{-2}\tmin\ans{k,\# s}$ for all $k$ and 
$s\sq I_k,$ and then
\dm
\psi(X)=\sum_{k=0}^\iy\psi_k(X\cap I_k)
\quad\hbox{and}\quad
\cI=\Fin_\psi\;;
\dm
it turns out that $\cI$ is an $\Fs$ P-ideal, but not 
summable. 
To show that $\cI$ distincts from any $\sui{r_n},$ 
Farah notes that there is a set $X$ 
(which depends on $\sis{r_n}{}$) 
such that the differences 
$|\mrn(X\cap I_k)-\psi_k(X\cap I_k)|,\msur$ 
$k=0,1,2,...\,,$ 
are unbounded. 

Further entry: Farah~\cite{f-tsir,f-bas,f-co} 
on {\it Tsirelson ideals\/}.

\punk{Erd\"os -- Ulam and density ideals}
\label{eu}

These are other types of Borel P-ideals.
Any sequence $\sis{r_n}{n\in\dN}$ of positive reals $r_n$ 
with $\sum r_n=\piy$ defines the ideal 
\dm
\eui{r_n}=\left\{x\sq\dN:\tlim_{n\to\piy}
\frac{\sum_{i\in x\cap\ir0n}r_i}
{\sum_{i\in \ir0n}r_i}=0\right\}.
\dm 
These ideals are called 
{\it Erd\"os -- Ulam\/ {\rm(or: EU)} ideals\/}. 
Examples: $\zo=\eui1$ and $\cZ_{\tt log}=\eui{1/n}$. 

This definition can be generalized. 
Let $\supp\mu=\ens{n}{\mu(\ans n)>0},$ for any measure 
$\mu$ on $\dN.$ 
Measures $\mu,\,\nu$ are {\it orthogonal\/} if we have 
$\supp\mu\cap\supp\nu=\pu.$  
Now suppose that $\vmu=\sis{\mu_n}{n\in\dN}$ is a sequence 
of pairwise orthogonal measures on $\dN,$ with finite 
sets $\supp\mu_i.$ 
Define $\vpi_\vmu(X)=\tsup_n\mu_n(X):$ this is a \lsc\ 
submeasure on $\dN.$ 
Let finally 
$\cD_\vmu=\Exh(\vpi_\vmu)=\ens{X}{\hv X{\vpi_\mu}=0}.$ 
Ideals of this form are called {\it density ideals\/} 
by Farah~\cite[\pff 1.13]{aq}. 
This class includes all EU ideals 
(although this is not immediately transparent), 
and some other ideals: for instance, $\ofi$ is a 
density but non-EU ideal.  
Generally density ideals are more complicated than 
summables. 
We obtain an even wider class if the requirement, 
that the sets $\supp\mu_n$ are finite, is dropped: 
this wider family includes all summmable ideals, too.

References \cite{jk84}, \cite[\pff1.13]{aq}. 

Further entries: 
1) Farah: structure of density ideals under $\obe,$ 
2) Farah: \dd{c_0}equalities, 
3) Relation to Banach spaces: Hjorth, SuGao. 

Which ideals are both summable and density ?

\punk{Some transfinite sequences of Borel ideals}
\label{non-P}

We consider three interesting families of Borel 
ideals (mainly, non-P-ideals), united by their 
relation to countable ordinals.  
Note that the underlying sets of the ideals below 
are countable sets different from $\dN$.\vom

{\bfit Fr\'echet ideals.\/} 
This family consists of ideals $\frt_\xi,\msur$ 
$\xi<\omi,$ obtained by inductive 
construction using Fubini products. 
We put $\frt_1=\ifi$ and 
$\frt_{\xi+1}=\fpd\ifi{\frt_\xi}$ for all $\xi.$ 
Limit steps cause a certain problem. 
The most natural idea would be to define 
$\frt_\la=\fps{\ifi_\la}{\frt_\xi}{\xi<\la}$ for any limit 
$\la,$ where $\ifi_\la$ is the ideal of all finite subsets 
of $\la,$ or perhaps 
$\frt_\la=\fps{\ibo_\la}{\frt_\xi}{\xi<\la},$ 
where $\ibo_\la$ is the ideal of all bounted subsets 
of $\la,$ or even $\frt_\la=\fps{0}{\frt_\xi}{\xi<\la},$ 
where $0$ is the ideal containing only the empty set, 
yet this appears not to be fully satisfactory in \cite{jkl}, 
where they define 
$\frt_\la=\fps{\ifi}{\frt_{\xi_n}}{n\in\dN},$ where 
$\sis{\xi_n}{}$ is a once and for all fixed cofinal increasing 
sequence of ordinals below $\la,$ with understanding that 
the result is independent of the choice of $\xi_n,$ modulo 
a certain equivalence.\vom

{\bfit Indecomposable ideals.\/} 
Let $\otp X$ be the order type of $X\sq\Ord.$ 
For any ordinals $\xi,\,\vt<\omi$ define: 
\dm
\iz\vt\xi \,=\, \ens{A\sq \vt}{\otp A<\om^\xi} 
\quad\hbox{(nontrivial only if $\vt\ge\om^\xi$)}\,.
\dm
To see that the sets $\iz\vt\xi$ are really ideals 
note that ordinals of the form $\om^\xi$ and only 
those ordinals are {\it indecomposable\/}, \ie, 
are not sums of a pair of smaller ordinals, hence, 
the set $\ens{A\sq \vt}{\otp A<\ga}$ is 
an ideal iff $\ga=\om^\xi$ for some $\xi.$

\vyk{
Note that any $\jz\xip$ is isomorphic to 
the Fubini product $\ifi\ti\jz\xi:$ 
for example, $\jz1=\ifi$ and $\jz2\cong\ifi\ti\ifi,$ 
yet the recursion at limit steps is not that 
easy. 
}

\vom

{\bfit Weiss ideals.\/} 
Let $\rkb X$ be the {\it Cantor-Bendixson rank\/} 
of $X\sq\Ord,$ 
\ie, the least ordinal $\al$ such that $\kb X\al=\pu.$ 
Here $\kb X\al$ is defined by induction on $\al:$ 
$\kb X0=X,\msur$ 
$\kb X\la=\bigcap_{\al<\la}\kb X\al$ at limit steps 
$\la,$ and finally $\kb X{\al+1}=(\kb X\al)',$ where  
$A',$ the Cantor-Bendixson 
derivative, is the set of all ordinals $\ga\in x$ which 
are limit points of $X$ in the interval topology.  
For any ordinals $\xi,\,\vt<\omi$ define: 
\dm
\iw\vt\xi \,=\, \ens{A\sq \vt}{\rkb A<\om^\xi}  
\quad\text{(nontrivial only if $\vt\ge\omm\xi$)}\,.
\dm
It is less transparent that all $\iw\vt\xi$ are 
ideals 
(Weiss, see Farah~\cite[\pff1.14]{aq}) 
while $\ens{A\sq \vt}{\rkb A<\ga}$ is not an ideal 
if $\ga$ is not of the form $\om^\xi$.

\punk{``Other'' ideals}
\label{o}

This title intends to include those interesting ideals 
which have not yet been subject of comprehensive study. 
A common method to obtain interesting ideals is to 
consider a countable set bearing a nontrivial structure, 
as the underlying set. 
In principle, there is no difference between different 
countable set as which of them is taken as the 
underlying set for the ideals considered. 
Yet if the set bears a nontrivial structure 
(\ie, more than just countability) 
then this 
gives additional insights as which ideals are 
meaningful. 
This is already transparent for the ideals defined 
in \prf{non-P}.

We give two examples. \vom

{\bfit Ideals on finite sequences.\/} 
The set $\nse$ of all finite sequences of natural 
numbers is countable, yet its own order structure is 
quite different from that of $\dN.$ 
We can exploit this in several ways, for instance, 
with ideals of sets $X\sq\nse$ which intersect every 
branch in $\nse$ by a set which belongs to a given 
ideal~on~$\dN$.

\parf{Introduction to \eqr s}
\label{ers}

The structure of Borel and analytic \er s under $\reb$ 
includes key \er s which play distinguished role. 
The plan of this section is to define some of them and 
outline their properties, then introduce some classes 
of \er s.

\punk{Basic equivalence relations}
\label{someq}

Equalities can be considered as the most elementary type
of \er s.
Let $\rav X$ denote the equality on a set $X,$ 
\index{equivalence relation, ER!DN@$\rav\dN=\alo$}%
\index{equivalence relation, ER!D2N@$\rav\dn=\cont$}%
\index{zzDN@$\rav\dN=\alo$}%
\index{zzD2N@$\rav\dn=\cont$}%
considered as an \eqr\ on $X.$
\index{equivalence relation, ER!equality@equality $\rav X$}%
\index{zzDX@$\rav X$}%

A much more diverse family is made of \eqr s generated
by ideals.
Recall that for any ideal $\cI$ on a set $A,$ $\rei$ is 
an \er\ on $\pws A,$ 
defined so that $X\rei Y$ iff ${X\sd Y}\in\cI.$  
\index{zzei@$\rei$}%
Equivalently, $\rei$ can be considered as an \er\ on
$2^A$ defined so that $f\rei g$ iff ${f\sd g}\in\cI,$
where $f\sd g=\ens{a\in A}{f(a)\ne g(a)}.$
Note that $\rei$ is Borel provided so is $\cI.$

This leads us to the following all-important \er s:
\bit
\item\msur
$\Eo=\rE_{\ifi},$ thus, $\Eo$ is a \er\ on $\pn$ 
\index{equivalence relation, ER!E0@$\Eo$}%
\index{zzE0@$\Eo$}%
and $x\Eo y$ iff $x\sd y\in\ifi$.

\item\msur
$\Ei=\rE_{\Ii},$ thus, $\Ei$ is a \er\ on 
$\pws{\dN\ti\dN}$ 
\index{equivalence relation, ER!E1@$\Ei$}%
\index{zzE1@$\Ei$}%
and $x\Eo y$ iff $\seq xk=\seq yk$ for all but finite 
$k,$ where, we recall, $\seq xk=\ens{n}{\ang{k,n}\in x}$ 
for $x\sq \dN\ti\dN.$ 

\item\msur
$\Ed=\rE_{\Id},$ thus, $\Ed$ is a \er\ on $\pn$ 
\index{equivalence relation, ER!E2@$\Ed$}%
\index{zzE2@$\Ed$}%
and $x\Ed y$ iff $\sum_{k\in x\sd y}k\obr<\iy.$ 

\item\msur
$\Et=\rE_{\It},$ thus, $\Ei$ is a \er\ on 
$\pws{\dN\ti\dN}$ 
\index{equivalence relation, ER!E3@$\Et$}%
\index{zzE3@$\Et$}%
and $x\Et y$ iff $\seq xk\Eo\seq yk\zd\kaz k.$
\eit
Alternatively, $\Eo$ can be viewed as an \eqr\ on 
$\dn$ defined as $a\Ei b$ iff $a(k)=b(k)$ 
for all but finite $k.$ 
Similarly, $\Ei$ can be viewed as a \er\ on 
$\pnqn,$ or even on $\dnqn,$
defined as $x\Ei y$ iff $x(k)=y(k)$ 
for all but finite $k,$ for all $x,y\in \pn^\dN,$
while $\Et$ can be viewed as a \er\ on 
$\pnqn,$ or on $\dnqn,$ defined as
$x\Et y$ iff $x(k)\Eo y(k)$ for all $k$.

Relations of the form $\rei$ are special case of a wider
family of \er s induced by group actions, see \prf{groa}
below.

The main structure relation between Borel \eqr s is
$\reb,$ Borel reducibility.
Some variations (see \prf{bored}) are involved in special
cases.


\bdf 
\label{typesER}
A Borel \eqr\ $\rE$ on a space $\dX$ is:\vtm 

\noi --
{\it countable\/}, \ 
if every \dde class
$\eke x=\ens{y\in\dX}{x\rE y}\zT x\in\dX,$ 
\index{equivalence relation, ER!countable}%
is countable;\vtm
\imar{Is any ctble $\fs11$ \er\ actually Borel ?}%

\noi --
{\it essentially countable\/}, \ 
if $\rE\reb\rF,$ where $\rF$  
\index{equivalence relation, ER!essentially countable}%
is a countable Borel \er;\vtm

\noi --
{\it finite\/}, \
if every \dde class $\eke x=\ens{y\in\dX}{x\rE y}\zT x\in\dX,$ 
\index{equivalence relation, ER!finite}%
is finite;\vtm

\noi --
{\it hyperfinite\/}, \
if $\rE=\bigcup_n\rE_n$ for an increasing sequence of Borel 
\index{equivalence relation, ER!hyperfinite}%
finite \er s $\rE_n$;\vtm

\noi --
{\it smooth\/}, \
if $\rE\reb\rav\dn$ --- then $\rE$ is obviously Borel;\vtm 
\index{equivalence relation, ER!smooth}%

\noi --
{\it hypersmooth\/}, \
if $\rE=\bigcup_n\rE_n$ for an increasing sequence of  
\index{equivalence relation, ER!hyperfinite}%
smooth \er s $\rE_n$.\qed
\eDf

Countable \eqr s form a widely studied family.
\bit
\item\msur
$\Ey$ is the \dd\reb largest, or
{\it universal countable Borel \er\/}. 
\index{equivalence relation, ER!Einfty@$\Ey$}%
\index{zzEinfty@$\Ey$}%
\eit
See Theorem~\ref{femo} on the existence and exact definition
of $\Ey$.

\vyk{%
The easiest way to define $\Ey$ is to employ
$F_2,$ the free group of two generators, 
which consists of (finite) words composed of the four 
symbols $a,\,b,\,a\obr,\,b\obr,$ including the empty 
word (the neutral element), and the group operation of 
concatenation. 
The {\it shift action\/} of $F_2$ on the compact space
$2^{F_2}$ is defined so that 
$(w\app x)(u)=x(w\obr u)$ for all $w,\,u\in F_2$ and 
$x\in 2^{F_2}.$ 
Now, for $x,\,y\in2^{F_2},$ put $x\Ey y$ iff 
$x=w\app y$ for some $w\in F_2.$
Despite the variety of different ways how a countable 
Borel \er\  can be defined, all of them are $\reb$ 
than $\Ey,$ this is why the latter is called 
{\it universal\/} (countable Borel \er). 
\imar{reference}%
This observation is marked by the framebox 
\framebox{\sf ctble} on the diagram.
}%

The next group includes \eqr s induced by actions of
(the additive groups of) some Banach spaces, in particular 
the following ones well known from textbooks:
\dm
\bay{rclclcl} 
\index{Banach space!lp@$\bel p$}%
\index{zzlp@$\bel p$}%
\index{Banach space!linfty@$\bel\iy$}%
\index{zzlinfty@$\bel\iy$}%
\index{Banach space!c@$\fc$}%
\index{zzc@$\fc$}%
\index{Banach space!c0@$\fco$}%
\index{zzc0@$\fco$}%
\bel p &=& \ens{x\in\rn}{\sum_n|x_n|^p < \iy}\;\;(p\ge1);&& 
\nr x p &=& (\sum_n |x_n|^p)^{\frac1p}\,;\\[0.8\dxii]
\index{norm!lpnorm@\dd{\bel p}norm $\nr\cdot p$}%
\index{zzpx@$\nr\cdot p$}%

{\bel\iy} &=& \ens{x\in\dR^\dN}{\tsup_n|x_n| < \iy};&&
\nr x\iy &=& \tsup_n |x_n|\,;\\[0.8\dxii]
\index{norm!linfnorm@\dd{\bel\iy}norm $\nr\cdot\iy$}%
\index{zzinf@$\nr\cdot\iy$}%

\fc &=& \ens{x\in\dR^\dN}{\tlim_n x_n<\iy\,\text{ exists}};&&
\nor x &=& \tsup_n |x_n|\,;\\[0.8\dxii]

\fco &=& \ens{x\in\dR^\dN}{\tlim_n x_n=0};&&
\nor x &=& \tsup_n |x_n|\,.\\[0.8\dxii]
\vyk{
\fC[a,b] &=& \ens{f\in\dR^{[a,b]}}{f\,\text{ is continuous}};&&
\nor f &=& \tmax_{a\le t\le b} |f(t)|\,;\\[0.8\dxii]

\fC(K) &=& \ens{f\in\dR^K}{f\,\text{ is continuous}};&&
\nor f &=& \tmax_{t\in K} |f(t)|\,;\\[0.8\dxii]

\beL p[a,b] &=&  
\ens{f\in\dR^{[a,b]}}{\int_a^b|f(t)|^p dt<\iy};&&
\nor f &=& (\int_a^b|f(t)|^p dt)^{\frac1p}.
}%
\eay
\dm
Note that ${\bel p}\yd{\fc}\yd{\fco}$ are separable
while ${\bel\iy}$ is non-separable.
The domain of each of the four spaces 
consists of infinite sequences 
$x=\sis{x_n}{n\in\dN}$ of reals, and is a subgroup of 
the group $\rn$ (with the componentwise addition). 
The latter can be naturally equipped with the Polish 
product topology, so that 
${\bel p}\yd{\bel\iy}\yd{\fc}\yd{\fco}$ are Borel subgroups 
of $\rn.$ 
(But not topological subgroups since the distances are different. 
The metric definitions as in ${\bel p}$ or ${\bel\iy}$ 
do not work for $\rn$.)

Each of the four mentioned Banach spaces defines
an orbit equivalence --- 
a Borel \eqr\ on $\rn$ also denoted by, resp.,  
${\bel p}\yd{\bel\iy}\yd{\fc}\yd{\fco}.$  
\index{equivalence relation, ER!lp@$\bel p$}%
\index{zzlp@$\bel p$}%
\index{equivalence relation, ER!linfty@$\bel\iy$}%
\index{zzlinfty@$\bel\iy$}%
\index{equivalence relation, ER!c@$\fc$}%
\index{zzc@$\fc$}%
\index{equivalence relation, ER!c0@$\fco$}%
\index{zzc0@$\fco$}%
For instance, $x\bel p y$ if and only if 
$\sum_k|x_k-y_k|^p<\piy$ (for all $x\yi y\in\rn$).
It is known (see Section~\ref{BAN})
that $\bel1\eqb\Ed$ and
${\bel p}\rebs{\bel q}$ whenever $1\le p<q,$ in 
particular, ${\bel1}\eqb\Ed\rebs{\bel q}$ for any $q>1.$
On the other hand, $\fco\eqb\rzo,$ where
$\rzo$ is the ``density 0'' \eqr:
\bit
\item\msur
$\rzo=\rE_{\Zo},$ thus, for $x,\,y$ in $\pn,$ 
$x\rzo y$ iff 
$\tlim_{n\to\iy}\frac{\#(x\sd y)}n=0.$ 
\eit

Another important \er\ is 
%
\bit
\item\msur
$\rtd,$ often called 
\index{equivalence relation, ER!t2@$\rtd$}%
\index{zzt2@$\rtd$}%
``the equality of countable sets of reals''. 
\eit
%
There is no reasonable way to turn $\pwc\dnn,$ the set 
of all at most countable subsets of $\dnn,$ into a 
Polish space, in order to directly define the equality
of countable sets of reals in terms of $\rav{\cdot}.$ 
However, nonempty members of $\pwc\dnn$ can be
identified with equivalence classes in $\nnn/\rtd,$ where 
$g\rtd h$ iff $\ran g=\ran h:$ for $g,\,h\in\nnn.$ 
(See below in Section~\ref{groas} on \eqr s
$\rT_\al$ for all $\al<\omi$.)

In addition to the families of \eqr s introduced by
Definition~\ref{typesER}, some more complicated families
will be considered below, including \er s induced by
Polish group actions, turbulent \er s,
\er s classifiable by countable structures, pinned \er s,
and some more.

\punk{Borel reducibility of basic \eqr s}
\label{brb}

The diagram
on page~\pageref{p-p}
begins, at the low end, 
with cardinals $1\le n\in\dN,\;\alo,\;\cont,$ 
which denote the \er s 
of equality on resp.\ finite, countable, uncountable 
Polish spaces.
As all uncountable Polish spaces are Borel isomorphic, 
the \eqr s $\rav\dX,$ $\dX$ a Polish space, are 
characterized, modulo $\reb,$
or even modulo Borel isomorphism between the domains,
by the cardinality of the domain,
which can be any finite $1\le n<\om,$ or 
$\alo,$ or $\cont=2^\alo.$


\begin{figure}[h]
\vspace*{1mm}

\setlength{\unitlength}{1mm} 
\begin{picture}(130,110)(0,0) 

\label{p-p}

\put(80,0){\tob}
\put(82,-1){$1$}
\put(80,0){\line(0,1){8}}
\put(80,8){\tob}
\put(82,7){$2=\rav{\ans{1,2}}$}
\put(80,9) {\toq}
\put(80,10){\toq}
\put(80,11){\toq}
\put(80,12){\toq}
\put(80,13){\toq}
\put(80,14){\toq}
\put(80,15){\toq}
\put(80,16){\tob}
\put(82,15){$n=\rav{\ans{1,2,...,n}}$}
\put(52,15){($1\le n<\alo$)}
\put(80,17){\toq}
\put(80,18){\toq}
\put(80,19){\toq}
\put(80,20){\toq}
\put(80,21){\toq}
\put(80,22){\toq}
\put(80,23){\toq}
\put(80,24){\tob}
\put(82,23){$\alo=\rav\dN$}
\put(80,24){\line(0,1){8}}
\put(80,32){\tob}
\put(82,31){$\cont=\rav{\dn}$}
\put(80,32){\line(0,1){8}}
\put(80,40){\tob}
\put(82,38){$\Eo$}

\put(80,40){\line(-3,1){30}}
\put(50,50){\line(-2,1){40}}
\put(10,70){\tob}
\put(10,70){\line(3,4){25}}
\put(4.5,69){$\Ei$}

\put(80,40){\line(-3,2){20.9}}
\put(55,56.7){\line(-3,2){20.3}}
\put(35,70){\tob}


\put(37,69){$\Ed\eqb{\bel1}$}
\put(55,55){\framebox{\sf ?}}

\put(80,40){\line(1,2){6.9}}
\put(95,70){\line(-1,-2){5.7}}
\put(95,70){\line(-1,1){19.5}}
\put(95,70){\tob}
\put(82.5,55){\framebox{\sf ctble}}
\put(96,69){$\Ey$}

\put(80,40){\line(3,2){45}}
\put(125,70){\tob}
\put(123,66){$\Et$}

\put(125,70){\line(0,1){14}}
\put(125,100){\line(0,-1){11}}
\put(125,100){\tob}
\put(117,102){$\rzo\eqb{\fco}$}
\put(118,86){\framebox{\sf \co eqs}}

\put(95,70){\line(0,1){30}}
\put(93,102){$\rtd$}
\put(95,100){\tob}
\put(95,100){\line(1,-1){30}}

\newlength{\borlen}
\settowidth{\borlen}{\small the {\sf non-P}\ }

\put(5,92){{\bmp{\borlen}{\small border of\\ the {\sf non-P}\\ domain}\emp}}

\put(35,70){\line(0,1){10}}
\put(35,85){\line(0,1){18}}
\put(35,103){\line(3,-1){40.5}}
\put(35,103){\tob}
\put(32.0,104){$\bel\iy$}
\put(31.5,81.2){\framebox{$\bel p$}}

\vyk{
\put(95,100){\mtir}
\put(91,100){\mtir}
\put(87,100){\mtir}
\put(80.5,99){$\ni$}
\put(79,100){\mtir}
\put(75,100){\mtir}
}


\linethickness{0.2mm}

\qbezier(9,68)(28,75)(39,109)

\put(17,90){\vector(1,-1){8.8}}

\end{picture}

\caption{Reducibility between some basic \er s}
%
\noindent
Connecting lines here indicate
Borel reducibility of lower ERs to upper ones.
\end{figure}
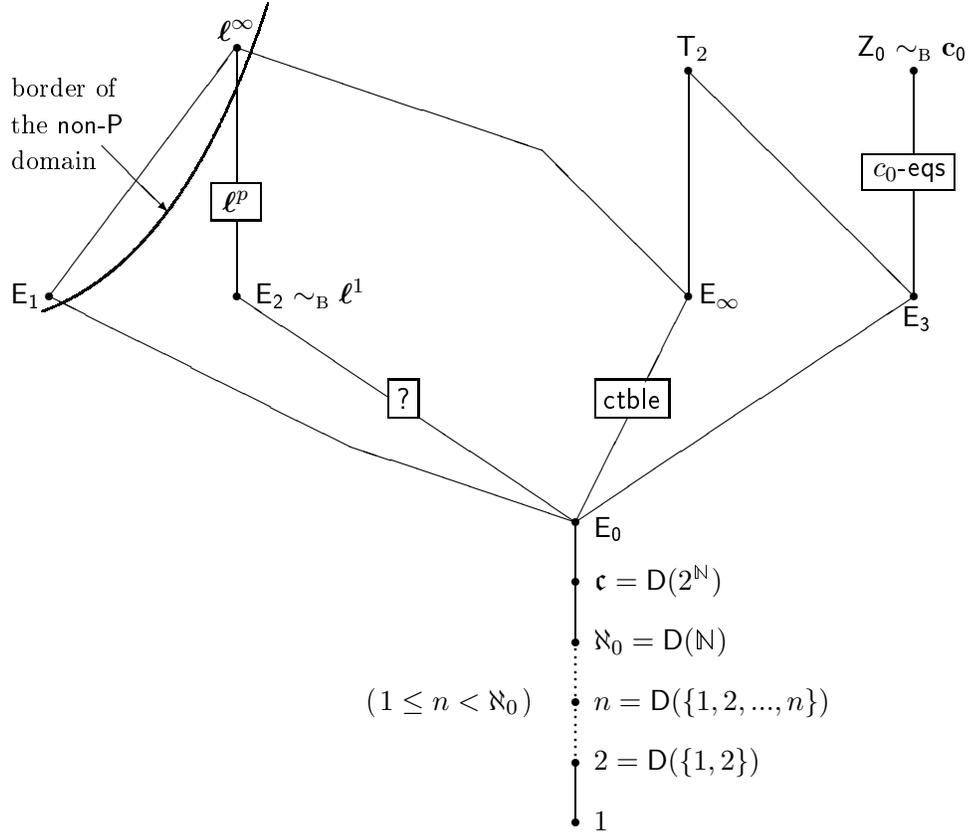
%


%
 
The $\Eo$ splitting is the key element of the diagram on
page~\pageref{p-p}.
That $\rav\dn\reb\Eo$ can be proved by a 
rather simple embedding while the strictness can be derived 
from an old result of Sierpi\'nski~\cite{sie}: 
any linear ordering of 
all \dd\Eo classes yields a Lebesgue non-measurable set of 
the same descriptive complexity as the ordering.
That every \er\ $\reb\Eo$ is $\eqb$ to some 
$n\ge 1\yt{\rav\dN}\yt{\rav\dn},$ or $\Eo$ itself, is 
witnessed by the following two classical results:

\bde
\itsep
\item[{\ubf 1st dichotomy\/}] 
(
Thm~\ref{1dih} below).  
{\it 
Any Borel, even any\/ $\fp11$ \er\/ $\rE$ \poq{either} 
has at most countably many equivalence classes, formally, 
$\rE\reb\alo=\rav\dN,$ \poq{or} satisfies\/
$\cont=\rav\dn\reb\rE$.\/} 

\item[{\ubf 2nd dichotomy\/}] 
(
Thm~\ref{2dih}).  
{\it 
Any Borel \er\/ $\rE$ satisfies either  
$\rE\reb\cont$ or\/ $\Eo\reb\rE.$\/} 
\ede
  
The linearity breaks above $\Eo:$ 
each one of the four \eqr s $\Ei\zd\Ed\zd\Et\zd\Ey$ of the 
next level is strictly \dd{\rebs}bigger than $\Eo,$ 
and they are pairwise \dd\reb incomparable 
with each other, see \prf{finr}.

\vyk{
\bpro
\label{eo1}
The \eqr s\/ ${\Ei}\zd{\Ed}\zd{\Et}\zd{\Ey}\zd$ are
pairwise\/ \dd\reb incomparable and satisfy\/
$\Eo\rebs\rE_i\yt i=1,2,3,\iy.$  
\epro
\bpf
Simple reductions witness $\Eo\rebs\rE_i$:
for instance, the map $\vt(x)=\ans0\ti x:\pn\to\pnn$
reduces $\Eo$ to $\Et.$
The incomparability claims will be proved below,
partially in this Section, partially as corollaries of
deep theorems in next sections. 
\epf
}

One naturally asks what is going on in the intervals
between $\Eo$ and these four \eqr s.
The following results provide some answers.  

\bde
\itsep
\item[{\ubf 3rd dichotomy\/}] 
(%
Thm~\ref{kelu}). 
{\it Any \er\/ $\rE\reb\Ei$ satisfies\/  
$\rE\reb\Eo$ or\/ $\rE\eqb\Ei.$} 

\item[{\ubf 4th dichotomy\/}] 
(
Thm~\ref{h}).  
{\it Any \er\/ $\rE\reb\Ed$ either is 
essentially countable or satisfies \/ $\rE\eqb\Ed$.}
\ede

See Definition~\ref{typesER} regarding  
essentially countable \er s in the 4th dichotomy. 
The ``either'' case there remains mysterious:
any countable Borel \er s $\rE\reb\Ed$ 
known so far are $\reb\Eo.$ 
It is a problem whether the ``either'' case can be
improved to $\reb\Eo.$  
This is marked by the framebox 
\framebox{\sf ?} on the diagram.

The fifth dichotomy theorem is a bit more special, it will
be addressed below.

\bde
\itsep
\item[{\ubf 6th dichotomy\/}] 
(
Thm~\ref{6dih}).  
{\it Any \er\/ $\rE\reb\Et$ satisfies\/ 
$\rE\reb\Eo$ or\/ $\rE\eqb\Et$.}

\item[{\ubf Adams--Kechris theorem\/}]
%
(not to be proved here).
{\it There is continuum many 
pairwise \dd\reb incomparable countable Borel \er s\/}.
\ede

The framebox \fbox{\sf\co eqs} denotes 
\co{\it equalities\/}, a family of Borel \er s introduced 
by Farah~\cite{f-co}, all of them are \dd\reb between 
$\Et$ and ${\fco}\eqb\rzo,$ and there is 
continuum-many \dd\reb incomparable among them.

The {\sf non-P} domain denotes the family of all 
\er s $\rei,$ where $\cI$ is a Borel ideal which is 
not a P-ideal. 
By Solecki~\cite{sol,sol'}, for a Borel ideal $\cI$ to be 
\imar{$\Ei$ and Polish grps action problem}%
not a P-ideal it is necessary and sufficient that 
$\Ii\reb\cI,$ or, equivalently,  
$\Ei\reb\rei.$%

\bqe
\label{min}
It there any reasonable ``basis'' of Borel \er s
above $\Eo$?
\eqe

It was once considered \cite{hk:nd} as a plausible
hypothesis that any Borel \er\ which is not $\reb\Ey,$ \ie, 
not an essentially countable \er, satisfies
$\rE_i\reb\rE$ for at least one $i=1,2,3.$ 
This turns out to be not the case: 
Farah~\cite{f-bas,f-tsir} and Velickovic~\cite{vel:tsir}
found an independent family of uncountable Borel \er s, 
based on {\it Tsirelson ideals\/}, \dd\reb incomparable 
with $\Ei\zi\Ed\zi\Et,$ see below.



It is the most interesting question whether 
the diagram on page~\pageref{p-p} is complete in the
sense that there is no \dd\reb connections betwen the
\eqr s mentioned in the diagram except for those
explicitly indicated by lines.
Basically, one may want to prove the following
non-reducibility claims:\vtm

$
\bay[t]{clcllll}
(1) & \Ei &\not\reb\;:& \Ed, & \rtd, & \fco;\\[0.7\dxii]

(2) & \bel\iy &\not\reb\;:& \Ei, & \Ed, & \rtd, & \fco;\\[0.7\dxii]

(3) & \Ed &\not\reb\;:& \Ei, & \rtd, & \fco;\\[0.7\dxii]

(4) & \Ey &\not\reb\;:& \Ei, & \Ed, & \fco;\\[0.7\dxii]

(5) & \Et &\not\reb\;:& \bel\iy;\\[0.7\dxii]

(6) & \rtd &\not\reb\;:& \bel\iy, & \fco;\\[0.7\dxii]

(7) & \fco &\not\reb\;:& \bel\iy, & \rtd.\\[0.7\dxii]
\eay
$\vtm

Beginning with (1), we note that $\Ei$ is not Borel
reducible to any \eqr\ induced by a Polish action
(of a Polish group) by Theorem~\ref{e1pga} below.
On the other hand, $\Ed\zd\rtd\zd {\fco}$ obviously
belong to this category of \er s.

(2) follows from (1) and (3) and can be omitted.

In (3), $\Ed\not\reb\Ei$ can be proved by an argument
rather similar to the proof of Theorem~\ref{t:>s}.
Alternatively, it will follow from Theorem~\ref{rig1}
that any Borel ideal $\cI$ with $\rE_\cI\reb\Ei$
is isomorphic, via a bijection between the underlying
sets, to $\Ii$ or to a trivial variation of $\ifi,$ 
but $\Id$ does not belong to this category.
The result $\Ed\not\reb{\fco}$ in (3) is
Theorem~\ref{t:>s}\ref{t>s2}.

The results $\Ed\not\reb\rtd$ and ${\fco}\not\reb\rtd$ in
(3) and (7) are proved below in Section~\ref{groat}
(Corollary~\ref{-t2}); this will involve the turbulence theory.

The result of (5) is Lemma~\ref{e3}.
It implies ${\fco}\not\reb{\bel\iy}$ in (7).

(6) will be established in Section~\ref{rtd}.

This leaves us with (4).
We don't know how to prove $\Ey\not\reb\Ei$
easily and directly.
The indirect way is to use Theorem~\ref{kelu} below,
according to which $\Ey\reb\Ei$ would imply 
either $\Ey\eqb\Ei$ --- impossible, see above, 
or $\Ey\le\Eo.$
The latter conclusion is also a contradiction since
$\Eo\rebs \Ey$ is known in the theory of countable Borel
\eqr s (see \cite[p.~210]{djk}).

\bqe 
\label{liy2co}
Is $\Ey$ Borel reducible to $\fco$? to $\bel1$ or any
other $\bel p$?
\eqe

\punk{Operations on \eqr s}
\label{opeer}

The following operations over \er s are in 
part parallel to the operations on ideals in \prf{opeid}.

\ben
\tenu{(o\arabic{enumi})} 
\itla{oe:beg}
\label{cun}
countable union 
(if it results in a \er) and countable intersection of 
\er s on one and the same space; 

\itla{cdun}
countable disjoint union $\rE=\bigvee_k\rE_k$ of 
\er s $\rE_k$ on Polish spaces $\spa_k,$  
that is, a \er\ on 
$\spa=\bigcup_k(\ans k\ti\spa_k)$ 
(with the topology generated by sets of the form 
$\ans k\ti U,$ where $U\sq \spa_k$ is open) 
defined as follows: $\ang{k,x}\rE \ang{l,y}$ iff 
$k=l$ and $x\rE_k y.$ 
(If $\spa_k$ are pairwise disjoint and open in 
$\spa'=\bigcup_k\spa_k$ then we can equivalently 
define $\rE=\bigvee_k\rE_k$ on $\spa'$ so that 
$x\rE y$ iff $x,\,y$ belong to the same $\spa_k$ 
and $x\rE_k y$.); 

\itla{prod} 
product $\rE=\prod_k\rE_k$ 
\index{product of ers@product of \er s}%
of \er s $\rE_k$ on spaces $\spa_k,$  
that is, the 
\er\ on the product space 
$\prod_k\spa_k$ defined by: 
$x\rE y$ iff $x_k\rE_k y_k$ for all $k$. 

\itla{fp}
the {\it Fubini product\/} (ultraproduct) 
\index{Fubini product!of ers@of \er s}%
\index{ultraproduct of ers@ultraproduct of \er s}%
$\fps{\cI}{\rE_{k}}{k\in\dN}$ 
of \er s $\rE_k$ on spaces $\spa_k,$ 
modulo an ideal $\cI$ on $\dN,$ 
that is, the \er\ on the product space 
$\prod_{k\in\dN}\spa_k$ defined as follows: 
$x\rE y$ iff $\ens{k}{x_k\nE_k y_k}\in\cI$; 

\itla{oe:end}
\label{cp}
{\it countable power\/} \er\ 
$\ei$ of a \er\ $\rE$ on a space $\spa$ is  
a \er\ on $\spn$ defined as follows: 
$x\ei y$ iff 
$\ens{[x_k]_{\rE}}{k\in\dN}=\ens{[y_k]_{\rE}}{k\in\dN},$ 
so that for any $k$ there 
is $l$ with $x_k\rE y_l$ and for any $l$ there 
is $k$ with $x_k\rE y_l$.
\een

These operations allow us to obtain a lot of interesting 
\er s starting just with very primitive ones. 
For instance, we can define the sequence of 
\er s $\rT_\xi,\msur$ $\xi<\omi,$ of 
H.~Friedman~\cite{frid} as follows~\footnote
{\ Hjorth~\cite{h} uses $\rF_\xi$ instead of $\rT_\xi$.}. 
Let $\rT_0=\rav\dN,$ the equality relation on $\dN.$ 
We put $\rT_{\xi+1}={\rT_\xi}^\iy.$ 
\index{zzTxi@$\rT_\xi$}%
If $\la<\omi$ is a limit ordinal, then put 
$\rT_\la=\bigvee_{\xi<\la}\rT_\xi$. 

In particular $\dom\rT_1=\dnn$ and 
$x\rT_1 y$ iff $\ran x=\ran y,$ for $x,\,y\in\dnn.$
Thus the map $\vt(x)=\ran x$ witnesses that 
$\rT_1\reb{\rav\pn}.$ 
To show the converse, define, 
for any infinite $u\sq\dN,$ $\ba(u)$ be the 
increasing bijection $\dN\onto u,$ while if 
$u=\ans{k_0,...,k_n}$ is finite, put 
$\ba(u)(i)=k_i$ for $i<n$ and $\ba(u)(i)=k_n$ for 
$i\ge n.$ 
Then $\ba$ witnesses ${\rav\pn}\reb{\rT_1},$ thus, 
${\rT_1}\eqb{\rav\pn}$. 
%
It easily follows that $\rT_2\eqb {{\rav\pn}}^\iy,$ 
in fact, $\rT_2\eqb {\rav\dX}^\iy$ for any uncountable 
Polish space $\dX$ as any such $\dX$ is Borel 
isomorphic to $\pn$ 
(or to $\dn,$ which is essentially the same). 
With $\dX=\dnn$ we obtain the definition of $\rT_2$ in 
\prf{someq}.

\punk{Orbit \eqr s of group actions}
\label{groa}

An {\it action\/} of a group $\dG$ on a space $\dX$ is 
\index{action}%
any map $a:{\dG\ti\dX}\to\dX,$ usually written as 
$a(g,x)=g\app x,$ such that 
\index{action!g.x@$g\app x$}%
\index{zzg.x@$g\app x$}%
$1)\msur$ $e\app x=x,$ and 
%
$2)\msur$ $g\app(h\app x)=(gh)\app x,$ --- 
then, for any $g\in\dG,$ the map $x\mapsto g\app x$ is 
a bijection $\dX$ onto $\dX$ with $x\mapsto g\obr\app x$ as the 
inverse map.
A {\it\dd\dG space\/} is a pair $\stk\dX a,$ where $a$ is an 
\index{Gspace@{}\dd\dG space}%
action of $\dG$ on $\dX\,;$ 
in this case $\dX$ itself is also called a \dd\dG space, 
and the {\it orbit \er\/},  
\index{equivalence relation, ER!orbit er@orbit \er\ $\ergx$}%
\index{equivalence relation, ER!induced@induced by an action}%
or {\it\er\ induced by the action\/}, 
\index{zzEGX@$\ergx$}%
$\aer a\dX=\ergx$ is defined on $\dX$ so that $x\ergx y$ iff 
there is $a\in\dG$ with $y=a\app x.$
{}\dd\ergx classes are the same as 
\dd\dG orbits, \ie,
\dm
[x]_\dG=[x]_{\ergx}=\ens{y}{\sus g\in\dG\;(g\app x=y)}\,.
\dm

A {\it homomorphism\/} (or \dd\dG homomorphism) 
of a \dd\dG space $\dX$ into 
a \dd\dG space $\dY$ is any map $F:\dX\to\dY$ compatible 
with the actions in the sense that $F(g\app x)=g\app F(x)$ 
for any $x\in\dX$ and $g\in\dG.$ 
A $1-1$ homomorphism is an {\it embedding\/}. 
\index{homomorphism!of actions}%
\index{embedding!of actions}%
\index{isomorphism!of actions}%
An embedding $\onto$ is an {\it isomorphism\/}. 
Note that a homomorphism $\stk\dX a\to\stk\dY b$ is a 
reduction of $\aer a\dX$ to $\aer b\dY,$ but not conversely.

A {\it Polish group\/} is a group whose underlying set is a 
Polish space and the operations are continuous; 
\index{group!Polish}%
\index{group!Borel}%
a {\it Borel group\/} is a group whose underlying set is a 
Borel set (in a Polish space) and the operations are Borel maps. 
%
\index{group!Polishable}%
A Borel group is {\it Polishable\/} if there is a Polish 
topology on the underlying set which 1) produces the same 
Borel sets as the 
original topology and 2) makes the group Polish.

\bit
\itsep
\item
If both $\dX$ and $\dG$ are Polish and the action continuous, 
then $\stk\dX a$ (and also $\dX$) is called a {\it Polish\/} 
\dd\dG space. 
\index{Gspace@{}\dd\dG space!Polish}%
\index{Gspace@{}\dd\dG space!Borel}%
If both $\dX$ and $\dG$ are Borel and the action is a 
Borel map, then $\stk\dX a$ (and also $\dX$) is called a 
{\it Borel\/} \dd\dG space.
\eit

\bex
\label{ex:ac}
(i)
Any ideal $\cI\sq\pn$ is a group with $\sd$ as the operation.
We cannot expect this group to be Polish in the product
topology inherited from $\pn$
(indeed, $\cI$ would have to be $\Gd$).
However if $\cI$ is a P-ideal then it is Polishable   
(see \prf{punk:poli}),
in other words, $\stk\cI\sd$
is a Polish group in an appropriate Polish topology
compatible with the Borel structure of $\cI.$
Given such a topology, the \dd\sd action of
\imar{correct?}%
(a P-ideal) $\cI$ on $\pn$ is Polish, too.

(ii)
Consider $\dG=\pwf\dN$ a countable subgroup of
$\stk\pn\sd.$
Define an action of $\dG$ on $\dn$ as follows:
$(w\app x)(n)=x(n)$ whenever $n\nin w$ and
$(w\app x)(n)=1-x(n)$ otherwise.
The orbit \eqr\ $\ergx$ of this action is obviously $\Eo$.
Note that this action is {\it free\/}:
$x=w\app x$ implies $w=\pu$ (the neutral element of $\dG$)
for any $x\in\dn.$

Now consider any Borel pairwise \dd\Eo inequivalent set
$T\sq\dn.$ 
Then $w\app T\cap T=\pu$ for any $w\ne\pu$ by the above. 
It easily follows that $T$ is meager in $\dn.$
(Otherwise $T$ is co-meager on a basic clopen set
$\bon s\dn=\ens{x\in\dn}{s\su x},$ where
$s\in\dln.$
Put $w=\ans{n},$ where $n=\lh s.$
Then $w\in\dG$ maps
$T\cap\bon{s\we0}\dn$ onto $T\cap\bon{s\we1}\dn.$
Thus $w\app T\cap T\ne\pu$ -- contradiction.)
We conclude that $\dG\app T=\bigcup_{w\in\dG}w\app T$
is still a meager subset of $\dn$ in this case, and hence
$T$ cannot be a full (Borel) transversal for $\Eo.$

(iii)
The {\it canonical\/ {\rm(or {\it shift\/})} action\/}
\index{action!canonicalofGonXG@canonical of $\dG$ on $X^\dG$}%
\index{action!shift}%
of a group $\dG$ on a set of the form $X^\dG$ ($X$ any set)
is defined as follows:
$g\app\sis{x_f}{f\in\dG}=\sis{x_{g\obr f}}{f\in\dG}$ for any
element $\sis{x_f}{f\in\dG}\in X^\dG$ and any $g\in\dG.$
This is easily a Polish action provided $\dG$ is countable,
$X$ a Polish space, and $X^\dG$ given the product topology.
The \eqr\ on $X^\dG$ induced by this action is denoted by
$\rE(\dG,X)$.
\index{equivalence relation, ER!EGX@$\rE(\dG,X)$}%
\index{zzEGX@$\rE(\dG,X)$}%
\index{equivalence relation, ER!induced by the shift action}%
\eex

The next theorem (rather difficult to be proved here) 
shows that the type of the group is the essential 
component in the difference between Polish and Borel actions: 
roughly, any Borel action of a Polish group $\dG$ is a Polish 
action of $\dG$. 

\bte[{{\rm\cite[5.2.1]{beke}}}]
\label{b2c}
Suppose that\/ $\dG$ is a Polish group and\/ $\stk\dX a$ is 
a Borel\/ \dd\dG space. 
Then\/ $\dX$ admits a Polish topology which 1) produces the 
same Borel sets as the original topology, and 2) makes 
the action to be Polish.\qeD
\ete

If $\stk\dX a$ is a Borel \dd\dG space 
(and $\dG$ is a Borel group) 
then $\ergx$ is easily a $\fs11$ \er\ on $\dX.$  
Sometimes $\ergx$ is even Borel: for instance, 
when $\dG$ is a countable group and the action is Borel, 
or if $\dG=\cI\sq\pn$ is a Borel ideal, considered as 
a group with $\sd$ as the operation, which acts on 
$\dX=\pn$ by $\sd,$ so that ${\aer\dG\pn}=\rE_{\cI}$ is 
Borel because $x \aer\dG\pn y$ iff $x\sd y\in\cI.$ 
Several much less trivial cases when $\ergx$ is Borel are 
described in \cite[Chapter 7]{beke}, 
for instance, if all \dd{\ergx}classes are Borel sets of 
bounded rank then $\ergx$ is Borel \cite[7.1.1]{beke}.
Yet rather surprisingly equivalence classes generated by
Borel actions are always Borel.

\bte
[{{\rm see \cite[15.14]{dst}}}]
\label{borcl} 
If\/ $\dG$ is a Polish group and\/ $\stk\dX a$ is 
a Borel\/ \dd\dG space  
then every equivalence class of\/ $\ergx$ is Borel.
\ete
\bpf
It can be assumed, by Theorem~\ref{b2c}, that the action 
is continuous.
Then for any $x\in\dX$ the {\it stabilizer\/} 
$\dG_x=\ens{g}{g\app x=x}$
is a closed subgroup of $\dG.$~\footnote 
{\ Kechris~\cite[9.17]{dst} gives an independent proof. 
Both $\dG_x$ and its topological closure, say, $G'$ are 
subgroups, moreover, $G'$ is a closed subgroup, hence, we can 
assume that $G'=\dG,$ in other words, that $\dG_x$ is dense 
in $\dG,$ and the aim is to prove that $\dG_x=\dG.$ 
By a simple argument, $\dG_x$ is 
either comeager or meager in $\dG.$ 
But a comeager subgroup easily coincides with the whole group, 
hence, assume that $\dG_x$ is meager (and dense) in 
$\dG$ and draw a contradiction. 

Let $\sis{V_n}{n\in\dN}$ be a basis of the topology of $\dX,$ 
and $A_n=\ens{g\in\dG}{g\app x\in V_n}.$ 
Easily $A_nh=A_n$ for any $h\in\dG_x.$ 
It follows, because $\dG_x$ is dense, that every $A_n$ is 
either meager or comeager. 
Now, if $g\in\dG$ then $\ans g=\bigcap_{n\in N(g)}A_n,$ 
where $N(g)=\ens{n}{g\app x\in V_n},$ thus, at least one of 
sets $A_n$ containing $g$ is meager. 
It follows that $\dG$ is meager, contradiction.} 
\imar{Quotient spaces~?}%
We can consider $\dG_x$ as continuously acting on $\dG$
by $g\app h=gh$ for all $g,h\in\dG.$
Let $\rF$ denote the associated orbit \er. 
Then every \ddf class $\ek g\rF=g\,\dG_x$
is a shift of $\dG_x,$ 
hence, $\ek g\rF$ is closed.
On the other hand, the saturation $\ek \cO\rF$
of any open set $\cO\sq\dG$ is obviously open. 
Therefore, by Lemma~\ref{transv}\ref{transv4} below,
$\rF$ admits a Borel transversal $S\sq\dG.$  
Yet $g\longmapsto g\app x$ is a Borel $1-1$ map 
of a Borel set $S$ onto $\ek x\rE,$ hence,
$\ek x\rE$ is Borel by \Cpro.
\epf

It follows that \poq{not} all $\fs11$ \er s are orbit \er s 
of Borel actions of Polish groups: indeed, take a non-Borel 
$\fs11$ set $X\sq\dnn,$ define $x\rE y$ if either $x=y$ or 
$x,\,y\in X,$ this is a $\fs11$ \er\ with a non-Borel class 
$X$.

\imar{$\fs11$ or Borel \er\ not induced by Borel grp~?}%

\imar{Borel \er\ not induced by Polish grp~?}%

\punk{Forcings associated with pairs of  \eqr s}
\label{f2}

The range of applications of this comparably new topic
is not yet clear, but at least it offers interesting
technicalities.

\bdf[{{\rm Zapletal~\cite{zap}}}]
\label{df:f2}
Suppose that $\rE$ is a Borel \eqr\ on a Polish space
$\dX,$ and $\rF\rebs\rE$ is another Borel \eqr.

$\zidef$ is the collection of all Borel
\index{ideal!ief@$\zidef$}%
\index{zzief@$\zidef$}%
sets $X\sq\dX$ such that ${\rE\res X}\reb\rF.$
Clearly $\zidef$ is an ideal in the algebra of all
Borel subsets of $\dX.$
The associated forcing $\zfoef$  
\index{forcing!pef@$\zfoef$}%
\index{zzpef@$\zfoef$}%
consists of all Borel sets $X\sq\dX\zt X\nin\zidef$.
\edf

For instance, the ideal $\zid{\rav\dn}{\rav\dN}$  consists
of all countable Borel sets $X\sq\dn,$ therefore 
$\zfo{\rav\dn}{\rav\dN}$ contains all
\poq{un}countable Borel sets $X\sq\dn$ and is equal to the
Sacks forcing.
The ideal $\zid\Eo{\rav\dn}$  consists
of all Borel sets $X\sq\dn$
such that $\Eo\res X$ is non-smooth
(since smoothness is equivalent to being $\reb\rav\dn$).
See \prf{zeo} on the associated forcing $\zfo\Eo{\rav\dn}$.

\vyk{
 
$C^n[a,b]$ funk
iz [a,b] v R, nepreryvn proizvodn do poryadka n vkl
s normoj  \sum_{k=0}^n \max_x |f^n(x)|,  a\le x\le b

$C^n(I_m)$ funk
iz m-mernogo kuba v R, nepreryvn proizvodn do poryadka n vkl
s ravnomernoj normoj po vsem proizvodn

$L_p[a,b]$ funk
izmerimy na [a,b], sunmmiruemy so stepen'yu p
norma kak l_p

This group includes, for any $p\ge 1,$ 
the \er\ $\bel p$ on $\rtn$ defined 
\imar{Are tsirelsons Banach ?\\[1ex]
Is any ctble Borel $\rE\reb{\fco}$ necessarily 
$\reb\Eo$?\\[1ex] 
Are $\bel p$ and $\fco$ comparable with tsirelsons~?}%
as follows: $x\bel p y$ iff $\sum_k|x_k-y_k|^p<\piy$ 
for all sequences 
$x=\sis{x_k}{k\in\dN}$ and $y=\sis{y_k}{k}$ in $\rtn.$ 
Clearly $\bel p$ is the orbit relation of the 
(additive group of the) Banach 
space $\ell^p=\ens{x\in\rtn}{\sum_k|x_k|^p<\piy}$ 
acting on $\rtn$ by componentwise addition. 
}

\parf{``Elementary'' stuff}
\label{BAN}

This Section gathers proofs of some
reducibility/irreducibility results related to the diagram
on page~\pageref{p-p}, elementary in the sense that they do not
involve any special concepts.
Some of them are really simple, some other quite tricky.

\punk{$\Et$ and $\rtd$: outcasts}
\label{e3t2}

These \eqr s, together with ${\fco}\eqb\rzo,$ are the only
non-$\fs02$ equivalences explicitly mentioned on the
diagram.

\ble
\label{e3}
$\Et$ is Borel irreducible to\/ ${\bel\iy}.$ 
\ele
\bpf
Suppose towards the contrary that $\vt:\dntn\to\rtn$
is a Borel reduction of $\Et$ to\/ ${\bel\iy}.$~\snos
{Recall that, for $x,y\in\dntn,$ $x\Et y$ means  
$\seq xk\Eo\seq yk\zd\kaz k,$ where
$\seq xk\in\dn$ is defined by $\seq xk(n)=x(k,n)$
for all $n$ while $a\Eo b$ means that
$a\sd b=\ens m{a(m)\ne b(m)}$ is finite.}
Since obviously ${\bel\iy}\eqb{\bel\iy}\ti{\bel\iy},$ 
Lemma~\ref{l:bc}
reduces the general case to the case of continuous $\vt.$
Define $\fo,\fr\in\dn$ by 
$\fo(n)=0\zd\fr(n)=1\zd\kaz n.$ 
Define $\bbo\in\dntn$  by
$\fo(k,n)=0$ for all $k,n,$ thus $\seq\bbo k=\fo\zd\kaz k.$
Finally, for any $k$ define $\fz_k\in\dn$ by
$\fz_k(n)=1$ for $n<k$ and $\fz_k(n)=0$ for $n\ge k$.

We claim that there are increasing sequences
of natural numbers $\sis{k_m}{}$ and $\sis{j_m}{}$ 
such that $|\vt(x)(j_m)-\vt(\bbo)(j_m)|>m$
for any $m$ and any $x\in\dntn$ satisfying
\dm
\seq x{k}=
\left\{
\bay{rl}
\fz_{k_i} & \text{whenever }\,i<m\,\text{ and }\,
k=k_i\\[0.7\dxii]

\fo & \text{for all }\,k<k_m\,
\text{not of the form }\,k_i.
\eay
\right.
\dm
To see that this implies contradiction
define $x\in\dntn$ so that
$\seq x{k_i}=\fz_{k_i}\zd\kaz i$ and
$\seq x{k}=\fo$ whenever $k$ does
not have the form $k_i.$
Then obviously $x\Et\bbo,$ but
$|\vt(x)(j_m)-\vt(\bbo)(j_m)|>m$ for all $m,$ hence
$\vt(x)\bel\iy\vt(\bbo)$ fails, as required.

We put $k_0=0.$ 
To define $j_0$ and $k_1,$ consider $x_0\in\dntn$ defined
by $\seq{x_0}0=\fr$ but $\seq{x_0}k=\fo$ for all $k\ge1.$
Then $x_0\Et \bbo$ fails, and hence
$\vt(x_0)\bel\iy\vt(\bbo)$ fails either.
Take any $j_0$ with $|\vt(x_0)(j_0)-\vt(\bbo)(j_0)|>0.$
As $\vt$ is continuous, there is a number $k_1>0$ such that
$|\vt(x)(j_0)-\vt(\bbo)(j_0)|>0$
holds for any $x\in\dntn$ with $\seq x0=\fz_{k_1}$
and $\seq xk=\fo$ for all $0<k<k_1$. 

To define $j_1$ and $k_2,$ consider $x_1\in\dntn$
defined so that $\seq{x_1}0=\fz_{k_1},$   
$\seq{x_1}k=\fo$ whenever $0<k<k_1,$ and
$\seq{x_1}{k_1}=\fr.$
Once again there is a number $j_1$ with
$|\vt(x_1)(j_1)-\vt(\bbo)(j_1)|>1,$
and a number $k_2>k_1$ such that
$|\vt(x)(j_1)-\vt(\bbo)(j_1)|>1$
for any $x\in\dntn$ with $\seq x0=\fz_{k_1},$
$\seq x{k_1}=\fz_{k_1},$
and $\seq xk=\fo$ for all $0<k<k_1$ and $k_1<k<k_2$.

Et cetera.
\epf

\ble
\label{e3'}
$\Et$ is Borel reducible to both\/ $\rtd$ and\/ $\fco$.
\ele
\bpf
(1)
If $a\in\dn$ and $s\in\dln$ then define $sx\in\dn$ by
$(sx)(k)=x(k)+_2s(k)$ for $k<\lh s$ and
$(sx)(k)=x(k)$ for $k\ge\lh s.$
If $m\in\dN$ then $m\we x\in\dn$ denotes the concatenation.
In these terms, if $x,y\in\dntn$ then obviously
\dm
x\Et y\leqv
\ens{m\we{(s\seq xm)}}{s\in\dln\yt m\in\dN}=
\ens{m\we{(s\seq ym)}}{s\in\dln\yt m\in\dN}.
\dm
Now any bijection $\dln\ti\dN\onto\dN$ yields a Borel
reduction of $\Et$ to $\rtd$.

(2)
To reduce $\Et$ to $\fco$ consider a Borel map
$\vt:\dntn\to\rtn$ such that
$\vt(x)(2^n(2k+1)-1)=n\obr \seq xn(k)$.
\epf

\ble
\label{e3"}
Any countable Borel ER is Borel reducible to\/ $\rtd$.
\ele
\bpf
Let $\rE$ be a countable Borel ER on $\dn.$
It follows from \Cenu\ that there is a Borel map
$f:\dn\ti\dN\to\dn$ such that 
$[a]_{\rE}=\ens{f(a,n)}{n\in\dN}$ for all $a\in\dn.$ 
The map $\vt$ sending any $a\in\dn$ to $x=\vt(a)\in\dntn$
such that $\seq xn=f(a,n)\zd\kaz n,$ is a reduction required.
\epf

See further study on $\rtd$ in Section~\ref{rtd}, where
it will be shown that $\rtd$ is not Borel reducible to
a big family of \eqr s that includes
${\fco}\zi{\bel p}\zi{\bel\iy}\zi{\Ei}\zi{\Ed}\zi{\Et}\zi{\Ey}.$
On the other hand, the \eqr s in this list,
with the exception of ${\Et}\zi{\Ey},$ 
are not Borel reducible to $\rtd$ --- this follows from the
turbulence theory presented in
Section~\ref{groat}.

\punk{Discretization and generation by ideals}
\label{e2i}

Some \eqr s on the diagram on page~\pageref{p-p} are
explicitly generated by ideals, like
$\rE_i\yt i=0,1,2,3.$
Some other \er s are defined differently.
It will be shown below (Section~\ref{rosen})
that {\ubf any} Borel \er\ $\rE$ is Borel reducible to a
\er\ of the form $\rei,$ $\cI$ a Borel ideal.
On the other hand, $\fco\zi\bel1\zi\bel\iy$ turn out to
be Borel equivalent to some meaningful Borel ideals. 
Moreover, these \eqr s admit ``discretization'' by means
of restriction to certain subsets of $\rtn.$

\bdf
\label{dx}
We define 
$\dX=\prod_{n\in\dN}X_n=\ens{x\in\rtn}{\kaz n\:(x(n)\in X_n)},$
where
$X_n=\ans{\frac0{2^n},\frac1{2^n},\dots,\frac{2^n}{2^n}}$.
\index{zzXXn@$\dX=\prod_nX_n$}%
\edf

\ble
\label{redx}
${\fco}\reb {{\fco}\res{\dX}}$ and\/
${\bel p}\reb {{\bel p}\res{\dX}}$ for any\/ $1\le p<\iy.$

On the other hand, ${\bel\iy}\reb {{\bel\iy}\res{\ztn}}$.
\ele
\bpf
We first show that ${\fco} \reb {{\fco}\res{\oin}}.$
Let $\pi$ be any bijection of $\dN\ti\dZ$ onto $\dN.$ 
For $x\in\rtn,$ define $\vt(x)\in\oin$ as follows. 
Suppose that $k=\pi(n,\eta)$ ($\eta\in\dZ$). 
If $\eta\le x(n)<\eta+1$ then let $\vt(x)(k)=x(n).$ 
If $x(n)\ge \eta+1$ then put $\vt(x)(k)=1.$ 
If $x(n)<\eta$ then put $\vt(x)(k)=0.$ 
Then $\vt$ is a Borel reduction of
${\fco}$ to ${\fco}\res{\oin}.$
Now we prove that ${{\fco}\res{\oin}}\reb {{\fco}\res{\dX}}.$
For $x\in\oin$ define $\psi(x)\in\dX$ so that $\psi(x)(n)$
the largest number of the form $\frac i{2^n}\zd 0\le i\le{2^n}$
smaller than $x(n).$
Then obviously $x\fco {\psi(x)}$ holds for any $x\in\oin,$
and hence $\psi$ is a Borel reduction of ${\fco}\res{\oin}$ to
${\fco}\res{\dX}$.

Thus ${\fco}\reb{{\fco}\res{\dX}},$ 
and hence in fact ${\fco} \eqb {{\fco}\res\dX}.$

The argument for $\bel1$ is pretty similar.
The result for $\bel\iy$ is obvious:
given $x\in\rtn,$ replace any $x(n)$ by the largest integer
value $\le x(n)$.

The version for $\bel p\yi1<p<\iy,$ needs some comments in
the first part (reduction to $\oin$).
Note that if $\eta\in\dZ$ and
$\eta-1\le x(n)<\eta<\za\le y(n)<\za+1$
then the value $(y(n)-x(n))^p$ in the distance 
$\nor{y-x}_p=(\sum_n|y(n)-x(n)|^p)^{\frac1p}$
is replaced by $(\za-\eta)+(\eta-x(n))^p+(y(n)-\za)^p$
in $\nor{\vt(y)-\vt(x)}_p.$
Thus if this happens infinitely many times then both
distances are infinite, while otherwise this case can be
neglected.
Further, if $\eta-1\le x(n)<\eta\le y(n)<\eta+1$
then $(y(n)-x(n))^p$ in $\nor{y-x}_p$
is replaced by $(\eta-x(n))^p+(y(n)-\eta)^p$
in $\nor{\vt(y)-\vt(x)}_p.$
However
$(\eta-x(n))^p+(y(n)-\eta)^p\le (y(n)-x(n))^p\le
2^{p-1} ((\eta-x(n))^p+(y(n)-\eta)^p),$
and hence these parts of the sums in $\nor{y-x}_p$
and $\nor{\vt(y)-\vt(x)}_p$ differ from each other by
a factor between $1$ and $2^{p-1}.$
Finally, if $\eta\le x(n)\zi y(n)<\eta+1$ for one and the
same $\eta\in\dZ$ then the term $(y(n)-x(n))^p$ in
$\nor{y-x}_p$ appears unchanged in $\nor{\vt(y)-\vt(x)}_p.$
Thus totally $\nor{y-x}_p$ is finite iff so is
$\nor{\vt(y)-\vt(x)}_p$.
\epf

\ble[{{\rm Oliver \cite{odis}}}]
\label{co=d}
${\fco}$ is\/ $\eqb$ to the \er\/ 
$\rzo=\rE_{\zo}$.
\ele
\bpf
Prove that ${\fco}\reb \rzo.$
It suffices,
by Lemma~\ref{redx}, to define a Borel reduction 
${{\fco}\res\dX}\to {\rzo},$
\ie, a Borel map $\vt:\dX\to\pn$ such 
that ${x\fco y}\leqv{\vt(x)\sd\vt(y)\in\zo}$ 
for all $x,\,y\in\dX.$ 
Let $x\in\dX.$ 
Then, for any $n,$ we have $x(n)=\fras{k(n)}{2^n}$ 
for some natural $k(n)\le 2^n.$ 
The value of $k(n)$ determines the intersection 
$\vt(x)\cap\ir{2^n}{2^{n+1}}:$
for each $j<2^n,$ we define $2^n+j\in\vt(x)$ iff $j<k(n).$
Then   
$x(n)=\frac{\#(\vt(x)\cap\ir{2^n}{2^{n+1}})}{2^n}$ 
for any $n,$ and moreover  
$|y(n)-x(n)|=
\fras{\#([\vt(x)\sd\vt(y)]\cap\ir{2^n}{2^{n+1}})}{2^n}$ 
для всех $x,y\in\dX$ и $n.$
This easily implies that $\vt$ is as required.

To prove $\rzo\reb {\fco},$ 
we have to define a Borel map $\vt:\pn\to\rtn$ such 
that ${X\sd Y\in\zo}\leqv{\vt(X)\fco\vt(Y)}.$ 
Most elementary ideas like $\vt(X)(n)=\frac{\#(X\cap\ir0n)}n$ 
do not work, the right way is based on the following 
observation: for any sets $s,\,t\sq\ir0n$ to satisfy 
$\#(s\sd t)\le k$ it is necessary and sufficient that 
$|\#(s\sd z)-\#(t\sd z)|\le k$ 
for any $z\sq\ir0n.$ 
To make use of this fact, let us fix an enumeration 
(with repetitions) $\sis{z_j}{j\in\dN}$ 
of all finite subsets of $\dN$ such that 
\dm
\ens{z_j}{2^n\le j<2^{n+1}}
\quad=\quad\hbox{all subsets of $\ir0n$}
\dm
for every $n.$ 
Define, for any $X\in\pn$ and $2^n\le j<2^{n+1},$ 
$\vt(X)(j)=\frac{\#(X\cap z_j)}n.$ 
Then $\vt:\pn\to[0,1]^\dN$ is a required reduction. 
\epf

Recall that for any sequence of reals $r_n\ge0,$
$\ern$ is an \eqr\ on $\pn$ generated by
the ideal $\srn=\ens{x\sq\dN}{\sum_{n\in x}r_n<\piy}$.

\ble
[{{\rm Attributed to Kechris in \cite[2.4]{h-ban}}}]
\label{l1=s}
If\/ $r_n\ge 0\yt r_n\to0\yt\sum_nr_n=\piy$ then\/
$\ern\eqb {\bel1}.$
In particular,\/ $\Ed=\eun$ satisfies\/
$\Ed\eqb{\bel1}.$
\ele
\bpf
To prove $\ern\reb {\bel1},$ 
define $\vt(x)\in\rtn$ for any $x\in\pn$ as 
follows: $\vt(x)(n)=r_n$ for any $n\in x,$ and
$\vt(x)(n)=0$ for any other $n.$ 
Then ${x\sd y\in\sui{r_n}}\leqv{\vt(x)\bel1\vt(y)},$ 
as required.

To prove the other direction, it suffices
to define a Borel reduction of 
${\bel1\res\dX}$ to $\ern.$
We can associate a (generally, infinite) set $s_{nk}\sq\dN$
with any pair of $n$ and $k<2^n,$ so that 
the sets $s_{nk}$ are pairwise disjoint and 
$\sum_{j\in s_{nk}}r_j=2^{-n}.$ 
The map   
$\vt(x)=\bigcup_n\bigcup_{k<2^nx(n)}s_{nk}\yt x\in\dX,$
is the reduction required.
\epf

\punk{Summables irreducible to density-0}
\label{summ}

The \dd\reb independence of $\bel1$ and $\fco,$ two
best known ``Banach'' \eqr s, is quite important.
In one direction it is provided by \ref{t>s2} of
the next theorem.
The other direction actually follows from Lemma~\ref{e3}.


Is there any example of Borel ideals $\cI\reb\cJ$ 
which do not satisfy $\cI\odl\cJ$? \  
Typically the reductions found 
to witness $\cI\reb\cJ$ are \dd\sd homomorphisms, and 
even better maps. 
The following lemma proves that Borel reduction yields 
\dd\orbpp reduction in quite a representative case. 
Let us say that $\cI\orbpp\cJ$ {\it holds exponentially\/} 
if there is a map $i\mapsto w_i$ withessing
$\cI\orbpp\cJ$~\snos
{Thus we have pairwise disjoint finite non-empty sets
$w_k\sq\dN$ (assuming $\cI,\cJ$ are ideals over $\dN$)
such that $A\in\cI\leqv w_A=\bigcup_{k\in A}w_k\in\cJ,$
and $\tmax w_k<\tmin w_{k+1}$.}
and in addition a sequence of natural numbers $k_i$ with 
$w_i\sq\ir{k_i}{k_{i+1}}$ and $k_{i+1}\ge 2k_i$. 

\bte
\label{t:>s}
Suppose that\/ $r_n\ge 0\yt r_n\to0\yt\sum_nr_n=\piy.$
Then
\ben
\renu
\itsep
\itla{t>s1}
{\rm(Farah~\cite[2.1]{f-tsir})} \ 
If\/ $\cJ$ is a Borel P-ideal and\/ $\sui{r_n}\reb\cJ$
then we have\/ $\sui{r_n}\orbpp\cJ$ exponentially$;$

\itla{t>s2}
{\rm(Hjorth~\cite{h-ban})} \
$\sui{r_n}$ is not 
Borel-reducible to\/ $\zo$.
\een
\ete
\bpf
\ref{t>s1}
Let a Borel $\vt:\pn\to\pn$ witness $\sui{r_n}\reb\cJ.$ 
Let, according to Theorem~\ref{sol}, $\nu$ be a \lsc\ 
submeasure on $\dN$ with $\cJ=\Exh_\nu.$ 
\vyk{
We are going to ``continualize'' $\vt$ in certain way: 
the new, continuous reduction will be a superposition of 
two continuous maps, $\ga$ and $\xi,$ where $\ga$ will 
reduce $\pn/\sui{r_n}$ to a certain ``generic'' part of 
itself, while $\xi$ will reduce the latter to $\pn/\cJ.$ 
}%
The construction makes use of stabilizers. 
%
Suppose that $n\in\dN.$ 
If $u,\,v\sq\ir0n$ then $(u\cup X)\sd(v\cup X)\in\sui{r_n}$ 
for any $X\sq\ir n\piy,$ hence, 
$\vt(u\cup X)\sd\vt(v\cup X)\in\cJ.$ 
It follows, by the choice of the submeasure $\nu,$
that for any $\ve>0$ 
there are numbers $n'>k>n$ and a set $s\sq\ir n{n'}$ 
such that    
\dm
\nu\skl(\vt(u\cup s\cup X)\sd\vt(v\cup s\cup X))
\cap\ir k\iy\skp\,<\,\ve
\dm
holds for all $u,\,v\sq\ir0n$ 
and all generic~\snos
{In the course of the proof, ``generic'' means Cohen-generic
over a 
sufficiently large countable model of a big enough fragment 
of $\ZFC$.}
$X\sq\ir{n'}\iy$. 

This allows us to define an increasing sequence of 
natural numbers 
$0=k_0=a_0<b_0<k_1<a_1<b_1<k_2<...$ and, for any 
$i,$ a set $s_i\sq\ir{b_i}{a_{i+1}}$ such that, for all 
generic $X,\,Y\sq\ir {a_{i+1}}\iy$ and all $u,\,v\sq\ir0{b_i},$
we have 
\ben
\tenu{(\arabic{enumi})}
\itla{stab1}\msur
$\nu\skl(\vt(u\cup s_i\cup X)\sd\vt(v\cup s_i\cup X))
\cap\ir{k_{i+1}}\iy\skp<2^{-i}$; 

\itla{stab2}\msur
$\skl\vt(u\cup s_i\cup X)\sd\vt(u\cup s_i\cup Y)\skp
\cap\ir0{k_{i+1}}=\pu$;

\itla{stabg}
any $Z\sq\dN,$ satisfying $Z\cap\ir{b_i}{a_{i+1}}=s_i$ 
for infinitely many $i,$ is generic;

\itla{stabk}\msur
$k_{i+1}\ge 2k_i$ for all $i$;
\een
and in addition, under the assumptions on $\sis{r_n}{}$, 
\ben
\tenu{(\arabic{enumi})}
\addtocounter{enumi}4
\itla{stab3}
there is a set $g_i\sq\ir{a_i}{b_i}$ such that 
$|r_i-\sum_{n\in g_i}r_n|<2^{-i}$.
\een
It follows from \ref{stab3} that 
$A\mapsto g_A=\bigcup_{i\in A}g_i$ 
is a reduction of $\sui{r_n}$ to $\sui{r_n}\res N,$ where 
$N=\bigcup_i\ir {a_i}{b_i}.$
Let $S=\bigcup_i s_i;$ note that $S\cap N=\pu.$

Put $\xi(Z)=\vt(Z\cup S)\sd \vt(S)$ for any $Z\sq N.$ 
\imar{why $\sd\vt(S)$ added?}%
Then, for any sets $X,\,Y\sq N,$ 
\dm
{X\sd Y \in\sui{r_n}}\eqv
{\vt(X\cup S)\sd\vt(Y\cup S)\in\cJ}\eqv
{\xi(X)\sd\xi(Y)\in\cJ,}
\dm
thus $\xi$ reduces $\sui{r_n}\res N$ to $\cJ.$ 
Now put $w_i=\xi(g_i)\cap\ir{k_{i}}{k_{i+1}}$ 
and $w_A=\bigcup_{i\in A} w_i.$ 
We assert that the map $i\mapsto w_i$ proves 
$\sui{r_n}\orbpp\cJ.$ 
In view of the above, 
it remains to show that $\xi(g_A)\sd w_A\in\cJ$ for 
any $A\in\pn$.

As $\cJ=\Exh_\nu,$ it suffices to demonstrate that 
$\nu\skl w_i\sd(\xi(g_A)\cap\ir{k_{i}}{k_{i+1}})\skp
<2^{-i}$ 
for all $i\in A$ while 
$\nu(\xi(g_A)\cap\ir{k_{i}}{k_{i+1}})<2^{-i}$ 
for $i\nin A.$ 
After dropping the common term $\vt(S),$    
it suffices to check that 
\ben
\tenu{(\alph{enumi})}
\itla{111}
$\nu\skl (\vt(g_i\cup S)\sd\vt(g_A\cup S))
\cap\ir{k_{i}}{k_{i+1}}\skp<2^{-i}$ 
for all $i\in A$ while 

\itla{222}
$\nu\skl(\vt(S)\sd\vt(g_A\cup S))\cap
\ir{k_{i}}{k_{i+1}}\skp<2^{-i}$ 
for $i\nin A.$ 
\een
Note that, as any set of the form $X\cup S,$ where $S\sq N,$ 
is generic by \ref{stabg}. 
It follows, by \ref{stab2}, that we can assume, in \ref{111} 
and \ref{222}, that $A\sq[0,i],$ \ie, resp.\ $\tmax A=i$ 
and $\tmax A<i.$ 
We can finally apply \ref{stab1}, with 
$u=A\cup\bigcup_{j<i}s_j,\msur$ $X=\bigcup_{j>i}s_j,$ and 
$v=u_i \cup\bigcup_{j<i}s_j$ if $i\in A$ while 
$v=\bigcup_{j<i}s_j$ if $i\nin A$.

\ref{t>s2}
Otherwise $\sui{r_n}\orbpp\zo$ exponentially by \ref{t>s1}. 
Let this be witnessed by $i\mapsto w_i$ and a sequence of 
numbers $k_i,$ so that $k_{i+1}\ge 2k_i$ and 
$w_i\sq\ir{k_i}{k_{i+1}}.$ 
If $d_i=\frac{\#(w_i)}{k_{i+1}}\to0$ 
then easily $\bigcup_iw_i\in\zo$ by the choice of $\sis{k_i}{}.$ 
Otherwise there is a set $A\in\sui{r_n}$ such that 
$d_i>\ve$ for all $i\in A$ and one and the same $\ve>0,$ 
so that $w_A=\bigcup_{i\in A}w_i\nin\zo.$ 
In both cases we have a contradiction with the assumption 
that the map $i\mapsto w_i$ witnesses $\sui{r_n}\orbpp\zo$.
\epf

\bqu
\label{far-s}
Farah~\cite{f-tsir} points out that
Theorem~\ref{t:>s}\ref{t>s1} also holds for $\ofi$
(instead of $\sui{r_n}$) and asks for which other ideals 
it is true.
\equ

\vyk{

Extension: full picture of summables. 

\bqu 
What about $\fco/\sui{1/n}$?
\equ

\bre
\label{goodm}
Claim \ref{assu4} in the proof of the theorem shows that
the map $\vt'(x)=\bigcup_n t^{x(n)}_n$ is still a
(continuous) reduction of $\ern$ to ${\fco}\res\dX:$
indeed, $|\vt(x)(k)-\vt'(x)(k)|\le 2^{-n}$ for any $k$
in $[\nu_n,\nu_{n+1})$.
\ere
}

\vyk{
The \dd\reb independence of $\bel1$ and $\fco,$ two
most elementary ``Banach'' \eqr s, is quite important.
In one direction it is provided by the next theorem.
The other direction actually follows from Lemma~\ref{e3}.

\bte[{{\rm Hjorth~\cite{h-ban}}}]
\label{lico}
In the assumptions of Lemma~\ref{l1=s}, 
$\ern$ is Borel irreducible to\/ $\fco.$
In particular, the ER\/ $\Ed\eqb {\bel1}$
is Borel irreducible to\/ $\fco.$~\snos
{It is worth to note that the corresponding ideals
(in the sense of \ref{l1=s} and \ref{co=d} $\sui{1/n}$
and $\Zo$ are somewhat different: $\sui{1/n}\sneq\Zo.$
To prove $\sui{1/n}\sq \zo,$ suppose that $X\nin\zo.$
Thus there is $\ve>0$ such that 
$\#(X\cap\ir0n)>2n\ve$ for infinitely many numbers $n.$ 
Take any such $n,$ big enough \vrt\ $\ve\obr.$ 
Then clearly $\#(X\cap\ir{n\ve}n)>n\ve.$ 
It follows that the sum 
$\sum_{k\in X,\;k>n\ve}k\obr$ is not smaller than 
$\frac{n\ve}n=\ve.$ 
As this holds for infinitely many $n,$ the set $X$ cannot 
belong to $\sui{1/n}.$
To define now a set in $\zo\dif\sui{1/n},$ 
let $m_j$ be the entire part of $j\log j.$ 
Then $X=\ens{m_j}{j\ge 2}\nin\sui{1/n}$ because the 
integral $\int_2^\iy(x\log x)\obr dx$ diverges. 
On the other hand, easily $X\in\zo.$
Generally, $\zo\not\sq \srn$ holds whenever $r_n\ge 0$
and $r_n\to 0.$ 
}
\ete
\bpf
Otherwise there is a Borel reduction $\vt$ of
$\ern$ to ${{\fco}\res\dX}.$
It will be more convenient to view $\ern$ as an ER
on $\dn$ defined so that $x\ern y$ iff $x\sd y\in \srn,$
where, in this case, $x\sd y=\ens{k}{x(k)\ne y(k)}.$
Contradiction will be obtained in the course of
several simplifying transformations of $\vt$ that
involve the stabilizers technique.

{\it Step 1\/}:
making $\vt$ and some other functions continuous.
To define those other functions, note the following.
Given $n\in\dN$ and $z\in 2^{[n,\piy)},$ we have
${(u\cup z)}\ern{(v\cup z)},$ and hence
$\tlim_{j\to\iy}|\vt(u\cup z)(j)-\vt(v\cup z)(j)|=0,$
for any two finite sequences $u,v\in 2^n.$
Let, for any $\ve>0,$ $\mu_{uv\ve}(z)$ be equal to the
least index $m>n$ with  
$|\vt(u\cup z)(j)-\vt(v\cup z)(j)|<\ve$
for all $j\ge m.$
Obviously all maps $\mu_{uv\ve}$ are Borel together with $\vt$.

There is a dense $\Gd$ set $D=\bigcap_nD_n\sq\dn,$
each $D_n$ being open dense in $\pn$ and $D_n\sq D_{n+1},$
such that both $\vt$ and all maps $\mu_{uv\ve},$ $\ve>0$
rational, are continuous on $D.$
Following the stabilizers argument (Lemma~\ref{l:bc}),
we define an increasing sequence of numbers 
$0=j_0<j_1<j_2<...$ and, for any $n,$ a sequence 
$s_n\in 2^{\ij{n}{n+1}},$ such that 
$x\res \ij{n}{n+1}=s_n\imp x\in D_n.$   
We obtain a partition $\dN=N_1\cup N_2$ onto the sets  
$N_1=\bigcup_n \ij{2n}{2n+1}$ and  
$N_2=\bigcup_n \ij{2n+1}{2n+2},$ and functions
$\sg_1=\bigcup_n s_{2n}\in 2^{N_1}$ and  
$\sg_2=\bigcup_n s_{2n+1}\in 2^{N_2}$
such that the sets
$X_i=\ens{x\in\dn}{\sg_i=x\res N_i}\yt i=1,2,$ 
are subsets of $D.$

{\it Step 2\/}:
contraction.
For at least one of $i=1,2,$ $\sum_{n\in N_i}r_n$ is
still infinite ---
suppose that $\sum_{n\in N_1}r_n=\piy.$ 
The increasing bijection $h:\dN\onto N_1$
induces a homeomorphism $H:\dn\onto X_2$
by $H(x)(h(n))=x(n)\zd\kaz n,$ and $H(x)\res N_2=\sg_2.$
Then the \dd Hpreimage of $\ern\res X_2$ is equal to
$\erpn,$ where $r'_n=r_{h(n)}\zd\kaz n.$
Moreover, $\vt'(x)=\vt(H(x))$ is a reduction of
$\erpn$ to ${{\fco}\res\dX}$ --- and $\vt'$ is
continuous since so is $\vt\res X_2.$
Finally, the functions $\mu'_{uv\ve},$ defined as above
but starting from $\erpn$ and $\vt',$ are similar
derivatives of the maps $\mu_{uv}\res X_2$
(for slightly different $u,v$) ---
and hence continuous as well for any rational $\ve>0.$
Thus it can be \noo\ assumed
from the very beginning that
\ben
\cenu
\itla{assu1}
The reduction
$\vt:\dn\to \dX$ and the maps $\mu_{uv\ve},$
$\ve$ rational,
are {\ubf continuous}.%
\een

{\it Step 3\/}:
restricting the values.
In the assumption \ref{assu1},
for any   
$j\in\dN$ and any rational $\ve>0$ there exist a number
$\nu=\nu(j,\ve)>n$ and
$s\in 2^{[j,\nu)}$ such that
\ben
\tenu{(\roman{enumi})}
\itla{s32}
for any $y\in Y_s$
the value of $\vt(y)\res{[0,\nu)}$ is determined by
$y\res{[0,j)},$
and

\itla{s31}
for any $k\ge \nu$ and
$y,y'\in Y_s=\ens{y\in\dn}{y\res{[j,\nu)}=s}$ with
$y\res{[\nu,\piy)}=y'\res{[\nu,\piy)},$ we have 
$|\vt(y)(k)-\vt(y')(k)|< \ve$.
\een 
To fulfill \ref{s31}, define
$\nu=\nu(j,\ve)=
\tsup_{u,v\in2^a,\;z\in 2^{[j,\piy)}}\mu_{uv\ve}(z)$
and choose any $s\in 2^{[j,\nu)}.$
The supremum is finite since  
$\mu_{uv\ve}$ are continuous functions.
Then to fulfill \ref{s32} extend $s,$ if necessarily, 
using the fact that $\vt$ is continuous.

This observation allows us to repeat the stabilizers
construction,
as in Step~1, in such a way that $j_{n+1}=\nu(j_n,2^{-n}).$ 
Subsequent contraction, as on Step~2, \noo\ reduces us to
the case when, in addition to \ref{assu1}, there is a
sequence of numbers $\nu_0<\nu_1<\dots$ with $n<\nu_n,$
satisfying the following:
\ben
\cenu
\itla{assu3}
For any $n$ and $y\in\dn$ the value of
$\vt(y)\res{[0,\nu_n)}$
is determined by $y\res{[0,n)}$.
\enuci

\itla{assu2}
The inequality $|\vt(y)(k)-\vt(y')(k)|< 2^{-n}$
holds for any $n,$ any $k\ge \nu_n,$ and any
$y,y'\in \dn$ with
$y\res{[n,\piy)}=y'\res{[n,\piy)}$.
\een

{\it Step 4\/}:
getting contradiction.
For any $n,$ define $\hn\in\dn$ by $\hn(n)=1$
and $\hn(k)=0$ for $k\ne n.$
Let $\fo\in\dn$ be the constant $0.$
Put $t^1_n=\vt(\hn)\res{[\nu_n,\nu_{n+1})},$ and
$t^0_n=\vt(\fo)\res{[\nu_n,\nu_{n+1})};$
thus $t^i_n\in \prod_{k=\nu_n}^{\nu_{n+1}-1}X_k.$
We claim the following:
\ben
\cenu
\itla{assu4}
{\it For any\/ $x\in\dn$ and\/ $k\in[\nu_n,\nu_{n+1}),$
$|\vt(x)(k)-t^{x(n)}_n(k)|\le 2^{-n}.$\/}
Thus $\vt(x)\res{[\nu_n,\nu_{n+1})}$
depends only on $x(n)$ modulo an addendum $2^{-n}\to0$.
\een
To prove \ref{assu4}, let $y\in\rtn$ satisfy $y(k)=x(k)$
for $k\ne n$ and $y(k)=0$ for all bigger $k.$
Then
$\vt(x)\res{[0,\nu_{n+1})}=\vt(y)\res{[0,\nu_{n+1})}$
by \ref{assu3}.
On the other hand, $y\res{[n,\piy)}$ coincides with either
$\hn\res{[n,\piy)}$ or $\fo\res{[n,\piy)}$
in cases, resp., $x(n)=1,0.$
Thus $|\vt(y)(k)-t^{x(n)}_n(k)|<2^{-n}$ for any
$k\ge\nu_n$ by \ref{assu2}.
We conclude that $|\vt(x)(k)-t^{x(n)}_n(k)|\le 2^{-n}$
for any $\nu_n\le k<\nu_{n+1},$ as required.


{\it Case 1\/}: 
$w_n=\tmax_{\nu_n\le k<\nu_{n+1}}|t^1_n(k)-t^0_n(k)|\to 0.$
Let $\fr\in\dn$ be the constant $1,$ thus $x\bel1 \fo$
fails. 
However, by \ref{assu4},
$|\vt(x)(k)-\vt(\fo)(k)|\le 2^{-n}+w_n$
whenever $\nu_n\le k<\nu_{n+1},$ and hence
$\vt(x)\fco \vt(\fo)$ holds, contradiction.

{\it Case 2\/}: 
$W(\ve)=\ens{n}{w_n>\ve}$ is infinite for some $\ve>0.$
Choose a
set $W\sq W(\ve)$ such that $\sum_{n\in W}r_n<\piy$ and
define $x(n)=1$ for $n\in W$ and $x(n)=0$ otherwise.
Then $x\ern \fo.$
On the other hand, by \ref{assu4},
$|\vt(x)(k)-\vt(\fo)(k)|\ge \ve-2^{-n}$
for all $n$ and all $k$ in $[\nu_n,\nu_{n+1}),$ and hence
$\vt(x)\fco \vt(\fo)$ fails, contradiction.
\epf
}

\punk
{The family $\bel p$}
\label{hd1}

It follows from the next theorem that Borel reducibility
between \eqr s $\bel p\yt 1\le p<\iy,$ is fully determined
by the value of $p$.

\bte[{{\rm Dougherty -- Hjorth \cite{dh}}}]
\label{hd-t}
If\/ $1\le p<q<\iy$ then\/ ${\bel p}\rebs{\bel q}$.
\ete
\bpf
{\ubf Part 1\/}:  
show that ${\bel q}\not\reb{\bel p}.$  

By Lemma~\ref{redx}, it suffices to prove that 
${\bel q\res\dX}\not\reb{\bel p\res \dX}.$ 
Suppose, on the contrary, that $\vt:\dX\to\dX$
is a Borel reduction of
${\bel q\res\dX}$ to ${\bel p\res \dX}.$
Arguing as in the proof of Theorem~\ref{t:>s},
we can reduce the general case to the case when
there exist increasing sequences of numbers
$0=j(0)<j(1)<j(2)<\dots$ and $0=a_0<a_1<a_2<\dots$ 
and a map $\tau:\dY\to\dX,$ where
$\dY=\prod_{n=0}^\iy X_{j(n)},$
which reduces ${\bel q\res\dY}$ to $\bel p\res \dX$ 
and has the form
$\tau(x)=\bigcup_{n\in\dN}t^{x(n)}_n,$
where $t^r_n\in \prod_{k=a_n}^{a_{n+1}-1}X_k$ for any
$r\in X_{j_n}.$
(See Definition~\ref{dx}.)

{\it Case 1\/}: 
there are\/ $c>0$ and a number\/ $N$ such that\/ 
$\nr{\tau^1_n-\tau^0_n}p\ge c$ for all\/ $n\ge N.$ 
Since $p<q,$ there is a non-decreasing sequence of natural 
numbers $i_n\le {j_n}\yt n=0,1,2,\dots,$ 
such that $\sum_n2^{p(i_n-j_n)}$ 
diverges but $\sum_n2^{q(i_n-j_n)}$ converges. 
({\it Hint\/}: $i_n\approx j_n-p\obr\log_2n$.) 

Now consider any $n\ge N.$ 
As $\nr{\tau^1_n-\tau^0_n}p\ge c$ and because $\nr{...}p$ 
is a norm, there exists a pair of rationals $u(n)<v(n)$ in 
$X_{j_n}$ with $v(n)-u(n)=2^{i_n-j_n}$ and 
$\nr{\tau^{v(n)}_n-\tau^{u(n)}_n}p\ge c\,2^{i_n-j_n}.$ 
In addition, put $u(n)=v(n)=0$ for $n<N.$ 
Then the \dd{\bel q}distance between the infinite
sequences $u=\sis{u(n)}{n\in\dN}$ 
and $v=\sis{v(n)}{n\in\dN}$ is equal to 
$\sum_{n=N}^\iy 2^{q(i_n-j_n)}<+\iy,$
while the \dd{\bel p}distance between $\tau(u)$ and
$\tau(v)$ is non-smaller than
$\sum_{n=N}^\iy c^p\,2^{p(i_n-j_n)}=\iy.$ 
But this contradicts the assumption that $\tau$
is a reduction.

{\it Case 2\/}: 
otherwise.
Then there is a strictly  increasing sequence 
$n_0<n_1<n_2<\dots$ with 
$\nr{\tau^1_{n_k}-\tau^0_{n_k}}p\le 2^{-k}$ for all $k.$ 
Let now $x\in \dY$ be the constant $0$ while $y\in\dY$ 
be defined by $y(n_k)=1\zd\kaz k$ and $y(n)=0$ 
for all other $n.$ 
Then $x\bel q y$ fails
($|y(n)-x(n)|\not\to 0$)
but $\tau(x)\bel p \tau(y)$ 
holds, contradiction.

{\ubf Part 2\/}: 
show that ${\bel p}\reb{\bel q}.$ 

It suffices to prove that ${\bel p\res\oin}\reb{\bel q}$
(Lemma~\ref{redx}).
We \noo\ assume that $q<2p$:
any bigger $q$ can be approached in several steps. 
For $\vec x=\ang{x,y}\in\dR^2,$ let 
$\nr{\vec x}h=(x^h+y^h)^{1/h}.$ 

\ble
\label{hd+1}
For any\/ $\frac12<\al<1$ there is a continuous 
map\/ $K_\al:[0,1]\to[0,1]^2$ and positive real numbers\/ 
$m<M$ such that for all\/ $x<y$ in $[0,1]$ we have\/
$m(y-x)^\al\le\nr{K_\al(y)-K_\al(x)}2\le M(y-x)^\al$. 
\ele 
\bpf[{{\rm Lemma}}]
The construction of such a map $K$ can be easier described
in terms of fractal geometry rather than by an analytic
expression. 
Let $r=4^{-\al},$ so that $\frac14<r<\frac12$ and 
$\al=-\log_4r.$ 
Starting with the segment $[(0,0)\,,\,(1,0)]$
of the horisontal axis of the cartesian plane, we replace it
by four smaller segments of length $r$ each
(thin lines on Fig.~2, left).
Each of them we replace by four segments of
length $r^2$ (thin lines on Fig.~2, right). 
And so on, infinitely many steps.
The resulting curve $K$ is parametrized by giving the
vertices of the polygons values equal to multiples of $4^{-n},$
$n$ being the number of the polygon.
For instance, the vertices of the left polygon on Fig.~2 are
given values $0,\frac14,\frac12,\frac34,1.$
%

\begin{figure}[h]
\label{fig2}


\setlength{\unitlength}{0.8mm} 
\begin{picture}(130,35)(2,0) 

\label{f-f}

\put(5,4) {(0,0)}
\put(65,4) {(1,0)}
\put(95,4) {(0,0)}
\put(155,4) {(1,0)}

\put(10,10){\tob}
\put(30,10){\tob}
\put(50,10){\tob}
\put(70,10){\tob}
\put(40,27.32){\tob}

\put(10,10){\line(1,0){20}}
\put(50,10){\line(1,0){20}}

\put(100,10){\tob}
\put(120,10){\tob}
\put(140,10){\tob}
\put(160,10){\tob}

\put(106.67,10){\tob}
\put(113.33,10){\tob}
\put(146.67,10){\tob}
\put(153.33,10){\tob}

\put(130,27.32){\tob}
\put(123.33,15.77){\tob}
\put(126.67,21.55){\tob}
\put(133.33,21.55){\tob}
\put(136.67,15.77){\tob}

\put(110,15.77){\tob}
\put(150,15.77){\tob}

\put(120,21.55){\tob}
\put(140,21.55){\tob}

\put(100,10){\line(1,0){6.67}}
\put(140,10){\line(1,0){6.67}}
\put(120,10){\line(-1,0){6.67}}
\put(160,10){\line(-1,0){6.67}}


\qbezier(30,10)(35,18.66)(40,27.32)
\qbezier(50,10)(45,18.66)(40,27.32)

\qbezier(106.67,10)(108.33,12.88)(110,15.77)
\qbezier(113.33,10)(111.67,12.88)(110,15.77)
\qbezier(146.67,10)(148.33,12.88)(150,15.77)
\qbezier(153.33,10)(151.67,12.88)(150,15.77)

\qbezier(120,10)(121.67,12.88)(123.33,15.77)
\qbezier(140,10)(138.33,12.88)(136.67,15.77)

\qbezier(120,21.55)(121.67,18.66)(123.33,15.77)
\qbezier(140,21.55)(138.33,18.66)(136.67,15.77)

\qbezier(120,21.55)(123.33,21.55)(126.67,21.55)
\qbezier(140,21.55)(136.67,21.55)(133.33,21.55)

\qbezier(133.33,21.55)(131.67,24.44)(130,27.32)
\qbezier(126.67,21.55)(128.33,24.44)(130,27.32)

\linethickness{0.5mm}



\linethickness{0.4mm}

\qbezier[20](10,10)(25,18.66)(40,27.32)
\qbezier[20](70,10)(55,18.66)(40,27.32)
\qbezier[30](10,9.3)(40,9.3)(70,9.3)

\qbezier[10](100,9.3)(110,9.3)(120,9.3)
\qbezier[6](100,10)(105,12.88)(110,15.77)
\qbezier[6](120,10)(115,12.88)(110,15.77)

\linethickness{0.2mm}

\qbezier[15](120,10)(120,15)(120,21.55)
\qbezier[15](120,21.55)(125,24.44)(130,27.32)
\qbezier[20](120.1,9.3)(125.1,17.96)(130.1,26.62)

\end{picture}

\caption{$r=\frac13$, \ left: step 1, \ right: step 2}

\end{figure}
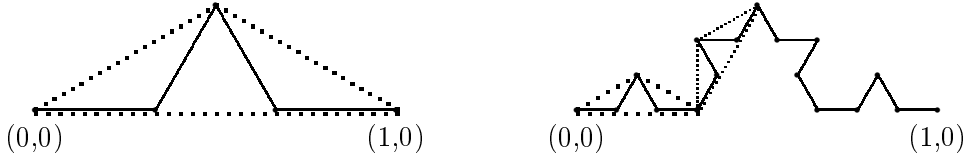

%
Note that the curve $K:[0,1]\to[0,1]^2,$
approximated by the polygons, is
bounded by certain triangles built on the sides of the
polygons.
For instance, the whole curve lies inside the triangle
bounded by dotted lines in Fig.~2, left.
(The dotted line that follows the basic side
$[(0,0)\,,\,(1,0)]$ of the triangle is drawn  
slightly below its true position.)
Further, the parts $0\le t\le\frac14$ and
$\frac14\le t\le\frac12$ of the curve 
lie inside the triangles bounded by (slightly different)
dotted lines in Fig.~2, right.
And so on.
Let us call those triangles {\it bounding triangles\/}.
\vyk{
It is quite obvious geometrically that for any $t$
in $[\frac{i-1}{4^n},\frac{i}{4^n}]$ ($1\le i\le4^n$)
the euclidean distance $\nr{K_r(t)-C^i_n}2$ between
$K(t)$ and the center $C^i_n$ of the straight interval
$[K(\frac{i-1}{4^n})\,,\,K(\frac{i}{4^n})]$
(inside the unit square)
satisfies $\nr{K_r(t)-C^i_n}2\le \frac{r^n}2.$
}

To prove the inequality of the lemma, consider any pair of
reals $x<y\in[0,1].$
Let $n$ be the least number such that $x,y$ belong to
non-adjacent intervals, resp.,
$[\fras{i-1}{4^n},\fras{i}{4^n}]$ and
$[\fras{j-1}{4^n},\fras{j}{4^n}],$ with $j>i+1.$
Then $4^{-n}\le |y-x|\le 8\cdot 4^{-n}.$

The points $K(x)$ and $K(y)$ then belong to one and the
same side or adjacent sides of the \dd{n-1}th polygon.
Let $C$ be a common vertice of these sides.
It is quite clear geometrically that the euclidean distances
from $K(x)$ and $K(y)$ to $C$ do not exceed $r^{n-1}$
(the length of the side), thus 
$\nr{K(x)-K(y)}2\le 2\,r^{n-1}.$

Estimation from below needs more work.
The points $K(x)\zd K(y)$ belong to the bounding triangles
built on the segments, resp.,
$[K(\frac{i-1}{4^n})\,,\,K(\frac{i}{4^n})]$ and
$[K(\frac{j-1}{4^n})\,,\,K(\frac{j}{4^n})],$ and
obviously $i+1<j\le i+8,$ so that there exist at most six
bounding triangles between these two.
Note that adjacent bounding triangles meet each other at
only two possible angles (that depend on $r$ but not on $n$),
and taking it as geometrically evident that non-adjacent
bounding triangles are disjoint, we conclude that there is
a constant $c>0$ (that depends on $r$ but not on $n$)
such that the distance between two non-adjacent
bounding triangles of rank $n,$ having at most $6$
bounding triangles of rank $n$ between them, does not exceed
$c\app r^n.$
In particular, $\nr{K(x)-K(y)}2\ge c\app r^{n}.$
Combining this with the inequalities above, we conclude that
$m(y-x)^\al\le\nr{K(y)-K(x)}2\le M(y-x)^\al,$
where $m=\frac c{8^\al}$ and $M=\frac2r$
(and $\al=-\log_4r$). 
\epF{Lemma}

Coming back to the theorem, let $\al=p/q,$ 
and let $K_\al$ be as in the lemma. 
Let $x=\ang{x_0,x_1,x_2,...}\in\oin.$ 
Then $K_\al(x_i)=\ang{x'_i,x''_i}\in[0,1]^2.$ 
We put $\vt(x)=\ang{x'_0,x''_0,x'_1,x''_1,x'_2,x''_2,...}.$ 
Prove that $\vt$ reduces $\bel p\res\oin$ to $\bel q$.

Let $x=\sis{x_i}{i\in\dN}$ and $y=\sis{y_i}{i\in\dN}$ 
belong to $\oin;$ we have to prove that $x\bel p y$ iff 
$\vt(x)\bel q\vt(y).$ 
To simplify the picture note the following: 
\dm
2^{-1/2}\nr w 2\le \tmax\ans{w',w''}\le \nr w q\le \nr w 1\le
2\nr w2
\dm 
for any $w=\ang{w',w''}\in\dR^2.$ 
The task takes the following form:
\dm
\sum_i(x_i-y_i)^p<\iy
\leqv 
\sum_i {\nr{K_\al(x_i)-K_\al(y_i)}2}^q<\iy\,.
\dm
Furthermore, by the choice of $K_\al,$ this converts to 
\dm
\sum_i(x_i-y_i)^p<\iy
\leqv 
\sum_i (x_i-y_i)^{\al q}<\iy\,,
\dm
which holds because $\al q=p.$ 
\epF{Theorem~\ref{hd-t}}

\punk{$\bel\iy$: maximal $\Ks$}
\label{liy}

Recall that $\Ks$ denotes the class of all \dd\fsg compact
sets in Polish spaces.
Easy computations show that this class contains, among
others, the \eqr s ${\Ei}\zd{\Ey}\zd{\bel p}\zd1\le p\le\iy,$
considered as sets of pairs in corresponding Polish spaces.
Note that if $\rE$ a $\Ks$ equivalence on a Polish space $\dX$
then $\dX$ is $\Ks$ as well 
since projections of compact sets are compact.
Thus $\Ks$ \er s on Polish spaces is one and the same as 
$\fs02$ \er s on $\Ks$ Polish spaces.

\bte
\label{eiy}
Any\/ $\Ks$ \eqr\ on a Polish space, in particular,
${\Ei}\zd{\Ey}\zd{\bel p},$
is Borel reducible to\/ ${\bel\iy}.$~\snos
{The result for $\bel p$ is due to Su Gao \cite{sudis}.
He defines
$d_p(x,s)=(\sum_{k=0}^{\lh s-1} |x(k)-s(k)|^p)^{\frac1p}$
for any $x\in\rtn$ and $s\in\dQ\lom$
(a finite sequence of rationals). 
Easily the \dd{\bel p}distance
$(\sum_{k=0}^{\iy} |x(k)-y(k)|^p)^{\frac1p}$
between any pair of $x,y\in\rtn$ is finite iff 
there is a constant $C$ such that
$|d_p(x,s)-d_p(y,s)|<C$ for all $s\in\dQ\lom.$
This yields a reduction required.}
\ete
\bpf[{{\rm from Rosendal~\cite{roz}}}]
Let $\dA$ be the set of all \dd\sq increasing sequences
$A=\sis{A_n}{n\in\dN}$ of subsets of $\dN$ --- a closed
subset of the Polish space $\cP(\dN)^\dN.$ 
Define an \er\ $\rH$ on $\dA$ by
\dm
\sis{A_n}{}\rH\sis{B_n}{}
\quad\text{iff}\quad
\sus N\:\kaz m\:
({A_m\sq B_{N+m}}\land {B_m\sq A_{N+m}}).
\dm

{\it Claim 1\/}:
$\rH\reb{\bel\iy}.$
This is easy.
Given a sequence $A=\sis{A_n}{n\in\dN},$ define
$\vt(A)\in \dN^{\dN\ti\dN}$
by $\vt(A)(n,k)$ to be the least $j\le k$ such that $n\in A_j,$
or $\vt(A)(n,k)=k$ whenever $n\nin A_k.$
Then $\sis{A_n}{}\rH\sis{B_n}{}$ iff there is $N$ such that 
$|\vt(A)(n,k)-\vt(B)(n,k)|\le N$ for all $n\zi k$.

{\it Claim 2\/}:
any $\Ks$ equivalence $\rE$ on a Polish space $\dX$ is Borel
reducible to $\rH.$ 
As a $\Ks$ set, $\rE$ has the form $\rE=\bigcup_nE_n,$ where
each $E_n$ is a compact subset of $\dX\ti\dX$
(not necessarily an \er)
and $E_n\sq E_{n+1}.$
We can \noo\ assume that each $E_n$ is reflexive and
symmetric on its domain $D_n=\dom{E_n}=\ran{E_n}$
(a compact set), in particular,
$x\in D_n\imp{\ang{x,x}\in E_n}.$
\vyk{
Furthermore, both $\dX$ and
$D=\rav\dX=\ens{\ang{x,x}}{x\in\dX}$ are $\Ks;$
let $D=\bigcup_nD_n,$ where $D_n$ are compact 
and $D_n\sq D_{n+1}.$
}%
Define $P_0=E_0$ and 
\dm
P_{n+1}=P_n\cup E_{n+1}\cup {P_n^{(2)}},
\;\text{ where }\,
P_n^{(2)}=\ens{\ang{x,y}}{\sus z\:
({\ang{x,z}\in P_n}\land{\ang{z,y}\in P_n})},
\dm
by induction.
Thus all $P_n$ are still compact subsets of $\dX\ti\dX,$
moreover, of $\rE$ since $\rE$ is an \eqr, and
$E_n\sq P_n\sq P_{n+1},$ therefore $\rE=\bigcup_nP_n.$

Let $\ens{U_k}{k\in\dN}$ be a basis for the topology of
$\dX.$
Put, for any $x\in\dX,$
$\vt_n(x)=\ens{k}{U_k\cap R_n(x)\ne\pu},$
where $R_n(x)=\ens{y}{\ang{x,y}\in R_n}.$
Then obviously $\vt_n(x)\sq\vt_{n+1}(x),$
and hence $\vt(x)=\sis{\vt_n(x)}{n\in\dN}\in\dA.$
Then $\vt$ reduces $\rE$ to $\rH.$

Indeed if $x\rE y$
then $\ang{y,x}\in P_n$ for some $n,$ and
for all $m$ and $z\in\dX$ we have 
${\ang{x,z}\in R_m}\imp{\ang{y,z}\in R_{1+\tmax\ans{m,n}}}.$
In other words, $R_m(x)\sq R_{1+\tmax\ans{m,n}}(y)$
and hence
$\vt_m(x)\sq \vt_{1+\tmax\ans{m,n}}(y)$
hold for all $m.$
Similarly, for some $n'$ we have
$\vt_m(y)\sq \vt_{1+\tmax\ans{m,n'}}(y)\zd \kaz m.$
Thus ${\vt(x)}\rH{\vt(y)}.$

Conversely, suppose that ${\vt(x)}\rH{\vt(y)},$ thus,
for some $N,$ we have 
$R_m(x)\sq R_{N+m}(y)$ and $R_m(y)\sq R_{N+m}(x)$
for all $m$ and $y.$
Taking $m$ big enough for $P_m$ to contain
$\ang{x,x},$ we obtain $x\in R_{N+m}(y),$ so that
immediately $x\rE y$.
\epf

\newpage

\parf{Smooth \er s and the first dichotomy}
\label{parf:1st}

This Section is mainly related to the node $\cont=\rav{\dn}$ 
in the diagram on page \pageref{p-p}. 
After a few rather simple results on smooth \er s which 
admit a Borel transversal, we show that countable, and 
\imar{where continual?}%
sometimes even continual unions of smooth \er s are smooth.
In the end, we prove the 1st dichotomy theorem.

\punk{Smooth and below}
\label{<smooth}
  
An important subspecies of smooth \er s consists of 
those having a Borel {\it transversal\/}: 
a set with 
\index{transversal}%
exactly one element in every equivalence class. 

\ble
\label{transv}
\ben
\renu
\itsep
\itla{transv1}
Any Borel \er\ that has a Borel transversal is smooth$;$

\itla{transv2}
any Borel finite (with finite classes) \er\ admits a
Borel transversal$;$

\itla{transv3}
any Borel countable smooth \er\ admits a Borel transversal;

\itla{transv4}
any Borel \er\ $\rE$ on a Polish space\/ $\dX,$ such that 
every\/ \dde class is closed and the saturation\/ 
$\ek\cO\rE$ of every open set\/ $\cO\sq\dX$ is Borel, 
admits a Borel transversal, hence, is smooth$.$~\snos
{Srivastava~\cite{sri} proved the result for 
\er s with $\Gd$ classes,  
which is the best possible as $\Eo$ is a Borel \er, whose 
classes are $\Fs$ and saturations of open sets are even open, 
but without any Borel transversal. 
See also \cite[18.20~iv)]{dst}.}

\itla{transv5}
$\Eo$ is not smooth.
\een
\ele
\bpf
\ref{transv1}
Let $T$ be a Borel transversal for $\rE.$   
The map $\vt(x)=$
``the only element of $T$ \dde equivalent to $x$''
reduces $\rE$ to $\rav T.$~\snos
{To see that a smooth \er\ does not necessarily have a 
Borel transversal  
take a closed set $P\sq\dnn\ti\dnn$ with $\dom P=\dnn,$ 
not uniformizable by a Borel set, and let  
$\ang{x,y}\rE\ang{x',y'}$ iff both $\ang{x,y}$ and 
$\ang{x',y'}$ belong to $P$ and $x=x'$.}

\ref{transv2}
Consider the set of the \dd<least elements of \dde classes, 
where $<$ is a fixed Borel linear order on the 
domain of $\rE$.

\ref{transv3}
Use \Cuni.

\ref{transv4}
Since any uncountable Polish space is a continuous image of 
$\dnn,$ 
we can assume that $\rE$ is a \er\  on $\dnn.$ 
Then, for any $x\in\dnn,$ $\ek x\rE$ is a closed subset of 
$\dnn,$ naturally identified with a tree, say, $T_x\sq\dN\lom.$ 
Let $\vt(x)$ denote the leftmost branch of $T_x.$ 
Then $x\rE\vt(x)$ and ${x\rE y} \imp{\vt(x)=\vt(y)},$ so that 
it remains to show that $Z=\ens{\vt(x)}{x\in\dnn}$ is Borel. 
Note that 
\dm
{z\in Z}\,\leqv\,
\kaz m\;\kaz s,\,t\in \dN^m\;
\skl
{s<_{\text{\tt lex}} t}\land {z\in\cO_t}\limp 
{\ek z\rE\cap\cO_t=\pu}
\skp,
\dm
where $<_{\text{\tt lex}}$ is the lexicographical order on 
$\dN^m$ 
and $\cO_s=\ens{x\in\dnn}{s\su x}.$ 
However $\ek x\rE\cap\cO_t=\pu$ iff $x\nin \ek{\cO_t}\rE$ 
and $\ek{\cO_t}\rE$ is Borel for any $t$.

\ref{transv5}
Otherwise $\Eo$ has a Borel transversal $T$ by \ref{transv3},
which is a contradiction, see Example~\ref{ex:ac}(ii).
\epf

\punk{Assembling smooth \eqr s}
\label{cud}

If $\rE$ and $\rF$ are smooth \er s on disjoint sets, resp.,  
$X$ and $Y,$ then easily $\rE\cup\rF$ is a smooth \er\ on 
$X\cup Y.$ 
The question becomes less clear when we have a Borel \er\ 
$\rE$ on a Polish space $X\cup Y$ such that both $\rE\res X$
and $\rE\res Y$ are smooth but the sets $X,Y$ not necessarily
\dd\rE invariant in $X\cup Y$ if even disjoint; is $\rE$ smooth? 
We answer this in the positive, even in the case of countable 
unions.

\bte
\label{cud1}
Let\/ $\rE$ be a Borel \er\ on a Borel set\/ 
$X=\bigcup_kX_k,$ with all\/ $X_k$ also Borel. 
Suppose that each\/ $\rE\res X_k$ is smooth. 
Then\/ $\rE$ is smooth.
\ete
\bpf\footnote
{\ The shortest proof is to note that otherwise 
$\Eo\reb\rE$ by the 2-nd dichotomy, easily leading to 
contradiction by a Baire category argument.
Yet we prefer to give a direct proof.
Note that even in the case when the sets $X_k$ are
pairwise disjoint, most obvious ideas like
``to define $\vt(x)$ take the least $k$ such that
$X_k$ intersects $\eke x$ and apply $\vt_k$'' do not
work.}
First consider the case of a union $X=Y\cup Z$ 
of just two Borel sets, so that a Borel \er\ $\rE$ is 
smooth on both $Y$ and $Z.$ 
We can assume that $Y\cap Z=\pu.$ 
Let the smoothness be witnessed by Borel reductions
$f:Y\to Q$ and $g:Z\to R,$ with $Q,\,R$
being disjoint Borel sets.
The set 
\dm
F=\ens{\ang{q,r}}{\sus y\in Y\;\sus z\in Z\;
\skl{f(y)=q}\land{g(z)=r}\land{y\rE z}\skp}
\sq Q\ti R
\dm
is a partial $\fs11$ map $Q\to R.$ 
Let $G:Q\to R$ be any Borel map with $F\sq G,$ 
and $H:R\to Q$ be any Borel map with $F\obr\sq H.$ 
Then $\Phi=G\cap H\obr$ is a 
$1-1$ Borel partial map $P\to Q$ with $F\sq\Phi.$ 
Now the $\fp11$ set 
\dm
P=\ens{\ang{q,r}\in \Phi}{\kaz y\in Y\;\kaz z\in Z\;
\skl f(y)=q\land g(z)=r\imp y\rE z\skp}\,,
\dm
satisfies $F\sq P\sq \Phi,$ hence, there is 
a Borel function $\Psi$ with $F\sq\Psi\sq P.$  
The sets $A=\dom\Psi$ and $B=\ran\Psi$ are Borel subsets 
of resp.\ $Q,R,$  
and it follows from the construction that  
$\Psi\cap(\dom F\ti \ran F)=F.$ 
Finally, put
\dm
D=\Psi\cup\ens{\ang{q,q}}{q\in Q\dif A}
\cup\ens{\ang{r,r}}{r\in R\dif B}\,,
\dm
then, for any $y\in Y$ there is unique 
$h(y)=\ang{q,r}\in D$ with $q=f(y),$ 
correspondingly, for any $z\in Z$ there is unique 
$h(z)=\ang{q,r}\in D$ with $r=g(z),$ and if $y\rE z$ 
then $h(y)=h(z)=\ang{f(y),g(z)},$ hence, $h$ 
witnesses that $\rE$ is smooth.

As for the general case, we can now assume 
that $X_k\sq X_{k+1}$ for all $k.$ 
Then there are disjoint Borel sets $W_k$ and Borel 
maps $f_k:X_k\to W_k$ which witness that $\rE\res X_k$ 
are smooth \er s. 
Let $R_k=\ran f_k$ (a $\fs11$ set) and 
\dm
F_k=\ens{\ang{a,b}\in R_k\ti R_{k+1}}
{\sus x\in X_k\:(f_k(x)=a\land f_{k+1}(x)=b)}\,,
\dm
this is a $\fs11$ set and a $1-1$ map $R_k\to R_{k+1}.$ 
For each $k$ there is a  
Borel $1-1$ map $G_k$ with $F_k\sq G_k.$  
Let $A_k=\dom G_k$ and $\ran G_k=B_k:$ these are  
Borel sets with $R_k\sq A_k.$ 
We can assume that $B_k\sq A_{k+1}.$ 
\imar{nuzhno li eto assume}
(Otherwise $G_k$ can be reduced in a certain iterative 
manner to achieve this property.) 
Then, for any $k$ and $b\in A_k$ there is the least 
$n=n(b)\le k$ such that the application 
\dm
h(b)=G\obr_n(G\obr_{n+1}(G\obr_{n+2}
(...G\obr_{k-1}(b)...)))
\dm
is possible, for instance, $n(b)=k$ and $h(b)=b$ 
whenever $b\in A_k\dif B_{k-1}.$ 
Then, $h(f_k(x))=h(f_{k+1}(x))$ holds for any 
$x\in X_k$ because $F_k\sq G_k,$ so that the map 
$g(x)=h(f_k(x))$ for $x\in X_k\dif X_{k-1}$ 
witnesses the smoothness of $\rE$. 
\epf

\punk{The 1st dichotomy theorem.}
\label{>smooth}
  
The following result is known as the 1st dichotomy theorem.

\bte[{{\rm Silver~\cite{sil}}}]
\label{1dih}
Any\/ $\fp11$ \er\/ $\rE$ on\/ $\dnn$ either has at 
most countably many equivalence classes or admits a 
perfect set of pairwise\/ \dde inequivalent reals, in 
other words, either\/ $\rE\reb\rav\dN$ or\/  
$\rav\dn\reb\rE$.
\ete
\bpf\footnote
{\ We present a forcing proof of Miller~\cite{am}, with 
some simplifications. 
See \cite{mw} for another proof, based on the Gandy -- 
Harrington topology.
In fact both proofs involve essentially the same
combinatorics.}
As usual, we can suppose that $\rE$ is a lightface
$\ip11$ relation.\vom

{\sl Case 1\/}: any $x\in\dnn$ belongs to a $\id11$ 
\dde equivalent set $X$ 
(\ie, all elements of $X$ are \dde equivalent to each other, 
in other words, 
the saturation $\eke X$ is an equivalence class). 
Then $\rE$  has at most countably many equivalence classes.\vom

{\sl Case 2\/}: otherwise. 
Then the set $H$ of all $x,$ which do \poq{not} belong to a 
$\id11$ pairwise \dde equivalent set 
({\it the domain of nontriviality\/}), is non-empty. 

\bct
\label{dih1:1}
$H$ is\/ $\is11.$ 
Any\/ $\is11$ set\/ $\pu\ne X\sq H$ 
is not pairwise\/ \dde equivalent.
\ect
\bpf
$x\in H$ iff for any $e\in\dN:$ 
\poq{if} $e$ codes a $\id11$ set, say, $W_e\sq\dnn$ and 
$x\in W_e$ 
\poq{then} $W_e$ is not \dde equivalent. 
The ``if'' part of this characterization is $\ip11$ while 
the ``then'' part is $\is11,$ by \Penu\ (see \prf{eff}).

If $X\ne\pu$ is a pairwise \dde equivalent $\is11$ set then 
$B=\bigcap_{x\in X}\eke x$ is a $\ip11$ \dde equivalence class  
and $X\sq B.$  
By \Sepa, there is a $\id11$ set $C$ with $X\sq C\sq B.$ 
Then, if $X\sq H$ then $C\sq H$ is a $\id11$ pairwise 
\dde equivalent set, a 
contradiction to the definition of $H$.
\epF{Claim}

Let us fix a countable transitive model 
\index{modelm@model $\mm$}%
$\mm$ of a big enought fragment of $\ZFC,$ and an elementary 
submodel of the universe \vrt\ all analytic formulas~\footnote
{\label{suf}\ 
For instance, $\mm$ models $\ZC$ and, in addition, Replacement 
for $\Sg_{100}$ \mem formulas and the first one million of 
instances of Replacement overall. 
Being an elementary submodel is useful to guarantee that 
relations like the inclusion orders of $\px$ and $\pg$ are 
absolute for $\mm$.}. 
Consider
$\dP=\ens{X\sq\dnn}{X\,\text{ is non-empty and }\,\is11}$ 
\index{zzP@$\dP$}%
as a forcing to extend $\mm$ 
(smaller sets are stronger conditions), 
the {\it Gandy -- Harrington forcing\/}. 
\index{Gandy -- Harrington!forcing}%
We have $\dP\nin$ and $\not\sq\mm,$ of course, but clearly 
$\dP$ can be adequately coded in $\mm,$ say, via a 
universal $\is11$ set.

\bcot[{{\rm from Theorem~\ref{s11pol}}}]
\label{d12}
If\/ $G\sq\dP$ is a\/ \dd\dP generic, over\/ $\mm,$ set, 
then\/ $\bigcap G$ contains a single real, denoted\/ 
$x_G$.\qeD
\ecot

Reals of the form $x_G,$ $G$ as in the Corollary, 
are called \dd\dP{\it generic\/} (over $\mm$).
Let $\dox$ be the name for $x_G.$ 
Then any\/ $A\in\dP$ forces that\/ $\dox\in A$.

Let $\dP^2$ consist of all ``rectangles'' $X\ti Y,$ with 
\index{zzpu2@$\dP^2$}%
$X,\,Y\in\dP.$  
It follows from the above by the product forcing lemmas 
that any \dd{\dP^2}generic, over $\mm,$ set $G\sq\dP^2$ 
produces a pair of reals 
(a \dd{\dP^2}{\it generic pair\/}),  
say, $x^G\ul$ and $x^G\ur,$ so that  
$\ang{x^G\ul,x^G\ur}\in W$ for any $W\in G.$  
Let $\doxl$ and $\doxr$ be their names. 

\blt
\label{d14}
$H\ti H$ \dd{\dP^2}forces\/ $\doxl \nE\doxr$.
\elt
\bpf
Otherwise a ``condition'' $X\ti Y\in\dP^2$ 
with $X\cup Y\sq H$ \dd{\dP^2}forces $\doxl\rE\doxr,$ 
so that any \dd{\dP^2}generic pair 
$\ang{x,y}\in X\ti Y$ satisfies $x\rE y.$ 
By the product forcing lemmas for any pair of \dd\dP generic 
$x',\,x''\in X$ there is $y\in Y$ such that both 
$\ang{x,y}$ and $\ang{x',y}$ are \dd{\dP^2}generic pairs, 
hence, we have 


\bit
\item[$(\ast)$] 
If\/  $x',\,x''\in X$ are\/ \dd\dP generic over\/ $\mm$ then\/ 
$x'\rE x''$.
\eit


\vyk{
The $\is11$ set $P=X^2\dif\rE$ is non-empty.
(Otherwise $X\ne\pu$ is a \dde equivalent $\is11$ 
subset of $H.$ 
Then $B=\bigcap_{x\in X}\eke x$ is a $\ip11$ set and 
an \dde equivalence class, and $X\sq B\sq H.$ 
By \Sepa, there is a $\id11$ set $C$ with $X\sq C\sq B.$ 
Then $C\sq H$ is a $\id11$ \dde equivalent set, a 
contradiction to the definition of $H.$)
}

The set $\dP_2$ of all non-empty $\is11$ 
\index{zzpd2@$\dP_2$}%
subsets of $\dnn\ti\dnn$ is just a copy of $\dP$ 
(not of $\dP^2$!) as a forcing, in particular, 
if $G\sq\dP_2$ is \dd{\dP_2}generic over $\mm$ 
then there is a unique pair of reals 
(\dd{\dP_2}{\it generic pair\/}) 
$\ang{x^G\ul,x^G\ur}$ which belongs to every $W$ in $G,$ 
and in this case, 
{\it both\/ $x^G\ul$ and $x^G\ur$ are\/ \dd\dP generic\/},  
because if $G\sq\dP_2$ is \dd{\dP_2}generic then 
the sets $G'$ and $G''$ of all projections of sets $W\in G$ 
to resp.\ 1st and 2nd co-ordinate, are easily 
\dd\dP generic.
Now let $G\sq\dP_2$ be a \dd{\dP_2}generic set, over $\mm,$  
containing the $\is11$ set $P=X^2\dif\rE.$  
(Note that $P\ne\pu$ by Lemma~\ref{dih1:1}.)
Then $\ang{x^G\ul,x^G\ur}\in P,$ hence, $x^G\ul\nE x^G\ur,$ 
however, as we observed, both $x^G\ul$ and $x^G\ur$ are 
\dd\dP generic elements of $X$ 
(because $P\sq X\ti X$), 
which contradicts $(\ast)$. 
\epF{Lemma~\ref{d14}}

Fix 
enumerations 
$\sis{\cD(n)}{n\in\dN}$ and $\sis{\cD^2(n)}{n\in\dN}$ 
of all dense subsets of resp.\ $\dP$ and $\dP^2$ which 
are coded in $\mm.$ 
Then there is a system 
$\sis{X_u}{u\in 2\lom}$ of sets $X_u,$ satisfying 
\ben
\tenu{(\roman{enumi})}
\itla{d1i}
$X_u\in\dP,$ moreover, $X_\La\sq H$ and $X_u\in\cD(n)$ 
whenever $u\in 2^n;$

\itla{d1d}
$X_{u\we i}\sq X_u$ for all $u\in2\lom$ and $i=0,1$;

\itla{d1t}
if $u\ne v\in 2^n$ then $X_u\ti X_v\in\cD^2(n)$.
\een
It follows from \ref{d1i} that, for any $a\in \dn,$ the 
set $\ens{X_{a\res m}}{m\in\dN}$ is \dd\dP generic over $\mm,$ 
hence, $\bigcap_mX_{a\res m}$ is a singleton, say, $x_a,$ by 
Corollary~\ref{d12}. 
Moreover the map $a\mapsto x_a$ is continuous as diameters 
of $X_u$ converge to $0$ uniformly with $\lh u\to 0,$ by 
\ref{d1i}. 
In addition, by \ref{d1t} and Lemma~\ref{d14},  
$x_a\nE x_b$ whenever $a\ne b,$ in particular, 
$x_a\ne x_b,$ hence,  
we have a perfect \dde inequivalent set 
$Y=\ens{x_a}{a\in\dn}$.\vom

\epF{Theorem~\ref{1dih}}

\parf{Hyperfinite and countable \er s}
\label{hyct}

This Section is mainly devoted to the node $\Eo$ in the 
diagram on page \pageref{p-p}. 
Together with the 2nd dichotomy theorem, we present some 
other properties of $\Eo,$ the ideal $\ifi,$ and 
hyperfinite (Borel) \eqr s.
This class of \eqr s is a very interesting object of 
study even aside of pure descriptive set theory. 
Papers \cite{djk,jkl} give a comprehensive account of  
most basic results, with further references.

After a rather simple theorem which shows that $\ifi$ is the 
least ideal in the sense of $\orbpp,\:\orb,\;\reb,$ we prove 
the ``Glimm--Effros'', or second, dichotomy which asserts that 
$\Eo=\rE_{\ifi}$ is the \dd\reb least among all non-smooth 
Borel \er s.
Finally, we present a characterization, in terms of the 
existence of transversals, of those Borel sets $X$ for which
\imar{where is this?}
$\Eo\res X$ is smooth.

\punk{$\ifi$ is the least !}

The proof of the following useful result is based on
a short argument involved in many other results.
A somewhat more pedestrian version of the argument
was used in
several proofs in Section~\ref{BAN}.

\bte
\label{jnmt}
\ben
\renu
\itsep
\itla{jnmt1}
{\rm\cite{jn,mat75,tal}}
If\/ $\cI$ is a\/ {\rm(nontrivial)} ideal on\/ $\dN,$ 
with the Baire property in the topology of $\pn,$ 
then\/ $\ifi\orbpp$ and\/ $\orb\cI\,;$ 

\itla{jnmt2}
however\/ $\rav\dn\rebs\Eo$ strictly, thus\/ $\rav\dn$
is not\/ \dd\eqb equivalent to an \eqr\ of the form\/
$\rei\,;$

\itla{jnmt3}
if\/ $\cI\orbp\cJ$ are Borel ideals, and there is an 
infinite set\/ $Z\sq\dom\cI$ such that\/ 
$\cI\res Z=\pwf Z,$ 
then\/ $\cI\orb\cJ$.
\een
\ete
\bpf
\ref
{jnmt1}
First of all $\cI$ must be meager in $\pn.$ 
(Otherwise $\cI$ would be comeager somewhere, 
easily leading to contradiction.) 
Thus, all $X\sq\dN$ ``generic'' 
(over a certain countable family of dense open subsets
of $\pn$) 
do not belong to $\cI.$ 
Now it suffices to define non-empty finite sets $w_i\sq\dN$ 
with $\tmax w_i<\tmin w_{i+1}$ such that any union of 
infinitely many of them is ``generic''. 
Clearly the following observation yields the result: 
if $D$ is an open dense subset of $\pn$ and $n\in\dN$ then 
there is $m>n$ and a set $u\sq[n,m]$ with $m,\,n\in u$ 
such that any $x\in\pn$ satisfying $x\cap[n,m]=u$ belongs 
to $D$.

Thus we have $\ifi\orbpp\cI.$ 
To derive $\ifi\orb\cI$ cover each $w_k$ by a 
finite set $u_k$ such that $\bigcup_{k\in \dN}u_k=\dN$ 
and still $u_k\cap u_l=\pu$ for $k\ne l$.

\ref{jnmt2}
That $\rav\dn\reb\Eo$ is witnessed by any perfect set 
$X\sq\dn$ which is a
{\it partial\/}
\index{transversal!partial}%
transversal for $\Eo$ 
(\ie, any $x\ne y$ in $X$ are \dd{\Eo}inequivalent). 
On the other hand, $\rav\dn$ is smooth but $\Eo$
is non-smooth by Lemma~\ref{transv}\ref{transv5}.

\ref{jnmt3}
Assume \noo\ that $\cI,\cJ$ are ideals over $\dN.$
Let pairwise disjoint finite sets $w_k\sq\dN$ 
witness $\cI\orbp\cJ.$ 
Put $Z'=\dN\dif Z\yt X=\bigcup_{k\in Z}w_k,$ 
and $Y=\bigcup_{k\in Z'}w_k.$ 
The reduction via $\sis{w_k}{}$ reduces $\pwf Z$ to 
$\cJ\res X$ and $\cI\res Z'$ to $\cJ\res Y.$ 
Keeping the latter, replace the former by a 
\dd\orb like reduction of $\pwf z$ to $\cJ\res Y',$ where 
$Y'=\dN\dif Y,$ which exists by Theorem~\ref{jnmt}.
\epf

Despite of Theorem~\ref{jnmt}, $\Eo=\rE_{\ifi}$ is 
\poq{not} the \dd\reb least among Borel \er s. 

Thus, $\rav\dn$ is not a 
\er\ generated by a Borel ideal, even modulo $\eqb$.

\punk{Countable \eqr s}
\label{cer}


This class of \eqr s, essentially bigger than hyperfinite
(modulo $\reb$),
is a subject of ongoing intence study. 
Yet we can only present here the following important
theorem and a few more results below, leaving
\cite{jkl,gab,kemi} as basic references in this domain.

\bte
[{{\rm \cite[Thm 1]{femo1}, \cite[1.8]{djk}}}]
\label{femo}
Any Borel countable \er\/ $\rE$ on a Polish space\/ $\dX\,{:}$
\ben
\tenu{{\rm(\roman{enumi})}}
\itsep
\itla{femo1}
is induced by a Polish action of a countable group\/
$\dG$ on\/ $\dX\,;$   

\itla{femo2}
satisfies\/ $\rE\reb \Ey=\rE(F_2,2),$ where\/
$F_2$ is the free group with two generators and\/
$\rE(F_2,2)$ is the \er\ induced by the 
shift action of\/ $F_2$ on\/ $2^{F_2}.$
\een
\ete
\bpf
\ref{femo1}
We \noo\ assume that $\dX=\dn.$
According to \Cenu\
(in a relativized version, if necessary, see
Remark~\ref{relrel}),
there is a sequence of Borel maps $f_n:\dn\to\dn$
such that $\eke a=\ens{f_n(a)}{n\in\dN}$ for each
$a\in\dn.$
Put $\Ga'_n=\ens{\ang{a,f_n(a)}}{a\in\dN}$
(the graph of $f_n$)
and $\Ga_n=\Ga'_n\dif\bigcup_{k<n}\Ga'_k.$
The sets $P_{nk}=\Ga_n\cap{\Ga_k}\obr$ form a partition 
of (the graph of) $\rE$ onto countably many Borel
injective sets.
Further define $\Da=\ens{\ang{a,a}}{a\in\dn}$ and let
$\sis{D_m}{m\in\dN}$ be an enumeration of all
non-empty sets of the form $P_{nk}\dif\Da.$
Intersecting the sets $D_m$ with the rectangles of the
form
\dm
R_s=\ens{\ang{a,b}\in\dn\ti\dn}{s\we0\su a\land s\we1\su b}
\quad\text{and}\quad
{R_s}\obr,
\dm
we reduce the general case to the case when
$\dom{D_m}\cap\ran{D_m}=\pu\zd\kaz m.$

Now, for any $m$ define $h_m(a)=b$ whenever either
$\ang{a,b}\in D_m$ or $\ang{a,b}\in {D_m}\obr,$ or
$a=b\nin \dom{D_m}\cup\ran{D_m}.$
Clearly $h_m$ is a Borel bijection $\dn\onto\dn.$
Thus $\sis{h_m}{m\in\dN}$ is a family of Borel
automorphisms of $\dn$ such that
$\eke a=\ens{h_m(a)}{m\in\dN}.$
It does not take much effort to expand this system to a Borel
action of $F_\om,$ the free group with $\alo$ generators,
on $\dn,$ whose induced \eqr\ is $\rE$.

\ref{femo2}
First of all, by \ref{femo1}, $\rE\reb\rR,$ where $\rR$
is induced by a Borel action $\app$ of $F_\om$ on $\dn.$
The map $\vt(a)=\sis{g\obr\app a}{g\in F_\om}\yt a\in\dn,$
is a Borel reduction of $\rR$ to $\rE(F_\om,\dn).$
If now $F_\om$ is a subgroup of a countable group $H$ then
$\rE(F_\om,\dn)\reb\rE(H,\dn)$ by means of the map sending
any $\sis{a_g}{g\in F_\om}$ to $\sis{b_h}{h\in H},$ where
$b_g=a_g$ for $g\in F_\om$ and $b_h$ equal to any fixed
$b'\in\dn$ for $h\in H\dif F_\om.$
As $F_\om$ admits a homomorphism into $F_2$~\snos
{Why ?.}
we conclude that $\rE\reb\rE(F_2,\dn).$

It remains to transform $\rE(F_2,\dn)$ to $\rE(F_2,2).$
The inequality $\rE(F_2,\dn)\reb\rE(F_2,2^{\dZ\dif\ans0})$
is clear.
Further $\rE(F_2,2^{\dZ\dif\ans0})\reb\rE(F_2\ti\dZ,3),$
by means of the map sending any
$\sis{a_g}{g\in F_2}\linebreak[0]\;\,(a_g\in2^{\dZ\dif\ans0})$ to
$\sis{b_{gj}}{g\in F_2,\;j\in\dZ},$ where
$b_{gj}= a_g(j)$ for $j\ne0$ and $b_{g0}= 2.$
Further, for any $G,$ $\rE(G,3)\reb\rE(G\ti \dZ_2,2)$
by means of the map sending any
$\sis{a_g}{g\in G}\;\,(a_g=0,1,2)$ to
$\sis{b_{gi}}{g\in G,\;i\in\dZ_2},$ where
\dm
b_{gi}=
\left\{
\bay{rclcl}
0,&\text{if}& a_g=0 &\text{\ubf or}&
a_g=1\,\text{ and }\, i=0,\\[0.7\dxii]

1,&\text{if}& a_g=2 &\text{\ubf or}&
a_g=1\,\text{ and }\, i=1.
\eay
\right.
\dm
Thus $\rE(F_2,\dn)\reb\rE(F_2\ti\dZ\ti\dZ_2,2).$
However, $F_2\ti\dZ\ti\dZ_2$ admits a homomorphism into
$F_\om,$ and then into $F_2$ (see above), so that
$\rE(F_2,\dn)\reb\rE(F_2,2),$ as required.
\epf

\punk{Hyperfinite \eqr s}
\label{hfer}

All Borel finite \er s are smooth (see \prf{<smooth}), 
accordingly, all hyperfinite \er s are hypersmooth. 
On the other hand, any finite or hyperfinite \eqr\ is 
countable, of course. 
It follows from the next theorem that, conversely, 
every {\it hypersmooth\/} countable \er\ is hyperfinite.
(But there exist countable non-hypersmooth \er s, 
for instance, $\Ey,$ which are not hyperfinite.) 

The theorem also shows that $\Eo$ is a universal 
hyperfinite \er.   
(To see that $\Eo$ is hyperfinite, let $x\rF_n y$ 
iff $x\sd y\sq[0,n)$ for $x,\,y\sq\dN.$) 

\bte
[{{\rm Theorems 5.1 and, partially, 7.1 in \cite{djk} 
and 12.1(ii) in \cite{jkl}}}]
\label{thf}
The following are equivalent for a Borel \er\/ $\rE$ on\/ 
a Polish space\/ $\dX:$
\ben
\tenu{{\rm(\roman{enumi})}}
\itla{thf1}\msur
$\rE\reb\Eo$ and\/ $\rE$ is countable$;$

\itla{thf2}\msur
$\rE$ is hyperfinite$;$

\itla{thf3}\msur
$\rE$ is  hypersmooth and countable$;$

\itla{thf4}
there is a Borel set\/ $X\sq\pnd$ such that\/ $\Ei\res X$ 
is a countable \er\ and\/ $\rE$ is isomorphic, 
via a Borel bijection of\/ $\dX$ onto\/ $X,$  
to\/ $\Ei\res X\;;$
         
\itla{thf5}\msur
$\rE$ is induced by a Borel action of\/ $\dZ,$ the additive 
group of the integers.

\itla{thf6} 
there exists a pair of Borel \er s\/ $\rF,\:\rR$ of type\/ 
$2$ such that\/ $\rE=\rF\lor\rR.$~\snos 
{An \eqr\ $\rF$ is {\it of type $n$\/} if any \ddf class 
contains at most $n$ elements. 
$\rF\lor\rR$ denotes the least \er\ which includes 
$\rF\cup\rR$.}
\een
\ete
\bpf
$\ref{thf2}\imp\ref{thf3}$ and $\ref{thf1}\imp\ref{thf3}$ 
are rather easy. 

$\ref{thf3}\imp\ref{thf4}.$ 
Let $\rE=\bigcup_n\rF_n$ be a countable and hypersmooth \er\
on a space $\dX,$ all $\rF_n$ being smooth (and countable), 
and $\rF_{n}\sq\rF_{n+1},\:\kaz n.$ 
We may assume that $\dX=\pn$ and $\rF_0=\rav{\pn}.$ 
Let $T_n\sq\dX$ be a Borel transversal for $\rF_n$ 
(recall Lemma~\ref{transv}\ref{transv3}).
Now let $\vt_n(x)$ be the only element of $T_n$ with 
$v\rF_n{\vt_n(x)}.$ 
Then $x\mapsto\sis{\vt_n(x)}{n\in\dN}$ is a $1-1$ Borel map 
$\dX\to \pnd$ and ${x\rE y}\eqv{\vt(x)\Ei\vt(y)}.$ 
Take $X$ to be the image of $\dX$.

$\ref{thf4}\imp\ref{thf5}.$
Let $X$ be as indicated. 
For any \dd\dN sequence $x$ and $n\in\dN,$ let 
$x\qc n=x\res{(n,\iy)}.$ 
It follows from (the relativized version of) \Cpro\ and \Cenu\ 
that for any $n$ the set
$X\qc n=\ens{x\qc n}{x\in X}$ is Borel 
and there is a countable family of Borel functions 
$g^n_i:X\qc n\to X\zT i\in\dN,$ such that the set  
$X_\xi=\ens{x\in X}{x\qc n=\xi}$ is equal to 
$\ens{g^n_i(\xi)}{i\in\dN}$ for any $\xi\in X\qc n,$ 
hence,  
$\ens{g^n_i(\xi)(n)}{i\in\dN}=\ens{x(n)}{x\in X_\xi}.$ 
        
For any $x\in\pnd$ let $\vpi(x)=\sis{\vpi_n(x)}{n\in\dN},$ 
where $\vpi_n(x)$ is the least number $i$ such that 
$x(n)=f^n_i(x)(n);$ thus, $\vpi(x)\in\dnn.$
Let $\mu(x)$ be the sequence
\dm
\vpi_0(x),\vpi'_0(x),
\vpi_1(x)+1,\vpi'_1(x)+1,\dots,
\vpi_n(x)+n,\vpi'_n(x)+n,\dots,
\dm
where $\vpi'_n(x)=\tmax_{k\le n}\vpi_k(x).$ 
Easily if $x\ne y\in X$ satisfy $x \Ei y,$ \ie, 
$x\qc n=y\qc n$ for some $n,$ then 
$\vpi(x)\qc n= \vpi(y)\qc n$ but 
$\vpi(x)\ne \vpi(y)\zT \mu(x)\ne\mu(y),$ 
and $\mu(x)\qc m= \mu(y)\qc m$ for some $m\ge n.$ 
  
Let $\alex$ be the anti-lexicographical partial order on 
$\dnn,$ \ie, $a\alex b$ iff there is $n$ such that 
$a\qc n=b\qc n$ and $a(n)<b(n).$  
For $x,\,y\in X$ define $x<_0 y$ iff $\mu(x)\alex \mu(y).$ 
It follows from the above that $<_0$ linearly orders every 
\dd\Ei class $\ek x\Ei\cap X$ of $x\in X.$ 
Moreover, it follows from the definition of $\mu(x)$ that 
any \dd\alex interval between some $\mu(x)\alex \mu(y)$ 
contains only finitely many elements of the form $\mu(z).$ 
(For $\vpi$ this would not be true.)
We conclude that any class $\ek x\Ei\cap X\zT x\in X,$ is 
linearly ordered by $<_0$ similarly to a subset of $\dZ,$ 
the integers. 
That $<_0$ can be converted to a required Borel action of 
$\dZ$ on $X$ is rather easy 
(however the \dd\Ei classes in $X$ ordered similarly 
to $\dN,$ the inverse of $\dN,$ or finite, should be 
treated separately). 

$\ref{thf5}\imp\ref{thf2}.$
Assume \noo\ that $\dX=\dn.$
An increasing sequence of \er s $\rF_n$ whose union 
is $\rE$ is defined separately on each \dd\rE class 
$C;$ they ``integrate'' 
into Borel \er s $\rF_n$ defined on the whole of $\dn$ 
because the action allows 
to replace quantifiers over a \dde class $C$ by quantifiers 
over $\dZ.$

Let $C$ be any \dde class of $x\in X.$ 
Note that if an element $x_C\in C$ can be chosen in some 
Borel-definable way then we can define $x\rF_n y$ iff there 
exist integers $j,\,k\in\dZ$ with $|j|\le n\zT |k|\le n,$ 
and $x=j\app x_C\zT y=k\app x_C.$  
This applies, for instance, when $C$ is finite, thus, we 
can assume that $C$ is infinite. 
Let $\lex$ be the lexicographical ordering of $\dn,$ 
and $\act$ be the partial order induced by the action, 
\ie, $x\act y$ iff $y=j\app x\zT j>0.$ 
By the same reason we can assume that neither of  
$a=\tinf_{\lex}C$ and $b=\tsup_{\lex}C$ belongs to $C.$ 
Let $C_n$ be the set of all $x\in C$ with 
$x\res n\ne a\res n$ and $x\res n\ne b\res n.$ 
Define $x\rF_n y$ iff $x,\,y$ belong to one and the same 
\dd\lex interval in $C$ lying entirely within $C_n,$ or 
just $x=y.$ 
In our assumptions, any $\rF_n$ has finite classes, and for 
any two $x,\,y\in C$ there is $n$ with $x\rF_n y$. 

$\ref{thf5}\imp\ref{thf1}.$
This is more complicated. 
A preliminary step is to show that $\rE\reb\rE(\dZ,\dn),$ 
where $\rE(\dZ,\dn)$ is the orbit equivalence induced by the 
shift action of $\dZ$ on $(\dn)^\dZ$: 
$k\app\sis{x_j}{j\in\dZ}=\sis{x_{j-k}}{j\in\dZ}$ for 
$k\in\dZ.$ 
Assuming \noo\ that $\rE$ is a \er\ on $\dn,$ we obtain a 
Borel reduction of $\rE$ to $\rE(\dZ,\dn)$ by 
$\vt(x)=\sis{j\app x}{j\in\dZ},$ where $\app$ is a Borel  
action of $\dZ$ on $\dn$ which induces $\rE.$ 
Then Theorem~7.1 in \cite{djk} proves that 
$\rE(\dZ,\dn)\reb\Eo$.

$\ref{thf6}\imp\ref{thf5}.$
\imar{wrong}
Let $\rE=\rF\lor\rR,$ where $\rF,\:\rR$ are of type $2.$ 
For any $x\in\dX$ (the domain of $\rE$), if $\ekf x$ contains 
another element $y\ne x$ then call $y$ the {\it left\/}, resp., 
{\it right neighbour\/} of $x$ if $y<x,$ resp., $y>x,$ where 
$<$ is a fixed Borel linear ordering of $\dX.$  
If the class $\ek x\rR$ also contains another element, say, 
$z,$ call it the neighbour of $x$ of the opposite side \vrt\ 
$y.$ 
The neighbour relation linearly orders any \dd\rE class 
similarly to a subset of $\dZ,$ which easily leads to \ref{thf5}. 

$\ref{thf5}\imp\ref{thf6}.$
The authors of \cite{jkl} present a short proof which 
refers to several difficult theorems on hyperfinite \er s. 
Here we give an elementary proof.

Let $\rE$ be induced by a Borel action of $\dZ.$ 
We are going to define $\rF$ and $\rR$ on any \dde class 
$C=\eke x.$
If we can choose an element $x_C\in C$ in some uniform 
Borel-definable way then a rather easy construction is 
possible, which we leave to the reader.
This applies, for instance, when $C$ is finite, hence, 
let us assume that $C$ is infinite. 
Let $\act$ be the linear order on $C,$ induced by the action 
of $\dZ;$ it is similar to $\dZ.$ 
Let $\lex$ be the lexicographical ordering of $\dn=\dom\rE.$ 

Our goal is to define $\rF$ on $C$ so that every \ddf class 
contains exactly two (distinct) elements. 
The ensuing definition of $\rR$ is then rather simple. 
(First, order pairs $\ans{x,y}$ of elements of $C$ in 
accordance with the \dd\act lexicographical ordering of 
pairs $\ang{\tmax_{\act}\ans{x,y},\tmin_{\act}\ans{x,y}},$ 
this is still similar to $\dZ.$  
Now, if $\ans{x,y}$ and $\ans{x',y'}$ are two \ddf classes, 
the latter being the next to the former in the sense just 
defined, and $x\act y\zT x'\act y',$ then define $y\rR x'.$)

Suppose that $W\sq C.$ 
An element $z\in W$ iz {\it lmin\/} (locally minimal)
in $W$ if it is \dd\lex smaller than both of its 
\dd\act neighbours in $W.$
Put $W_{\text{lmin}}=\ens{z\in W}{z\,\text{ is lmin in }\,W}.$
If $C_{\text{lmin}}$ is \poq{not} unbounded in $C$ in both 
directions 
then an appropriate choice of $x_C\in C$ is possible. 
(Take the \dd\act least or \dd\act largest point in 
$C_{\text{lmin}},$ 
or if $C_{\text{lmin}}=\pu,$ so that, for instance, 
$\act$ and $\lex$ coincide on $C,$ we can choose something 
like a \dd\lex middest element of $C.$)
Thus, we can assume that $C_{\text{lmin}}$ is unbounded in 
$C$ in both directions. 

Let a {\it lmin-interval\/} be any \dd\act semi-interval 
$[x,x')$ between two consecutive elements $x\act x'$ 
of $C_{\text{lmin}}.$ 
Let $[x,x')=\ans{x_0,x_1,...,x_{m-1}}$ be the enumeration 
in the \dd\act increasing order ($x_0=x$). 
Define $x_{2k}\rF x_{2k+1}$ whenever $2k+1<m.$ 
If $m$ is odd then $x_{m-1}$ remains unmatched. 
Let $C^1$ be the set of all unmatched elements. 
Now, the nontrivial case is when $C^1$ is unbounded in $C$ in 
both directions. 
We define $C^1_{\text{lmin}},$ as above, and repeat the same 
construction, extending $\rF$ to a part of $C^1,$ with, perhaps, 
a remainder $C^2\sq C^1$ where $\rF$ remains indefined. 
{\sl Et cetera\/}.

Thus, we define a decreasing sequence 
$C=C^0\supseteq C^1 \supseteq C^2 \supseteq \dots$ of subsets 
of $C,$ and the equivalence relation $\rF$ on each difference 
$C^n\dif C^{n+1}$ whose classes contain exactly two points each, 
and the nontrivial case is when every $C^n$ is \dd\act unbounded 
in $C$ in both directions. 
(Otherwise there is an appropriate choice of $x_C\in C.$) 
If $C^\iy=\bigcap_nC^n=\pu$ then $\rF$ is defined on $C$ and we 
are done. 
If $C^\iy=\ans{x}$ is a singleton then $x_C=x$ chooses an 
element in $C.$ 
Finally, $C^\iy$ cannot contain two different elements as 
otherwise one of $C^n$ would contain two \dd\act neighbours 
$x\act y$ which survive in $C^{n+1},$ which is easily impossible.
\epf

\punk{Non-hyperfinite countable equivalence relations}
\label{no-hypf}

It follows from Theorem~\ref{thf}\ref{thf1},\ref{thf2}
that hyperfinite \eqr s form an initial segment,
in the sense of $\reb,$ among all countable \eqr s.
Let us show that not all countable \eqr s
are hyperfinite.

\bte
The \eqr\ $\Ey$ is not hyperfinite.
\ete
\bpf
A clean elementary proof is given in \cite{slst}.
\vyk{
The argument is based the fact that all hyperfinite \eqr s
satisfy a certain property, amenability, which
fails for $\Ey.$
We are not going to discuss amenability as a general
notion
(see \cite{k-amen,jkl} and references there on this topic),
but rather take just as much of it as is necessary for
the proof.

A group $\dG$ is {\it amenable\/}
\index{amenable}%
iff there is a finitely additive probability measure
(\fap\ in brief)
$\mu:\pws\dG\to{[0,1]},$
left-invariant in the sense that
$\mu(X)=\mu(aX)$ for any $X\sq\dG$ and $a\in G,$
where $aX=\ens{ax}{x\in X}.$

For instance the additive group $\dZ$ of all integers
is amenable.
Indeed fix a nontrivial
(\ie, not containing singletons)
ultrafilter $U$ over $\pws\dN.$
Take any set $A\sq\dZ.$
To define $\mu(A),$ put
$\mu_n(A)=\frac{\card{A\cap{[-n,n]}}}{2n+1}$ for any $n.$
For any $p\in{[0,1]},$ one of the two complementary sets
\dm
N_p^+(A)=\ens{n}{p>\mu_n(A)}\,,\quad 
N_p^-(A)=\ens{n}{p\le\mu_n(A)}\,,
\dm
belongs to $U.$
The sets $P^+(A)=\ens{p}{N_p^+(A)\in U},$ 
$P^-(A)=\ens{p}{N_p^-(A)\in U}$ form a Dedekind cut in
$[0,1].$
Let $\mu(A)$ be either the least real in $P^+(A)$ or the
largest real in $P^-(A).$~\snos
{In other words, $\mu(A)$ is the limit over $U$ of the
sequence of numbers $\mu_n=\mu_n(A)$.}
The left-invariance of $\mu$ follows from the fact
that for any ${A\sq\dZ}\yt{z\in\dZ}\zi{\ve>0}$ we have
$|\mu_n(A)-\mu_n(z+A)|<\ve$ for all sufficiently large
natural $n.$

On the other hand, $F_2,$
the free group with two generators, say, $a\zi b,$ 
is not amenable.
Indeed, suppose that $\nu:\pws{F_2}\to{[0,1]}$
is a left-invariant \fap.
Let, for $x\in\ans{a,b,a\obr,b\obr},$
$W_x$ denote the set of all words in $F_2$ beginning
with $x.$
As obviously $F_2=W_a\cup aW_{a\obr},$ we have
$\nu(W_a)+\nu(W_{a\obr})=1.$
Similarly, $\nu(W_b)+\nu(W_{b\obr})=1.$
On the other hand,
$F_2=\ans1\cup W_a\cup aW_{a\obr}\cup W_b\cup aW_{b\obr},$
thus $\nu(W_a)+\nu(W_{a\obr})+\nu(W_b)+\nu(W_{b\obr})\le1,$
contradiction.

Now suppose towards the contrary that the \eqr\
$\rE=\rE(F_2,2)$ induced by the 
shift action of $F_2$ on $2^{F_2},$ is hyperfinite.
Then by Theorem~\ref{thf} $\rE$ is induced by a
Borel action $\app$ of $\dZ$ on $2^{F_2}.$
Let $\mu:\pws\dZ\to{[0,1]}$ witness the amenability of
$\dZ,$ as above.
A similar measure can be defined on each
\dd{\rE}class
$C=\eke f=\ens{g\in2^{F_2}}{f\rE g},$ where $f\in2^{F_2},$
as follows:
$\vpi_C(X)=\mu(\ens{z\in\dZ}{z\app f\in X})$ 
for any $X\sq C.$
It follows from the invariance of $\mu$ that $\vpi_C$ in
does not depend on the choice of $f\in C.$
}%
\epf

\punk{Assembling hyperfinite equivalence relations}
\label{assem:hypf}

The following theorem shows that, similarly to the case of
smooths \er s (Thm~\ref{cud1}),
hyperfinite ones possess a certain 
form of countable additivity.

\bte
\label{cud2}
Let\/ $\rE$ be a Borel \er\ on a Borel set\/ 
$X=\bigcup_kX_k,$ with all\/ $X_k$ also Borel. 
Suppose that\/ $\rE\res X_k\reb\Eo$ for each\/ $k.$ 
Then\/ $\rE\reb\Eo$.
\ete
\bpf
We consider only the case when $X_k\sq X_{k+1}$ for all $k$ 
\imar{not the best prf}%
(the result will be used below only for this particular case), 
the general case needs to consider separately the two--sets 
case, as in Theorem~\ref{cud1}, which we leave to the reader. 

There are disjoint Borel sets $B_k\sq \pn$ and Borel 
maps $f_k:X_k\to B_k$ which witness that 
$\rE\res X_k\reb\Eo.$ 
We shall assume that the sets $B_k$ are \deo incompatible 
in the sense that if $k\ne n$ then $a\Eo b$ does not hold 
for any $a\in B_k$ and $b\in B_n.$ 
Let $R_k=\ran f_k$ (a $\fs11$ subset of $B_k$). 
Then  
\dm
F_k=\ens{\ang{a,b}\in R_k\ti R_{k+1}}
{\sus x\in X_k\:
\skl{f_k(x)=a}\land {f_{k+1}(x)=b}\skp}\,,
\dm
is a $\fs11$ set, 
$1-1$ modulo $\Eo$ in the sense that if $\ang{a,b}$ and 
$\ang{a',b'}$ belong to $F_k$ then $a\Eo a'\eqv b\Eo b'.$  
As ``to be $1-1$ modulo $\Eo$'' is a $\ip11$ property 
in the codes (of $\is11$ subsets of $\pn^2$), there 
is, by \Refl, a $\id11$ set $F'_k$ with 
$F_k\sq F'_k\sq B_k\ti B_{k+1}$ and still $1-1$ modulo 
$\Eo.$ 
The following $\id11$ set 
\dm
G_k=\ens{\ang{a',b'}}
{\sus\ang{a,b}\in F'_k\;({a\Eo a'}\land{b\Eo b'})}
\dm
is still $1-1$ modulo $\Eo,$ hence, both ``vertical'' and 
``horisontal'' cross-sections of $G_k$ are countable, thus, 
$A_k=\dom G_k$ and $B_k=\ran G_k$ are \deo invariant 
Borel sets (and $R_k=\dom F_k\sq A_k$), 
and there are Borel maps $h_k:B_k\to A_k$ such that 
$\ang{h_k(b),b}\in G_k$ whenever $b\in B_k.$ 
It follows still from the ``$1-1$ modulo $\Eo$'' property 
that if $b\in B_{k}$ and $b'\Eo b$ then $b'\in B_{k}$ 
and $h_k(b)\Eo h_k(b')$. 

We can assume that $B_{k+1}\sq A_k$ for all $k.$ 
Then, for any $k$ and $b\in A_k$ there is the least 
$n=n(b)\le k$ such that the application 
\dm
h(b)=h_n(h_{n+1}(h_{n+2}(...h_{k-1}(b)...)))
\dm
is possible, for instance, $n(b)=k$ and $h(b)=b$ 
whenever $b\in A_k\dif B_{k-1}.$ 
As in the proof of Theorem~\ref{cud1},  
the map $g(x)=h(f_k(x))$ for $x\in X_k\dif X_{k-1}$ 
witnesses $\rE\reb\Eo$.
\epf

\parf{The 2nd dichotomy}
\label{parf2d}

The following result is known as 2nd, or 
``Glimm--Effros'', dichotomy. 

\bte[{{\rm Harrington, Kechris, Louveau \cite{hkl}}}]
\label{2dih}
If\/ $\rE$ is a Borel \er\ then either\/ 
$\rE$ is smooth or\/ $\Eo\emn\rE$.
\ete

\punk{The Gandy -- Harrington closure}
\label{ghc}

Beginning the proof of Theorem~\ref{2dih}
(it will be completed in \prf{d2splitk}),
we suppose, as usual, that $\rE$ is a lightface $\id11$ 
\er\ on $\dnn.$ 
Consider an auxiliary \er\ $x\Ec y$ iff $x,\,y\in\dnn$ 
belong to the same \dde invariant $\id11$ sets. 
(A set $X$ is \dde{\it invariant\/} iff $X=\eke X.$) 
\index{set!einvariant@\dde invariant}%
Easily $\rE\sq\Ec.$  
To see that $\Ec$ is the closure of $\rE$ in the Gandy -- 
Harrington topology, prove


\blt
\label{d20}
If\/ $\rF$ is a\/ $\is11$ \er\ on\/ $\dnn,$ and 
$X,Y\sq\dnn$ are disjoint\/ \ddf in\-var\-iant\/ $\is11$ 
sets, then there is an\/ \ddf invariant\/ $\id11$ 
set\/ $X'$ separating\/ $X$ from\/~$Y.$ 
\elt
\bpf
By \Sepa, for any $\is11$ set $A$ with $A\cap Y=\pu$ there is 
a $\id11$ set $A'$ with $A\sq A'$ and $A'\cap Y=\pu$ --- note 
that then $\ekf{A'}\cap Y=\pu$ because $Y$ is \ddf invariant. 
It follows that that there is a sequence 
$X=A_0\sq A'_0\sq A_1\sq A'_1\sq ...,$  
where 
$A'_i$ are $\id11$ sets, accordingly, $A_{i+1}=\ekf{A'_i}$ 
are $\is11$ sets, and $A_i\cap Y=\pu.$ 
Then $X'=\bigcup_nA_n=\bigcup_nA'_n$ and 
is an \ddf invariant Borel set which separates $X$ from $Y.$  
To make $X'$ $\id11$ we have to maintain the choice of sets 
$A_n$ effectively. 

Let $U\sq \dN\ti\dnn$ be a ``good'' universal $\is11$ set
(see \prf{eff}).
Then there is a recursive $h:\dN\to\dN$ such that
$\ekf{U_n}=U_{h(n)}$ for each $n.$
Moreover, applying Lemma~\ref{effred}
(to the complement of $U$ as a ``good'' universal $\ip11$ set,
and with a code for $Y$ fixed),
we obtain a pair of recursive functions $f\zi g:\dN\to\dN$
such that for any $n,$ if $U_n\cap Y=\pu$ then
$U_{f(n)}\zi U_{g(n)}$ are complementary sets
(hence, either of them is $\id11$)
containing, resp., $U_n$ and $Y.$ 
A suitable iteration of $h$ and $f,g$ allows us
to define a sequence $X=A_0\sq A'_0\sq A_1\sq A'_1\sq ...$
as above effectively enough for the union of those sets
to be $\id11$. 
\epF{Lemma}

\blt
\label{d21}
$\Ec$ is a \/ $\is11$ relation.
\elt
\bpf
Let $C\sq\dN$ and $W,\,W'\sq\dN\ti\dnn$ be as in \Penu\ 
(\prf{eff}). 
The formula $\inva(e)$ saying that $e\in C$ and 
$W_e=W'_e$ is \dde invariant, \ie,
\dm
e\in C\lland 
\kaz a,\,b\:
\skl {a\in W_e\land b\nin W'_e}\imp {a\nE b}\skp
\dm
is obviously $\ip11,$ however $x\Ec y$ iff
\dm
\kaz e\; 
\skl
\inva(e)\limp 
(x\in W_e\imp y\in W'_e) \land (y\in W_e\imp x\in W'_e)
\skp
\eqno\square\:(\hbox{\sl Lemma})
\dm
\ePf

Let us return to the proof of the theorem. 
We have two cases.\vom

{\bfsl Case 1\/}:
$\rE=\Ec,$ \ie, $\rE$ is Gandy -- Harrington closed. 
\vyk{
Let $\sis{D_n}{n\in\dN}$ be an arbitrary enumeration of all 
\dde invariant $\id11$ sets $D\sq\dnn.$ 
Then the map 
$\vt(x)=\ens{n}{x\in D_n}$ is a Borel reduction of $\rE$ 
to $\rav{\dn},$ as required.

Although it is of no use for Theorem~\ref{2dih}, let us 
prove that the smoothness in Case~1 can be achieved in 
effective manner.
}

\blt
\label{d2eff}
If\/ $\rE=\Ec$ then there is a\/ $\id11$ reduction of\/ 
$\rE$ to\/ $\rav\dn$.
\elt
\bpf
Let $C\sq\dN$ and $W,\,W'\sq\dN\ti\dnn$ be as in the \Penu\ 
of \prf{eff}. 
By \Kres\ there is a $\id11$ function $\vpi:X^2\to C$ such 
that $W_{\vpi(x,y)}=W'_{\vpi(x,y)}$ is a \dde invariant 
$\id11$ set containing $x$ but not $y$ whenever $x,\,y\in X$ 
are \dde inequivalent. 
Then $R=\ran \vpi$ is a $\is11$ subset of $C,$ hence, by \Sepa, 
there is a $\id11$ set $N$ with $R\sq N\sq C.$ 
The map $\vt(x)=\ens{n\in N}{x\in D_n}$ is a $\id11$ reduction 
of $\rE$ to $\rav{\dn}$. 
\epF{Lemma and Case 1}

{\bfsl Case 2\/}: $\rE\sneq\Ec.$ 
Then the $\is11$ set 
$H=\ens{x}{\eke x\sneq\ekec x}$ 
(the union of all \ddec classes containing more than one 
\dde class) is non-empty.

\blt
\label{e=oe}
If\/ $X\sq H$ is a\/ $\is11$ set then\/ $\rE\sneq\Ec$ on\/ $X$.
\elt
\bpf
Suppose that ${\rE\res X}={\Ec\res X}.$ 
Then $\rE=\Ec$ on $Y=\eke X$ as well. 
(If $y,\,y'\in Y$ then there are $x,\,x'\in X$ such 
that $x\rE y$ and $x'\rE y',$ so that if $y\Ec y'$ then 
$x\Ec x'$ by transitivity, hence, $x\rE x',$ and $y\rE y'$ 
again by transitivity.) 
It follows that $\rE=\Ec$ on an even bigger set, 
$Z=\ekco X.$ 
(Otherwise the $\is11$ set 
$Y'=Z\dif Y=\ens{z}{\sus x\in X\:(x\Ec y\land{x\nE y})}$ 
is non-empty and \dde invariant, together with $Y,$ hence
by Lemma~\ref{d20} 
there is a \dde invariant $\id11$ set $B$ with $Y\sq B$ 
and $Y'\cap B=\pu$, which implies that 
no point in $Y$ is \dd{\Ec}equivalent to a point in $Y',$ 
contradiction.) 
Then by definition $Z\cap H=\pu$.
\epF{Lemma}

\blt
\label{d27}
If\/ $A,B\sq H$ are non-empty\/ $\is11$ sets with\/ 
$A\rE B$ then there exist non-empty \poq{disjoint}\/ $\is11$ 
sets\/ $A'\sq A$ and\/ $B'\sq B$ still satisfying\/  
$A'\rE B'$.
\elt
\bpf
We assert that there are points $a\in A$ and $b\in B$ 
with $a\ne b$ and $a\rE b.$ 

(Otherwise $\rE$ is the equality on $X=A\cup B.$ 
Prove that then $\rE=\Ec$ on $X,$   
a contradiction to Lemma~\ref{e=oe}. 
Take any $x\ne y$ in $X.$ 
Let $U$ be a clopen set containing $x$ but not $y.$ 
Then $A=\eke {U\cap X}$ and $C=\eke {X\dif U}$ are two 
disjoint \dde invariant $\is11$ sets containing resp.\ $x,\,y.$ 
Then $x\Ec y$ fails by Lemma~\ref{d20}.) 

Thus let $a,\,b$ be as indicated. 
Let $U$ be a clopen set containing $a$ but not $b.$ 
Put $A'=A\cap U\cap\eke{\doP U}$ and 
$B'=B\cap{\doP U}\cap\eke U$.
\epF{Lemma}

\punk{Restricted product forcing}
\label{d2forc}

Recall that forcing notions $\dP$ and $\dP_2$ were
introduced in \prf{>smooth}.
In continuation of the proof of Theorem~\ref{2dih} (Case 2), 
let $\dpe$ be the collection of all sets of the form 
\index{zzp2e@$\dpe$}%
$X\ti Y,$ where $X,\,Y\sq\dnn$ are non-empty $\is11$ sets and 
$X\rE Y$ (which means here that $\eke X=\eke Y$).  
\index{zzXEY@$X\rE Y$}%
Easily $\dP_2\sq\dpe\sq\dP^2$. 
The forcing~\snos
{Over a countable 
model $\mm$ chosen in accordance with the requirements in 
Footnote~\ref{suf}.}
$\dpe$ is not really a product, yet if $X\ti Z\in\dpe$ and 
$\pu\ne X'\sq X$ is $\is11$ then $Z'=Z\cap\eke{X'}$ is $\is11$ 
and $X'\ti Z'\in\dpe.$ 
It follows that any \dd\dpe generic set $G\sq\dpe$ produces 
a pair of \dd\dP generic sets 
$G\ul=\ens{\dom P}{P\in G}$ and 
$G\ur=\ens{\ran P}{P\in G},$ hence, produces a pair 
of \dd\dP generic reals $x^G\ul$ and $x^G\ur,$ whose names 
will be $\doxl$ and $\doxr$.

\blt
\label{d23}
In the sense of the forcing\/ $\dpe,$ 
any\/ $P=X\ti Z\in\dpe$ forces\/ $\ang{\doxl,\doxr}\in P$ 
and forces\/ $\doxl\Ec\doxr,$ but\/  
$H\ti H$ forces\/ $\doxl\nE\doxr$.
\elt
\bpf
To see that $\doxl\Ec\doxr$ is forced suppose otherwise. 
Then, by the definition of $\Ec,$ 
there is a condition $P=X\ti Z\in\dpe$ and an  
\dde invariant $\id11$ set $B$
such that $P$ forces  $\doxl\in B$ but $\doxr\nin B.$ 
Then easily $X\sq B$ but $Z\cap B=\pu,$ a contradiction with 
$\eke X=\eke Z$.

To see that $H\ti H$ forces $\doxl\nE\doxr$ suppose towards 
the contrary that some $P=X\ti Z\in\dpe$ with $X\cup Z\sq H$ 
forces $\doxl\rE\doxr,$ thus,

\ben
\tenu{(\arabic{enumi})}
\itla{d23-}
$x\rE z$ holds 
for every \dd\dpe generic pair $\ang{x,z}\in P$. 
\een

\bct
\label{d24}
If\/ $x,y\in X$ are\/ \dd\dP generic over\/ $\mm,$ 
and\/ $x\Ec y,$ then\/ $x\rE y$.
\ect
\bpf
We assert that 
\ben
\tenu{(\arabic{enumi})}
\addtocounter{enumi}1
\itla{d23a}
${x\in A}\eqv {y\in A}$ holds
for each \dde invariant $\is11$ set $A$. 
\een
Indeed, if, say, $x\in A$ but $y\nin A$ then by the 
genericity of $y$ there is a $\is11$ set $C$ with 
$y\in C$ and $A\cap C=\pu.$ 
As $A$ is \dde invariant, Lemma~\ref{d20} yields an 
\dde invariant $\id11$ set $B$ such that $C\sq B$ but 
$A\cap B=\pu.$ 
Then $x\nin B$ but $y\in B,$ a contradiction to 
$x\Ec y$.

Let 
$\sis{\cD_n}{n\in\dN}$ be an enumeration 
of all dense subsets of $\dpe$ which are coded in $\mm.$ 
We define 
two sequences $P_0\qs P_1\qs...$ and 
$Q_0\qs Q_1\qs...$ of conditions $P_n=X_n\ti Z_n$ and 
$Q_n=Y_n\ti Z_n$ in $\dpe,$ so that $P_0=Q_0=P,$ 
$x\in X_n$ and $y\in Y_n$ for any $n,$ and finally 
$P_n,\,Q_n\in\cD_{n-1}$ for $n\ge1.$ 
If this is done then we have a real $z$ 
(the only element of $\bigcap_nZ_n$) 
such that both $\ang{x,z}$ and $\ang{y,z}$ are \dd\dpe generic, 
hence, $x\rE z$ and $y\rE z$ by \ref{d23-}, hence, 
$x\rE y$. 

Suppose that $P_n$ and $Q_n$ have been defined.  
As $x$ is generic, there is (we leave details for the reader) 
a condition $P'=A\ti C\in \cD_{n}$ and $\sq P_n$ such that  
$x\in A.$ 
Let $B=Y_n\cap\eke A:$ then $y\in B$ by \ref{d23a}, and 
easily $\eke B=\eke C=\eke A$ 
(as $\eke{X_n}=\eke{Z_n}=\eke{Y_n}$), thus, $B\ti C\in\dpe,$ 
so there is a condition $Q'=V\ti W\in\cD_n$ and 
$\sq B\ti C\sq Q_n$ such that $y\in V.$ 
Put $Y_{n+1}=V,\msur$ $Z_{n+1}=W,$ and $X_{n+1}=A\cap\eke W$.
\epF{Claim}

It follows that $\rE=\Ec$ on\/ $X.$  
(Otherwise $S=\ens{\ang{x,y}\in X^2}{x\Ec y\land x\nE y}$ 
is a non-empty $\is11$ set,  
and any \dd{\dP_2}generic pair $\ang{x,y}\in S$ implies a 
contradiction to Claim~\ref{d24}.
Recall that $\dP_2=$ all non-empty $\is11$ subsets of 
$\dnnp2.$) 
But this implies $X\cap H=\pu$ by Lemma~\ref{e=oe}, 
contradiction.
\vyk{
To accomplish the proof of Lemma~\ref{d23}, note that 
\dm
X\sq C'=\ens{c\in\dnn}{\kaz x\in X\:(x\Ec c\imp x\rE c)}\,,
\quad\text{and}\quad C'\:\hbox{ is }\,\ip11\,,
\dm 
hence, there is a $\id11$ set $B'$ with $X\sq B'\sq C'.$ 
Further, 
\dm
X\sq C=\ens{c\in B'}{\kaz b\in B'\:(b\Ec c\imp b\rE c)}\,,
\dm 
hence, there is a $\id11$ set $B$ with $X\sq B\sq C,$ 
and by definition ${\rE\res B}={\Ec\res B}.$ 
It easily follows that $\rE,\,\Ec$ coincide also on $W=\eke B,$ 
hence, $W$ cannot intersect $H.$ 
However $\pu\ne X\sq B\cap H,$ contradiction.
}
\epF{Lemma~\ref{d23}}

\vyk{
\blt
\label{d27}
$H\ti H$ forces, in the sense of\/ $\dP_2,$ that\/ 
$\doxl\ne\doxr$.
\elt
\bpf
Otherwise there is a condition in $\dP_2,$ of the form 
$X\ti X,$ where $\pu\ne X\sq H$ is a non-empty set, such 
that ${\rE}\res X$ is the equality. 
Let us prove that then $\rE=\rF$ on $X:$ then we will have 
a contradiction as in the end of the proof of Lemma~\ref{d23}. 
Take any $x\ne y$ in $X.$ and
To show $x\nF y$ let $U$ be a clopen set containing $x$ but 
not $y.$ 
Then $A=\eke {U\cap X}$ and $C=\eke {X\dif U}$ are two 
disjoint \dde invariant $\is11$ sets containing resp.\ $x,\,y.$ 
Then by Lemma~\ref{d20} there is a \dde invariant $\is11$ set 
$B$ containing $x$ but not $y,$ as required.
\epF{Lemma}
}

\punk{Splitting system}
\label{d2split}

Let us fix 
enumerations 
$\sis{\cD(n)}{n\in\dN},\msur$ $\sis{\cD_2(n)}{n\in\dN},\msur$  
$\sis{\cD^2(n)}{n\in\dN}$ 
of all dense subsets of resp.\ $\dP,\:\dP_2,\:\dpe,$ which 
belong to $\mm\,;$ we assume that 
$\cD(n+1)\sq\cD(n),\msur$ $\cD_2(n+1)\sq\cD_2(n),$ and 
$\cD^2(n+1)\sq\cD^2(n).$ 
If $u,\,v\in 2^m$ (binary sequences of length $m$) have the 
form $u=0^k\we0\we w$ and $v=0^k\we1\we w$ for some $k<m$ and 
$w\in 2^{m-k-1}$ then we call $\ang{u,v}$ a {\it crucial pair\/}. 
It can be proved, \eg, by induction on $m,$ that $2^m$ is 
a connected tree 
(\ie, a connected graph without cycles) of crucial pairs, with 
sequences beginning with $1$ as the endpoints of the graph. 
We define 
a system of sets $X_u$ ($u\in2\lom$) and 
$\rR_{uv}\,,\msur$ $\ang{u,v}$ being a crucial pair, so that 
the following conditions are satisfied:

\ben
\tenu{(\roman{enumi})}
\itla{d2i1}
$X_u\in\dP,$ moreover, $X_\La\sq H,$ and 
$X_u\in \cD(n)$ for any $u\in2^n$; 

\itla{d2i2}
$X_{u\we i}\sq X_u$ for all $u$ and $i$;

\itla{d2i3}
$\rR_{uv}\in\dP_2,$ moreover, $\rR_{uv}\in \cD_2(n)$
for any crucial pair  $\ang{u,v}$ in $2^n$; 

\itla{d2i4}
$\rR_{uv}\sq \rE$ and $X_u \rR_{uv} X_v$ 
for any crucial pair $\ang{u,v}$ in $2^n$; 

\itla{d2i5} 
$\rR_{u\we i\,,\,v\we i}\sq\rR_{uv}$;

\itla{d2i6}
if $u,v\in2^n$ and $u(n-1)\ne v(n-1)$ then 
$X_u\ti X_v\in \cD^2(n)$ and also $X_u\cap X_v=\pu$.
\een
Note that \ref{d2i4} implies that $X_u\rE X_v$ for any crucial 
pair $\ang{u,v},$ hence, also for any pair in 
$2^n$ because any $u,v\in2^n$ are connected by a 
unique chain of crucial pairs. 
It follows that $X_u\ti X_v\in\dpe$ for any pair of 
$u,v\in 2^n,$ for~any~$n$.

Assume that such a system has been defined.  
Then for any $a\in\dn$ the sequence 
$\sis{X_{a\res n}}{n\in\dN}$ is \dd\dP generic over $\mm,$ 
hence, $\bigcap_n X_{a\res n}=\ans{x_a},$ where $x_a$ is 
\dd\dP generic,  
and the map $a\mapsto x_a$ is continuous since diameters 
of $X_u$ converge to $0$ uniformly with $\lh u\to 0$ by 
\ref{d2i1}, and is $1-1$ by the last condition~of~\ref{d2i6}. 

Let $a,b\in\dn.$ 
If ${a\neo b}$ 
then, by \ref{d2i6}, $\ang{x_a,x_b}$ is a \dd\dpe generic 
pair, hence, $x_a\nE x_b$ by Lemma~\ref{d23}.
Now suppose that $a\Eo b,$ prove that then $x_a\rE x_b.$
We can suppose that $a=w\we 0\we c$ and $b=w\we0\we c,$ 
where $w\in2\lom$ and $c\in\dn$ 
(indeed if $a\Eo b$ then $a,\,b$ can be connected by a finite 
chain of such special pairs). 
Then $\ang{x_a,x_b}$ is \dd{\dP_2}generic, actually, the only 
member of the intersection 
$\bigcap_n \rR_{w\we0\we(c\res n)\,,\,w\we1\we(c\res n)}$ 
by \ref{d2i3} and \ref{d2i4}, in particular, $x_a\rE x_b$ 
because we have $R_{uv}\sq{\rE}$ for all $u,\,v$.

Thus we have 
a continuous $1-1$ reduction of $\Eo$ to $\rE.$ 
\vom

\qeDD{Case 2 in Theorem~\ref{2dih} modulo the construction}

\punk{Construction of a splitting system}
\label{d2splitk}

Let $X_\La$ be any member of $\cD(0)$ satisfying $X_\La\sq H.$ 
Now suppose that $X_s$ and $\rR_{st}$ have been defined for all 
$s\in2^n$ and all crucial pairs in $2^n,$ and extend the 
construction on $2^{n+1}.$ 
Temporarily, define $X_{s\we i}=X_s$ and 
$\rR_{s\we i\,,\,t\we i}=\rR_{st}:$ this leaves 
$\rR_{0^n\we0\,,\,0^n\we 1}$ still undefined, so we put 
$\rR_{0^n\we0\,,\,0^n\we 1}={\rE}\cap X_{0^n}\ti X_{0^n}.$ 
Note that the such defined system of sets $X_u$ and 
relations $\rR_{uv}$ at level $n+1$ 
satisfies all requirements of \ref{d2i1} -- \ref{d2i6} 
except for the requirement of membership in the dense sets 
-- say in this case that the system is ``coherent''. 
It remains to produce a still ``coherent'' system of smaller 
sets and relations which also satisfies the 
membership in the dense sets. 
This will be achieved in several steps.\vom

{\sl Step 1\/}: achieve that $X_u\in\cD(n+1)$ for any 
$u\in 2^{n+1}.$ 
Take any particular $u_0\in 2^{n+1}.$ 
There is, by the density, 
$X'_{u_0}\in\cD(n+1)$ and $\sq X_{u_0}.$ 
Suppose that $\ang{u_0,v}$ is a crucial pair. 
Put
$\rR'_{u_0,v}=\ens{\ang{x,y}\in\rR_{u_0,v}}{x\in X'_{u_0}}$ 
and $X'_v=\ran \rR'_{u_0,v}.$ 
This shows how the change spreads along the whole set 
$2^{n+1}$ viewed as the tree of crucial pairs. 
Finally we obtain a coherent system with the additional 
requirement that $X'_{u_0}\in\cD(n+1).$ 
Do this consecutively for all $u_0\in 2^{n+1}.$ 
The total result -- we re-denote it as still $X_u$ and $\rR_{uv}$ 
-- is a ``coherent'' system with $X_u\in\cD(n+1)$ for all $u.$ 
Note that still $X_{0^n\we 0}=X_{0^n\we 1}$ and 
\dm
\rR_{0^n\we0\,,\,0^n\we 1}=
{\rE}\cap (X_{0^n\we 0}\ti X_{0^n\we 1})\,.
\eqno(\ast)
\dm

{\sl Step 2\/}: achieve that 
$X_{s\we0}\ti X_{t\we1}\in\cD^2(n+1)$ for all $s,\,t\in 2^{n+1}.$ 
Consider a pair of $u_0=s_0\we0$ and $v_0=t_0\we 1$ in 
$2^{n+1}.$ 
By the density there is a set 
$X'_{u_0}\ti X'_{v_0}\in\cD^2(n+1)$ and $\sq X_{u_0}\ti X_{v_0}.$ 
By definition we have $X'_{u_0}\rE X'_{v_0},$ but, due to 
Lemma~\ref{d27} we can maintain that $X'_{u_0}\cap X'_{v_0}=\pu.$
The two ``shockwaves'', from the changes at $u_0$ and $v_0,$ as 
in Step 1, meet only at the pair $0^m\we0,\,0^m\we1,$ where the 
new sets satisfy $X'_{0^m\we0}\rE X'_{0^m\we1}$ just because 
\dde equivalence is everywhere kept and preserved though the 
changes. 
Now, in view of $(\ast),$ we can define 
$\rR'_{0^n\we0\,,\,0^n\we 1}=
{\rE}\cap (X'_{0^n\we 0}\ti X'_{0^n\we 1}),$ 
preserving $(\ast)$ as well. 
All pairs considered, we will be left with a coherent system 
of sets and relations, re-denoted as $X_u$ and $\rR_{uv},$ 
which satisfies the \dd{\cD(n+1)}requirements in \ref{d2i1} 
and \ref{d2i6}.\vom

{\sl Step 3\/}: achieve that $\rR_{uv}\in\cD_2(n+1)$ for any 
crucial pair at level $n+1,$ and also that 
$X'_{0^n\we 0}\cap X'_{0^n\we 1}=\pu.$ 
Consider any crucial pair $\ang{u_0,v_0}.$ 
If this is not $\ang{0^n\we 0,0^n\we 1p}$ then let 
$\rR'_{u_0v_0}\sq\rR_{u_0v_0}$ be any set in $\cD_2(n+1).$ 
If this is $u_0=0^n\we 0$ and $v_0=0^n\we 1$ then first we 
choose (Lemma~\ref{d27}) disjoint non-empty $\is11$ sets 
$U\sq X_{0^n\we 0}$ and $V\sq X_{0^n\we 1}$ still with 
$U\rE V,$ and only then a set 
$\rR'_{u_0v_0}\sq\rE\cap(U\ti V)$ which belongs to 
$\in\cD_2(n+1).$ 
In both cases, put $X'_{u_0}=\dom \rR'_{u_0v_0}$ and 
$X'_{v_0}=\ran \rR'_{u_0v_0}.$ 
It remains to spread the changes, along the chain of crucial 
pairs, to the left of $u_0$ and  to the right of $v_0,$ exactly 
as in Case 1. 
Executing such a reduction for all crucial pairs $\ang{u_0,v_0}$
at level $n+1$ one by one, we end up with a system of sets 
fully satisfying \ref{d2i1} -- \ref{d2i6}.\vtm

\qeDD{Theorem~\ref{2dih}}

\punk{A forcing notion associated with $\Eo$}
\label{zeo}

We here consider the forcing notion $\zfo{\Eo}{\rav\dn}$
(see \prf{f2}), 
that will be denoted by $\peo$ below.
\index{forcing!peo@$\peo$}%
\index{zzpeo@$\peo$}%
Thus by definition $\peo$ consists of all Borel sets
$X\sq\dn$ such that $\Eo\res X$ is non-smooth while 
the related ideal $\ieo=\zid{\Eo}{\rav\dn}$ consists
\index{ideal!ieo@$\ieo$}%
\index{zzieo@$\ieo$}%
of all Borel sets $X\sq\dn$ such that $\Eo\res X$ is
smooth.

\ble
\label{bay}
\ben
\tenu{{\rm(\roman{enumi})}}
\itla{bay1}
$\ieo$ is a\/ \dd\fsg additive ideal.
Let\/  $X\sq\dn$ be a Borel set. 

\itla{bay2}
$X$ belongs to\/ $\peo$ iff\/
${\Eo}\emn{\Eo\res X}$ (by a continuous injection).

\itla{bay3}
$X$ belongs to\/ $\ieo$ iff\/
${{\Eo}\res X}$ admits a Borel transversal. 
\een
\ele
\bpf
\ref{bay1} immediately follows from Theorem~\ref{cud1}.
In \ref{bay2},
if $X\in\peo$ then ${\Eo}\emn{\Eo\res X}$ by 
Theorem~\ref{2dih}, while if ${\Eo}\emn{\Eo\res X}$ 
then $\Eo\res X$ is not smooth since $\Eo$ itself is 
not smooth by Lemma~\ref{transv}\ref{transv5}.
In \ref{bay3},
if ${{\Eo}\res X}$ admits a Borel transversal then it
is smooth by Lemma~\ref{transv}\ref{transv1} and hence
$X$ belongs to $\ieo.$
To prove the converse apply Lemma~\ref{transv}\ref{transv3}.
\epf

Note that any $X\in\peo$ contains a closed subset $Y\sq X$
also in $\peo$ by Theorem~\ref{2dih}.
(Apply the theorem for $\rE=\Eo\res X.$
As $\Eo\res X$ is not smooth, we have ${\Eo}\emn{\Eo\res X},$
by  a continuous reduction $\vt.$
Take as $Y$ the full image of $\vt.$
$Y$ is compact, hence closed.)
Such sets $Y$ can be chosen in a special family.

\bdf
[{{\rm Zapletal~\cite{zap}}}]
\label{df:zap}
Suppose that two binary sequences $u^0_n\ne u^1_n\in \dln$
of equal length $\lh{u^0_n}=\lh{u^1_n}\ge1$
are chosen for each $n,$ together with one more
sequence $u_0\in\dln.$
Define
$\vt(a)=
u_0\we u^{a(0)}_0\we u^{a(1)}_1\we\dots$
for any $a\in\dn.$
Easily $\vt$ is a continuous injection $\dn\to\dn,$
$Y=\ran \vt$ is a closed set in $\dn,$ $\vt$ witnesses
$\Eo\emn{\Eo\res Y},$ and hence $Y\in\peo.$

Let $\pep$ denote the collection of all sets $Y$
definable in such a form.
\edf

\bte
[{{\rm Zapletal~\cite{zap}}}]
\label{jin}
$\pep$ is a dense subset of\/ $\peo:$
for any\/ $X\in\pep$ there exists\/ $Y\in\peo\yt Y\sq X.$
In addition, $\peo$ forces that the ``old'' continuum\/
$\gc$ remains uncountable.
\ete
\bpf
The proof employs splitting technique
for the forcing $\peo.$
This technique somewhat differs from
the splittings used in the proof of Theorem~\ref{2dih}.
First of all, as mentioned above, we can consider only
closed sets in $\peo,$ that enables us to replace the
Gandy -- Harrington stuff by a simple compactness
argument.
Second, the equivalence relation considered has the form
$\Eo\res X.$

For any sequences $r,w\in\dln$ with $\lh r\le\lh w,$
define $rw\in\dln$ (the \dd rshift of $w$)
so that $\lh{rw}=\lh w$ and 
$(rw)(k)=1-w(k)$ whenever $k<\lh r$ and $r(k)=1,$ and
$(rw)(k)=w(k)$ otherwise.
Clearly $r(rw)=w.$
Similarly define $ra\in\dn$ for $a\in\dn,$ and
$rX=\ens{ra}{a\in X}$ for any set $X\sq\dn.$

We are going to define sequences $u\in\dln$ and
$u^0_n\ne u^1_n\in\dln\;(n\in\dN)$ such that
$\lh{u^0_n}=\lh{u^1_n},$ as in
Definition~\ref{df:zap}, and also a system of closed sets
$X_s\in\peo\;\,(s\in\dln)$ satisfying the following:
\ben
\tenu{(\roman{enumi})}
\itla{jin1}\msur
$X_\La\sq X$ and $X_{s\we i}\sq X_s$;

\itla{jin2}\msur
$X_s\sq\cO_{w_s},$ where
$w_s=u_0\we u^{s(0)}_0\we u^{s(1)}_1\we\dots
\we u^{s(k-1)}_{k-1}\in\dln,$
$k=\lh s,$ and $\cO_w=\ens{a\in\dn}{w\su a}$ for $w\in\dln;$

\itla{jin3} 
if $s,t\in 2^n$ for some $n$ then $X_t=w_t w_s X_s$.
\een
Then define the map $\vt$ as in Definition~\ref{df:zap}.
The set $Y=\ran\vt=\bigcap_n\bigcup_{s\in2^n}X_s\sq X$
belongs to $\pep,$ proving the density claim of the theorem.

{\it Step 0\/}.
We put $X_\La=X$ and let $u_0\in\dln$
be the largest sequence such that $X_\La\sq \cO_{u_0}.$
Let $\ell_0=\lh{u_0}$.

{\it Step 1\/}.
Here we define $u^i_0$ and $X_{\ang i}$ for $i=0,1.$ 
Let $R$ be the set of all sequences $r\in\dln$ containing
at least one term equal to $1$
(and hence $ra\ne a$ for any $a$).
Consider the union $Z=\bigcup_{r\in R}Z_r$ of all sets
$Z_r=\ens{a\in X_\La}{ra\in X_\La};$ each $Z_r$ is closed.
The difference $D=X_\La\dif Z$ is pairwise
\dd\Eo inequivalent, hence $D\in\ieo$ by Lemma~\ref{bay}.
Thus at least
one of $Z_r\yt r\in R,$ belongs to $\peo$ by Lemma~\ref{bay}.
Let $r_1$ be any $r\in R$ of this sort.
Put $\ell_1=\lh r_1;$
clearly $\lh{u_0}=\ell_0<\ell_1$ and $r_1\res\ell_0$
consists only of terms equal to $0$.

There is a sequence $w_{\ang0}\in\dln$ such that
$\lh{w_{\ang0}}=\ell_1$
and the set $X_{\ang 0}=Z_{r_1}\cap\cO_{w_{\ang0}}$
still belongs to $\peo.$
Put $w_{\ang1}=r_1w_{\ang0}.$
Then the set
$X_{\ang 1}=r_1X_{\ang 0}=\ens{r_1a}{a\in X_{\ang 0}}=
Z_{r_1}\cap\cO_{w_{\ang1}}$
belongs to $\peo$ together with $X_{\ang 0}.$
%
Note that $u_0\su w_{\ang i},$ and hence there exist sequences
$u^0_0\ne u^1_0\in\dln$ of length $\ell_1-\ell_0$ such that
$w_{\ang0}=u_0\we u^0_0$ and $w_{\ang1}=u_0\we u^1_0.$
It follows from the construction that
$w_{\ang0}w_{\ang1}=r_1,$ therefore
$X_{\ang 1}=w_{\ang0}w_{\ang1}X_{\ang 1},$ and \ref{jin3}
holds.

{\it Step 2\/}.
Here we define $u^i_1$ for $i=0,1$ and $X_s$ for $s\in\dln$
with $\lh s=2.$
Once again there is a sequence $r_2\in R$ such that
the (closed) set
$Z_{r_2}=\ens{a\in X_{\ang 0}}{ra\in X_{\ang 0}}$
still belongs to $\peo.$
Put $\ell_2=\lh {r_2};$
then $\lh{r_{1}}=\ell_1<\ell_2$ and
$r_2\res{\ell_1}$ consists
only of terms equal to $0.$
Once again there is a sequence $w_{\ang{0,0}}\in\dln$ such that
$\lh{w_{\ang{0,0}}}=\ell_2$
and the set $X_{\ang{0,0}}=Z_{r_2}\cap\cO_{w_{\ang{0,0}}}$
belongs to $\peo.$
Put $w_{\ang{0,1}}=r_2w_{\ang{0,0}}.$
Then the set
$X_{\ang{0,1}}=r_2X_{\ang{0,0}}=Z_{r_2}\cap\cO_{w_{\ang{0,1}}}$
belongs to $\peo$ together with $X_{\ang{0,0}}.$ 
Also, put $w_{\ang{1,i}}=r_1w_{\ang{0,i}}$ and
$X_{\ang{1,i}}=r_1X_{\ang{0,i}}=Z_{r_2}\cap\cO_{w_{\ang{1,i}}}$
for $i=0,1$ --- these sets also belong to $\peo.$
As for \ref{jin3} at this level, take, for instance,
$s=\ang{0,1}$ and $t=\ang{1,0}.$
By definition
$X_{\ang{1,0}}=r_1X_{\ang{0,0}}=r_2r_1X_{\ang{0,1}},$
on the other hand, $w_{\ang{1,0}}=r_2r_1w_{\ang{0,1}},$ too.

Finally, there exist sequences
$u^0_1\ne u^1_1\in\dln$ of length $\ell_2-\ell_1$ such that
$w_{\ang{i,j}}=u_0\we u^i_0\we u^j_1$ for $i,j=0,1.$

{\it Steps $\ge3$\/}.
Et cetera.
The construction results in a system of sets and sequences
satisfying requirements \ref{jin1}, \ref{jin2}, \ref{jin3},
as required.

To prove the additional claim of the theorem, the splitting
construction has to be modified so that for any $n$ the sets
$X_s\yt s\in 2^n,$ belong to the \dd nth dense subset of
$\peo,$ in the sense of a given countable sequence of dense
sets.
\epf

We observe that $\peo$ as a forcing is somewhat closer to
Silver rather than Sacks forcing.
The property of minimality of the generic real, common to
both Sacks and Silver, holds for $\peo$ as well, the proof
resembles known arguments, but in addition the following
is applied: if $X\in\peo$ and $f:X\to\dn$ is a Borel
\dd\Eo invariant map (that is, ${x\Eo y}\imp{f(x)=f(y)}$)
then $f$ is constant on a set $Y\in\peo\yt Y\sq X.$~\snos
{Suppose, for the sake of brevity, that $X=\dn.$
For any $n,$ the set $Y^0_n=\ens{a}{f(a)(n)=0}$ is Borel and
\dd\Eo invariant.
It follows that $Y^0_n$ is either meager or comeager.
Put $b(n)=0$ iff $Y^0_n$ is comeager.
Then $D=\ens{a}{f(a)=b}$ is comeager.
A splitting construction as in the proof of Theorem~\ref{jin}
yields a set $Y\in\peo\yt Y\sq D$.}

\vyk{
It follows from the 2nd dichotomy that if $X\sq \dn$ is a 
Borel set then either $\Eo\res X$ is smooth or 
${\Eo\res X}\eqb \Eo.$ 
Let is show that the existence of a full Borel transversal 
is necessary and sufficient for the smoothness. 
Recall that a \dde{\it transversal\/} is a set 
\index{transversal}%
\index{transversal!partial}%
having exactly one element in common with every \dde class. 
Let a {\it partial transversal\/} be a set 
having $\le1$ element in common with every \dde class, in 
other words, it is required that the set is 
pairwise \dde inequivalent.

\bte
\label{baby}
Let\/ $X\sq\dn$ be\/ $\is11.$ 
Then exactly one of the following two conditions is 
satisfied$:$ 
\ben
\tenu{{\rm(\Roman{enumi})}}
\itla{bab1}
there is a\/ $\is11$ transversal\/ $S$ for\/ $\Eo\res X\;;$ 

\itla{bab2}
there is a closed set\/ $F\sq X$ such that\/ 
$\Eo\eqb{\Eo\res F}$.
\een
\ete 
\bpf
To see that \ref{bab1} and \ref{bab2} are incompatible, prove

\blt
\label{transd}
Any\/ $\is11$ partial transversal\/ $S$ for\/ $\Eo$ can be 
covered by a\/ $\id11$ partial transversal\/ $D$.
\elt
\bpf
$P=\ens{x\in\dn}{x\nin\ekeo{S\dif\ans{x}}}$ 
is a $\ip11$ set and $S\sq P:$ 
$P$ is the union of $S$ and all \dd\Eo classes 
disjoint from $S.$  
If $U$ is a $\id11$ set with $S\sq U\sq P$ then 
$D=\ens{x\in U}{\card(\ekeo{x}\cap U)=1}$ 
is as required,  
that $D$ remains $\id11$ follows from the countability  
of \dd\Eo classes.
\epF{Lemma}

Now, if $F$ is as in \ref{bab2} then 
$\Eo\res F$ has a $\id11$ transversal by the lemma, 
hence, $\Eo$ itself has a Borel transversal, then 
$\Eo$ is smooth, contradiction.

To see that one of $\ref{bab1}$ and $\ref{bab2}$ holds, 
define ${x\Eco y}$ iff 
$x,\,y$ belong to the same \dd\Eo invariant $\id11$ sets. 
The rest depends on the relationships between the given 
$X$ and the $\is11$ set  
$H_0=\ens{x}{\ekeo x\sneq\ekco x}$.\vom

{\sl Case 1\/}: 
$X\cap H_0=\pu.$ 
Following the proof of Lemma~\ref{d2eff}, we obtain a 
$\id11$ set $Y$ with $X\sq Y\sq \doP{H_0}$ and a $\id11$ 
reduction $\vt:Y\to\dn$ of $\Eo\res Y$ to $\rav\dn.$ 
The set $P=\ens{\ang{a,y}}{y\in Y\land \vt(y)=a}$ is then a 
$\id11$ set whose cross-sections $P_a=\vt\obr(a)$ are just 
classes of $\Eo\res Y.$ 
Then by \Cenu\ of \prf{eff} there is $\id11$ set 
$T\sq\dN\ti Y$ such that each $T_n=\ens{y}{T(n,y)}$ is a 
transversal for $\Eo\res Y$ and $Y=\bigcup_n T_n.$ 
Now it takes some efforts, which we leave for the reader 
(in particular \Kres\ is applied) 
to assemble a $\is11$ transversal for $\Eo\res X$ from 
parts of the transversals $T_n$.\vom

{\sl Case 2\/}: 
$Z=X\cap H_0\ne\pu.$ 
Then the argument of Case 2 of the proof of Theorem~\ref{2dih} 
(begin with $Z$ rather than with the whole $H$!) shows that 
there is a perfect set $F\sq Z$ such that 
$\Eo\reb{\Eo\res F}$.\vom

\epF{Theorem~\ref{baby}}
}

\parf{Ideal $\Ii$ and P-ideals}
\label{idI1}

By definition the ideal $\fio=\Ii$ consists of all sets 
$x\sq\pnn$ such that all, except for finitely many, 
cross-sections 
$\seq xn=\ens{k}{\ang{n,k}\in x}$ are empty. 
%

\punk{Ideals below $\Ii$}
\label{<i1}

It turns out that there exist only three different ideals 
Borel reducible to $\Ii,$ they are $\ifi,$ the disjoint sum 
$\ifi\oplus\pn,$ and $\Ii$ itself.

\bdf
An ideal $\cI$ is a {\it trivial variation\/} of $\cJ$
\index{trivial variation}%
if there is an infinite set $D$ such that 
$I\res D\cong\cJ$~\snos
{Recall that $\cI\cong\cJ$ means isomorphism 
via a bijection between the underlying sets.}
while $\cI\res\dop D=\cP(\dop D).$ 
(The last condition is equivalent to 
$\cI=\ens{x}{x\cap D\in\cI\res D}$.)
\edf

\bte[{{\rm Kechris~\cite{rig}}}]
\label{rig1}
If\/ $\cI\reb\Ii$ is a Borel 
(nontrivial) ideal on\/ $\dN$ then either\/ 
$\cI\cong\Ii$ or\/ $\cI$ is a trivial variation 
of\/ $\ifi.$
\ete

\bupt
Prove that any trivial variation of $\cI_1$ is isomorphic 
to $\cI_1$ 
while any trivial variation of $\ifi$ is isomorphic 
either to $\ifi$ or to the disjoint sum 
$\ifi\oplus\pn,$ \eg, realized in the form of 
$\ens{x\sq\dN}{x\cap{\textsc{odd}}\in\ifi}$. 
\eupt

\bpf[{{\rm Theorem}}]
We begin with another version of the method 
used in the proof of Theorem~\ref{jnmt}. 
Suppose that $\sis{\cB_k}{k\in\dN}$ is a fixed 
system of Borel subsets of $\pn.$ 
(It will be specified later.) 
Then there exists an increasing sequence of integers 
$0=n_0<n_1<n_2<...$ and sets $s_k\sq\il k{k+1}$ such 
that 
\ben
\tenu{(\arabic{enumi})}
\itla{ri1}
any $x\sq\dN$ with $\kai k\:(x\cap\il k{k+1}=s_k)$ 
is ``generic''~\footnote
{\label{genmean}\ 
We mean, Cohen generic over a certain 
fixed countable transitive model $\mm$ of a big enough 
fragment of $\ZFC,$ which contains Borel codes for all 
sets $\cB_k$.}~;

\itla{ri2}
if $k'\ge k$ and $u\sq\ilo{k'}$ then $u\cup s_{k'}$ 
decides $\cB_k$ in the sense that either any 
``generic'' $x\in\pn$ with $x\cap\ilo{k'+1}=u\cup s_{k'}$ 
belongs to $\cB_k$ or any ``generic'' $x$ with 
$x\cap\ilo{k'+1}=u\cup s_{k'}$ does not belong to $\cB_k$.
\een
Now put $\cD_0=\ens{x\cup S_1}{x\sq Z_0}$ and  
$\cD_1=\ens{x\cup S_0}{x\sq Z_1},$ where
\dm
\textstyle
S_0=\bigcup_ks_{2k}\sq Z_0=\bigcup_k\il{2k}{2k+1}\,,
\quad 
S_1=\bigcup_ks_{2k+1}\sq Z_1=\bigcup_k\il{2k+1}{2k+2}.
\dm
Clearly any $x\in\cD_0\cup\cD_1$ is ``generic'' by 
\ref{ri1}, hence, by \ref{ri2}, 
\ben
\tenu{(\arabic{enumi})}
\addtocounter{enumi}2
\itla{ri3}
each $\cB_k$ is clopen on both $\cD_0$ and $\cD_1$.
\een

As $\cI\reb\Ii,$ 
it follows from Lemma~\ref{l:bc} 
(and the trivial fact that $\cI_1\oplus\cI_1\cong\cI_1$) 
that there exists a {\it continuous\/} reduction 
$\vt:\pn\to\cP(\dN\ti\dN)$ of $\cI$ to $\Ii.$
Thus $\rE_\cI$ is 
the union of an increasing sequence of (topologically) 
closed \er s $\rR_m\sq\rE_\cI$ just because $\cI_1$ 
admits such a form. 
We now require that $\sis{\cB_k}{}$ includes  
all sets 
$B^m_l=\ens{x\in\pn}{\kaz s\sq\ir0l\;x\rR_m(x\sd s)}.$ 
Then by \ref{ri3} and the compactness of $\cD_i$ 
for any $l$ there is $m(l)\ge l$ satisfying 
\ben
\tenu{(\arabic{enumi})}
\addtocounter{enumi}3
\itla{ri4}
$\kaz x\in\cD_0\cup\cD_1\:\kaz s\sq\ir0l\;
\skl x\rR_{m(l)}(x\sd s)\skp$.
\een

To prove the theorem it suffices to obtain a sequence 
$x_0\sq x_1\sq x_2\sq...$ of sets $x_k\in\cI$ with 
$\cI=\bigcup_n\pws{x_n}:$ 
that in this case $\cI$ is as required is an easy exercise. 
As any topologically closed ideal is easily $\pws x$ for 
some $x\sq\dN,$ it suffices to show that $\cI$ is a union 
of a countable sequence of closed subideals. 
It suffices to demonstrate this fact separately for 
$\cI\res Z_0$ and $\cI\res Z_1.$ 
Prove that 
$\cI\res Z_0$ is a countable union of closed subideals, 
ending the proof of the theorem.

If $m\in \dN$ and $s\sq u\sq Z_0$ are finite then let 
\dm
I^m_{us}=\ens{A\sq Z_0}{\kaz x\in\cD_0\:
(x\cap u=s\imp {\skl x\cup(A\dif u)\skp\rR_m x})}\,.
\dm

\blt
\label{rig1l}
Sets\/ $I^m_{us}$ are closed topologically and 
under\/ $\cup,$ and\/ $I^m_{us}\sq\cI$.
\elt
\bpf
$I^m_{us}$ are topologically closed because so are $\rR_m$. 

Suppose that $A,\,B\in I^m_{us}.$  
To prove that $A\cup B\in I^m_{us},$ let 
$x\in\cD_0$ satisfy $x\cap u=s.$ 
Then $x'=x\cup(A\dif u)\in\cD_0$ satisfies $x'\cap u=s,$ 
too, hence, as $B\in I^m_{us},$ we have 
$\skl{x'\cup(B\dif u)}\skp\rR_m x',$ thus, 
$\skl{x\cup((A\cup B)\dif u)}\skp\rR_m x'.$ 
However $x'\rR_m x$ just because $A\in I^m_{us}.$ 
It remains to recall that $\rR_m$ is a \er.

To prove that any $A\in I^m_{us}$ belongs to 
$\cI$ take $x=s\cup S_1.$ 
Then we have ${x\cup(A\dif u)}\rR_m x,$ thus, 
$A\in\cI$ as $s$ is finite and $\rR_m\sq\rE_\cI$.
\epF{Lemma}

\blt
\label{rig1L}
$\cI\res Z_0=\bigcup_{m,\,u,\,s}I^m_{us}$.
\elt
\bpf
Let $A\in\cI,\msur$ $A\sq Z_0.$ 
The sets 
$Q_m=\ens{x\in\cD_0}{\skl x\cup A\skp\rR_m x}$ 
are closed and satisfy $\cD_0=\bigcup_m Q_m.$ 
It follows that one of them has a non-empty interior 
in $\cD_0,$ thus, there exist finite sets 
$s\sq u\sq Z_0$ and some $m_0$ with 
\dm
\kaz x\in\cD_0\:
(x\cap u=s\imp {\skl x\cup A\skp\rR_{m_0} x})\,.
\dm
This is not exactly what we need, however, 
by \ref{ri4}, there exists a number 
$m=\tmax\ans{m_0,m(\tsup u)}$ 
big enough for 
\dm
\kaz x\in\cD_0:
\skl x\cup A\skp\rR_{m} \skl x\cup(A\dif u)\skp\,.
\dm
It follows that $A\in I^m_{su},$ as required.
\epF{Lemma}

Let $J^m_{su}$ be the hereditary hull of $I^m_{su}$ 
(all subsets of sets in $I^m_{su}$). 
It follows from Lemma~\ref{rig1l} that any $J^m_{su}$ 
is a topologically closed subideal of $\cI\res Z_0,$ 
however, $\cI\res Z_0$ is the union of those 
ideals by Lemma~\ref{rig1L}, as required.
\epf

\vyk{
\bdf
\label{api}
An ideal $\cI$ is a {\it P-ideal\/} if for any sequence of 
sets $x_n\in\cI$ there is a set $x\in\cI$ such that 
$x_n\sqa x$ (\ie, $x_n\dif x\in\ifi$) for all $n$.
\edf

For instance, the ideals $\ifi,\,\Id,\,\It,\,\Zo$ 
(but not $\Ii$!) are P-ideals. 

Recall that a {\it submeasure\/} on a set $A$ is any map 
$\vpi:\cP(A)\to[0,\piy],$ satisfying $\vpi(\pu)=0,$ 
$\vpi(\ans a)<\piy$ for any $a\in A,$
and $\vpi(x)\le \vpi(x\cup y)\le \vpi(x)+\vpi(y).$  
To be a {\it measure\/}, a sumbeasure $\vpi$ has to 
satisfy, in addition, that 
$\vpi(x\cup y)=\vpi(x)+\vpi(y)$ whenever $x,\,y$ 
are disjoint. 

Suppose that $\vpi$ is a submeasure on $\dN.$ 
Define the {\it tailsubmeasure\/} 
$\vpy(x)=\hv x\vpi=\tinf_n(x\cap\iry n).$ 
The following ideals are considered:
\dm
\bay{rcllll}
\Fin_\vpi &=& \ens{x\in\pn}{\vpi(x)<\piy} &&&;\\[1ex]

\Nul_\vpi &=& \ens{x\in\pn}{\vpi(x)=0} &&&;\\[1ex]

\Exh_\vpi &=& \ens{x\in\pn}{\vpy(x)=0}&=& \Nul_{\vpy}&.\\[0.5ex]
\eay
\dm
A submeasure $\vpi$ on $\dN$ is lover semicontinuous, or 
\lsc\ for brevity, if we have 
$\vpi(x)=\tsup_n\vpi(x\cap\ir0n)$ for all $x\in\pn.$ 
Note that if $\vpi$ is \lsc\ then $\vpy$ is not necessarily 
\lsc\ itself.
Any \dds additive submeasure is \lsc.

\bex
\label{exh:e}
$\ifi=\Exh_\vpi=\Nul_\vpi,$ where $\vpi(x)=1$ for any 
$x\ne \pu.$ 
We also have $\ofi=\Exh_\psi,$ where 
$\psi(x)=\sum_k\,2^{-k}\,\vpi(\ens{l}{\ang{k,l}\in x})$ 
is \lsc.
\eex

\bdf
\label{polI}
An ideal $\cI$ is {\it polishable\/} if 
there is a Polish group topology $\tau$ on $\cI$ which 
produces the same Borel subsets of $\cI$ as the product 
topology of $\dn,$ restricted on $\cI,$ does. 
($\cI$ is considered as a group with $\sd$ as the 
operation.) 
\edf
}

\vyk{
Let $T$ be the ordinary Polish product topology on $\pn$. 

\ble
\label{sol:?}
Suppose that an ideal\/ $\cI\sq\pn$ is polishable. 
Then there is only one Polish group topology\/ 
$\tau$ on\/ $\cI.$ 
This topology\/ refines\/ $T\res\cI$ 
and is metrizable by a\/ \dd\sd invariant metric. 
If\/ $Z\in\cI$ then\/ $\tau\res\cP(Z)$ coincides 
with\/ $T\res\cP(Z).$ 
In addition, $\cI$ itself is\/ \dd TBorel.
\ele
\bpf
Let $\tau$ witness that $\cI$ is polishable. 
The identity map 
$f(x)=x:\stk\cI\tau\to\stk\pn T$ 
is a \dd\sd homomorphism and is Borel-measurable 
because all \dd{(T\res\cI)}open sets are \dd\tau Borel, 
hence, by the Pettis theorem (Kechris~\cite[??]{dst}), 
$f$ is continuous. 
It follows that all \dd{(T\res\cI)}open subsets of 
$\cI$ are \dd\tau open, and that $\cI$ is 
\dd TBorel in $\pn$ because $1-1$ continuous 
images of Borel sets are Borel. 

A similar ``identity map'' argument shows that 
$\tau$ is unique if exists. 

It is known (Kechris~\cite[]{dst}) 
that any Polish group topology admits a 
left-invariant compatible metric, which, in this case, 
is right-invariant as well since $\sd$ is an 
abelian operation.

Let finally $Z\in\pn.$ 
Then $\cP(Z)$ is \dd Tclosed, hence, 
\dd\tau closed by the above, subgroup of $\cI,$ and 
$\tau\res\cP(Z)$ is a Polish group topology on $\cP(Z).$ 
Yet $T\res\cP(Z)$ is another Polish group topology on 
$\cP(Z),$ with the same Borel sets. 
It follows, by the same ``identity map'' argument, 
that $T$ and $\tau$ coincide on $\cP(Z)$.
\epf

\bex
\label{fio}
$\cI_1=\fio$ is not polishable. 
Indeed we have $\fio=\bigcup_n W_n,$ where  
$W_n=\ens{x}{x\sq\ans{0,1,...,n}\ti\dN}.$  
Let, on the contrary, $\tau$ be a Polish group topology 
on $\cI_1.$ 
Then $\tau$ and the ordinary topology $T$ coincide on 
each set $W_n$ by the lemma, in particular, 
each $W_n$ remains \dd\tau nowhere dense 
in $W_{n+1},$ hence, in $\cI_1,$ a contradiction with 
the Baire category theorem for $\tau$.
\eex

On the contrary, $\ifi$ and $\ofi$ are polishable, that 
follows from the results of Example~\ref{exh:e} 
and Theorem~\ref{sol} below, saying that P-ideals, 
polishable ideals, and those of the form $\Exh_\vpi,$ 
where $\vpi$ is \lsc, is one and the same.
}

\punk{$\Ii$ and P-ideals}
\label{sol: prf}

Thus $\Ii$ is a \dd\reb minimal ideal over $\ifi:$ we have 
$\ifi\rebs\Ii$ and the \dd\rebs interval $(\ifi,\Ii)$ is 
empty. 
Although $\Ii$ is not {\it the least\/} over $\ifi,$ still 
it turns out that $\Ii$ is the least among all Borel 
ideals which are not P-ideals.  


The next theorem is of great importance for the whole 
theory of Borel ideals.

\bte[{{\rm Solecki~\cite{sol,sol'}}}]
\label{sol}
The following families of ideals on\/ $\dN$ coincide$:$
\ben
\tenu{{\rm(\roman{enumi})}}
\itla{s2}
ideals of the form\/ $\Exh_\vpi,$ where\/ $\vpi$ is a 
\lsc\ submeasure on\/ $\dN\;;$

\itla{s3} 
polishable ideals. 

\itla{s1}
analytic P-ideals$;$

\itla{s4} 
analytic ideals $\cI$ with $\cI_1\not\orb \cI\;;$ 

\itla{s5} 
analytic ideals $\cI$ such that all countable unions of\/ 
\ddi small sets are \ddi small,
{\rm where a set $X\sq\pn$ is \ddi{\it small\/} if there
is $A\in\cI$ such that 
$X\res A=\ens{x\cap A}{x\in X}\sq\cP(A)$ 
is meager in $\cP(A)$\/}.
\een
\ete

It follows that all analytic P-ideals actually belong to 
$\fp03,$ just because any ideal of type \ref{s2} is 
easily $\fp03$.

\bpf 
The formal scheme of the proof is: 
$\ref{s2}\imp\ref{s3}\imp\ref{s1}\imp\ref{s4}
\imp\ref{s5}\imp\ref{s2}.$ 
The hard part will be $\ref{s5}\imp\ref{s2},$ the 
rest is rather elementary but tricky in some points. 
The elementary part of the proof is organized so that 
the proofs that $\ref{s2}\eqv\ref{s3}$ and 
$\ref{s1}\eqv\ref{s4}\eqv\ref{s5}$ and that the 
first group implies the second, are obtained 
independently of the hard part.\vom 
\imar{\bfit Give corollaries of Thm~\ref{sol}\/} 

$\ref{s2}\imp\ref{s3}$ 
If $\vpi(\ans n)>0$ for all $n$ then the required metric 
on $\cI=\Exh_\vpi$ can be defined by 
$\dpi(x,y)=\vpi(x\sd y).$ 
Then any set $U\sq\cI$ open in the sense of the ordinary 
topology (the one inherited from $\pn$) is \dd\dpi open, 
while any \dd\dpi open set is Borel in the ordinary sense.
In the general case we assemble the required metric of 
$\dpi$ on the domain $\ens{n}{\vpi(\ans n)>0}$ and the 
ordinary Polish metric on $\pn$ on the complementary 
domain.\vom

$\ref{s3}\imp\ref{s2}$ 
Let $\tau$ be a Polish group topology on $\cI,$ 
generated by a \dd\sd invariant compatible metric $d.$ 
It can be shown (Solecki~\cite[p.\ 60]{sol'}) that 
$\vpi(x)=\tsup_{y\in\cI,\:y\sq x}\,d(\pu,x)$ is 
a \lsc\ submeasure with $\cI=\Exh_\vpi.$ 
The key observation is that for any $x\in\cI$ the 
sequence $\sis{x\cap\ir0n}{n\in\dN}$ \dd dconverges 
to $x$ by the last statement of Lemma~\ref{sol:?}, which 
implies both that $\vpi$ is \lsc\ 
(because the supremum above can be restricted to 
finite sets $y$) 
and that $\cI=\Exh_\vpi$ 
(where the inclusion $\supseteq$ needs another 
``identity map'' argument).\vom

$\ref{s2}\imp\ref{s1}$ 
That any $\cI=\Exh_\vpi,$ $\vpi$ being \lsc, is 
a P-ideal, is an easy exercise:  
if $x_1,\,x_2,\,x_3,\,...\in\cI$ then define 
an increasing sequence of numbers $n_i\in x_i$ with 
$\vpi(x\cap\iry{n_i})\le2^{-n}$ and put 
$x=\bigcup_i(x\cap\iry{n_i})$.\vom

$\hbox{Any of \ref{s1},\,\ref{s2},\,\ref{s3},\,\ref{s5}}
\limp \ref{s4}$ \ 
This is because $\cI_1$ easily does not satisfy any of 
the four properties indicated. 
For the formal purpose to complete the proof of 
Theorem~\ref{sol}, we need here only the implication 
$\ref{s1}\imp\ref{s4}$.\vom 

$\ref{s4}\imp\ref{s5}$ \ 
Suppose that sets $X_n\sq\pn$ are \ddi small, so that 
$X_n\res A_n$ is meager in $\cP(A_n)$ for some $A_n\in\cI,$  
but $X=\bigcup_n X_n$ is not \ddi small, and prove 
$\cI_1\orb \cI.$ 
Arguing as in the proof of Theorem~\ref{jnmt}, we use 
the meagerness to find, for any $n,$ a sequence of 
pairwise disjoint non-empty 
finite sets $w^n_k\sq x_n,\msur$ $k\in\dN,$ and 
subsets $u^n_k\sq w^n_k,$ such that 
\ben
\tenu{(\alph{enumi})}
\itla{sa1}
if $x\sq \dN$ and $\exi k\:(x\cap w^n_k=u^n_k)$ 
then $x\nin X_n$. 
\een
Dropping some sets $w^n_k$ away and reenumerating the 
rest, we can strengthen the disjointness to the following: 
$w^n_k\cap w^m_l=\pu$ unless both $n=m$ and $k=l$. 

Now put $w^n_{ij}=w^n_{2^i(2j+1)-1}.$  
The sets $\ovw_{ij}=\bigcup_{n\le i} w^n_{ij}$ 
are still pairwise disjoint, and satisfy the following 
two properties: 

\ben
\tenu{(\alph{enumi})}
\addtocounter{enumi}1
\itla{sa2}\msur
$\bigcup_j\ovw_{ij}\sq x_n,$ hence, $\in \cI,$ 
for any $i$;

\itla{sa3}
if a set $Z\sq\dN\ti\dN$ does not belong to $\cI_1,$ 
\ie, $\exi i\:\sus j\:(\ang{i,j}\in Z),$ then 
$\kaz n\:\exi k\:(w^n_k\sq\ovw_Z),$ 
where $\ovw_Z=\bigcup_{\ang{i,j}\in K}\ovw_{ij})$.
\een
We assert that the map $\ang{i,j}\mapsto\ovw_{ij}$ 
witnesses $\cI_1\orbp\cI.$ 
(Then a simple argument, as in the proof of 
Theorem~\ref{jnmt}, gives $\cI_1\orb\cI.$) 

Indeed if $Z\sq\dN\ti\dN$ belongs to $\cI_1$ then 
$\ovw_Z\in\cI$ by \ref{sa2}. 
Suppose that $Z\nin\cI_1.$ 
It suffices to show that $X_n\res \ovw_Z$ is meager in 
$\cP(\ovw_Z)$ for any $n.$ 
Note that by \ref{sa3} the set 
$K=\ens{k}{w^n_k\sq\ovw_Z}$ is infinite and in fact 
$\ovw_Z\cap x_n=\bigcup_{k\in K}w^n_k.$ 
Therefore, any $x\sq \ovw_Z$ satisfying 
$x\cap w^n_k=u^n_k$ for infinitely many $k\in K,$ 
does not belong to $X_n$ by \ref{sa1}. 
Now the meagerness of $X_n\res \ovw_Z$ is clear.\vom

$\ref{s5}\imp\ref{s1}$ \ 
This also is quite easy: if a sequence of sets $Z_n\in\cI$ 
witnesses that $\cI$ is not a P-ideal, then the union of 
\ddi small sets $\cP(Z_n)$ is not \ddi small.

\punk{The hard part}

We prove $\ref{s5}\imp\ref{s2},$ the hard part of 
Theorem~\ref{sol}. 
A couple of definitions before the key lemma.  
\bit
\item
Let $C(\cI)$ be the collection of all hereditary 
(\ie, ${y\sq x\in K}\imp {y\in K}$) compact 
\ddi large sets $K\sq\pn.$ 

\item
Given sets $A,\,B\sq\pn,$ let 
$A+B=\ens{x\cup y}{x\in A\land y\in B}$.
\eit

\ble
\label{sl2}
Assuming that\/ $\cI$ is of type \ref{s5}, there is 
a countable sequence of sets\/ $K_m\in C(\cI)$ 
such that\/ 
for any set\/ $K\in C(\cI)$ there are\/ $m,\,n$ with\/ 
$K_m+K_n\sq K$. 
\ele
\bpf
Fix a continuous map $f:\bn\onto\cI.$ 
For any $s\in\nse,$ we define 
\dm
N_s=\ens{a\in\bn}{s\su a}\quad\hbox{and}\quad
B_s=f\ima N_s \quad\text{(the \dd fimage of $N_s$)}\,. 
\dm
Consider the set 
$T=\ens{s}{B_s\,\hbox{ is \ddi large}}.$ 
As $\cI$ itself is clearly \ddi large, $\La\in T.$ 
On the other hand, the assumption \ref{s5} easily implies 
that $T$ has no endpoints and no isolated branches, hence, 
$P=\ens{a\in\bn}{\kaz n\:(a\res n\in T)}$ is a perfect 
set. 
Moreover, $A_s=f\ima(P\cap N_s)$ is \ddi large for 
any $s\in T$ because $B_s\dif A_s$ is a countable 
union of \ddi small sets.

Now consider any set $K\in C(\cI).$ 
By definition, if $x,\,y\in\cI$ then $Z=x\cup y\in\cI,$ 
thus, $K\res Z$ is not meager in $\cP(Z),$ 
hence, by the compactness, 
$K\res Z$ includes a basic nbhd of $\cP(Z),$ hence, 
by the hereditarity, there is a number $n$ such that 
$Z\cap\iry n\in K.$ 
We conclude that $P^2=\bigcup_n Q_n,$ where each 
$Q_n=\ens{\ang{a,b}\in P^2}
{(f(a)\cup f(b))\cap\iry n\in K}$ 
is closed in $P$ because so is $K$ and $f$ is continuous. 
Thus, there are $s,\,t\in T$ such that 
$P^2\cap (N_s\ti N_t)\sq Q_n,$ in other words, 
$(A_s+A_t)\res\iry n\sq K,$ hence, 
$(\nad{A_s}+\nad{A_t})\res\iry n\sq K,$ where 
$\nad{\vphantom|...}$ denotes 
the topological closure of the hereditary hull.
Thus we can take, as $\sis{K_m}{},$ all sets of the 
form $K_{sn}=\nad{A_s}\res n$. 
\epf

Using the fact that $C(\cI)$ is a filter 
(as easy exercise which makes main use if the 
hereditarity), we can define 
(still in the assumption that $\cI$ is of type \ref{s5}) 
a \dd\sq decreasing sequence of sets $K_n\in C(\cI)$ 
such that
\ben
\tenu{(\arabic{enumi})}
\itla{ck1} 
for any $K\in C(\cI)$ there is $n$ with $K_n\sq K$,
\een
and $K_{n+1}+K_{n+1}\sq K_n$ for any $n.$ 
Taking any other term of the sequence, we can sharpen 
the latter requirement to
\ben
\tenu{(\arabic{enumi})}
\addtocounter{enumi}1
\itla{ck2}
for any $n:$ $K_{n+1}+K_{n+1}+K_{n+1}\sq K_n$.
\een

This is the starting point for the construction of a 
\lsc\ submeasure $\vpi$ with $\cI=\Exh_\vpi.$ 
Assuming that, in addition, $K_0=\pn,$ let, for any 
$x\in\pwf\dN$, 
\dm
\bay{lcll}
\vpi_1(x) &=& \tinf\left\{\,2^{-n}:x\in K_n\,\right\}&,
\hbox{ and}\\[1ex]

\vpi_2(x) &=& \tinf\left\{\,\textstyle
\sum_{i=1}^{m}\vpi_1(x_i):m\ge 1\land 
x_i\in\pwf\dN\land x\sq\bigcup_{i=1}^{m}x_i\,
\right\}&.
\eay
\dm
Then set $\vpi(x)=\tsup_n\vpi_2(x\cap\ir0n)$ for 
any $x\sq\dN.$ 
A routine verification shows that $\vpi$ submeasure 
and that $\cI=\Exh_\vpi.$ 
(See Solecki~\cite{sol'}. 
To check that any $x\in\Exh_\vpi$ belongs to $\cI$ 
we use the following observation: $x\in\cI$ iff for 
any $K\in C(\cI)$ there is $n$ such that  
$x\cap\ir0n\in K$.)

\epF{Theorem~\ref{sol}}

\bcor
\label{<pideal}
Suppose that\/ $\cJ$ is an analytic P-ideal. 
Then any ideal\/ $\cI\reb\cJ$ is an analytic P-ideal, too.
\ecor
\bpf
Use equivalence $\ref{s4}\eqv\ref{s1}$ of the theorem. 
(The result can be obtained via a more direct argument, 
of course.)
\epf

\parf{Equivalence relation $\Ei$}
\label{Ei:er}

The ideal $\Ii$ naturally defines the \er\ $\Ei=\rE_{\Ii}$ 
on $\pnn$ so that  
$x\Ei y$ iff $x\sd y\in\Ii.$ 
\index{equivalence relation, ER!E1@$\Ei$}%
\index{zzE1@$\Ei$}%
We can as well consider $\Ei$ as an \er\ on $\dX^\dN$
for any uncountable Polish space $\dX,$ 
defined as $x\Ei y$ iff $x(k)=y(k)$ for all but finite $k.$

\punk{$\Ei$ and hypersmoothness}
\label{einc}

The following notation will be rather useful in our study 
of subsets of $\pnd$ or $(\dn)^\dN.$
If $x$ is a function defined on $\dN$ then, for any $n,$ let 
\dm
x\rmq n= x\res{[0,n)}\,,\;\; 
x\rme n= x\res{[0,n]}\,,\;\;
x\qc n = x\res{(n,\iy)}\,,\;\; 
x\qec n= x\res{[n,\iy)}\,.
\dm 
\index{zzx<n@$x\rmq n$}%
\index{zzx<=n@$x\rme n$}%
\index{zzx>n@$x\qc n$}%
\index{zzx>=n@$x\qec n$}%
For any set $X$ of \dd\dN sequences, let 
$X\rmq n=\ens{x\rmq n}{x\in X},$ and similarly for 
${\le}\zd{>}\zd{\ge}.$
If $\xi\in X\qc n$ then let 
$\srez X\xi=\ens{x(n)}{x\in X\land x\qc n=\xi}.$
\index{zzSxxi@$\srez X\xi$}%

Recall that a hypersmooth \er\ is a countable increasing 
union of Borel smooth \er s. 
The following lemma shows that $\Ei$ is universal in this 
class.

\ble
\label{hs<e1}
For a Borel \er\/ $\rE$ to be hypersmooth it is necessary 
and sufficient that\/ $\rE\reb\Ei$.
\ele
\bpf
Let $\dX$ be the domain of $\rE.$ 
Assume that $\rE$ is hypersmooth, \ie, $\rE=\bigcup_n\rE_n,$ 
where $x\rE_n y$ iff $\vt_n(x)=\vt_n(y),$ each 
$\vt_n:\dX\to\pn$ is Borel, and $\rE_n\sq\rE_{n+1},\:\kaz n.$ 
Then $\vt(x)=\sis{\vt_n(x)}{n\in\dN}$ witnesses 
$\rE\reb\Ei.$ 
Conversely, if $\vt:\dX\to\pnd$ is a Borel reduction 
of $\rE$ to $\Ei$
then the sequence of \er s $x\rE_n y$ iff 
$\vt(x)\qec n=\vt(y)\qec n$ witnesses that $\rE$ 
is hypersmooth.
\epf

This Subsection contains a couple of results which 
describe the relationships between hypersmooth and countable 
\er s.
The following result is given in \cite{hypsm} with a 
reference to earlier papers.


\ble
\label{eoei'}
\ben
\tenu{{\rm(\roman{enumi})}}
\itsep
\itla{eoei'1}
$\Ei$ is not essentially countable, \ie, there is no Borel 
countable\/ {\rm(that is, with at most countable classes)} 
\er\ $\rE$ such that\/ $\Ei\reb\rE$.

\itla{eoei'2}
$\Eo\rebs\Ei,$ in other words,\/ $\ifi\rebs\Ii$.
\een
\ele
\bpf
\ref{eoei'1}
(A version of the argument in \cite{hypsm}, 1.4 and 1.5.) 
Let $\dX$ be the domain of $\rE,$ and 
$\vt:\pnd\to\dX$ a Borel map satisfying  
${x\Ei y}\limp{\vt(x)\rF\vt(y)}.$  
Then $\vt$ is continuous on a dense $\Gd$ set $D\sq\pnd.$ 
We begin with a few definitions. 
Let ``generic'' mean Cohen generic over a certain 
fixed countable transitive model $\mm$ of a big enough 
fragment of $\ZFC,$ which contains codes for 
$D\zd {\vt\res D}\zd \dX$. 

We are going to define, for any $k,$ a pair of 
$x_k\ne y_k\in\pn,$ a number $\ell(k)$ and a tuple 
$\za_k\in\pn^{\ell(k)}$ such that 
\ben
\tenu{(\arabic{enumi})}
\itla{^1}
both $x=\ang{x_0}\we\za_0\we \ang{x_1}\we\za_1\we...$ and 
$y=\ang{y_0}\we\za_0\we \ang{y_1}\we\za_1\we...$ are ``generic'' 
elements of $\pnd$;

\itla{^2}
for any $k,$ 
$\za_{\le k}=\ang{x_0,y_0}\we\za_0\we 
\ang{x_1,y_1}\we\za_1\we...\we\ang{x_k,y_k}\we\za_k$ 
is ``generic'', hence, so are 
$\xi_{\le k}=\ang{x_0}\we\za_0...\we\ang{x_k}\we\za_k$ and 
$\eta_{\le k}=\ang{y_0}\we\za_0...\we\ang{y_k}\we\za_k$; 

\itla{^3}
for any $k$ and any $z\in\pnd$ such that $\za_{\le k}\we z$ 
is ``generic'' we have 
$\vt(\xi_{\le k}\we z)=\vt(\eta_{\le k}\we z)$.
\een

If this is done then we can choose, using \ref{^2}, a point 
$z_{(k)}\in\pnd$ for any $k$ so that 
$\za_{\le k}\we z_{(k)}\in\pnd$ is ``generic'', hence, by 
\ref{^3}, for $x_{(k)}=\xi_{\le k}\we z_{(k)})$ and 
$y_{(k)}=\eta_{\le k}\we z_{(k)})$ we have 
$\vt(x_{(k)})=\vt(y_{(k)}).$ 
Note that $x_{(k)}\to x$ and $y_{(k)}\to y,$ and on the other 
hand, all of $x_{(k)},\ x,\,y_{(k)},\,y$ belong to $D$ because 
all are ``generic''. 
It follows that $\vt(x)=\vt(y)$ by the choice of $D.$ 
However obviously $\neg\;{x\Ei y},$ so that $\vt$ is not a 
reduction, as required.

To define $x_0,\,y_0,\,\za_0$ note that, by an ordinary 
splitting argument, there is a set $X\sq\pn$ of cardinality 
$\cont$ and $z\in\pnd$ such that $\ang{a,b}\we z$ is 
``generic'' for any two $a\ne b\in X.$ 
In particular, all $\ang a\we z\zT a\in X,$ are ``generic''. 
But all of them are pairwise \dd\Ei equivalent, hence, $\vt$ 
sends all of them into one and the same \ddf class, which 
is a countable set by the choice of $\rF.$ 
It follows that there is a pair of $a\ne b$ in $X$ such that 
$\vt(\ang a\we z)\ne\vt(\ang b\we z).$  
This equality is a property of the ``generic'' object 
$\ang{a,b}\we z,$ hence, it is forced in the sense that 
there is a number $\ell$ such that  
$\vt(\ang a\we z')\ne\vt(\ang b\we z')$ whenever 
$\ang{a,b}\we z'$ is ``generic'' with $z'\res\ell=z\res\ell.$ 
Put $x_0=a\zT y_0=b\zT \za_0=z\res\ell.$ 

The induction step is carried out by the same argument.

\ref{eoei'2}
That $\Eo\reb\Ei$ is witnessed by the map 
$f(x)=\ens{\ang{0,n}}{n\in x}.$ 
\epf

\vyk{
\bcot
\label{eoei"}
$\Eo\rebs\Ei,$ in other words,\/ $\ifi\rebs\Ii$.
\imar{why $\Ei$ and $\Et$ incomparable?}%
\ecot
\bpf
That $\Eo\reb\Ei$ is witnessed by the map 
$f(x)=\ens{\ang{0,n}}{n\in x}.$ 
\epf
}

While $\Ei$ is not countable, the conjunction of 
hypersmootheness and countability characterizes the 
essentially more primitive class of hyperfinite \er s. 

\vyk{
The following lemma is a necessary step towards the full 
result below.

\blt[{{\rm a part of Theorem~5.1 in \cite{djk}}}]
\label{hf&c=hs}
Suppose that\/ $X\sq\pnd$ is a Borel set and\/ $\Ei\res X$ 
is a countable \er. 
Then\/ $\Ei\res X$ is hyperfinite.
\elt
\bpf

Now, that $\Ei\res X$ is hyperfinite follows from 
Theorem~\ref{thf}.
\epf
}

\punk{The 3rd dichotomy}
\label{hypersm}

The following major result is called the 3rd dichotomy theorem.
 
\bte[{{\rm Kechris and Louveau \cite{hypsm}}}]
\label{kelu}
Suppose that\/ $\rE$ is a Borel\/ \er\ on some Polish 
space, and\/ $\rE\reb\Ei.$ 
Then either\/ $\rE\reb\Eo$ or\/ $\Ei\reb\rE$.
\ete
\bpf
Starting the proof, we may assume that $\rE$ is a 
$\id11$ \er\ on $\pn,$ and that there is a reduction 
$\rho$ of $\rE$ to $\Ei,$ of class $\id11.$ 
Then $R=\ran\rho$ is a $\is11$ subset of $\pnd.$ 
The idea behind the proof is to show that the set $R$ 
is either small enough for $\Ei\res R$ to be 
Borel reducible to $\Eo,$ or otherwise it is big enough 
to contain a closed subset $X$ such that $\Ei\res X$ 
is Borel isomorphic to $\Ei$.

Relations $\cl$ and $\cle$ will denote the inverse 
order relations on $\dN,$ \ie, $m\cle n$ iff $n\le m,$ 
and $m\cl n$ iff $n<m.$  
If $x\in\pnd$ then $x\rec n$ denotes the restriction 
of $x$ (a function defined on $\dN$) 
on the domain ${\cle n},$ \ie, $\iry n.$
If $X\sq\pnd$ then let $X\rec n=\ens{x\rec n}{x\in X}.$ 
Define $x\rc n$ and $X\rc n$ similarly.
In particular, $\pnd\rec n=\pn^{\cle n}=\pn^{\iry n}.$  

For a sequence $x\in\pn^{\cle n},$ let $\glu x$ 
(the {\it depth\/} of $x$) 
be the number (finite or $\iy$) of elements of the set 
$\nab x=\ens{j\cle n}{x(j)\nin \id11(x\rc{j})}.$   
The formula $\glu x\ge d$ 
(of two variables, $d$ running over $\dN\cup\ans\iy$) 
is obviously $\is11$. 

We have two cases:\vtm

{\bfit Case 1\bf:} \ 
all $x\in R=\ran \rho$ satisfy $\glu x<\iy$.\vtm

{\bfit Case 2\bf:} \ 
there exist $x\in R$ with $\glu x=\iy$.

Case 1 is the easier case. 
First of all we observe that $R,$ a $\is11$ set, is a 
subset of the $\ip11$ set $Z=\ens{x}{\glu x<\iy},$ 
hence, there is a $\id11$ set $Y$ with 
$\ran\rho\sq Y\sq Z.$ 
The following lemma ends the argument.

\blt
\label{dep3}
Suppose that\/ $X\sq\pnd$ is a\/ $\id11$ set and any\/ 
$x\in X$ satisfies\/ $\glu x<\iy.$ 
Then\/ $\rE_1\res X\reb\rE_0$.
\elt
\bpf
By the choice of $X$ for any $x\in X$ there is a number 
$n$ such that 
${\kaz m\cle n\:\skl x(m)\in\id11(x\rc{m})\skp}.$ 
As the relation between $x$ and $n$ here is clearly $\ip11,$ 
the ``Kreisel selection'' theorem yields a $\id11$ map 
$\nu:X\to\dN$ such that $x(m)\in\id11(x\rc{n})$ 
holds whenever $x\in X$ and $m\cle\nu(x).$ 
Now define, for each $x\in X,$ $\vt(x)\in\pnd$ as follows: 
$\vt(x)\rec{\nu(x)}=x\rec{\nu(x)},$ but $\vt(x)(j)=\pu$ 
for all $j<\nu(x).$ 
Note that $x\rE_1\vt(x)$ for any $x\in X$.

The other important thing is that 
$\ran\vt\sq Z=\ens{x\in\pnd}{\glu x=0},$ 
where $Z$ is a $\ip11$ set, hence, 
there is a $\id11$ set $Y$ with $\ran\vt\sq Y\sq Z.$ 
In particular $\vt$ reduces $\rE_1\res X$ to 
$\rE_1\res Y.$ 
We observe that $\rE_1\res Y$ is a countable \er: 
any \dd{\rE_1}class in $\pnd$ intersects $Y$ by an 
at most countable set 
(as so is the property of $Z,$ a bigger set). 
Thus, $\Ei\res Y$ is hyperfinite by Theorem~\ref{thf}. 
\vyk{
Thus $\rE_1\res Y$ is both countable and hypersmooth. 
It follows, by a theorem of \cite{djk}, that 
\imar{check \cite{djk} !!!}%
$\rE_1\res Y\reb\rE_0,$ as required.
}%
\epf

\punk{Case 2}

Since $\glu x=\iy$ is a $\is11$ formula, it suffices 
to show that for any non-empty $\is11$ set $R\sq\pnd$ 
with $\kaz x\in R\:(\glu x=\iy)$ we have a $\id11$ 
subset $X\sq R$ with $\rE_1\reb\rE_1\res X.$ 
Fix a set $R,$ as indicated, for the course of the proof.
The subset $X$ of $R$ will be defined with the help of 
a splitting construction developed in \cite{nwf} for 
the study of ``ill''founded Sacks iterations. 

We shall define a map $\vpi:\dN\to\dN,$ which assumes 
infinitely many values and assumes each its value 
infinitely many times (but $\ran\vpi$ may be a proper 
subset of $\dN$), and, for each $u\in\bse,$ a non-empty 
$\is11$ subset $X_u\sq R,$ which satisfy a quite long 
list of properties. 
First of all, if $\vpi$ is already defined at least on 
$\ir0n$ and $u\ne v\in\bse$ then let 
$
\nu_\vpi[u,v]=
\tmin_\cle\ens{\vpi(k)}{k<n\land u(k)\ne v(k)}.
$
(Note that the minimum is taken in the sense of $\cle,$ 
hence, 
it is $\tmax$ in the sense of $\le,$ the usual order). 
Separately, put $\vpi[u,u]=-1$ for any $u$.

Now we give the list of requirements.

\ben
\tenu{{\rm(\roman{enumi})}}
\itla{z1}
if $\vpi(n)\nin\ens{\vpi(k)}{k<n}$ then 
$\vpi(n)\cl \vpi(k)$ for any $k<n$;

\itla{z2}
every $X_u$ is a non-empty $\is11$ subset of $R$;

\itla{z3}
if $u\in2^n,\msur$ $x\in X_u,$ and $k<n,$ then 
$\vpi(k)\in \nab x$;

\itla{z4}
if $u,\,v\in2^n$ then 
$X_u\rc{\npi[u,v]}=X_v\rc{\npi[u,v]}$;

\itla{z5}
if $u,\,v\in2^n$ then 
$X_u\rec{\npi[u,v]}\cap X_v\rec{\npi[u,v]}=\pu$;

\itla{z6}\msur
$X_{u\we i}\sq X_u$ for all $u\in\bse$ and $i=0,1$;

\itla{z7}\msur
$\tmax_{u\in2^n}\dia X_u\to0$ as $n\to\iy$ 
(a reasonable 
Polish metric on $\pnd$ is assumed to be fixed);

\itla{z8}
a certain condition, in terms of the Choquet game, which 
connects each $X_{u\we i}$ with $X_u$ so that, as a 
consequence, $\bigcap_nX_{a\res n}\ne\pu$ 
for any $a\in\dn$. 
\een

Let us demonstrate how such a system of sets and a function 
$\vpi$ accomplish Case 2. 
According to \ref{z7} and \ref{z8}, for any $a\in\dn$ 
the intersection $\bigcap_nX_{a\res n}$ contains a single 
point, let it be $F(a),$ and $F$ is continuous and 
$1-1$.

Put $J=\ran\vpi=\ens{j_m}{m\in\dN},$
in the \dd<increasing order; $J\sq\dN$ is infinite. 
Let $n\in\dN.$ 
Then $\vpi(n)=j_m$ for some (unique) $m:$ we put 
$\psi(n)=m.$ 
Thus $\psi:\dN\onto\dN$ and the preimage 
$\psi\obr(m)=\vpi\obr(j_m)$ is an infinite subset of 
$\dN$ for any $m.$ 
This allows us to define a parallel system of 
sets $Y_u,\msur$ $u\in\bse,$ as follows. 
Put $Y_\La=\pnd.$ 
Suppose that $Y_u$ has been defined, $u\in2^n.$ 
Put $j=\vpi(n)=j_{\psi(n)}.$ 
Let $K$ be the number of all indices $k<n$ still 
satisfying $\vpi(k)=j,$ perhaps $K=0.$ 
Put $Y_{u\we i}=\ens{x\in Y_u}{x(j)(K)=i}$ for 
$i=0,1$. 

Each of $Y_u$ is clearly a basic clopen set in $\pnd,$ 
and one easily verifies that conditions 
\ref{z1} -- \ref{z7}, except for \ref{z3},  
are satisfied for the sets $Y_u$ (instead of $X_u$) and 
the map $\psi$ (instead of $\vpi$), in particular, for 
any $a\in\dn,$  $\bigcap_nY_{a\res n}=\ans{G(a)}$ 
is a singleton, and the map $G$ is continuous and $1-1.$ 
(We can, of course, define $G$ explicitly: 
$G(a)(m)(l)=a(n),$ where $n\in\dN$ is chosen so that 
$\psi(n)=m$ and there is exactly $l$ numbers $k<n$ 
with $\psi(k)=m$.) 
Note finally that $\ens{G(a)}{a\in\dn}=\pnd$ since 
by definition $Y_{u\we 1}\cup Y_{u\we 0}=Y_u$ for 
all $u$.

We conclude that the map $\vt(x)=F(G\obr(x))$ is a 
continuous bijection 
(hence, in this case, a homeomorphism by compactness) 
$\pnd\onto X.$ 
We further assert that $\vt$ satisfying the following: 
for each $y,\,y'\in \pnd$ and $m$,
\dm
y\rec m= y'\rec m\quad\hbox{iff}\quad
\vt(y)\rec{j_m}=\vt(y')\rec{j_m}
\,.
\eqno(\ast)
\dm
Indeed, let $y=G(a)$ and $x=F(a)=\vt(y),$ and 
similarly $y'=G(a')$ and $x'=F(a')=\vt(y'),$ where 
$a,\,a'\in\dn.$ 
Suppose that $y\rec m= y'\rec m.$ 
According to \ref{z5} for $\psi$ and the sets $Y_u,$ 
we then have $m\cl\nsi[a\res n,a'\res n]$ for any $n.$ 
It follows, by the definition of $\psi,$ that 
$j_m\cl\npi[a\res n,a'\res n]$ for any $n,$ hence, 
$X_{a\res n}\rec{j_m}=X_{a\res n}\rec{j_m}$ 
for any $n$ by \ref{z4}. 
Assuming now that Polish metrics on all spaces 
$\pn^{\cle j}$ are chosen so that 
$\dia Z\ge\dia {(Z\rec j)}$ for all $Z\sq\pn$ and $j,$ 
we easily obtain that $x\rec{j_m}=x'\rec{j_m},$ \ie, 
the right-hand side of $(\ast).$ 
The inverse implication in $(\ast)$ is proved similarly. 

Thus we have $(\ast),$ but this means that $\vt$ is a 
continuous reduction of $\Ei$ to $\Ei\res X,$ thus, 
$\Ei\reb {\Ei\res X},$ as required.\vtm

\qeDD{Theorem~\ref{kelu} modulo the construction\/ 
\ref{z1} -- \ref{z8}}

\punk{The construction}

Recall that $R\sq\pnd$ is a fixed non-empty $\is11$ set 
such that $\glu x=\iy$ for each $x\in R.$ 
Set $X_\La=R$.

Now suppose that the sets $X_u\sq R$ with $u\in2^n$ 
have been defined and satisfy the applicable part of 
\ref{z1} -- \ref{z8}.\vom 

{\sl Step 1\/}. 
Our 1st task is to choose $\vpi(n).$  
Let $\ans{j_1<...<j_m}=\ens{\vpi(k)}{k<n}.$ 
For any $1\le p\le m,$ let $N_p$ be the number of all 
$k<n$ with $\vpi(k)= j_p.$\vom 

{\sl Case 1a\/}.  
If some numbers $N_p$ are $<m$ then  
choose $\vpi(n)$ among $j_p$ with the least $N_p,$ and 
among them the least one.\vom

{\sl Case 1b\/}:  
$N_p\ge m$ (then actually $N_p=m$) for all $p\le m.$ 
It follows from our assumptions, in particular \ref{z4}, 
that $X_u\rc{j_m}=X_v\rc{j_m}$ for all $u,\,v\in2^n.$ 
Let $Y=X_u\rc{j_m}$ for any such $u.$ 
Take any $y\in Y.$ 
Then $\nab y$ is infinite, hence, there is some 
$j\in\nab y$ with $j\cl j_m.$ 
Put $\vpi(n)=j$. 

We have something else to do in this case. 
Let $X'_u=\ens{x\in X_u}{j\in\nab y}$ for any $u\in2^m.$ 
Then we easily have 
$X'_u=\ens{x\in X_u}{x\rc{j_m}\in Y'},$ where 
$Y'=\ens{y\in Y}{j\in\nab y}$ is a non-empty $\is11$ set, 
so that the sets $X'_u\sq X_u$ are non-empty $\is11.$ 
Moreover, as $j_m$ is the \dd\cle least in 
$\ens{\vpi(k)}{k<n},$ we can easily show that the system 
of sets $X'_u$ still satisfies \ref{z4}.
This allows us to assume, without any loss of generality, 
that, in Case~1b, $X'_u=X_u$ for all $u,$ or, 
in other words, that any $x\in X_u$ for any $u\in2^n$ 
satisfies $j=\vpi(n)\in\nab x.$ 
(This is true in Case~1a, of course, because then 
$\vpi(n)=\vpi(k)$ for some $k<n$.)\vom

Note that this manner to choose $\vpi(n)$ implies 
\ref{z1} and also implies that $\vpi$ takes 
infinitely many values and takes each its value 
infinitely many times.

The continuation of the construction requires the following

\blt
\label{suz}
If\/ $u_0\in2^n$ and\/ $X'\sq X_{u_0}$ is a non-empty\/ 
$\is11$ set then there is a system of\/ $\is11$ sets\/ 
$\pu\ne X'_u\sq X_u$ with\/ $X'_{u_0}=X',$ which still 
satisfies\/ \ref{z4}. 
\elt
\bpf
For any $u\in2^n,$ let 
$X'_u=\ens{x\in X_u}{x\rc{n(u)}\in X'\rc{n(u)}},$  
where $n(u)=\npi[u,u_0].$ 
In particular, this gives $X'_{u_0}=X',$ because 
$\npi[u_0,u_0]=-1.$ 
The sets $X'_u$ are as required, via a routine verification.
\epF{Lemma}

{\it Step 2\/}. 
First of all put $j=\vpi(n)$ and $Y_u=X_u\rc j.$  
(All $Y_u$ are equal to $Y$ in Case~1b, but the argument 
pretends to make no difference between 1a and 1b). 
Take any $u_1\in2^n.$ 
By the construction any element $x\in X_{u_1}$ satisfies 
$j\in\nab x,$ so that $x(j)\nin\id11(x\rc j).$ 
As $X_{u_1}$ is a $\is11$ set, it follows that  
$\ens{x'(j)}{x'\in X_{u_1}\land x'\rc j=x\rc j}$ 
is not a singleton, in fact is uncountable. 
It follows that there is a number $l_{u_1}$ having the 
property that the $\is11$ set 
\dm
Y'_{u_1}\:=\:
\ens{y\in Y_{u_1}}{\sus x,x'\in X_{u_1}\:
\skl x'\rc j=x\rc j=y\land 
l_{u_1}\in x(j)\land l_{u_1}\nin x'(j)\skp}
\dm
is non-empty. 
We now put $X'=\ens{x\in X_{u_1}}{x\rc j\in Y'_{u_1}}$ 
and define $\is11$ sets $\pu\ne X'_u\sq X_u$ as in the 
lemma, in particular, $X'_{u_1}=X',\msur$ 
$X'_{u_1}\rc j=Y'_{u_1},$ still \ref{z4} is satisfied, 
and in addition 
\dm
\kaz y\in X'_{u_1}\rc j\;
\sus x,x'\in X'_{u_1}\:
\skl x'\rc j=x\rc j=y\land 
l_{u_1}\in x(j)\land l_{u_1}\nin x'(j)\skp
\eqno(1)
\dm

Now take some other $u_2\in2^n.$ 
Let $\nu=\npi[u_1,u_2].$ 
If $j\cl \nu$ then $X_{u_1}\rc j=X_{u_2}\rc j,$ so 
that we already have, for $l_{u_2}=l_{u_1},$ that 
\dm
\kaz y\in X'_{u_2}\rc j\;
\sus x,x'\in X'_{u_2}\:
\skl x'\rc j=x\rc j=y\land 
l_{u_2}\in x(j)\land l_{u_2}\nin x'(j)\skp\,,
\eqno(2)
\dm
and can pass to some $u_3\in2^n.$ 
Suppose that $\nu\cle j.$ 
Now things are somewhat nastier. 
As above there is a number $l_{u_2}$ such that 
\dm
Y'_{u_2}\:=\:
\ens{y\in Y_{u_2}}{\sus x,x'\in X_{u_2}\:
\skl x'\rc j=x\rc j=y\land 
l_{u_2}\in x(j)\land l_{u_2}\nin x'(j)\skp}
\dm
is a non-empty $\is11$ set, thus, we can define 
$X''=\ens{x\in X_{u_1}}{x\rc j\in Y'_{u_1}}$ and 
maintain the construction of Lemma~\ref{suz}, getting 
non-empty $\is11$ sets $X''_u\sq X'_u$ still satisfying 
\ref{z4} and $X''_{u_2}=X'',$ therefore, 
we still have $(2)$ for the set $X''_{u_2}.$

Yet it is most important in this case that $(1)$ 
is preserved, \ie, 
it still holds for the set $X''_{u_1}$ instead of 
$X'_{u_1}$! \ 
Why is this ? \ 
Indeed, according to the construction in the proof of 
Lemma~\ref{suz}, we have 
$X''_{u_1}=\ens{x\in X'_{u_1}}{x\rc\nu\in X''\rc \nu}.$ 
Thus, although, in principle, $X''_{u_1}$ is smaller than 
$X'_{u_1},$ for any $y\in X''_{u_1}\rc j$ we have 
\dm
\ens{x\in X''_{u_1}}{x\rc j=y}\;=\;
\ens{x\in X'_{u_1}}{x\rc j=y}\,,
\dm
simply because now we assume that $\nu\cle j.$ 
This implies that $(1)$ still holds.

Iterating this construction so that each $u\in2^n$ is 
eventually encountered, we obtain, in the end, a system 
of non-empty $\is11$ sets, let us call them ``new'' 
$X_u,$ but they are subsets of the ``original'' $X_u,$ 
still satisfying \ref{z4}, still satisfying that 
$\vpi(n)\in\nab x$ for each $x\in\bigcap_{u\in2^n}X_u,$ 
and, in addition, for any $u\in2^n$ there is a number 
$l_u$ such that ${j\cl \npi[u,v]}\limp {l_u=l_v}$ and 
\dm
\kaz y\in X_{u}\rc j\;
\sus x,x'\in X_{u}\:
\skl x'\rc j=x\rc j=y\land 
l_{u}\in x(j)\land l_{u}\nin x'(j)\skp\,.
\eqno(\ast)
\dm

{\it Step 3\/}. 
We define the \dd{(n+1)}th level of sets by 
$X_{u\we 0}=\ens{x\in X_u}{l_u\in x(j)}$ and 
$X_{u\we 1}=\ens{x\in X_u}{l_u\nin x(j)}$ for 
all $u\in2^n,$ where still $j=\vpi(n).$ 
It follows from $(\ast)$ that all these $\is11$ sets 
are non-empty. 

\blt
\label{n+1}
The just defined system of sets\/ 
$X_s,\msur$ $s\in2^{n+1},$ satisfies\/ \ref{z4}, 
\ref{z5}.
\elt
\bpf 
Let $s=u\we i$ and $t=v\we i'$ belong to $2^{n+1},$ 
so that $u,\,v\in2^n$ and $i,i'\in\ans{0,1}.$ 
Let $\nu=\npi[u,v]$ and $\nu'=\npi[s,t]$.
\vom

{\it Case 3a\/}: $\nu\cle j=\vpi(n).$ 
Then easily $\nu=\nu',$ 
so that \ref{z5} immediately follows from \ref{z5} 
at level $n$ for $X_u$ and $X_v.$
As for \ref{z4}, we have $X_s\rc\nu=X_u\rc\nu$ 
(because by definition $X_s\rc j=X_u\rc j$), and 
similarly $X_t\rc\nu=X_v\rc\nu,$  
therefore, $X_t\rc{\nu'}=X_s\rc{\nu'}$ since 
$X_u\rc{\nu}=X_v\rc{\nu}$ by \ref{z4} at level $n.$ 
\vom

{\it Case 3b\/}: $j\cl\nu$ and $i=i'.$ 
Then still $\nu=\nu',$ thus we have \ref{z5}.
Further, $X_u\rc\nu=X_v\rc\nu$ by \ref{z4} at 
level $n,$ hence, $X_u\rec j=X_v\rec j,$ hence, 
$l_u=l_v$ (see above). 
Now, assuming that, say, $i=i'=1$ and $l_u=l_v=l,$ 
we conclude that
\dm
X_s\rc{\nu'}=\ens{y\in X_u\rc\nu}{l\in y(j)}
=\ens{y\in X_v\rc\nu}{l\in y(j)}=X_t\rc{\nu'}\,.
\dm

{\it Case 3c\/}: $j\cl\nu$ and $i\ne i',$ say, 
$i=0$ and $i'=1.$ 
Now $\nu'=j.$ 
Yet by definition $X_s\rc j=X_u\rc j$ and 
$X_t\rc j=X_v\rc j,$ so it remains to apply \ref{z4} 
for level $n.$ 
As for \ref{z5}, 
note that by definition 
$l\nin x(j)$ for any $x\in X_s=X_{u\we 0}$ while 
$l\in x(j)$ for any $x\in X_t=X_{v\we 1},$ where 
$l=l_u=l_v$.
\epF{Lemma}

{\it Step 4\/}.
In addition to \ref{z4} and \ref{z5}, we already have 
\ref{z1}, \ref{z2}, \ref{z3}, \ref{z6} at level $n+1.$ 
To achieve the remaining properties \ref{z7} and \ref{z8}, 
it suffices to consider, one by one, all elements 
$s\in2^{n+1},$ 
finding, at each such a substep, a non-empty $\is11$ 
subset of $X_s$ which is consistent with the requirements 
of \ref{z7} and \ref{z8} 
(for instance, for \ref{z7}, just take it so the diameter 
is $\le 2^{-n}$), 
and then reducing all other sets $X_t$ by Lemma~\ref{suz} 
at level $n+1$.
\vtm

\epF{Construction and Theorem~\ref{kelu}}

\punk{Above $\Ei$}
\label{abovei1}

Recall that an embedding is a $1-1$ reduction, and 
an invariant embedding is an embedding $\vt$ such that 
its range is an invariant set, see Subsection~\ref{reli} 
above.

\bte[{{\rm Kechris and Louveau \cite{hypsm}}}]
\label{abi1}
Suppose that\/ $\Ei\reb\rF,$ where\/ $\rF$ is an 
analytic \er\ on a Polish space\/ $\dY.$ 
Then both\/ $\Ei\emn\rF$ 
and\/ $\Ei\embi\rF$.
\ete
\bpf
To prove the first statement, 
let $\cle$ be the inverted order on $\dN,$ \ie, 
$m\cle n$ iff $n\le m.$ 
Let $\gP$ be the collection of all sets $P\sq\pn^\dN$ 
such that there is a continuous $1-1$ map 
$\eta:\pn^\dN\onto P$ such that we have 
\dm
{x\rec{n}=y\rec{n}}\leqv
{\eta(x)\rec{n}=\eta(y)\rec{n}}
\dm
for all $n$ and $x,\,y\in\pn^\dN,$ where 
$x\rec{n}=\sis{x_i}{i\cle n}$ for any 
$x=\sis{x_i}{}\in\pn^\dN.$ 
Clearly any such a map is a continuous embedding of 
$\Ei$ into itself.

This set $\gP$ is a forcing notion to extend the 
universe by a sequence of reals $x_i$ so that each 
$x_n$ is Sacks--generic over $\sis{x_i}{i\curle n},$ 
an example of iterated Sacks extensions with an 
ill-founded ``skeleton'' of iteration, which we 
defined in~\cite{nwf}. 
Here, the ``skeleton'' is $\dN$ with the 
inverted order $\cle$.

The method of \cite{nwf} contains a study of 
continuous and Borel functions on 
sets in $\gP.$ 
In particular it is shown there that Borel maps 
admit the following {\it cofinal classification\/} 
on sets in $\gP:$ 
if $\dY$ is Polish, $P'\in\gP,$ and 
$\vt:P'\to\dY$ is Borel then there is a set 
$P\in\gP,\msur$ $P\sq P',$ on which $\vt$ is 
continuous, and either a constant or, for some $n,$  
$1-1$ on $P\rec{n}$ in the sense that,   
\dm
\hbox{for all }\;x,\,y\in P:\quad
{x\rec{n}=y\rec{n}}\leqv{\vt(x)=\vt(y)}\,.
\eqno(\ast)
\dm
We apply this to a Borel map 
$\vt:\pn^\dN\to\dY$ which reduces $\Ei$ to 
$\rF.$ 
We begin with $P'=\pn^\dN$ and find a set $P\in\gP$ 
as indicated. 
Since $\vt$ cannot be a constant on $P$ 
(indeed, any $P\in\gP$ contains many pairwise 
\dd\Ei inequivalent elements), 
we have $(\ast)$ for some $n.$ 
In other words, there is a $1-1$ continuous map 
$f:P\rec n\to \dY$ 
(where ${P\rec n}=\ens{x\rec n}{x\in P}$) 
such that $\vt(x)=f(x\rec n)$ for all $x\in P.$ 
Now, let $x=\sis{x_i}{i\in\dN}\in\pn^\dN.$ 
Define $\zeta(x)=z=\sis{z_i}{i\in\dN}$ so that 
$z_i=\pu$ for $i<n$ and $z_{n+i}=x_i$ for all $i.$ 
Finally set $\vt'(x)=f(\eta(\zeta(x))\resic n)$ 
for all $x\in\pn^\dN:$ this is a continuous embedding 
of $\Ei$ in $\rF$. 

Now we prove the second claim. 
We can assume that $\dY=\pn$ and that $\vt:\pnd\to\pn$ 
is already a continuous embedding $\Ei$ into $\rF.$ 
Let $Y=\ran\vt$ and $Z=[Y]_{\rF}.$ 
Normally $Y,\,Z$ are analytic, but in this case they are 
even Borel. 
Indeed $Z$ is the projection of 
$P=\ens{\ang{z,x}}{z\rF\vt(x)},$ a Borel subset of 
$\pn\ti\pnd$ whose all cross-sections are 
\dd\Ei equivalence classes, \ie, \dds compact sets. 
It is known 
\imar{reference}%
that in this case $Z$ is Borel and, moreover, there is 
a Borel map $f:Z\to\pnd$ such that $f(z)\Ei x$ 
whenever $z\rF\vt(x)$.

We can convert $f$ to a $1-1$ map $g:\pn\to\pnd$ 
with the same properties: 
$g(z)_n=f(z)_n$ for $n\ge 1,$ but 
$g(z)_0=z.$ 
Then $f:\pnd\to Z\sq\pn$ and $g:Z\to\pnd$ 
are Borel $1-1$ maps 
($\vt$ is even continuous, but this does not matter now), 
and, for any $x\in\pnd,$ $\vt$ maps $[x]_{\Ei}$ into 
$[\vt(x)]_{\rF}\sq Z,$ and $g$ maps $[\vt(x)]_{\rF}$ 
back into $[x]_{\Ei}.$ 
It remains to apply the construction from the Cantor -- 
Bendixson theorem, to get a Borel embedding, say, $F$ 
of $\Ei$ into $\rF$ with $\ran F=Z,$ \ie, an 
invariant embedding.
\epf

The following theorem shows that orbit \eqr s of  
Polish group actions cannot reduce $\Ei$. 
 
\bte[{{\rm Kechris and Louveau \cite{hypsm}}}]
\label{e1pga}
Suppose that\/ $\dG$ is a Polish group and\/ $\dX$ 
is a Borel\/ \dd\dG space. 
Then\/ $\Ei$ is \poq{not} Borel reducible to\/ $\ergx$.
\ete
\bpf
Towards the contrary, let $\vt:\pnd\to\dX$ 
be a Borel reduction of $\Ei$ to $\rE.$ 
We can assume, by Theorem~\ref{abi1}, that $\vt$ is 
in fact an invariant embedding, \ie, $1-1$ and 
$Y=\ran\vt$ is an \dde invariant set. 
Define, for $g\in\dG$ and $x\in\pnd,$ 
$g\ac x=\vt\obr(g\ac\vt(x)).$ 
Then this is a Borel action of $\dG$ on $\pnd$ such 
that the induced relation $\aer\dG\pnd$ 
coincides with $\Ei$. 

Let us fix $x\in\pnd.$ 

Consider any $y=\sis{y_n}n\in[x]_{\Ei}.$ 
Then $[x]_{\Ei}=\bigcup_nC_n(y),$ where each set 
$C_n(y)=\ens{y'\in\pnd}{\kaz m\ge n\:(y_n=y'_n)}$  
is Borel (even compact). 
It follows that $\dG=\bigcup_nG_n(y),$ where each 
$G_n(y)=\ens{g\in\dG}{g(x)\in C_n(y)}$ is Borel. 
Thus, as $\dG$ is Polish, there is a number 
$n$ such that $G_n(y)$ is not meager in $\dG$ 
(then this will hold for all $n'\ge n,$ of course). 
Let $n(y)$ be the least such an $n$. 

We assert that for any $n$ the set 
$Y_n(x)=\ens{y\res{\iry n}}{y\in[x]_{\Ei}\land n(x)=n}$ 
is at most countable. 
Indeed suppose that $Y_n(x)$ is not countable. 
Note that if $y_1$ and $y_2$ in $[x]_{\Ei}$ have 
different restrictions $y_i\res{\iry n}$ then the sets 
$C_n(y_1)$ and $C_n(y_2)$ are disjoint, therefore, 
the sets $G_n(y_1)$ and $G_n(y_2)$ are disjoint, 
so we would have uncountably many pairwise disjoint 
non-meager sets in $\dG,$ contradiction. 
Thus all sets $Y_n(x)$ are countable. 

It is most important that $Y_n(x)$ depends on $[x]_{\Ei}$ 
rather than $x$ itself, more exactly, if $x'\in[x]_{\Ei}$ 
then $Y_n(x)=Y_n(x'):$ this is because any set $G_n(y)$ 
in the sense of $x'$ is just a shift, within $\dG,$ of 
$G_n(y)$ in the sense of $x.$ 
Therefore, putting  
$Y(x)=\bigcup_n\ens{\bar u}{u\in Y_n(x)},$ 
where, for $u\in\pn^{\iry n},$ $\bar u\in\pnd$ is 
defined by $\bar u\res{\iry n}=u$ and $\bar u(k)=\pu$ 
for $k<n,$ we have the set $Y=\bigcup_{x\in\pnd}Y(x)$ 
with the property that $Y\cap[x]_{\Ei}$ is non-empty 
and at most countable for any $x\in\pnd$. 

The other important fact is that the relation $y\in Y(x)$ 
is Borel: this is because it is assembled from Borel 
relations via the Vaught quantifier ``there exists 
nonmeager-many'', known to preserve the Borelness. 
\imar{reference}%
It follows that
\dm
Y=\ens{y}{\sus x\:(y\in Y_x)}=
\ens{y}{\kaz x\:(x\in[y]_{\Ei}\imp y\in Y(x)}
\dm
is a Borel subset of $\pnd.$ 
By the uniformization theorem for 
\imar{reference}%
Borel sets with countable sections, there is a Borel 
map $f$ defined on $\pnd$ so that $f(x)\in Y(x)$ for 
any $x,$ which implies $\Ei\reb{\Ei\res Y}.$ 
On the other hand, $\Ei\res Y$ is a countable \er\  
by the above, which is a contradiction to Lemma~\ref{eoei'}.
\epf

\parf{Actions of the infinite symmetric group}
\label{groas}

This Section is connected with the next one (on turbulence). 
We concentrate on a main result in this area, due to 
Hjorth, that turbulent \er s are not reducible to 
those induced by actions of $\isg.$ 
In particular, we shall prove the following:

\ben
\tenu{\Roman{enumi}}
\def\labelenumi{\theenumi.}
\itla{isg1}
Lopez-Escobar: any invariant Borel set of countable models 
is the truth domain of a formula of $\lww$.

\itla{isg2}
Any orbit \er\ of a Polish action of a closed subgroup of 
$\isg$ is classifiable by countable structures (up to isomorphism). 

\itla{isg3}
Any \er, classifiable by countable structures, is Borel reducible 
to isomorphism of countable ordered graphs.

\itla{isg4}
Any \poq{Borel} \er, classifiable by countable structures, 
is Borel reducible to one of \er s $\rT_\xi$.

\itla{isg5}
Any \er, classifiable by countable structures and induced by a 
Polish action (of a Polish group), 
is Borel reducible to one of \er s $\rT_\xi$ on a comeager set.

\itla{isg6}
Any ``turbulent'' \er\ $\rE$ is generically 
\dd{\rT_\xi}ergodic for any $\xi<\omi,$ in particular, $\rE$ 
is not Borel reducible to $\rT_\xi$.

\itla{isg7}
Any ``turbulent'' \er\ is not classifiable by countable 
structures: a corollary of \ref{isg6} and \ref{isg5}. 

\itla{isg8}
A generalization of \ref{isg7}: 
any ``turbulent'' \er\ is not Borel reducible to a \er\ which 
can be obtained from $\rav\dN$ using operations 
defined in \prf{opeer}.
\een

Scott's analysis, involved in proofs of \ref{isg4} and \ref{isg5}, 
appears only in a rather mild and self-contained version.

\punk{Infinite symmetric group $\isg$}
\label{loac}

Let $\isg$ be the group of all permutations 
(\ie, 1--1 maps $\dN\onto\dN$) of $\dN,$ with the 
superposition as the group operation. 
Clearly $\isg$ is a $\Gd$ subset of $\dnn,$ hence, a 
Polish group. 
A compatible complete metric on $\isg$ can be defined by 
$D(x,y)=d(x,y)+d(x\obr,y\obr),$ where $d$ is the ordinary 
complete metric of $\dnn,$ \ie, $d(x,y)=2^{-m-1},$ where 
$m$ is the least such that $x(m)\ne y(m).$ 
Yet $\isg$ admits no compatible left-invariant 
complete metric \cite[1.5]{beke}. 
\imar{Proof of $\isg$ not \cli~?}%

For instance isomorphism relations of various kinds of 
countable structures are orbit \er s induced by $\isg.$ 
Indeed, suppose that $\cL=\sis{R_i}{i\in I}$ is a countable 
\index{language!countable relational}%
relational language, \ie, $0<\card I\le\alo$ and each $R_i$ 
is an \dd{m_i}ary relational symbol. 
We put~\footnote
{\ $X_\cL$ is often used to denote $\mox\cL$.}   
$\mox\cL=\prod_{i\in I}\cP(\dN^{m_i}),$ 
\index{zzModL@$\mox\cL$}%
\index{action!logic@logic action $\loa\cL$ of $\isg$}%
\index{zzjL@$\loa\cL$}%
the space of (coded) \dd\cL{\it structures\/} on $\dN.$ 
The {\it logic action\/} $\loa\cL$ of $\isg$ on $\mox\cL$ 
is defined as follows: 
if $x=\sis{x_i}{i\in I}\in\mox\cL$ and $g\in\isg$ 
then $y=\loa\cL(g,x)=g\app x=\sis{y_i}{i\in I}\in\mox\cL,$ 
where we have 
\dm
{\ang{k_1,...,k_{m_i}}\in x_i} \leqv 
{\ang{g(k_1),...,g(k_{m_i})}\in y_i}  
\dm
for all $i\in I$ and $\ang{k_1,...,k_{m_i}}\in\dN^{m_i}.$ 
Then $\stk{\mox\cL}{\loa\cL}$ is a Polish \dd\isg space  
and \dd{\loa\cL}orbits in $\mox\cL$ are exactly the 
isomorphism classes of \dd\cL structures, which is a reason 
to denote the associated equivalence relation 
\index{isomorphism!of st@of structures, $\ism\cL$}%
\index{zzconl@$\ism\cL$}%
$\aer{\loa\cL}{\mox\cL}$ as ${\ism\cL}$. 

If $G$ is a subgroup of $\isg$ then $\loa\cL$ restricted to 
$G$ is still an action of $G$ on $\mox\cL,$ whose orbit \er\ 
\index{isomorphism!of st@induced by $G,$ $\izm\cL G$}%
\index{zzconlg@$\izm\cL G$}%
will be denoted by $\izm\cL G,$ \ie, $x\izm\cL G y$ iff 
${\sus g\in G\:(g\app x=y)}$.

\punk{Borel invariant sets}
\label{BIS}

A set $M\sq\mox\cL$ is {\it invariant\/} if 
\index{set!invariant}%
$\ek M{\ism\cL}=M.$ 
There is a convenient characterization of {\it Borel\/} 
invariant sets, in terms of $\lww,$ an infinitary extension 
of $\cL=\sis{R_i}{i\in I}$ 
\index{language!Lww@$\lww$}%
\index{zzLww@$\lww$}%
by countable conjunctions and disjunctions. 
To be more exact, 
\ben
\tenu{\arabic{enumi})}
\item 
any $R_i(v_0,...,v_{m_i-1})$ is an atomic 
formula of $\lww$ (all $v_i$ being variables over $\dN$ 
and $m_i$ is the arity of $R_i$), and 
propositional connectives and quantifiers ${\sus},\,{\kaz}$ 
can be applied as usual; 

\item
if $\vpi_i,\;i\in\dN,$ are formulas of $\lww$ whose free 
variables 
are among a finite list $v_0,...,v_n$ then $\bigvee_i\vpi_i$ 
and $\bigwedge_i\vpi_i$ are formulas of $\lww$.
\een
If $x\in\mox\cL,\msur$ 
$\vpi(v_1,...,v_n)$ is a formula of $\lww,$ and 
$i_1,...,i_n\in\dN,$ then $x\mo\vpi(i_1,...,i_n)$ means that 
$\vpi(i_1,...,i_n)$ is satisfied on $x,$ in the 
usual sense that involves transfinite induction on the ``depth'' 
of $\vpi,$ see \cite[16.C]{dst}.

\bte[{{\rm Lopez-Escobar, see \cite[16.8]{dst}}}] 
\label{lopes}
A set\/$M\sq\mox\cL$ is invariant and Borel iff\/ 
$M=\ens{x\in\mox\cL}{x\mo\vpi}$ for a 
closed formula\/ $\vpi$ of\/ $\lww$.
\ete
\bpf
To prove the nontrivial direction let $M\sq\mox\cL$ 
be invariant and Borel. 
Put $B_s=\ens{g\in\isg}{s\su g}$ for any injective 
$s\in\dN\lom$ (\ie, $s_i\ne s_j$ for $i\ne j$), this is a 
clopen subset of $\isg$ 
(in the Polish topology of $\isg$ inherited from $\dnn$). 
If $A\sq\isg$ then let $s\forc A(\dog)$ mean that the set 
$B_s\cap A$ is co-meager in $B_s,$ \ie, $g\in A$ holds for \pv\ 
$g\in \isg$ with $s\su g.$ 
The proof consists of two parts:

\ben
\tenu{(\roman{enumi})}
\itla{pari}\msur
$M=\ens{x\in\mox\cL}{\La\forc \dog\app x\in M}$ 
(where $g\app x=\loa\cL(g,x),$ see above);

\itla{pard}
For any Borel $M\sq\mox\cL$ and any $n$ there is a formula 
$\vpi_{M}^n(v_0,...,v_{n-1})$ of $\lww$ such that we have, for 
every $x\in\mox\cL$ and every injective $s\in\dN^n:$  
$x\mo\vpi_{M}^n(s_0,...,s_{n-1})$ iff 
$s\forc\dog{}\obr\app x\in M$.
\een

\ref{pari} is clear: since $M$ is invariant, we have 
$g\app x\in M$ for all $x\in M$ and $g\in\isg,$ on the other 
hand, if $g\app x\in M$ for at least one $g\in\isg$ then 
$x\in M$. 

To prove \ref{pard} we argue by induction on the Borel 
complexity of $M.$ 
Suppose, for the sake of simplicity, that $\cL$ contains 
a single binary predicate, say, $R(\cdot,\cdot);$ then 
$\mox\cL=\cP(\dN^2).$ 
If $M=\ens{x\sq\dN^2}{\ang{k,l}\nin x}$ for some $k,\,l\in\dN$ 
then take
\dm
\TS
\kaz u_0\,...\,\kaz u_m\;
\skl
\bigwedge_{i<j\le m}(u_i\ne u_j)\,\land\,
\bigwedge_{i<n}(u_i=v_i)\,\imp\,\neg\:R(u_k,u_l)
\skp\,,
\dm
where $m=\tmax\ans{l,k,n},$ as $\vpi_{M}^n(v_0,...,v_{n-1}).$ 
Further, take 
\dm
\TS
\bay{rl}
\bigwedge_{k\ge n} 
\,\kaz u_0\,...\,\kaz u_{k-1}\;
\bigvee_{m\ge k}
\,\sus w_0\,...\,\sus w_{m-1}\;
\skl
\bigwedge_{i<j<k}(u_i\ne u_j)\,\land\,
\bigwedge_{i<n}(u_i=v_i)&\\[\dxii]

\limp 
\bigwedge_{i<j<m}(w_i\ne w_j)\,\land\,
\bigwedge_{i<k}(w_i=v_i)\,\land\,\vpi^m_{M}(w_0,...,w_{m-1})
\skp
\eay
\dm
as $\vpi_{\neg\,M}^n(v_0,...,v_{n-1}).$ 
Finally, if $M=\bigcap_jM_j$ then we take 
$\bigwedge_j\vpi^n_{M_j}(v_0,...,v_{n-1})$ 
as $\vpi_M^n(v_0,...,v_{n-1})$. 
\epF{Theorem~\ref{lopes}}

\punk{\er s classifiable by countable structures}
\label{erclass}

The classifiability by countable structures means that 
we can associate, in a Borel way, a countable \dd\cL structure, 
say, $\vt(x)$ with any point $x\in\dX=\dom\rE$ so that $x\rE y$ 
iff $\vt(x)$ and $\vt(y)$ are isomorphic. 

\bdf[{{\rm Hjorth~\cite[2.38]h}}]
\label{classif}
An \er\ $\rE$ is {\it classifiable by countable structures\/} 
\index{equivalence relation, ER!classifiable by countable structures}%
if there is a countable relational language $\cL$ such that 
$\rE\reb{\ism\cL}$.
\edf

\bre
\label{classifex}
Any $\rE$ classifiable by countable structures is $\fs11,$ 
of course, and many of them are Borel. 
The \eqr s $\rtd\yt{\Et},$ all countable Borel \er s
(see the diagram on page \pageref{p-p})
are classifiable by countable structures, but 
$\Ei\yt\Ed,$ Tsirelson \er s are not.
\ere

\bte[{{\rm Becker and Kechris \cite{beke}}}]  
\label{siy2cc}
Any orbit \er\ of a Polish action of a closed subgroup of 
$\isg$ is classifiable by countable structures. 
\ete

Thus all orbit \er s of Polish actions of $\isg$ and its 
closed subgroups are Borel reducible to a very special kind 
of actions of $\isg.$ 

\bpf
First show that any orbit \er\ of a Polish action of $\isg$ 
itself is classifiable by countable structures.  
Hjorth's simplified argument~\cite[6.19]h is as follows. 
Let $\dX$ be a Polish \dd\isg space with basis 
$\sis{U_l}{l\in\dN},$ and let $\cL$ be the language with 
relations $R_{lk}$ where each $R_{lk}$ has arity $k.$ 
If $x\in \dX$ then define $\vt(x)\in\mox\cL$ by 
stipulation that $\vt(x)\mo R_{lk}(s_0,...,s_{k-1})$ 
iff $1)\msur$ $s_i\ne s_j$ whenever $i<j<k,$ and 
$2)\msur$ $\kaz g\in B_s\:(g\obr\app x\in U_l),$ where 
\imar{Hjorth requires $\in\nad{U_l}.$ Why~? 
Also, it seems that $\forall^\ast g\in B_s$ extends the proof to 
Borel actions.}%
$B_s=\ens{g\in\isg}{s\su g}$ and 
$s=\ang{s_0,...,s_{k-1}}\in\dN^k.$  
Then $\vt$ reduces $\aer\isg\dX$ to $\ism\cL.$ 

To accomplish the proof of the theorem, it remains to 
apply the following result 
(an immediate corollary of Theorem~2.3.5b in \cite{beke}):

\bprt
\label{subg2g}
If\/ $\dG$ is a closed subgroup of a Polish group\/ $\dH$ 
and\/ $\dX$ is a Polish\/ \dd\dG space then there is a  
Polish\/ \dd\dH space\/ $\dY$ such that\/ 
$\ergx\reb{\aer\dH\dY}$.
\eprt
\bpf
Hjorth~\cite[7.18]h outlines a proof as follows. 
Let $Y=\dX\ti\dH\,;$ define $\ang{x,h}\approx\ang{x',h'}$ if 
$x'=g\app x$ and $h'=gh$ for some $g\in\dG,$ and consider the 
quotient space $\dY=Y/{\approx}$ with the topology induced by 
the Polish topology of $Y$ via the surjection 
$\ang{x,h}\mapsto\ek{\ang{x,h}}{\approx},$ on which $\dH$ 
acts by 
$h'\app\ek{\ang{x,h}}{\approx}=\ek{\ang{x,{hh'}\obr}}{\approx}.$ 
Obviously $\ergx\reb{\aer\dH\dY}$ via the map 
$x\mapsto\ek{\ang{x,1}}{\approx},$ hence, it remains to prove 
that $\dY$ is a Polish \dd\dH space, which is not really 
elementary --- we 
refer the reader to \cite[7.18]h or \cite[2.3.5b]{beke}. 
\epF{Proposition}

To bypass \ref{subg2g} in the proof of Theorem~\ref{siy2cc}, 
we can use a characterization of all closed subgroups of $\isg.$ 
Let $\cL$ be a language as above, and $x\in\mox\cL.$  
Define $\aut x=\ens{g\in\isg}{g\app x=x}:$  
the group of all automorphisms of $x.$ 

\bprt[{{\rm see \cite[1.5]{beke}}}] 
\label{closubg}
$G\sq\isg$ is a closed subgroup of\/ $\isg$ iff there is 
an\/ \dd\cL structure\/ $x\in\mox\cL$ of a countable 
language\/ $\cL,$ such that\/ $G=\aut x$. 
\eprt 
\bpf
For the nontrivial direction, 
let $G$ be a closed subgroup of $\isg.$ 
For any $n\ge 1,$ let $I_n$ be the set of all \dd Gorbits in 
$\dN^n,$ \ie, equivalence classes of the \er\ 
$s\sim t$ iff $\sus g\in G\:(t=g\circ s),$ thus, $I_n$ is an 
at most countable subset of $\cP(\dN^n).$ 
Let $I=\bigcup_nI_n,$ and, for any $i\in I_n,$ let 
$R_i$ be an \dd nary relational symbol, and 
$\cL=\sis{R_i}{i\in I}.$  
Let $x\in\mox\cL$ be defined as follows: 
if $i\in I_n$ then $x\mo R_i(k_0,...,k_{n-1})$ iff 
$\ang{k_0,...,k_{n-1}}\in i.$ 
Then $G=\aut x,$ actually, if $G$ is not necessarily closed 
subgroup then $\aut x=\nad G$.
\epF{Proposition}

Now come back to Theorem~\ref{siy2cc}. 
The same argument as in the beginning of the proof shows that 
any orbit \er\ of a Polish action of $G,$ a closed 
subgroup of $\isg,$ is $\reb{\izm\cL G}$ for an appropriate 
countable language $\cL.$ 
Yet, by \ref{closubg}, $G=\aut{y_0}$ where $y_0\in\mox{\cL'}$ 
and $\cL'$ is a countable language disjoint from $\cL.$ 
The map $x\longmapsto\ang{x,y_0}$ witnesses that 
${\izm\cL G}\reb{\ism{\cL\cup\cL'}}$.
 
\epF{Theorem~\ref{siy2cc}}

\punk{Reduction to countable graphs}
\label{gra}

It could be expected that the more complicated a language 
$\cL$ is accordingly the more complicated
isomorphism \eqr\ $\ism\cL$ it 
produces. 
However this is not the case. 
Let $\cG$ be the language of (oriented binary) graphs, \ie, 
\index{language!g@$\cG$ of graphs}%
$\cG$ contains a single binary predicate, say $R(\cdot,\cdot)$.

\bte
\label{grafs}
If\/ $\cL$ is a countable relational language then\/ 
${\ism\cL}\reb{\ism\cG}.$ 
Therefore, an \er\/ $\rE$ is classifiable by countable 
structures iff\/ $\rE\reb{\ism\cG}.$
In other words, a 
single binary relation can code structures of any countable 
language. 
\ete  

Becker and Kechris \cite[6.1.4]{beke} outline a proof based 
on coding in terms of lattices, unlike the following argument, 
yet it may in fact involve the same idea.

\bpf
Let $\hfn$ be the set of all hereditarily finite sets over the 
set $\dN$ considered as the set of atoms, and $\ve$ be the 
associated ``membership'' (any $n\in\dN$ has no 
\dd\ve elements, $\ans{0,1}$ is different from $2,$ \etc\/.). 
Let $\ihf$ be the $\hfn$ version of $\ism\cG,$ \ie, if 
$P,\,Q\sq\hfn^2$ then $P\ihf Q$ means that there is a 
bijection $b$ of $\hfn$ such that 
$Q=b\app P=\ens{\ang{b(s),b(t)}}{\ang{s,t}\in P}.$ 
Obviously $(\ism\cG)\eqb(\ihf),$ thus, we have to prove that 
${\ism\cL}\reb{\ihf}$ for any $\cL$.

An action of $\isg$ on $\hfn$ is defined as follows. 
If $g\in\isg$ then $g\circ n=g(n)$ for any $n\in\dN,$ 
and, by \dd\ve induction,  
$g\circ\ans{a_1,...,a_n}=\ans{g\circ a_1,...,g\circ a_n}$ 
for all $a_1,...,a_n\in\hfn.$
Clearly the map $a\mto g\circ a$ ($a\in\hfn$) 
is an \dd\ve isomorphism of $\hfn,$ for any fixed $g\in\isg.$ 

\blt
\label{grafs1}
Suppose that\/ $X,\,Y\sq\hfn$ are\/ \dd\ve transitive subsets 
of\/ $\hfn,$ 
the sets\/ $\dN\dif X$ and\/ $\dN\dif Y$ are infinite, and\/ 
${\ve\res X}\ihf{\ve\res Y}.$ 
Then there is\/ $f\in\isg$ such that\/ 
$Y=f\circ X=\ens{f\circ s}{s\in X}$.
\elt 
\bpf
It follows from the assumption ${\ve\res X}\ism\hfn{\ve\res Y}$ 
that there is an \dd\ve iso\-mor\-phism $\pi:X\onto Y.$  
Easily $\pi\res(X\cap\dN)$ is a bijection of $X_0=X\cap\dN$ 
onto $Y_0=Y\cap\dN,$ hence, there is $f\in\isg$ such that 
$f\res X_0=\pi\res X_0,$ and then we have $f\circ s=\pi(s)$ 
for any $s\in X$.
\epF{Lemma}

Coming back to the proof of Theorem~\ref{grafs}, we first show 
that ${\ism{\cG(m)}}\reb{\ihf}$ for any $m\ge 3,$ where 
$\cG(m)$ is the language with a single \dd mary predicate.  
Note that $\ang{i_1,...,i_m}\in\hfn$ whenever  
$i_1,...,i_m\in\dN.$ 

Put $\vT(x)=\ens{\vt(s)}{s\in x}$ for every element
$x\in\mox{\cG(m)}=\cP(\dN^m),$ where 
$\vt(s)=\tce{\ans{\ang{2i_1,...,2i_m}}}$ for each 
$s=\ang{i_1,...,i_m}\in \dN^m,$ and finally, for $X\sq\hfn,$ 
$\tce X$ is the least \dd\ve transitive set $T\sq\hfn$ with 
$X\sq T.$  
It easily follows from Lemma~\ref{grafs} that 
$x\ism{\cG(m)}y$ iff 
${\ve\res\vT(x)}\ihf{\ve\res\vT(y)}.$ 
This ends the 
proof of ${\ism{\cG(m)}}\reb{\ihf}$. 

It remains to show that ${\ism{\cL'}}\reb{\ihf},$ where 
$\cL'$ is the language with infinitely many binary 
predicates. 
In this case $\mox{\cL'}=\cP(\dN^2)^\dN,$ so that we can 
assume that every $x\in\mox{\cL'}$ has the form 
$x=\sis{x_n}{n\ge1},$ with 
$x_n\sq(\dN\dif\ans0)^2$ for all $n.$ 
Let $\vT(x)=\ens{s_n(k,l)}{n\ge1\land\ang{k,l}\in x_n}$ for 
any such $x,$ where 
\dm
s_n(k,l)=\tce{
\ans{\{...\{\ang{k,l}\}...\}\,,\,0}}
\,,\,\text{ with }\,n+2\,\text{ pairs of brackets }\,
\{\,,\,\}\,.
\dm
Then $\vT$ is a continuous reduction of $\ism{\cL'}$ to $\ihf$.
\epF{Theorem}

\punk{Borel countably classified \er s: reduction to $\rT_\xi$}
\label{2t}

Equivalence relations $\rT_\xi$ of \prf{opeer} offer a 
perfect calibration tool for those Borel \er s which admit 
classification by countable structures. 
First of all, 

\bpro
\label{taccs}
Every\/ $\rT_\xi$ admits classification by 
countable structures.
\epro
\bpf
$\rT_0,$ the equality on $\dN,$ is the orbit \er\ of the 
action of $\isg$ by $g\app x=x$ for all $g,\,x.$ 
The operation \ref{cdun} of \prf{opeer} 
(countable disjoint union) 
easily preserves the property of being Borel 
reducible to an orbit \er\ of continuous action of $\isg.$ 

Now consider operation \ref{cp} of countable power. 
Suppose that a \er\ $\rE$ on a Polish space $\dX$ 
is Borel reducible to $\rF,$ the orbit relation of 
a continuous action of $\isg$ on some Polish $\dY.$ 
Let $D$ be the set of all points 
$x=\sis{x_k}{k\in\dN}\in\dX^\dN$ such that either 
$x_k\nE x_l$ whenewer $k\ne l,$ or there is $m$ such 
that $x_k\rE x_l$ iff $m$ divides $|k-l|.$
Then $\rE^\iy\reb{(\rE^\iy\res D)}$ 
(via a Borel map $\vt:\dX^\dN\to D$ such that 
$x\rE^\iy \vt(x)$ for all $x$). 
On the other hand, obviously $(\rE^\iy\res D)\reb \rF',$ 
where, for $y,\,y'\in\dY^\dN,$ $y\rF'y'$ means that there 
is $f\in\isg$ such that $y_k\rF y'_{f(k)}$ for all $k.$ 
Finally, $\rF'$ is the orbit \er\ of a 
continuous action of $\isg\ti\isg^\dN,$ which can be realized 
as a closed subgroup of $\isg,$ so it remains to apply 
Theorem~\ref{subg2g}. 
\epf  

The relations $\rT_\al$ are known in different versions, 
which reflect 
the same idea of coding sets of \dd\al th cumulative level 
over $\dN,$ as, \eg, in \cite[\pff1]{sinf}, where results 
similar to Proposition~\ref{taccs} are obtained in much 
more precise form.  

\bte
\label{2tal}
If\/ $\rE$ is a Borel \er\ classifiable by countable 
structures then\/ $\rE\reb\rT_\xi$ for some\/ $\xi<\omi$. 
\ete
\bpf
The proof (a version of the proof in \cite{frid}) is based on 
Scott's analysis. 
Define, by induction on $\al<\omi,$ a family of Borel \er s
$\rrO\al$ on $\dN\lom\ti\cP(\dN^2):$ 
\bit
\item[$\mtho\ast$]\msur
$\rrq\al stAB$ means $\rro\al stAB$;
\eit  
thus, all $\rrt\al st$ ($s,\,t\in\dN\lom$) 
are binary relations on $\cP(\dN^2),$ and among them all 
relations $\rrt\al ss$ are \er s;

\bit
\item\msur
$\rrq0stAB$ \ iff \ 
$A(s_i,s_j)\eqv B(t_i,t_j)$ for all $i,\,j<\lh s=\lh t$;

\item\msur
$\rrq{\al+1}stAB$ \ iff \ 
$\kaz k\:\sus l\:(\rrQ\al{s\we k}{t\we l}AB)$ 
and $\kaz l\:\sus k\:(\rrQ\al{s\we k}{t\we l}AB)$;

\item 
if $\la<\omi$ is limit then: \ 
$\rrq\la stAB$ \ iff \ $\rrq\al stAB$ for all $\al<\la$.
\eit
Easily ${\rrO\ba}\sq{\rrO\al}$ whenever $\al<\ba$. 
%

Recall that, for $A,\,B\sq\dN^2,$ $A\ism\cG B$ means that there 
is $f\in\isg$ with $A(k,l)\eqv B(f(k),f(l))$ for all $k,\,l.$  
Then we have ${\ism\cG}\sq\bigcap_{\al<\omi}{\rrt\al\La\La}$ 
by induction on $\al$ 
(in fact $=$ rather than $\sq,$ see below), 
where $\La$ is the empty sequence. 
Call a set $P\sq \cP(\dN^2)\ti\cP(\dN^2)$ {\it unbounded\/} if 
$P\cap{\rrt\al\La\La}\ne\pu$ for all $\al<\omi$. 

\blt
\label{2tal1}
Any unbounded\/ $\fs11$ set\/ $P$ contains\/ 
$\ang{A,B}\in P$ with\/ $A\ism\cG B.$
\elt

It follows that $A\ism\cG B$ iff $\rrq\al\La\La AB$ for all 
$\al<\omi$ (take $P=\ans{\ang{A,B}}$). 

\bpf
Since $P$ is $\fs11,$ there is a continuous map $F:\dnn\onto P.$   
For $u\in\dN\lom,$ let $P_u=\ens{F(a)}{u\su a\in\dnn}.$ 
There is a number $n_0$ such that $P_{\ang{n_0}}$ is still 
unbounded. 
Let $k_0=0.$ 
By a simple cofinality argument, there is $l_0$ such that 
$P_{\ang{n_0}}$ is still unbounded 
{\it over\/} $\ang{k_0},\,\ang{l_0}$ 
in the sense that there is no ordinal $\al<\omi$ such that 
$P_{\ang{i_0}}\cap{\rrt\al{\ang{k_0}}{\ang{l_0}}}=\pu.$ 
Following this idea, we can define infinite sequences of 
numbers $n_m,\,k_m,\,l_m$ such that both $\sis{k_m}{m\in\dN}$ 
and $\sis{l_m}{m\in\dN}$ are permutations of $\dN$ and, for 
any $m,$ the set $P_{\ang{n_0,...,n_m}}$ is still 
unbounded over $\ang{k_0,...,k_m},\,\ang{l_0,...,l_m}$ in the 
same sense. 
Note that $a=\sis{n_m}{m\in\dN}\in\dN$ and 
$F(a)=\ang{A,B}\in P$ (both $A,\,B$ are subsets of $\dN^2$). 

Prove that the map $f(k_m)=l_m$ witnesses $A\ism\cG B,$ \ie, 
$A(k_j,k_i)$ iff $B(l_j,l_i)$ for all $j,\,i.$ 
Take $m>\tmax\ans{j,i}$ big enough for the 
following: if $\ang{A',B'}\in P_{\ang{i_0,...,i_m}}$ 
then $A(k_j,k_i)$ iff $A'(k_j,k_i),$ and 
similarly $B(l_j,l_i)$ iff $B'(l_j,l_i).$ 
By the construction, there is a pair 
$\ang{A',B'}\in P_{\ang{i_0,...,i_m}}$ with 
$\rrq0{\ang{k_0,...,k_m}}{\ang{l_0,...,l_m}}{A'}{B'},$ 
in particular, $A'(k_j,k_i)$ iff $B'(l_j,l_i),$ as required.
\epF{Lemma}

\bcot[{{\rm See, \eg, Friedman~\cite{frid}}}] 
\label{2tal3} 
If\/ $\rE$ is a Borel\/ \er\ and\/ $\rE\reb{\ism\cG}$ then\/ 
$\rE\reb{\rrt\al\La\La}$ for some\/ $\al<\omi$.
\ecot
\bpf
Let $\vt$ be a Borel reduction of $\rE$ to ${\ism\cG}.$ 
Then 
$\ens{\ang{\vt(x),\vt(y)}}{x\nE y}$ is a $\fs11$ subset of 
$\cP(\dN^2)\ti\cP(\dN^2)$ which does not intersect $\ism\cG,$ 
hence, it is bounded by Lemma~\ref{2tal1}. 
Take an ordinal $\al<\omi$ which witnesses the boundedness.
\epf

Now, if $\rE$ is a Borel \er\ classifiable by countable 
structures then $\rE\reb{\ism\cG}$ by Theorem~\ref{grafs}, 
hence, it remains to establish the following: 

\bprt
\label{r2t}
Any \er\/ $\rrO\al$ is Borel reducible to some\/ $\rT_\xi$.
\eprt
\bpf 
We have ${\rrO0}\reb{\rT_0}$ since $\rrO0$ has countably 
many equivalence classes, all of which are clopen sets. 
To carry out the step $\al\mapsto \al+1$ note that the map   
$\ang{s,A}\mapsto \sis{\ang{s\we k,A}}{k\in\dN}$ is a Borel 
reduction of $\rrO{\al+1}$ to $(\rrO\al)^\iy.$ 
To carry out the limit step, let 
$\la=\ens{\al_n}{n\in\dN}$ be a limit ordinal, and 
${\rR}=\bigvee_{n\in\dN}{\rrO{\al_n}},$ \ie, $\rR$ is 
a \er\ on $\dN\ti\dN\lom\ti\cP(\dN^2)$ defined so that 
$\ang{m,s,A}\rR\ang{n,t,B}$ iff $m=n$ and $\rrq{\al_m} stAB.$ 
However the map $\ang{s,A}\mapsto\sis{\ang{m,s,A}}{m\in\dN}$ 
is a Borel reduction of $\rrO\la$ to $\rR^\iy.$
\epF{Proposition}

\epF{Theorem~\ref{2tal}}

\parf{Turbulent group actions}
\label{groat}

This is an entirely different class of orbit \er s, 
disjoint with those which admit classification by countable 
structures.

\punk{Local orbits and turbulence}
\label{tu:def}

Suppose that a group $\dG$ acts on a space $\dX.$ 
If $G\sq\dG$ and $X\sq\dX$ then let 
\dm
\rr XG=\ens{\ang{x,y}\in X^2}{\sus g\in G\:(x=g\ap y)}
\index{zzrgx@$\rr XG$}%
\dm 
and let $\sym XG$ denote the \er-hull of $\rr XG,$ \ie, 
\index{zzsymgx@$\sym XG$}%
the \dd\sq least \er\ on $X$ such that 
${x\rr XG y}\imp {x\sym XG y}.$ 
In particular ${\sym\dX\dG}={\ergx},$ but generally 
we have ${\sym XG}\sneq{\ergx\res X}.$  
Finally, define 
$\lo XGx=\ek x{\sym XG}=\ens{y\in X}{x\sym XG y}$ 
for $x\in X$ -- the {\it local orbit\/} of $x.$ 
\index{orbit!local}%
\index{zzouxx@$\lo XGx$}%
In particular, $\ek x\dG=\ek x{\ergx}=\lo\dX\dG x,$ 
the full \dd\dG orbit of $x\in\dX$.

\bdf[{{\rm This particular version  taken from 
Kechris~\cite[\S~8]{apg}}}]  
\label{df:durb}
Suppose that $\dX$ is a Polish space and $\dG$ is 
a Polish group acting on $\dX$ continuously. 
\ben
\tenu{(t\arabic{enumi})}
\itla{T1}
A point $x\in\dX$ is {\it turbulent\/} if for any open
\index{turbulent!point}%
non-empty set $X\sq \dX$ containing $x$ and any nbhd  
$G\sq\dG$ (not necessarily a subgroup) of  
$\ong,$ the local orbit $\lo XGx$ is somewhere 
dense (\ie, not a nowhere dense set) in $\dX$.

\itla{T2} 
An orbit $\ek x\dG$ is {\it turbulent\/} if $x$ is such 
\index{turbulent!orbit}%
\index{orbit!turbulent}%
(then all $y\in\ek x\dG$ are turbulent).

\itla{T3}
The action (of $\dG$ on $\dX$) is 
{\it generically}~\footnote
{\ In this research direction, ``generically'', or, in our 
abbreviation, ``gen.'' (property) intends to 
mean that (property) holds on a comeager domain.}, 
or {\it\gen turbulent\/} and $\dX$ is a 
{\it \gen turbulent\/} Polish \dd\dG space, if 
\vyk{
there is a dense and 
\dd\dG invariant $\Gd$ set $D\sq\dX$ such that, for all 
$x\in D,$ the orbit $[x]_\dG$ is dense, turbulent, and meager. 
(For this it is necessary and sufficient \cite[8.7]{apg} 
that 
}%
the union of all dense, turbulent, and meager orbits 
$[x]_\dG$ is comeager.\qeD
\index{action!\gen turbulent}%
\index{space!\gen turbulent}%
\index{turbulent!\gen}%
\index{gengen@\gen= generically}%
\een
\eDf

Our proof of the following theorem, based on ideas in 
\cite[\pff3.2]h, \cite[\pff12]{apg}, \cite{frid}, 
is designed so that only 
quite common tools of descriptive set theory are involved. 
It will also be shown that ``turbulent'' \er s are 
not reducible actually to a much bigger family of \er s than 
orbit \er s of Polish actions of $\isg$.

\bte[{{\rm Hjorth~\cite h}}] 
\label{turb1}
Suppose that\/ $\dG$ is a Polish group, $\dX$ is a 
\gen turbulent Polish\/ \dd\dG space. 
Then\/ $\ergx$ is \poq{not} \bm\ reducible~\footnote
{\ Reducible via a Baire measurable function. 
This is weaker than Borel reducibility, of course.} 
to a Polish action of\/ $\isg,$ hence, \poq{not} 
classifiable by countable structures.
\ete

We begin the proof with two rather simple technical results. 

\blt
\label{deno}
In the assumptions of the theorem, suppose that 
$\pu\ne X\sq\dX$ is an open set, $G\sq\dG$ is a nbhd 
of\/ $\ong,$ and\/ $\lo XGx$ is dense in\/ $X$ for\/ 
\dd Xco\-meager many\/ $x\in X.$ 
Let\/ $U,\,U'\sq X$ be non-empty open and\/ $D\sq X$ 
comeager in\/ $X.$ 
Then there exist points\/ $x\in D\cap U$ and\/ 
$x'\in D\cap U'$ with\/ $x\sym XG x'$.
\elt
\bpf
Under our assumptions there exist points $x_0\in U$ and 
$x'_0\in U'$ with $x_0\sym XG x'_0,$ \ie, there are elements 
$g_1,...,g_n\in G\cup G\obr$ such that 
$x'_0=g_n\app g_{n-1}\app...\app g_1\app x_0$ and in addition 
$g_k\app ...\app g_1\app x_0\in X$ for all $k\le n.$ 
Since the action is continuous, there is a nbhd $U_0\sq U$ of 
$x_0$ such that $g_k\app ...\app g_1\app x\in X$ for all $k$ 
and $g_n\app g_{n-1}\app...\app g_1\app x\in U'$ for all  
$x\in U_0.$ 
Since $D$ is comeager, easily there is $x\in U_0\cap D$ such that  
$x'=g_n\app g_{n-1}\app...\app g_1\app x\in U'\cap D$. 
\epF{Lemma}

\blt
\label{kl}
In the assumptions of the theorem, 
for any open non-empty\/ $U\sq\dX$ and\/ 
$G\sq\dG$ with\/ $\ong\in G$ there is an open non-empty\/ 
$U'\sq U$ such that the local orbit\/ $\lo{U'}Gx$ is 
dense in\/ $U'$ for\/ \dd{U'}comeager many\/ $x\in U'$.
\elt
\bpf
Let $\incl X$ be the interior of the closure of $X.$ 
If $x\in U$ and $\lo UGx$ is somewhere dense (in $U$) 
then the set $U_x=U\cap\incl{\lo UGx}\sq U$ is open and\/ 
\dd{\sym UG}invariant 
(an observation made, \eg, in \cite[proof of 8.4]{apg}), 
moreover, ${\lo UGx}\sq U_x,$ hence, ${\lo UGx}={\lo{U_x}Gx}.$ 
It follows from the invariance that the sets $U_x$ are 
pairwise disjoint, and it follows from the turbulence that 
the union of them is dense in $U.$  
Take any non-empty $U_x$ as $U'.$ 
\epF{Lemma}

\punk{Ergodicity}
\label{ergo}

The non-reducibility in Theorem~\ref{turb1} will be established 
in a special stronger form. 
Let $\rE,\,\rF$ be \er s on Polish spaces 
resp.\ $\dX,\,\dY.$ 
A map $\vt:\dX\to\dY$ is 

\bit
\item\msur
\ef{\it invariant\/} \ if \ 
\index{map!efinvariant@\ef invariant}%
${x\rE y}\imp {\vt(x)\rF\vt(y)}$ for all 
$x,\,y\in\dX$; 

\item
{\it \gen}\ef invariant \ if \
${x\rE y}\imp {\vt(x)\rF\vt(y)}$ holds for all $x,\,y$ 
\index{map!genefinvariant@\gen\ef invariant}%
in a comeager subset of $\dX$;

\vyk{
\item\msur 
\ddf{\it constant\/} if we have 
\index{map!fconstant@\ddf constant}%
${\vt(x)\rF\vt(y)}$ for all $x,\,y\in\dX$; 
}

\item 
{\it \gen}\ddf constant \ if \ ${\vt(x)\rF\vt(y)}$ 
for all $x,\,y$ 
\index{map!genfconstant@\gen\ddf constant}%
in a comeager subset~of~$\dX.$
\eit 
Finally, following Hjorth and Kechris, say that $\rE$ is 
{\it \gen\ddf ergodic\/} if every \bm\  
\index{equivalence relation, ER!genferodic@\gen\ddf ergodic}%
\ef invariant map is \gen\ddf constant. 

\bprt
\label{nonpw}
$\rE$ is \gen\/\ddf ergodic if and only if every 
Borel \gen\/\ef invariant map is \gen\/\ddf constant.
\eprt
\bpf
Let $\rE,\,\rF$ live in resp.\ $\dX,\,\dY.$ 
Suppose that $\vt:\dX\to\dY$ is a Borel \gen   
\ef invariant map. 
There is a Borel comeager set $D\sq\dX$ on which 
$\vt$ is \ef invariant. 
Then we can extend $\vt\res D$ to a \bm\ map 
$\vt':\dX\to\dY$ which is still (everywhere) 
\ef invariant. 
This proves implication $\imp$ of the lemma. 
To prove the opposite implication, let $\vt:\dX\to\dY$ 
be a \bm\ \ef invariant map. 
Then $\vt\res D$ is Borel for a suitable comeager Borel 
set $D\sq\dX.$ 
Let $\vt'$ be any Borel extension of $\vt\res D$ to 
the whole $\dX$.
\epf

\bprt
\label{t2nonr}
Suppose that\/ $\rE$ is \gen\/\ddf ergodic and 
does not have a comeager equivalence class. 
Then\/ $\rE$ is not Borel reducible to\/ $\rF$.\qeD
\eprt
This is exactly how the non-reducibility is often 
established. 
\footnote
{\ Yet there are cases when $\rE$ is neither 
\ddf ergodic nor Borel reducible to $\rF,$ for instance, 
among the \er s of the form $\bel p.$} 
Our proof of Theorem~\ref{turb1} is of this type. 
It consists of two parts~\footnote
{\ There are slightly different ways to the same goal. 
Hjorth~\cite[3.18]h proves outright and with different 
technique, that any \gen turbulent 
\er\ is \gen ergodic \vrt\ any Polish action of $\isg.$ 
Kechris~\cite[\pff12]{apg} proves that 1) any 
\gen  \dd{\rT_2}ergodic \er\ is \gen ergodic 
\vrt\ any 
Polish action of $\isg,$ and 2) any turbulent \er\ is 
\gen \dd{\rT_2}ergodic.}: 

\blt
\label{tul1}
If\/ $\dG$ is a Polish group, $\dX$ a Polish\/ \dd\dG space, 
and\/ $\ergx$ is \bm\ reducible 
to a Polish action of $\isg,$ then there is a comeager\/ $\Gd$ 
set\/ $D\sq \dX$ such that\/ ${\ergx}\res D$ is Borel 
reducible to one of \er s $\rT_\xi$.
\elt
In other words, any \er, \bm\ reducible to a Polish action of 
$\isg,$ is ``generically'' Borel reducible to one of $\rT_\xi.$ 
Note that any \er\ Borel reducible, in proper sense, to one of 
$\rT_\xi,$ is Borel.
 
\blt
\label{tul2} 
Any \er\ induced by a \gen turbulent Polish action 
is \gen\dd{\rT_\xi}ergodic for every $\xi$. 
\elt

\qeDD{Theorem~\ref{turb1} modulo \ref{tul1} and \ref{tul2}}

\punk{``Generical'' reduction of countably classified 
\er s to $\rT_\xi$}
\label{lem3}

Here, we prove Lemma~\ref{tul1}. 
Suppose that $\dG$ is a Polish group, $\dX$ a Polish \dd\dG space, 
and the orbit \er\ $\rE=\ergx$ is \bm\ reducible 
to a Polish action of $\isg.$ 
Then, according to Theorems \ref{siy2cc} and \ref{grafs}, there 
is a \bm\ reduction $\rho:\dX\to\cP(\dN^2)$ of $\rE$ to 
$\ism\cG,$ the isomorphism of binary relations on $\dN.$ 
The remainder of the argument borrows notation from the 
proof of Theorem~\ref{2tal}. 

There is a dense $\Gd$ set $D_0\sq\dX$ such that 
$\vt={\rho\res D_0}$ is continuous on $D_0.$ 
%
By definition, we have ${x\rE y}\imp{\vt(x)\ism\cG\vt(y)}$ and 
${x\nE y}\imp{\vt(x)\not\ism\cG\vt(y)}$ for all $x,\,y\in D_0.$ 
We are mostly interested in the second implication, and the aim is 
to find a $\Gd$ dense set $D\sq D_0$ such that, 
for some $\al<\omi,$ we have 
\bit
\item[$\mtho(\ast)$] 
implication \ 
${x\nE y}\imp{\nrq\al\La\La{\vt(x)}{\vt(y)}}$ \ 
holds for all $x,\,y\in D$.
\eit
(Recall that $A\not\ism\cG B$ iff 
$\sus\al<\omi\:\nrq\al\La\La AB,$ see a remark 
after Lemma~\ref{2tal1}.)  

To find such an $\al$ we apply a Cohen forcing argument. 
\vyk{
Let $\px$ be the (countable) set of all basic open sets of 
$\dX,$ and let $\pg$ be the (countable) set of all basic 
open sets of $\dG.$ 
We consider $\px$ and $\pg$ as forcing notions, with smaller sets 
being stronger conditions, as usual ($\fC$ stands for ``Cohen''). 
}%
Let us fix a countable transitive model $\mm$ of $\zhc,$ 
\index{theory!zfhc@$\zhc$}%
\index{zfhc@$\zhc$}%
\index{zzzfhc@$\zhc$}%
\ie, $\ZFC$ minus the Power Set axiom but plus the axiom: 
``every set belongs to $\HC
=\ens{x}{x\,\text{ is hereditarily countable}}$''. 

We shall assume that
$\dX$ is coded in $\mm$ in the sense that there is a  
set $D_\dX\in\mm$ which is a dense (countable) subset of 
$\dX,$ and $d_\dX\res D_\dX$ 
(the distance function of $\dX$ restricted to $D_\dX$) 
also belongs to $\mm.$ 
Further, $\dG,$ the action, $D_0,\msur$ $\vt$ are also 
assumed to be coded in $\mm$ in a similar sense. 
%
In this assumption, in particular, the notion of a Cohen 
generic, over $\mm,$ point of $\dX,$ or of $\dG,$ makes 
sense, 
in particular, the set 
$D$ of all Cohen generic, over $\mm,$ points of $\dX$ is 
a dense $\Gd$ subset of $\dX$ and $D\sq D_0.$ 
We are going to prove that $D$ fulfills $(\ast)$.

Suppose that $x,\,y\in D,$ and $\ang{x,y}$ is a Cohen generic, 
pair over $\mm.$  
If $x\ergx y$ is false then we have 
$\vt(x)\not\ism\cD\vt(y),$ 
moreover, this fact holds in $\mm[x,y]$ by 
the Mostowski absoluteness, hence, arguing in $\mm[x,y]$ 
(which is still a model of $\zhc$) we find an ordinal 
$\al\in\Ord^\mm=\Ord^{\mm[x,y]}$ with 
$\nrq\al\La\La {\vt(x)}{\vt(y)}.$ 
Moreover, since the Cohen forcing satisfies \ccc, 
there is an ordinal $\al\in\mm$ such that we have 
$\nrq\al\La\La {\vt(x)}{\vt(y)}$ for \poq{every} 
Cohen generic, over $\mm,$ pair $\ang{x,y}\in D^2$ 
such that $x\ergx y$ is false.  
It remains to show that this also holds when $x,\,y\in D$ 
(are generic separately, but)  
do not form a pair, Cohen generic over $\mm$. 

Let $g\in\dG$ be Cohen generic over $\mm[x,y].$~\footnote
{\label{mxy}\ 
In this case, we cannot, generally speaking, define $\mm[x,y]$ 
as a generic extension of $\mm,$ hence, let $\mm[x,y]$ be any 
(countable transitive) model of $\zhc$ containing $x,\,y,$ and 
all sets in $\mm.$ 
It is not really harmful here that $\mm[x,y]$ can contain more 
ordinals than $\mm.$}
Then $x'=g\app x$ is easily Cohen generic over $\mm[x,y]$ 
(because the action is continuous), furthermore, $x'\ergx x,$ 
hence, $x'\ergx y$ \poq{fails}. 
Yet $y$ is generic over $\mm$ and $x'$ is generic over $\mm[y],$ 
thus, $\ang{x',y}$ is Cohen generic over $\mm,$ hence, we 
have $\nrq\al\La\La {\vt(x')}{\vt(y)}$ by the choice of 
$\al.$ 
On the other hand, $\rrq\al\La\La {\vt(x)}{\vt(x')}$ holds 
because $x'\ergx x,$ thus, we finally obtain 
$\nrq\al\La\La {\vt(x')}{\vt(y)},$ as required.\vtm

\qeDD{Lemma~\ref{tul1}}

\punk{Ergodicity of turbulent actions \vrt\ $\rT_\xi$}
\label{lem4}

Here, we prove Lemma~\ref{tul2}. 
The proof involves a somewhat stronger 
property than \gen ergodicity in \prf{ergo}. 
Suppose that $\rF$ is an \er\ on a Polish space $\dX$. 
\bit
\item
An action of $\dG$ on $\dX$ and the induced \eqr{} $\ergx$
are {\it hereditarily generically\/} 
({\it\hgp,} for brevity) \ddf ergodic  
\index{equivalence relation, ER!hgenferodic@\hg\ddf ergodic}%
if \er\ $\sym XG$ is generically \ddf ergodic  
whenever $X\sq \dX$ is a non-empty open set,  
$G\sq\dG$ is a non-empty open set containing $\ong,$ 
and the local orbit $\lo{X}Gx$ is 
dense in $X$ for comeager (in $X$) many $x\in X$. 
\eit
This obviously implies \gen \ddf ergodicity of $\ergx$  
provided the action is \gen turbulent. 
Therefore, Lemma~\ref{tul2} is a corollary of the following 
theorem:

\bte
\label{tut2}
Let\/ $\dX$ be a \gen turbulent Polish\/ \dd\dG space. 
Suppose that an \er\ $\rF$ belongs to\/ $\cF_0,$ 
the least collection of \er s containing\/ 
$\rav\dN$ 
{\rm(the equality on\/ $\dN$)} 
and closed under the operations\/ \ref{oe:beg} -- \ref{oe:end} 
of\/ \prf{opeer}.  
\index{zzFo@$\cF_0$}%
\index{equivalence relation, ER!zzFo@$\cF_0,$ a family of \er s}%
Then\/ $\ergx$ is \hg \ddf ergodic, in particular, is 
not Borel reducible to\/ $\rF$.
\ete
%

\bret
\label{tut:rem}
Due to the other creative operation, the Fubini product, $\cF_0$ 
contains a lot of \er s very different from $\rT_\xi\,,$ among  
them some Borel \er s which do not admit classification 
by countable structures, \eg, all \er s of the form $\rE_\cI,$ 
where $\cI$ is one of Fr\'echet ideals, indecomposable ideals, 
or Weiss ideals of \prf{non-P}. 
(In fact it is not so easy to show that ideals of the two last 
families produce \er s in $\cF_0$.) 
In particular, it follows that 
{\it no \gen turbulent \er\ is Borel reducible to a Fr\'echet, 
or indecomposable, or Weiss ideal\/}.
\eret

Our proof of Theorem~\ref{tut2} goes on by induction on the number 
of applications of the basic operations, in several following 
subsections. 

Right now, we begin with the initial step: prove that, under the 
assumptions of the theorem, $\ergx$ is \hg\dd{\rav\dN}ergodic. 
Suppose that $X\sq\dX$ and $G\sq\dG$ are non-empty open 
sets, $\ong\in G,$ and $\lo{X}Gx$ is dense in $X$ for 
\dd Xcomeager many $x\in X,$ and prove that $\sym XG$ is 
generically \dd{\rav\dN}ergodic. 

Consider, accordingly with Proposition~\ref{nonpw}, a Borel 
\gen\dd{\sym XG,\rav\dN}invariant map $\vt:\dX\to\dN.$ 
Suppose, on the contrary, that $\vt$ is not 
\gen\dd{\rav\dN}constant. 
Then there exist two open non-empty sets $U_1,\,U_2\sq X,$ 
two numbers $\ell_1\ne \ell_2,$ and a comeager set $D\sq X$ 
such that $\vt(x)=\ell_1$ for all $x\in D\cap U_1,$ 
$\vt(x)=\ell_2$ for all $x\in D\cap U_2,$ and $\vt\res D$ 
is ``strictly'' \dd{\sym XG,\rav\dN}invariant. 
Lemma~\ref{deno} yields a pair of points $x_1\in U_1\cap D$ and 
$x_2\in U_2\cap D$ with $x_1\sym XG x_2,$ contradiction. 

\punk{Inductive step of countable power}
\label{tut:iy}

To carry out this step in the proof of Theorem~\ref{tut2},
suppose that 

\bit
\item\msur 
$\dX$ is a \gen turbulent Polish \dd\dG space, 
$\rF$ is a Borel \er\ on a Polish space $\dY,$ 
and the action of $\dG$ on $\dX$ is \hg\ddf ergodic, 
\eit
and prove that the action is \hg\dd\rfy ergodic. 
Fix a nonempty open set $X_0\sq\dX$ and a nbhd $G_0$ of $\ong$ 
in $\dG,$ such that $\lo{X_0}{G_0}x$ is dense in $X_0$ for 
\dd{X_0}comeager many $x\in X_0.$  
Consider, accordingly to Proposition~\ref{nonpw}, 
a Borel function $\vt:X_0\to\dY^\dN,$  
\dd{\sym{X_0}{G_0},\rfy}invariant on a 
dense $\Gd$ set $D_0\sq X_0,$ so that 
\dm
x\sym{X_0}{G_0} x'\;\limp\;
\kaz k\;\sus l\;\skl\vt_k(x)\rF\vt_l(x')\skp
\quad:\qquad\hbox{for all }\;x,\,x'\in D_0\,,
\dm
where $\vt_k(x)=\vt(x)(k),$ $\vt_k: X_0\to\dY,$ 
and prove that $\vt$ is \gen \dd\rfy constant. 

Below, let $\px$ be the Cohen forcing for $\dX,$ 
which consists of rational balls with centers in  
a fixed dense countable subset of $\dX,$ and let $\pg$ 
be the Cohen forcing for $\dG$ defined similarly 
(the dense subset is assumed to be a subgroup). 
Smaller sets are stronger conditions.  
Let us fix a countable transitive model 
$\mm$ of $\zhc$ (see above),  
which contains all relevant objects or their codes, in 
particular, codes of the topologies of $\dX,\,\dG,\,\dY$ 
and the Borel map $\vt$.

\bct
\label{g*}
Suppose that\/ $\ang{x,g}\in \dX\ti\dG$ is\/ 
\dd{\px\ti\pg}generic over\/ $\mm.$  
Then\/ $g\ap x$ is\/ \dd\px generic over\/ $\mm$. \ 
{\rm(Because the action is continuous.)}\qeD
\ect

Coming back to the theorem, fix $k\in\dN.$ 
Consider an open non-empty $U\sq U_0.$ 
By the invariance of $\vt$ and Claim~\ref{g*} 
there are conditions $U'\in\px,\msur$ $U'\sq U,$ 
and $Q\in\pg,\msur$ $Q\sq G_0,$ and
a number $l,$ such that $\vt_k(x)\rF\vt_l(g\ap x)$ 
holds for any \dd{\px\ti\pg}generic over $\mm$ 
pair $\ang{x,g}\in U'\ti Q.$ 
As $Q$ is open, there is $g_0\in Q\cap\mm$ and a 
nbhd $G\sq G_0$ of $\ong$ such that $g_0\,G\sq Q$.

\bct[{{\rm The key point of the turbulence}}]
\label{41*}
If\/ $x,\,x'\in U'$ are\/ \dd\px generic over\/ $\mm$  
and\/ $x\sym{U'}{G}x'$ then we have\/ $\vt_k(x)\rF\vt_k(x')$.
\ect
\bpf 
We argue by induction on $n(x,x')=$ the least number $n$ 
such that there exist $g_1,...,g_n\in G$ satisfying
\bit
\item[\mtho$(\ast)$]
$x'=g_n\app g_{n-1}\app...\app g_1\app x,$ \ and \
$g_k\app...\app g_1\app x\in U'$ \ for all \ $k\le n$.
\eit
Suppose that $n(x,x')=1,$ thus, $x=h\ap x'$ for some $h\in G.$  
Take any \dd\pg generic, over $\mm[x,x']$ 
(see Footnote~\ref{mxy}) 
element $g\in Q\cup Q\obr,$ close enough to $g_0$ for 
$g'=gh\obr$ to belong to $Q.$ 
Then $g$ is \dd{\pg}generic over $\mm[x],$ hence, 
$\ang{x,g}$ is \dd{\px\ti\pg}generic over $\mm$ 
by the product forcing theorem. 
Therefore $\vt_k(x)\rF\vt_l(g\ap x).$ 
Moreover, $g'$ also is \dd{\pg}generic over $\mm[x'],$ 
so that $\vt_k(x')\rF\vt_l(g'\ap x')$ by the same argument. 
Yet we have $g'\ap x'=gh\obr\ap(h\ap x)={g\ap x}$. 

As for the inductive step, suppose that $(\ast)$ holds for some 
$n\ge 2.$ 
Take a \dd\pg generic, over $\mm[x],$ element $g'_1\in G$ close 
enough to $g_1$ for $g'_2= g_2\,g_1\,{g'_1}\obr$ to belong to $G$ 
and for $x^\ast=g'_1\app x$ to belong to $U'.$ 
Note that $x^\ast$ is \dd\px generic over $\mm$ (product forcing) 
and $n(x^\ast,x')\le n-1$ because 
$g'_2\app x^\ast=g_2\app g_1\app x$. 
\epF{Claim}

To summarize, we have shown that for any $k$ and any 
open $\pu\ne U\sq U_0$ there 
exist: an open set $\pu\ne U'\sq U,$ and an open 
$G\sq G_0$ with $\ong\in G,$ such that 
$\vt_k(x)$ is \gen\dd{\sym{U'}{G},\rF}invariant on $U'.$ 
We can also assume that the orbit $\lo{U'}Gx$ 
is dense in $U'$ for \dd{U'}comeager many $x\in U',$ 
by Lemma~\ref{kl}. 
Then, by the \hg\ddf ergodicity, $\vt_k$ is 
\gen\ddf constant on $U',$ that is, there 
is a dense $\Gd$ set $D'\sq U'$ and $y'\in \dY$ such 
that $\vt_k(x)\rF y'$ for all $x\in D'.$ 

We conclude that there exist: an \dd{U_0}comeager set 
$D\sq U_0,$ 
and a countable set $Y=\ens{y_j}{j\in\dN}\sq \dY$ such that, 
for any $k$ and for any $x\in D$ there is $j$ with  
$\vt_k(x)\rF y_j.$ 
Let $\eta(x)=\bigcup_{k\in\dN}\ens{j}{\vt_k(x)\rF y_j}.$ 
Then, for any pair $x,\,x'\in D,$ 
$\vt(x)\rfy \vt(x')$ iff $\eta(x)=\eta(x'),$ so that, by 
the invariance of $\vt,$ we have: 
\dm
{x\sym{U_0}{G_0}x'}\,\limp \,{\eta(x)=\eta(x')}
\quad:\qquad\hbox{for all }\;x,\,x'\in D\,.
\eqno(\ast)
\dm
It remains to show that $\eta$ is a constant on a comeager 
subset of $D$. 

Suppose, on the contrary, that there exist two non-empty 
open sets $U_1,\,U_2\sq U_0,$ a number $j\in\dN,$ and a 
comeager set $D'\sq D$ such that 
$j\in\eta(x_1)$ and $j\nin\eta(x_2)$ for all 
$x_1\in D'\cap U_1$ and $x_2\in D'\cap U_2.$ 
Now Lemma~\ref{deno} yields a contradiction to $(\ast),$
as in the end of \prf{lem4}.\vtm

\vyk{
Take any $x\in U_1,$ \dd\px generic over $\mm.$ 
Since comeager many orbits are dense, the orbit 
$\lo{U_0}{G_0}x$ intersects $U_2,$ in other words, there exist 
$g_1,...,g_n\in G_0$ such that all intermediate points 
$x_k=g_kg_{k-1}...g_1\ap x$ belong to $U_0$ and $x_n\in U_2.$ 
As $\mm\cap\dG$ is dense in $\dG,$ we can assume that all 
$g_i$ belong to $\mm,$ and subsequently all $x_i$ are 
\dd\px generic over $\mm,$ so that, in particular 
$y=x_n\in U_2$ is generic. 
It follows, by the choice of $U_1,\,U_2,$ that 
$j\in\eta(x)$ but $j\nin\eta(y).$ 
Yet we have $x\sym{U_0}{G_0} y$ by the construction, 
and $x,\,y\in D$ by the genericity. 
But this contradicts $(\ast),$ as required.\vtm
}

\qeDD{Inductive step of countable power in Theorem~\ref{tut2}}

\punk{Inductive step of the Fubini product}
\label{tut:Fu}

To carry out this step in the proof of Theorem~\ref{tut2},
suppose that 

\bit
\item\msur 
$\dX$ is a \gen turbulent Polish \dd\dG space, 
for any $k,$ $\rF_k$ be a Borel \er\ on a Polish space 
$\dY_k,$ the action of $\dG$ on $\dX$ is \hg\dd{\rF_k}ergodic 
for any $k,$ 
and $\rF=\fps{\ifi}{\rF_k}{k}$ is, accordingly, 
a Borel \er\ on $\dY=\prod_k\dY_k$, 
\eit
and prove that the action is \hg\dd\rF ergodic. 

Fix a nonempty open set $U_0\sq\dX$ and a nbhd $G_0$ of $\ong$ 
in $\dG,$ such that \dd{U_0}comeager many orbits 
$\lo{U_0}{G_0}x$ with $x\in U_0$ are dense in $U_0.$   
Consider a Borel function $\vt:U_0\to\dY,$ 
\dd{\sym{U_0}{G_0},\rF}invariant on a 
dense $\Gd$ set $D_0\sq U_0,$ \ie,%
\dm
x\sym{U_0}{G_0} y\;\limp\;
\sus k_0\;\kaz k\ge k_0\;\skl\vt_k(x)\rF_k\vt_k(y)\skp
\quad:\qquad\hbox{for all }\;x,\,y\in D_0\,,
\dm
where $\vt_k(x)=\vt(x)(k),$ and prove that $\vt$ is 
\gen\ddf constant.

Consider an open non-empty set $U\sq U_0.$ 
By the invariance of $\vt$ and Claim~\ref{g*}  
there are conditions $U'\in\px,\msur$ $U'\sq U,$ 
and $Q\in\pg,\msur$ $Q\sq G_0,$ and
a number $k_0,$ such that $\vt_k(x)\rF_k\vt_k(g\ap x)$ 
holds for all $k\ge k_0$ and 
for any \dd{\px\ti\pg}generic over $\mm$ 
pair $\ang{x,g}$ of $x\in U'$ and $g\in Q.$ 
As $Q$ is open, there is $g_0\in Q\cap\mm$ 
and a symmetric nbhd $G\sq G_0$ of $\ong$ such that 
$g_0\,G\sq Q$.

\bct
\label{42}
If\/ $k\ge k_0$ and points\/ 
$x,\,y\in U'$ are\/ \dd\px generic over\/ $\mm$  
and\/ $x\sym{U'}{G}y$ then\/ $\vt_k(x)\rF_k\vt_k(y)$. \ 
{\rm(Similarly to Claim~\ref{41*}.)}\qeD
\ect
\vyk{
\bpf
Similar to the proof of Claim~\ref{41*}.
\epf
}

Thus, for any open non-empty $U\sq U_0$ there exist: 
a number $k_0,$ an open non-empty $U'\sq U,$ and a nbhd  
$G\sq G_0$ of $\ong,$ such that $\vt_k(x)$ 
is \gen\dd{\sym{U'}{G},\rF_k}invariant on $U'$ 
for all $k\ge k_0.$ 
We can assume that \dd{U'}comeager many orbits $\lo{U'}Gx$ 
are dense in $U',$ by Lemma~\ref{kl}. 
Now, by the \hg\dd{\rF_k}ergodicity, any $\vt_k$ with 
$k\ge k_0$ is \gen\dd{\rF_k}constant on such a set $U',$ 
hence, $\vt$ itself is \gen\dd{\rF}constant on $U'$ since  
$\rF=\fps{\ifi}{\rF_k}{k}.$ 
It remains to show that these constants are \ddf equivalent  
to each other. 

Suppose, on the contrary, that there exist two non-empty 
open sets $U_1,\,U_2\sq U_0$ and a pair of $y\nF y'$ 
in $\dY$ such that $\vt(x)\rF y$ and $\vt(x')\rF y'$ for 
comeager many $x\in U_1$ and $x'\in U_2.$ 
Contradiction follows as in the end of \prf{tut:iy}.\vtm

\vyk{
Take any $x\in U,$ \dd\px generic over $\mm.$ 
Exactly as in the end of \prf P, there exists a 
\dd\px generic $x'\in U'$ with $x\sym{U'}{G} x'.$ 
Then we have $\vt(x)\rF y$ and $\vt(x')\rF y'$ by the 
genericity, and $\vt(x)\rF\vt(x')$ by the invariance of 
$\vt,$ which contradicts the assumption $y\nF y'$. 
}

\qeDD{Inductive step of Fubini product in Theorem~\ref{tut2}}

\punk{Other inductive steps}
\label{tut:oth}

Here, we accomplish the proof of Theorem~\ref{tut2}, by 
carrying out induction steps, related to operations 
\ref{cun}, \ref{cdun}, \ref{prod} of \prf{opeer}.\vtm

{\sl Countable union\/}. 
Suppose that $\rF_1,\,\rF_2,\,\rF_3,\,...$ are Borel \er s   
on a Polish space\/ $\dY,$ and $\rF=\bigcup_k{\rF_k}$ is 
still a \er, and the Polish and \gen turbulent action of 
$\dG$ on $\dX$ is \hg\dd{\rF_k}ergodic for any $k,$ 
and prove that it remains \hg\ddf ergodic. 

Fix a nonempty open set $U_0\sq\dX$ and a nbhd $G_0$ 
of $\ong$ in $\dG,$ such that \dd{U_0}comeager many orbits 
$\lo{U_0}{G_0}x$ with $x\in U_0$ are dense in $U_0.$   
Consider a Borel function $\vt:U_0\to\dY,$  
\dd{\sym{U_0}{G_0},\rF}invariant on a 
dense $\Gd$ set $D_0\sq U_0.$  
It follows from the invariance that for any open 
$\pu\ne U\sq U_0$ there exist: a number $k$ and open 
non-empty sets $U'\sq U$ and $Q\sq G_0$  
such that $\vt(x)\rF_k\vt(g\ap x)$ holds
for any \dd{\px\ti\pg}generic, over $\mm,$    
pair $\ang{x,g}\in U'\ti Q.$  
We can find, as above, $g_0\in Q\cap\mm$ 
and a nbhd $G\sq G_0$ of $\ong$ such that $g_0\,G\sq Q.$
Similarly to Claims~\ref{41*} and \ref{42}, we have 
$\vt(x)\rF_k\vt(x')$ for any pair of  
\dd\px generic, over $\mm,$ elements $x,\,x'\in U',$ 
satisfying $x\sym{U'}{G}x'.$
It follows, by the ergodicity, that 
$\vt$ is \dd{\rF_k}constant, hence, \ddf constant,  
on a comeager subset of $U'.$ 
It remains to show that these \ddf constants are 
\ddf equivalent to each other, which is demonstrated exactly 
as in the end of \prf{tut:iy}.\vtm

{\sl Disjoint union\/}. 
Let\/ $\rF_k$ be Borel \er s on Polish spaces\/ 
$\dY_k,\msur$ $k=0,1,2,...\;.$ 
By definition, $\bigvee_k{\rF_k}=\bigcup_k{\rF'_k},$
where each $\rF'_k$ is a Borel \er\ defined on the space
$\dY=\bigcup_k{\ans k\ti\dY_k}$ as follows: 
$\ang{l,y}\rF'_k \ang{l',y'}$ iff either $l=l'$ and 
$y=y'$ or $l=l'=k$ and $y\rF_k y'$.\vtm

{\sl Countable product\/}. 
Let\/ $\rF_k$ be \er s on a Polish spaces $\dY_k.$ 
Then $\rF=\prod_k{\rF_k}$ is a \er\ on the space 
$\dY=\prod_k\dY_k.$ 
For any map $\vt:\dX\to\dY,$ to be \gen\dd{\rE,\rF}invariant 
(where $\rE$ is any \er\ on $\dX$) 
it is necessary and sufficient that every 
co-ordinate map $\vt_k(x)=\vt(x)(k)$ 
is \gen\dd{\rE,\rF_k}invariant. 
This allows to easily accomplish this induction step.\vtm

\qeDD
{Theorem~\ref{tut2}, Lemma~\ref{tul2}, Theorem~\ref{turb1}}

\punk{An application to the shift actions of ideals}
\label{turbI}

Say that a Borel ideal\/ $\cZ\sq\pn$ is {\it special\/} 
if there is a sequence of reals\/ $r_n>0$ with\/ 
$\ans{r_n}\to0,$ such that\/ $\sui{r_n}\sq\cZ.$ 
{\it Nontrivial\/} in the next theorem means: containing 
\index{ideal!special}%
\index{ideal!nontrivial}%
no cofinite sets.

\bte
\label{speid}
Suppose that\/ $\cZ$ is a nontrivial Borel special ideal, 
and\/ $\rF$ belongs to the family\/ $\cF_0$ of Theorem~\ref{tut2}. 
Then\/ $\rez$ is generically \ddf ergodic, hence, is 
not Borel reducible to\/ $\rF$.
\ete
\bpf
The ``hence'' statement follows because by the nontriviality 
all \dd\rez equivalence classes are meager subsets of $\pn$.
 
As $\cZ$ is special, let $\sis{r_k}{}\to 0$ be a sequence of 
positive reals such that $\srn\sq\cZ.$ 
It obviously suffices to prove that $\ern=\rE_{\srn}$ is 
generically \ddf ergodic. 
Further, by Theorem~\ref{tut2}, it suffices to prove that the 
shift action of $\srn$ on $\pn$ is Polish and \gen turbulent. 

\vyk{
\epf

\ble
\label{ket}
If\/ $\cZ\zT\vpi\zT r_k$ are as above, and\/ 
$\sis{r_k}{}\to0,$ then the shift action of\/ $\cZ$ on\/ 
$\pn$ is \gen turbulent. 
\ele
\bpf
}

The ideal $\srn$ is easily a P-ideal, hence, 
a polishable group (with $\sd$ as the operation). 
For instance, $\srn$ is a Polish group in the topology 
generated by the metric $\dpr(a,b)=\vpr(a\sd b)$ on 
$\srn,$ where 
\bit
\item
$\vpr(x)=\sum_{n\in x}r_n$ for $x\in\pn,$ 
so that $\srn=\ens{x}{\vpr(x)<\piy}$.  
\eit
The shift action of $\srn$ by $x\app y=x\sd y$ on $\pn$ 
(considered in the product topology; $\pn$ is here 
identified with $\dn$) 
is then continuous. 
It remains to verify the turbulence. 

Let $x\in\pn.$  
The orbit $\ek x\srn=\srn\sd x$ is easily dense and meager, 
hence, it suffices to prove that $x$ is a turbulent point 
of the action. 
Consider an open set $X\sq\pn$ containing $x,$ and 
a \dd\dpr hbhd $G$ of $\pu$ (the neutral element of $\srn$); 
we may assume that, for some $k,$ 
$X=\ens{y\in\pn}{y\cap [0,k)=u},$ where $u=x\cap[0,k),$ and 
$G=\ens{g\in\srn}{\vpi(g)<\ve}$ for some $\ve>0.$ 
Prove that the local orbit $\lo XGx$ is somewhere dense in $X$. 

Let $l\ge k$ be big enough for $r_n<\ve$ for all $n\ge l.$ 
Put $v=x\cap [0,l)$ and prove that $\lo XGx$ is dense in 
$Y=\ens{y\in\pn}{y\cap [0,l)=v}.$ 
Consider an open set  
$Z=\ens{z\in Y}{z\cap [l,j)=w},$ where $j\ge l\zT w\sq[l,j).$ 
Let $z$ be the only element of $Z$ with 
$z\cap{[j,\piy)}=x\cap{[j,\piy)},$ thus, 
${x\sd z}=\ans{l_1,...,l_m}\sq[l,j).$ 
Each $g_i=\ans{l_i}$ belongs to $G$ by the choice of $l$ 
(indeed, $l_i\ge l$). 
Moreover, easily 
$x_i=g_i\sd g_{i-1}\sd...\sd g_1\sd x=
\ans{l_1,...,l_i}\sd x$ 
belongs to $X$ for any $i=1,...,m,$ and $x_m=z,$ 
thus, $z\in \lo XGx,$ as required.
\epf

The next corollary returns us to the discussion in
the end of \prf{brb}.

\bcor
\label{-t2}
The equivalence relations\/ $\fco$ and\/ $\Ed$ are
not Borel reducible to any ideal\/ $\rF$ in 
the family\/ $\cF_0$ of Theorem~\ref{tut2}, in particular,  
are not Borel reducible to\/ $\rtd$.
\ecor
\bpf
According to lemmas \ref{co=d} and \ref{l1=s}, it suffices to
prove that the ideals $\zo$ (density $0$) and $\sui{1/n}$
are special.
The latter is special by definition.
As for the former, see ??? (that $\sui{1/n}\sq\zo$).
\epf

\newpage

\parf{Ideal $\It$ and the \eqr\ $\Et$}
\label{It:id}

The ideal $\ofi$ is traditionally denoted by $\It.$   
It consists of all sets $x\sq\pnn$ such that all 
cross-sections $\seq xn=\ens{k}{\ang{n,k}\in x}$ are finite.
It defines the \er\ $\Et=\rE_{\It}$ on 
$\pnn$ by $x\Et y$ iff $x\sd y\in\It.$ 
But we rather consider $\Et$ as an \er\ on $\pnd$ defined 
\index{equivalence relation, ER!E3@$\Et$}%
\index{zzE3@$\Et$}%
by $x\Et y$ iff $x(n)\Eo y(n)$ for all $n:$ 
here $x\yo y$ belong to $\pnd$.

\punk{Ideals below $\It$}
\label{<it}

\ble
\label{III}
$\ifi\rebs \It.$
$\It$ and\/ $\Ii$ are\/ \dd\reb incomparable.
\ele
\bpf
To see that $\ifi\rebs \It$ take 
$\vt(x)=\ens{\ang{n,0}}{n\in x}.$ 
That $\It\not\reb\Ii$ can be shown as follows: otherwise 
by Theorem~\ref{rig1} $\It$ would be isomorphic either 
to one of $\ifi,\msur$ $\Ii,$ or to a trivial variation 
of $\ifi,$ which can be easily shown to be not the case. 
To see that $\Ii\not\reb\It$ recall that $\It=\ofi$ is 
of the form $\Exh_\psi$ for a \lsc\ submeasure $\psi$ 
(Example~\ref{exh:e}) and apply Theorem~\ref{sol}. 
\epf

The following theorem is analogous to Theorem~\ref{rig1}, 
yet the method of its proof is absolutely different.

\bte[{{\rm Kechris~\cite{rig}}}]
\label{rig2}
If\/ $\cI\reb\It$ is a Borel (nontrivial) ideal on\/ $\dN$ 
then either\/ $\cI\cong\It$ or\/ $\cI$ is a trivial 
variation of\/ $\ifi$. 
\ete
\bpf
First of all we make use of Theorem~\ref{sol}: 
$\Ii\not\reb \cI$ according to Lemma~\ref{III}, therefore, 
$\cI=\Exh_\vpi$ for a \lsc\ submeasure $\vpi$ on $\dN.$ 
We can \noo\ suppose that $\vpi(x)\le 1$ for any $x\in\pn.$ 
Now put $U_n=\ens{k}{\vpi(\ans k)\le\frac1n}$.

We assert that $\tlim_{n\to\iy}\vpi(U_n)=0.$ 
Indeed, otherwise $\vpi(U_n)>\ve$ for some $\ve>0$ and all $n.$  
As $\vpi$ is \lsc\ we can choose a sequence 
of numbers $n_1<n_2<n_3<...$ and for any $l$ a finite set 
$w_l\sq U_{n_l}\dif U_{n_{l+1}}$ with $\vpi(w_l)>\ve.$ 
Then $W=\bigcup_lw_l\nin\cI$ and obviously 
$\ans{\vpi(\ans k)}_{k\in W}\to 0.$ 
Note that the Borel ideal $\cZ=\cI\res W$ satisfies 
$\cZ\reb\cI$ (via the identity map), because $W\nin\cI.$ 
On the other hand, $\cZ$ is isomorphic to a special ideal 
(see \prf{turbI})
via the order preserving bijection of $W$ onto $\dN.$ 
It follows from Theorem~\ref{speid} that $\rez$ is not Borel 
reducible to any \eqr\ in $\cF_0,$ hence, neither is $\rei.$ 
But $\rE_{\It}=\Et$ obviously belongs to $\cF_0,$ which is 
a contradiction because $\cI\reb\It$.

Thus $\vpi(U_n)\to 0.$ 
Then clearly a set $x\in\pn$ belongs to $\cI$ iff 
$x\cap (U_n\dif U_{n+1})$ is finite for any $m,$ which 
easily implies that $\cI$ is as required. 
\imar{check the proof}
\epf

\punk{Assembling equivalence relations}
\label{asser}

The next theorem, similar to a couple of results above,
\imar{give ref}
will be used in the proof of a dichotomy theorem related to
$\Et$.

\bte
\label{t:ass} 
Suppose that\/ $\dX,\,\dY$ are 
Polish spaces, $P\sq\dX\ti\dY$ is a Borel set, 
$\rE$ is a Borel \er\ on\/ $P,$ and\/ 
$\dG$ is a countable group acting on\/ $\dX$ in a Borel 
way so that\/ $\ang{x,y}\rE\ang{x',y'}$ implies\/ 
$x\ergx x'.$ 
Finally, assume that\/ $\rE\res{P(x)}$ is smooth for 
each\/ $x\in \dX,$ where\/ $P(x)=\ens{\ang{x',y}\in P}{x'=x}.$
Then\/ $\rE$ is Borel-reducible to a Borel action of\/ 
$\dG$.
\ete
\bpf
We can assume that $\dX=\dY=\dn$ and both $P$ and $\rE$ 
are $\id11.$ 
We can also assume that the action of $\dG$ 
(a countable group) is $\id11.$ 
Then clearly ${x\ergx x'}\imp{\id11(x)=\id11(x')}.$ 
Define 
$\tP(x)=\bigcup_{a\in\dG}P(a\app x)$ for $x\in\dX$.

\bct
\label{ass1}
Suppose that\/ $\ang{x,y}$ and\/ $\ang{x',y'}$ belong 
to\/ $P$ and\/ $x\ergx x'.$ 
Then\/ $\ang{x,y}\rE\ang{x',y'}$ iff the equivalence\/ 
$\ang{x,y}\in U\eqv\ang{x',y'}\in U$ holds for any\/  
\dd{\rE\res{\tP(x)}}-invariant $\id11(x)$ 
set\/ $U\sq\tP(x)$.
\ect
\bpf
Note that $\rE\res{\tP(x)}$ is still smooth by 
Theorem~\ref{cud1} because $\dG$ is countable. 
In addition $\rE\res{\tP(x)}$ is $\id11(x).$ 
This observation yields the result, because otherwise, 
\ie, if the \er, defined om $\tP(x)$ by intersections 
with \dd{\rE\res{\tP(x)}}invariant $\id11(x)$ sets, 
is coarser than $\rE\res{\tP(x)},$ then it is known 
from the proof of the 2nd dichotomy theorem 
(Theorem~\ref{2dih}) 
that we would have $\Eo\reb{\rE\res{\tP(x)}},$ 
a contradiction with the smoothness.
\epF{Claim}

For any $x\in\dX$ let $E(x)$ be the set of all $e\in\dN$ 
which code a $\id11(x)$ subset of $P,$ and, for 
$e\in E(x),$ let $\dc ex$ be the $\id11(x)$ subset of $P$ 
coded by $e.$ 
(It is known that $\ens{\ang{x,e}}{e\in E(x)}$ is 
$\ip11.$) \ 
Let $\inva(x,e)$ be the formula 
\dm
x\in\dX\lland e\in E(x)\lland \dc ex\sq\tP(x)\lland
\hbox{$\dc ex$ is \dd{\rE\res{\tP(x)}}-invariant}\,. 
\dm

\bcot
\label{ass3}
Let\/ $\ang{x,y},\,\ang{x',y'}$ be as in 
Claim~\ref{ass1}. 
Then\/ $\ang{x,y}\rE\ang{x',y'}$ iff\/ 
$\ang{x,y}\in \dc ex\eqv\ang{x',y'}\in \dc ex$ 
holds for any\/ $e$ with\/ $\inva(x,e)$.\qeD
\ecot

Implication $\nmp$ of the ``iff'' in this Corollary can 
be considered as a property of the $\ip11$ set 
$C=\ens{\ang{x,e}}{\inva(x,e)},$ \ie, the property that 

\bit 
\item 
for all pairs $\ang{x,y}$ and $\ang{x',y'}$ in $P$ 
with $x\ergx x',$ we have:\\[0.7ex]
if $
\kaz\ang{x,e}\in C\;
(\ang{x,y}\in \dc ex\seqv\ang{x',y'}\in \dc ex)$ 
then $\ang{x,y}\rE\ang{x',y'}$.
\eit
This is easily a $\ip11$ property in the codes, hence, 
by the $\ip11$ Reflection, there is a $\id11$ set 
$B\sq C$ satisfying the same property, that is, we have

\vyk{
Unfortunately the formula $\inva$ is a $\ip11$ formula, 
which makes it difficult to use Corollary~\ref{ass3} 
directly if one wants to define a Borel reduction. 
Yet there is a typical trick which allows to extract a 
sufficient $\id11$ part. 
Note that  
\dm
Q=\ens{\ang{x,y,x',y'}}{\ang{x,y}\in P\lland
\ang{x',y'}\in P\lland {x\ergx x'} \lland 
{\ang{x,y}\nE\ang{x',y'}}} 
\dm
is a $\id11$ subset of $P\ti P$ while 
\dm
S=\ens{\ang{x,y,x',y',e}\in Q\ti\dN}
{\inva(x,e)\lland \neg\:
(\ang{x,y}\in \dc ex\eqv\ang{x',y'}\in \dc ex)}
\dm
is a $\ip11$ set whose projection on $Q$ coincides with 
$Q.$ 
Then there is a $\id11$ map $\ve:Q\to\dN$ 
such that $\ang{x,y,x',y',\ve(x,y,x',y')}\in S$ holds 
whenever $\ang{x,y,x',y'}\in Q.$ 
Now
$A=\ens{\ang{x,\ve(x,y,x',y')}}{\ang{x,y,x',y'}\in Q}$ 
is a $\is11$ subset of the $\ip11$ set 
$C=\ens{\ang{x,e}}{\inva(x,e)}\sq\dX\ti\dN,$ 
hence, there is a $\id11$ set $B$ with $A\sq B\sq C$. 
}%

\bcot
\label{ass4}
Let\/ $\ang{x,y},\,\ang{x',y'}$ be as in 
Claim~\ref{ass1}. 
Then\/ $\ang{x,y}\rE\ang{x',y'}$ iff\/ 
$\ang{x,y}\in \dc ex\eqv\ang{x',y'}\in \dc ex$ 
holds for any\/ $e$ with\/ $\ang{x,e}\in B$.\qeD
\ecot
\vyk{
\bpf
The direction $\imp$ of the ``iff'' follows from 
Corollary~\ref{ass3} because $B\sq C.$ 
The opposite implication holds by the definition of 
$\ve$.
\epF{Corollary}
}

To continue the proof of the theorem, define, for 
any $\ang{x,y}\in P$, 
\dm
D_{xy}=\ens{\ang{a,e}}{a\in \dG\lland 
\ang{a\app x,e}\in B
\lland \ang{x,y}\in\dc e{a\app x}}\,.
\dm
Clearly $\ang{x,y}\mapsto D_{x,y}$ is a $\id11$ map 
$P\to\cP(\dG\ti\dN)$.
 
If $D\sq\dG\ti\dN$ and $b\in\dG$ then put 
$b\circ D=\ens{\ang{ab\obr,e}}{\ang{a,e}\in D}$.

\bct
\label{ass2}
Suppose that\/ $\ang{x,y}$ and\/ $\ang{x',y'}$ belong 
to\/ $P,$ $b\in\dG,$ and\/ $x'=b\app x.$ 
Then\/ $\ang{x,y}\rE\ang{x',y'}$ iff\/ 
$b\circ D_{xy}=D_{x'y'}$.
\ect
\bpf
Assume that $b\circ D_{xy}=D_{x'y'}.$ 
According to Corollary~\ref{ass4}, to prove 
$\ang{x,y}\rE\ang{x',y'}$ it suffices to show that 
${\ang{x,y}\in\dc ex}\eqv{\ang{x',y'}\in\dc ex}$ 
holds whenever $\ang{x,e}\in B.$ 
We have 
\dm
{\ang{x,y}\in\dc ex}\seqv{\ang{\La,e}\in D_{xy}} 
\seqv{\ang{b\obr,e}\in D_{x'y'}}\seqv
{\ang{x',y'}\in\dc e{b\obr\app x'}=\dc ex}\,,
\dm
as required. 
Conversely, let $\ang{x,y}\rE\ang{x',y'}.$ 
If $\ang{a,e}\in D_{xy}$ then $\ang{a\app x,e}\in B$ 
and $\ang{x,y}\in\dc e{a\app x},$ hence, 
$\ang{x',y'}\in\dc e{a\app x},$ too, because 
the set $\dc e{a\app x}$ is invariant and 
$\ang{x,y}\rE\ang{x',y'}.$  
Yet $a\app x=ab\obr\app x',$ therefore, 
by definition, $\ang{ab\obr,e}\in D_{x'y'}.$ 
The same argument can be carried out in the opposite 
direction, so that $\ang{a,e}\in D_{xy}$ 
iff $\ang{ab\obr,e}\in D_{x'y'},$ that means 
$b\circ D_{xy}=D_{x'y'}$.
\epF{Claim}

To end the proof of the theorem, consider 
$\dZ=\dX\ti\cP(\dG\ti\dN),$ a Polish space.  
Define a Borel action   
$b\app\ang{x,D}=\ang{b\app x, b\circ D}$ 
of $\dG$ on $\dZ.$  
We assert that $\vt(x,y)=\ang{x,D_{xy}}$ 
is a Borel reduction of $\rE\res P$ to 
the action $\aer\dG\dZ.$ 
Indeed, let $\ang{x,y}$ and $\ang{x',y'}$ belong 
to $P.$ 
Suppose that $\ang{x,y}\rE\ang{x',y'}.$ 
Then $x\aer\dG\dX x',$ so that $x'=b\app x$ for some 
$b\in \dG.$ 
Moreover, $b\circ D_{xy}=D_{x'y'}$ by Claim~\ref{ass2}, 
hence, $\vt(x',y')=b\app\vt(x,y).$  
Let, conversely, $\vt(x',y')=b\app\vt(x,y),$ so that 
$x'=b\app x$ and $D_{x'y'}=b\circ D_{xy}.$ 
Then $\ang{x,y}\rE\ang{x',y'}$ by Claim~\ref{ass2}, 
as required. 
\epf

\punk{The 6th dichotomy}
\label{Et:er}

\bte[{{\rm Hjorth and Kechris \cite{hk:nd,hk:rd}}}] 
\label{6dih}
If\/ $\rE\reb\Et$ is a Borel \er\ then either\/ 
$\rE\reb\Eo$ or\/ $\rE\eqb \Et$.
\ete
\bpf[{{\rm a modification of the proof in \cite{hk:rd}}}]  
We may assume that $\rE$ is a $\id11$ \er\ on a 
recursively presented Polish space $\dX,$ and there is 
a $\id11$ reduction $\vt:\dX\to\pnd$ of $\rE$ to 
$\Et.$ 
Let $Q=\ran \vt,$ a $\is11$ subset of $\pnd$.
 
For $x,\,y\in\pnd$ and $n\in\dN,$ define 
$x\eqn y$ iff $x\Et y$ and $x\rmq n=y\rmq n$ 
(the latter requirement means $x_k=y_k$ for all $k<n$). 
For $n,k,p\in\dN$ put~\footnote
{\ Hjorth and Kechris~\cite{hk:rd} define $\cA_{nkp}$ 
with $\kaz x,y\in Q\cap A$ instead of $\kaz x,y\in A.$ 
Let us use $\cA'_{nkp}$ to denote their 
version, thus, $\cA_{nkp}\sq\cA'_{nkp}.$ 
However if Case 1 holds in the sense of $\cA'_{nkp}$ 
then it also holds in the sense of $\cA_{nkp}$ because 
$A\in \cA'_{nkp}$ iff $A\cap Q\in \cA_{nkp}$.}
\dm
\cA_{nkp}=\ens{A\sq\pnd}{A\,\hbox{ is }\,\is11\lland 
\kaz x,y\in A\;
\skl{x\eqn y}\imp
{x_k\sd y_k\sq\ir0p}
\skp}\,.
\dm 

\bct
\label{utv1}
If\/ $A\in\cA_{nkp}$ then there is a\/ $\id11$ set\/ 
$B\in\cA_{nkp}$ with\/ $A\sq B$.
\ect
\bpf
(Reflection)
\epF{Claim}

Put $A_{nkp}=\bigcup\ens{A}{A\in\cA_{nkp}}$ and 
$\wA=\bigcup_n\bigcap_{k\ge n}\bigcup_p\:A_{nkp}$
\vtm\vom

\vyk{
\bcot
\label{utv2}
$A_0$ is\/ $\ip11$.
\ecot
\bpf
A standard computation using Claim~\ref{utv1}.
\epF{Corollary}
}

\vyk{
$x\rmq n$ $x\rlq n$
}

{\bfsl Case 1\/\bf:} \ 
$Q\sq \wA$.\hfill
{\bfsl Case 2\/\bf:} \ \hspace*{\mathsurround}otherwise.
\hfill $\,$

\punk{Case 1}

We are going to prove that in this case $\rE\reb\Eo.$
\vyk{
First of all put 
\dm
\cB_{nkp}=\ens{B\sq\pnd}{B\,\hbox{ is }\,\id11\lland 
\kaz x,y\in B\;
\skl{x\eqn y}\imp
{x_k\sd y_k\sq\ir0p}
\skp}\,,
\dm 
note the difference with $\cA_{nkp}.$ 
Note also the following: 
if $B\in\cB_{nkp}$ then for any $x\in B$ the set 
$\ens{y_k}{y\in B\cap[x]_{\eqn}}$ is finite.

Let 
$B_0=\bigcup_n\bigcap_{k\ge n}\bigcup_p\:
\ens{B}{B\in\cB_{nkp}}$.

\bct
\label{..1}
$Q\sq B_0$. 
\ect
\bpf
It suffices to prove that for any $\id11$ set 
$A\in\cA_{nkp}$ there is a $\id11$ set $B\in\cB_{nkp}$ 
with $Q\cap A\sq B.$ 
Let
\dm
X=\ens{x\in A}{\kaz y\in A\cap Q\;
\skl
{x\eqn y}\imp{x_k\sd y_k\sq\ir0p}
\skp}\,,
\dm
this is a $\ip11$ set with $A\cap Q\sq X$ by he definition 
of $\cA_{nkp},$ hence, there is a $\id11$ set $Y$ with 
$A\cap Q\sq Y\sq X.$ 
Now the $\ip11$ set 
\dm
Z=\ens{x\in Y}{\kaz y\in Y\;
\skl
{x\eqn y}\imp{x_k\sd y_k\sq\ir0p}
\skp}
\dm
still satisfies $A\cap Q\sq Z,$ hence, 
there is a $\id11$ set $B$ with $A\cap Q\sq B\sq Z.$ 
However easily $B\in\cB_{nkp}$.
\epF{Claim}
}%

As easily $\wA$ is $\ip11$ by Claim~\ref{utv1} and 
a standard computation, 
there is a $\id11$ set $R$ such that $Q\sq R\sq \wA.$ 
Thus, for $\rE\reb\Eo$ it suffices now to prove

\ble
\label{.key} 
${\Et\res R}\reb\Eo$ 
for any $\id11$ set $R\sq \wA$. 
\ele
\bpf
By \Kres\ there exists a $\id11$ map 
$\nu:R\to\dN$ such that 
\dm
\kaz k\ge\nu(x)\;\sus p\;\sus B\in \cA_{\nu(x),k,p}\;
(x\in B\in \id11)
\dm
for any $x\in R.$ 
Let $R_n=\ens{x\in R}{\nu(x)\le n},$ these are increasing 
$\id11$ subsets of $R,$ and $R=\bigcup_nR_n.$ 
According to Theorem~\ref{cud2}, it suffices to prove 
that ${\Et\res R_n}\reb\Eo$ for any $n.$ 
Thus let us fix $n.$ 
By definition we have
\dm
\kaz x\in R_n\;\kaz k\ge n\;\sus p\;
\sus B\in \cA_{nkp}\;
(x\in B\in \id11)\,.
\eqno(\ast)
\dm 

Recall that $\fC$ is the least class of sets containing 
all open sets and closed under the A-operation and the 
complement. 
A map $f$ is called {\it \fC-measurable\/} iff all 
\dd fpreimages of open sets belong to $\fC$.

\bct
\label{.k1}
For any\/ $n$
there is a \dd\fC measurable map\/ $f:R_n\to\pnd$ such 
that\/ $f(x)=f(y)\eqn x$ whenever\/ $x,\,y\in R_n$ 
satisfy\/ $x\eqn y$.
\ect
\bpf
Let $E\sq\dN$ be the $\ip11$ set of all codes of 
$\id11$ subsets of $\pnd,$ and let 
$W_e\sq\pnd$ be the $\id11$ set coded by $e\in E.$ 
We have, by $(\ast)$,
\dm
\kaz x\in R_n\;
\kaz k\ge n\;\sus p\;\sus e\in E\;
(x\in W_e\in \cA_{nkp})\,,
\dm 
and an ordinary application of the Kreisel selection 
yields a pair of $\id11$ maps 
$\pi,\,\ve:R_n\ti\dN\to\dN$ such that  
$\ve(x,k)\in E$ and 
$x\in W_{\ve(x,k)}\in\cA_{n,k,\pi(x,k)}$ 
hold whenever $x\in R_n$ and $k\ge n.$ 
Let $\pit(x,k)$ and $\tve(x,k)$ 
to be the least, in the sense of any fixed recursive 
\dd\om long wellordering of $\dN\ti\dN,$ of all 
possible pairs $\pi(x',k)$ and $\ve(x',k)$ with 
$x'\in R_n\cap[x]_{\eqn}.$ 
Then $\pit$ and $\tve$ are \dd{\eqn}invariant in the 
1st argument. 
In addition, we have  
$W_{\tve(x,k)}\in\cA_{n,k,\pit(x,k)}$ and the set 
$Z_{xk}=R_n\cap[x]_{\eqn}\cap W_{\tve(x,k)}$ 
is nonempty, 
whenever $x\in R_n$ and $k\ge n$. 

Let $x\in R_n.$ 
For any $k\ge n,$ the set 
$Y_{xk}=\ens{y_k}{y\in Z_{xk}}\sq\pn$ is finite 
(and nonempty) by the definition of $\cA_{nkp}\,,$ 
thus, let $f_k(x)$ be the least member of $Y_{xk}$ 
in the sense of the lexicographical order of $\pn$. 
Define $f(x)\in\pnd$ so that $f(x)_k=x_k$ for $k<n$ 
and $f(x)_k=f_k(x)$ for $k\ge n$. 

That $f(x)=f(y)$ whenever $x\eqn y$ follows from the 
invariance of $\ve$ and $\pi.$ 
To see that $f(x)\eqn x$ note that by definition 
$f_k(x)\Eo x_k$ for $k\ge n:$ 
indeed, $f_k(x)=y_k$ for some $y\in[x]_{\eqn},$ but 
$x\eqn y$ implies $x_k\Eo y_k$ for all $k.$ 
Finally, the \dd\fC measurability needs a routine check.
\epF{Claim}

For any $u\in\pn^n$ let 
$R_n(u)=\ens{x\in R_n}{x\rmq n=u}$.

\bct
\label{utv4}
If\/ $u\in\pn^n$ then\/ $\Et\res R_{n}(u)$ is 
smooth.
\ect
\bpf 
As $\Et$ and $\eqn$ coincide on $R_n(u),$ the relation 
$\Et\res R_n(u)$ is smooth via a \dd\fC measurable, hence,
a Baire-measurable map. 
Suppose, towards the contrary, that it is not really 
smooth, \ie, via a Borel map. 
Then, by the 2-nd dichotomy theorem, we have 
$\Eo\reb{\Et\res R_n(u)},$ hence, $\Eo$ turns out to 
be smooth via a Baire-measurable map, which is easily 
impossible.
\epF{Claim}

To complete the proof of the lemma, let 
$\dG={\pwf\dN}^n,$ acting on $\dX=\pn^n$ componentwise 
and by $\sd$ at each of the $n$ co-ordinates, so that,  
for $u,\,v\in\dX,$ we have $u\ergx v$ iff $u_k\Eo v_k$ 
for all $k<n.$  
Let us apply Theorem~\ref{t:ass} with $\dG$ and $\dX$ 
as indicated, and $P=R_n$ and $\rE={\Et\res R_n},$  
Claim~\ref{utv4} witnesses the principal requirement. 
We obtain: $\Et\res R_n$ is Borel reducible to a \er\ 
induced by a Borel action of $\dG.$ 
Yet $\dG$ is the increasing union of a countable 
sequence of its finite subgroups, hence, any \er\ 
induced by a Borel action of $\dG$ is hyperfinite, hence, 
Borel reducible to $\Eo$.\vtm

\epF{Lemma~\ref{.key} and Case 1 in Theorem~\ref{6dih}}

\punk{Case 2}

Then the $\is11$ set $H=Q\dif \wA$ is non-empty. 
Our idea will be to define a Borel subset $X$ of $H$ such 
that $\Et\res X\eqb\Et,$
the ``or'' case of Theorem~\ref{6dih}.

By definition, 
$H=\bigcap_n\bigcup_{k>n} H_{nk},$ where
$H_{nk}=H\dif \bigcup_pA_{nkp}.$
Note that 
\dm
\bay{rcl}
H_{nk}=
%
\ens{x\in H}
{\kaz p\:\kaz A\in\id11\:(x\in A\imp A\nin\cA_{nkp})}
\eay
\dm
by Claim~\ref{utv1}, and hence $H_{nk}$ is $\is11$ by
rather elementary computation.

Let $b$ be any recursive bijection $\dN^2\onto\dN,$
increasing in each argument.
Put $L(n)=\tmax\ens{r}{b(r,0)\le n}$ -- thus for any
$\ell>L(n)$ we have $b(\ell,j)>n\zd \kaz j.$

The splitting system used here will contain non-empty
$\is11$ sets $X_s\sq \pnd\yt s\in\dln,$ numbers $k_m\yt m\in\dN,$
and elements $g_s\in\pnd\yt s\in\dln,$ satisfying the
following requirements \ref{3i} -- \ref{3last}:
\ben
\tenu{(\roman{enumi})}
\itla{3i}\msur
$X_\La\sq H\yt X_{s\we i}\sq X_s\yt \dia{X_s}\le2^{-\lh s},$ 
and a certain condition, in terms of the Choquet game,
holds, connecting each $X_{s\we i}$ with $X_s$ so that, as a 
consequence, $\bigcap_nX_{a\res n}\ne\pu$ 
for any $a\in\dn$.

\itla{3ii}\msur
$0<k_0<k_1<\dots$ and
$X_{0^{n+1}}\sq \bigcap_{r<L(n)}H_{r,k_r}$.~\snos
{Recall that $0^m$ is a sequence of $m$ zeros.}

\itla{3iii}
If $s\in2^{n+1}$ then 
$g_s(i)$ is finite for all $i$ and $=\pu$
for all $i>k_{L(n)};$
in addition, $g_{0^{n+1}}(i)=\pu$ for all $i$.

\itla{3iv}
For any $s\in2^{n+1},$ we have
$\kaz x\in X_{0^{n+1}}\:
\sus y\in X_s\:(y\equiv_{k_{L(n)}}g_s\app x)$;~\snos
{For $g,x\in\pnd,$ $g\app x=y\in\pnd$ is defined by
$y(n)=g(n)\sd x(n)\zd\kaz n$.}

\itla{3v}

\itla{3vi}
\label{3last}

\een

\vyk{
\dm
\kaz n\;\sus k>n\;\kaz p\;\kaz A\in\is11\:
(x\in A\imp A\nin \cA_{nkp}), \qquad
\text{that is,}
\dm
\bus
\label{++1}
\kaz n\;\sus k>n\;\kaz p\;\kaz A\in\is11\:
\skl
x\in A\imp \sus y,z\in A\;
({y\eqn z}\land
{y_k\sd z_k\not\sq\ir0p})\skp.
\eus
}

\epF{Theorem~\ref{6dih}}

\parf{Summable ideals}
\label{summI}

Farah~\cite[\pff 1.12]{aq} gives the following classification 
of summable ideals $\sui{r_n},$ 
based on the distribution of numbers $r_n$:
\ben
\tenu{(S\arabic{enumi})}
\itla{S1}
{\it Atomic\/} ideals: there is $\ve>0$ such that 
the set $A_{\ve}=\ens{n}{r_n\ge\ve}$ is infinite and 
satisfies $\mrn(\dop{A_{\ve}})<\piy.$ 
In this case $\sui{r_n}=\ans{X:X\cap A_{\ve}\in\ifi};$ 
Kechris~\cite{rig} called such ideals 
{\it trivial variations of\/ $\ifi$}.

\itla{S2}
{\it Dense\/} (summable) ideals: $r_n\to 0$.

\itla{S3}
There is a decreasing sequence of 
positive reals $\ve_n\to0$ sich that all sets 
$D_n=A_{\ve_{n+1}}\dif A_{\ve_n}$ are infinite.

\itla{S4}
Ideals of the form $\ifi\oplus\hbox{dense}:$ 
there is a real $\ve>0$ such that 
\imar{define $\oplus$ somewhere}%
the set $A_{\ve}$ is infinite, 
$\mrn(\dop{A_{\ve}})=\piy,$ and 
$\tlim_{n\to\iy\,,\;n\in\dop{A_{\ve}}}r_n=0$.
\een 

In the sense of $\reb,$ all ideals of types \ref{S2}, 
\ref{S3}, \ref{S4} are equivalent to each other, 
and all ideals of type \ref{S1} are equivalent to each other, 
so that we have just $2$ summable ideals modulo $\eeb,$ 
namely, $\ifi$ and $\sui{1/n}.$ 
The structure under $\orb$ or $\obe$ is much more 
complicated (Farah ?).

\vyk{

\punk{Above the summable ideals}
\label{>s}

Is there any example of Borel ideals $\cI\reb\cJ$ 
which do not satisfy $\cI\odl\cJ$? \  
Typically the reductions found 
to witness $\cI\reb\cJ$ are \dd\sd homomorphisms, and 
even better maps. 
The following lemma proves that Borel reduction yields 
\dd\orbpp reduction in quite a representative case. 
Let us say that $\cI\orbpp\cJ$ {\it holds exponentially\/} 
if there is a map $i\mapsto w_i$ withessing
$\cI\orbpp\cJ$~\snos
{Thus we have pairwise disjoint finite non-empty sets
$w_k\sq\dN$ (assuming $\cI,\cJ$ are ideals over $\dN$)
such that $A\in\cI\leqv w_A=\bigcup_{k\in A}w_k\in\cJ,$
and $\tmax w_k<\tmin w_{k+1}$.}
and in addition a sequence of natural numbers $k_i$ with 
$w_i\sq\ir{k_i}{k_{i+1}}$ and $k_{i+1}\ge 2k_i$. 

\bte
\label{t:>s}
Suppose that\/ $r_n\ge 0\yt r_n\to0\yt\sum_nr_n=\piy.$
Then
\ben
\renu
\itsep
\itla{t>s1}
{\rm(Farah~\cite[2.1]{f-tsir})} \ 
If\/ $\cJ$ is a Borel P-ideal and\/ $\sui{r_n}\reb\cJ$
then we have\/ $\sui{r_n}\orbpp\cJ$ exponentially$;$

\itla{t>s2}
{\rm(Hjorth~\cite{h-ban})} \
$\sui{r_n}$ is not 
Borel-reducible to\/ $\zo$.
\een
\ete
\bpf
\ref{t>s1}
Let a Borel $\vt:\pn\to\pn$ witness $\sui{r_n}\reb\cJ.$ 
Let, according to Theorem~\ref{sol}, $\nu$ be a \lsc\ 
submeasure on $\dN$ with $\cJ=\Exh_\nu.$ 
\vyk{
We are going to ``continualize'' $\vt$ in certain way: 
the new, continuous reduction will be a superposition of 
two continuous maps, $\ga$ and $\xi,$ where $\ga$ will 
reduce $\pn/\sui{r_n}$ to a certain ``generic'' part of 
itself, while $\xi$ will reduce the latter to $\pn/\cJ.$ 
}%
The construction makes use of stabilizers. 
%
Suppose that $n\in\dN.$ 
If $u,\,v\sq\ir0n$ then $(u\cup X)\sd(v\cup X)\in\sui{r_n}$ 
for any $X\sq\ir n\piy,$ hence, 
$\vt(u\cup X)\sd\vt(v\cup X)={u\sq v}\in\cJ.$ 
It follows, by the choice of the submeasure $\nu,$
that for any $\ve>0$ 
there are numbers $n'>k>n$ and a set $s\sq\ir n{n'}$ 
such that    
\dm
\nu\skl(\vt(u\cup s\cup X)\sd\vt(v\cup s\cup X))
\cap\ir k\iy\skp\,<\,\ve
\dm
holds for all $u,\,v\sq\ir0n$ 
and all generic~\snos
{In the course of the proof, ``generic'' means Cohen-generic
over a 
sufficiently large countable model of a big enough fragment 
of $\ZFC$.}
$X\sq\ir{n'}\iy$. 

This allows us to define an increasing sequence of 
natural numbers 
$0=k_0=a_0<b_0<k_1<a_1<b_1<k_2<...$ and, for any 
$i,$ a set $s_i\sq\ir{b_i}{a_{i+1}}$ such that, for all 
generic $X,\,Y\sq\ir {b_{i+1}}\iy$ and all $u,\,v\sq\ir0{a_i},$ 
we have 
\ben
\tenu{(\arabic{enumi})}
\itla{stab1}\msur
$\nu\skl(\vt(u\cup s_i\cup X)\sd\vt(v\cup s_i\cup X))
\cap\ir{k_{i+1}}\iy\skp<2^{-i}$; 

\itla{stab2}\msur
$\skl\vt(u\cup s_i\cup X)\sd\vt(u\cup s_i\cup Y)\skp
\cap\ir0{k_{i+1}}=\pu$;

\itla{stabg}
any $Z\sq\dN,$ satisfying $X\cap\ir{b_i}{a_{i+1}}=s_i$ 
for infinitely many $i,$ is generic;

\itla{stabk}\msur
$k_{i+1}\ge 2k_i$ for all $i$;
\een
and in addition, under the assumptions on $\sis{r_n}{}$, 
\ben
\tenu{(\arabic{enumi})}
\addtocounter{enumi}4
\itla{stab3}
there is a set $g_i\sq\ir{a_i}{b_i}$ such that 
$|r_i-\sum_{n\in g_i}r_n|<2^{-i}$.
\een
It follows from \ref{stab3} that 
$A\mapsto g_A=\bigcup_{i\in A}g_i$ 
is a reduction of $\sui{r_n}$ to $\sui{r_n}\res N,$ where 
$N=\bigcup_i\ir {a_i}{b_i}.$
Let $S=\bigcup_i s_i;$ note that $S\cap N=\pu.$

Put $\xi(Z)=\vt(Z\cup S)\sd \vt(S)$ for any $Z\sq N.$ 
\imar{why $\sd\vt(S)$ added?}%
Then, for any sets $X,\,Y\sq N,$ 
\dm
{X\sd Y \in\sui{r_n}}\eqv
{\vt(X\cup S)\sd\vt(Y\cup S)\in\cJ}\eqv
{\xi(X)\sd\xi(Y)\in\cJ,}
\dm
thus $\xi$ reduces $\sui{r_n}\res N$ to $\cJ.$ 
Now put $w_i=\xi(g_i)\cap\ir{k_{i}}{k_{i+1}}$ 
and $w_A=\bigcup_{i\in A} w_i.$ 
We assert that the map $i\mapsto w_i$ proves 
$\sui{r_n}\orbpp\cJ.$ 
In view of the above, 
it remains to show that $\xi(g_A)\sd w_A\in\cJ$ for 
any $A\in\pn$.

As $\cJ=\Exh_\nu,$ it suffices to demonstrate that 
$\nu\skl w_i\sd(\xi(g_A)\cap\ir{k_{i}}{k_{i+1}})\skp
<2^{-i}$ 
for all $i\in A$ while 
$\nu(\xi(g_A)\cap\ir{k_{i}}{k_{i+1}})<2^{-i}$ 
for $i\nin A.$ 
After dropping the common term $\vt(S),$    
it suffices to check that 
\ben
\tenu{(\alph{enumi})}
\itla{111}
$\nu\skl (\vt(g_i\cup S)\sd\vt(g_A\cup S))
\cap\ir{k_{i}}{k_{i+1}}\skp<2^{-i}$ 
for all $i\in A$ while 

\itla{222}
$\nu\skl(\vt(S)\sd\vt(g_A\cup S))\cap
\ir{k_{i}}{k_{i+1}}\skp<2^{-i}$ 
for $i\nin A.$ 
\een
Note that, as any set of the form $X\cup S,$ where $S\sq N,$ 
is generic by \ref{stabg}. 
It follows, by \ref{stab2}, that we can assume, in \ref{111} 
and \ref{222}, that $A\sq[0,i],$ \ie, resp.\ $\tmax A=i$ 
and $\tmax A<i.$ 
We can finally apply \ref{stab1}, with 
$u=A\cup\bigcup_{j<i}s_j,\msur$ $X=\bigcup_{j>i}s_j,$ and 
$v=u_i \cup\bigcup_{j<i}s_j$ if $i\in A$ while 
$v=\bigcup_{j<i}s_j$ if $i\nin A$.

\ref{t>s2}
Otherwise $\sui{r_n}\orbpp\zo$ exponentially by \ref{t>s1}. 
Let this be witnessed by $i\mapsto w_i$ and a sequence of 
numbers $k_i,$ so that $k_{i+1}\ge 2k_i$ and 
$w_i\sq\ir{k_i}{k_{i+1}}.$ 
If $d_i=\frac{\#(w_i)}{k_{i+1}}\to0$ 
then easily $\bigcup_iw_i\in\zo$ by the choice of $\sis{k_i}{}.$ 
Otherwise there is a set $A\in\sui{r_n}$ such that 
$d_i>\ve$ for all $i\in A$ and one and the same $\ve>0,$ 
so that $w_A=\bigcup_{i\in A}w_i\nin\zo.$ 
In both cases we have a contradiction with the assumption 
that the map $i\mapsto w_i$ witnesses $\sui{r_n}\orbpp\zo$.
\epf

\bqu[{{\rm Farah~\cite{f-tsir}}}]
For which ideals other than summables
Theorem~\ref{t:>s}\ref{t>s1} is true~? 
(Farah points out that it also holds for $\ofi$.)
\equ
 
\vyk{
\bcor[{{\rm Hjorth ?}}]
\label{s2d}
If\/ $r_n$ are as in the lemma then\/ $\sui{r_n}$ is not 
Borel-reducible to\/ $\zo$.
\ecor
\bpf

\epf
}

}

\punk{A useful lemma}

\vyk{
The two lemmas below, of quite general nature, 
will be useful below.

\ble
\label{is}
Let\/ $\rF$ be a\/ $\is11$ \er, $Y\sq\dom\rF$ be an\/ 
\ddf invariant\/ $\ip11$ set, and\/ $X\sq Y$ be a\/ 
$\is11$ set.  
Then there is a\/ $\id11$ \ddf invariant set\/ 
$Z$ with\/ $X\sq Z\sq Y$. 
\ele
\bpf
By Separation there is a $\id11$ set $Z_1$ with 
$X\sq Z_1\sq Y.$  
The $\is11$ set 
$X_1=[Z_1]_{\rE}=\ans{y:[y]_{\rE}\cap Z_1\ne\pu}$ 
(the \ddf{\it saturation\/} of $Z_1$) 
satisfies $Z_1\sq X_1\sq Y$ 
because $Y$ is \ddf invariant. 
Applying Separation again, we get a $\id11$ set 
$Z_2$ with $X_1\sq Z_2\sq Y.$ 
Let $X_2=[Z_2]_{\rF}.$ 
As above there is a $\id11$ set $Z_3$ with 
$X_2\sq Z_3\sq Y.$ 
{\sl Et cetera\/}. 
The set $Z=\bigcup_n Z_n=\bigcup_n X_n$ is 
\ddf invariant and $X\sq Z\sq Y$ by the construction. 
Moreover, the choice of the $\id11$ sets $Z_n$ can 
be made effective enough for $Z$ itself to be $\id11,$ 
not merely Borel.
\epf
}

\vyk{
\ble
\label{s2p}
Let\/ $\rF$ is a\/ $\ip11$ \er\/ on a\/ $\ip11$   
set\/ $Y,$ and\/ $\rE$ be a\/ $\is11$ \er\/ on a\/ 
$\is11$ set\/ $X\sq Y,$ and\/ $\rE={\rF\res X}.$ 
Then there is a\/ $\id11$ \er\/ $\rD\sq\rF$ on a\/ 
$\id11$ set\/ $Z$ such that\/ $X\sq Z\sq Y$ and 
still\/ $\rE={\rD\res X}$. 

If moreover\/ $\rE$ is a countable \er\ then it can 
be required that\/ $\rD$ is countable.
\ele
\bpf
By Separation there is a $\id11$ relation $P_1$ with 
$\rE\sq P_1\sq\rF.$ 
Then the \der{\it hull\/} $\rD_1$ of $P_1$ 
(\ie, the least \er\ includeing $P_0$) 
is $\is11,$ $P_1\sq\rD_1\sq\rF,$ 
and still $\rE={\rD_1\res X}.$ 
In the same manner, we choose a $\id11$ relation $P_2$ 
with $\rD_1\sq P_2\sq\rF,$ and let $\rD_2$ be the 
$\is11$ \der hull of $P_2.$ 
Then choose a $\id11$ $P_3$ with $\rD_2\sq P_3\sq\rF,$ 
and let $\rD_3$ be the \der hull of $P_3.$ 
{\sl Et cetera\/}. 
The relation $\rD=\bigcup_n\rD_n=\bigcup_n P_n$ 
is a \er\ as the union of an increasing sequence of 
\er s, and easily $\rE={\rD\res X}.$  
Moreover we can arrange the choice of $\id11$ sets $P_n$ 
effectively enough for $\rD=\bigcup_n P_n$ itself to be 
$\id11,$ not merely Borel. 

To prove the ``moreover'' part assume that $\rE$ is 
countable. 
By the above, both $\rF$ and $Y$ can be assumed 
already $\id11.$ 
As ``to be an index of a countable section'' (of a 
$\id11$ set) is a $\ip11$ notion, the set 
$Y'$ of all $y\in Y$ with at most countable $[y]_{\rF}$ 
is $\ip11,$ and $X\sq Y'$ because $\rE={\rF\res X}$ 
is countable. 
It remains to apply the main part of the lemma to 
the $\ip11$ \er\ ${\rF'}={\rF\res Y'}$.
\epf
}

\ble[{{\rm Attributed to Kechris in \cite{h-ban}}}]
\label{KT}
Suppose that\/ $A,\,X$ are Borel sets, 
$\rE$ is a\/ Borel \er\ on\/ $A,$ and\/ 
$\rho:A\to X$ is a Borel map satisfying the 
following$:$ first, the\/ \dd\rho image of any\/ 
\dde class is at most countable, 
secong,\/ \dd\rho images of any different\/ 
\dde classes are disjoint. 
Then\/ $\rE$ is an essentially countable \er.
\ele
\bpf
The relation: 
$x\rR y$ iff $x,\,y\in Y$ belong to the \dd\rho image 
of one and the same \dde class in $A,$ 
is a $\fs11$ \er\ on the set $Y=\ran \vt,$ 
moreover, 
\dm
\rR\sq P=\ans{\ang{x,y}:\neg\;\sus a,\,b\in A\;
(a\nE b\land x=\rho(a)\land y=\rho(b))}\,,
\dm
where $P$ is $\fp11,$ hence, there is a Borel set $U$ 
with $\rR\sq U\sq P,$ in particular, $U\cap(Y\ti Y)=\rR.$ 
As all \dd\rR equivalence classes are at most countable, 
we can assume that all cross-sections of $U$ are 
at most countable, too.   

Now it suffices to find a Borel \er\ $\rF$ with 
$\rR\sq\rF\sq U.$ 
%
Say that a set $Z\sq X$ is ``stable'' if $U\cap(Z\ti Z)$ 
is a \er, for example, $Y$ is ``stable''. 
We observe that the set 
$D_0=\ans{y:Y\cup\ans{y}\,\hbox{ is ``stable''}}$
is $\fp11$ and satisfies $Y\sq D_0,$ hence, there is a 
Borel set $Z_1$ with $Y\sq Z_1\sq D_0.$ 
Similarly, 
\dm
D_1=\ans{y'\in Z_1:Y\cup\ans{y,y'}\,
\hbox{ is ``stable'' for any }\,y\in Z_1}
\dm
is $\fp11$ and satisfies $Y\sq D_1$ by the definition 
of $Z_1,$ so that there is a Borel set $Z_2$ with 
$Y\sq Z_2\sq D_1.$ 
Generally, we define 
\dm
D_n=\ans{y'\in Z_n:Y\cup\ans{y_1,...,y_n,y'}\,
\hbox{ is ``stable'' for all }\,y_1,...,y_n\in Z_n}\,
\dm
find that $Y\sq D_n,$ and choose a Borel set $Z_n$ 
with $Y\sq Z_n\sq D_n.$ 
Then, by the construction, $Y\sq Z=\bigcap_nZ_n,$ 
and, for any finite 
$Z'\sq Z,$ the set $Y\cup Z'$ is ``stable'', so that 
$Z$ itself is ``stable'', and we can take 
$\rF=U\cap(Z\ti Z)$.
\epf

\punk{Under the summable ideal}
\label S

Subsets of $\dN$ will be systematically identified 
with their characteristic functions.

For $a,\,b\in \dn$ put $a\sd b=\ans{n:a(n)\ne b(n)}$ 
(identified with the function $c(n)=1$ iff $a(n)\ne b(n)$) 
and $\sm(a,b)=\sum_{n\in a\sd b}\frac1{n+1}.$ 
(This can be a nonnegative real or $\piy.$) 
Generally, we define 
$\sm_k^m(a,b)=\sum_{n\in a\sd b\,,\;k\le n\le m}\frac1{n+1},$
and accordingly
$\sm_k^\iy(a,b)=\sum_{n\in a\sd b\,,\;k\le n<\iy}\frac1{n+1}.$ 
Define $\sm(a)=\sum_{\ans{n:a(n)=1}}\frac1{n+1}$ 
and similarly $\sm^m_k(a)$ and $\sm^\iy_k(a)$. 


Recall that the {\it summable ideal\/} is defined as 
\dm
\sun=\ans{a\in\dn:\sm(a)<\piy}\,.
\dm
(The notation $\cI_2$ and $\cI_0$ is also used.) 
$\esn$ will denote the associated Borel \er\ on $\dn,$ 
\ie, $a\esn b$ iff $\sm(a,b)<\piy$. 

\bte
\label h
Let\/ $\rE$ be a Borel \er\ on a Polish space\/ $\dX,$ 
and\/ $\rE\reb \esn.$ 
Then either\/ $\rE\eqb\esn$ or\/ $\rE$ is essentially 
countable.
\ete 
\bpf
This is a long proof. 
Let $\vt:\dX\to\dn$ be a Borel reduction $\rE$ to 
$\esn.$ 
We can assume that $\vt$ is in fact continuous: indeed 
it is known that there is a stronger Polish topology 
on $\dX$ which makes $\vt$ continuous but does not add 
new Borel subsets of $\dX.$ 
Now, as any Polish $\dX$ is a $1-1$ continuous image 
of a closed subset of $\bn,$ we can assume that 
$\dX=\bn$. 

Finally, we can assume that $\vt$ is $\id11,$ not 
merely Borel.

If $a\in A\sq\dn$ and $q\in\dqp$ then let $\gal qAa$ 
be the set of all $b\in A$ such that there is a finite 
chain $a=a_0,\,a_1,\,...,\,a_n=b$ of reals $a_i\in A$ 
such that $\sm(a_i,a_{i+1})<q$ for all $i,$ the 
\dd q{\it galaxy of\/ $a$ in\/ $A$.\/} 

\bdt
\label{d:good}
A set $A\sq\dn$ is {\it\gra q\/}, where 
$q\in\dqp,$ iff $\sm(a,b)<1$ for all $a\in A$
and $b\in\gal qAa.$ 
A set $A$ is {\it\grap\/} if it is \gra q for some $q\in\dqp.$ 
(In other words it is required that the galaxies are 
rather small.) 
\edt


\bct
\label{cl1}
Any\/ \gra q $\is11$ set\/ $A\sq\dn$ is covered  
by a\/ \gra q $\id11$ set. 
\ect 
\bpf\footnote
{\label{rp}\ 
The result can be achieved as a routine application 
of a reflection principle, yet we would like to show 
how it works with a low level technique.} 
The set 
$D_0=\ans{b\in\dn:A\cup\ans b\,\hbox{ is \gra q}}$ 
is $\ip11$ and $A\sq D_0,$ hence, there is a $\id11$ set 
$B_1$ with $A\sq B_1\sq D_0.$ 
Note that $A\cup\ans a$ is \gra q for any $a\in B_1.$ 
It follows that the $\ip11$ set 
\dm
D_1=\ans{b\in B_1:A\cup\ans{a,b}\,
\hbox{ is \gra q for any }\,a\in B_1}
\dm
still contains $A,$ hence, there is a $\id11$ set 
$B_2$ with $A\sq B_2\sq D_1\sq B_1.$ 
Note that $A\cup\ans{a_1,a_2}$ is \gra q for any 
$a_1,\,a_2\in B_2.$ 
In general, as soon as we have got a $\id11$ set $B_n$ 
with $A\sq B_n$ and such that $A\cup\ans{a_1,...,a_n}$ is 
\gra q for any $a_1,...,a_n\in B_n,$ then 
the $\ip11$ set 
\dm
D_{n}=\ans{b\in B_n:A\cup\ans{a_1,...,a_n,b}\,
\hbox{ is \gra q for any }\,a_1,...,a_n\in B_n}
\dm
contains $A,$ hence, there is a $\id11$ set 
$B_{n+1}$ with $A\sq B_{n+1}\sq D_{n}\sq B_n$.

As usual in similar cases, the choice of the sets $B_n$ 
can be made effective enough for the set 
$B=\bigcap_nB_n$ to be still $\id11,$ not merely Borel. 
On the other hand, $A\sq B$ and $B$ is \gra q.
\epF{Claim}

Coming back to the proof of the theorem, 
let $C$ be the union of all \grap\ $\id11$ sets. 
An ordinary computation shows that $C$ is $\ip11.$ 
We have two cases.\vtm

{\bfit Case 1\/\bf:} \ \msur$\ran\vt\sq C.$ 
\hfill
{\bfit Case 2\/\bf:} \ otherwise.\hfill $\;$

\punk{Case 1}

We are going to prove that, in this case, $\rE$ is 
essentially countable.
First note that, 
by \Sepa, there is a $\id11$ set $\aH\sq\dn$ 
with $\ran \vt\sq \aH\sq C$. 

Fix a standard enumeration $\sis{W_e}{e\in E}$ of all 
$\id11$ subsets of $\dn,$ where, as usual, $E\sq\dN$ is 
a $\ip11$ set. 
By \Kres, there exist $\id11$ functions $a\longmapsto e(a)$ 
and $a\longmapsto q(a),$ defined on $\aH,$ 
such that for any $a\in\aH$ the $\id11$ set 
$W(a)=W_{e(a)}$ contains $a$ and is \gra{q(a)}. 
The final point of our argument will be an application 
of Lemma~\ref{KT}, where $\rho$ will be a derivate of 
the function $G(a)=\gal{q(a)}{W(a)}a.$ 
We prove 

\bct
\label{cl3}
If\/ $a\in\aH$  then\/ 
$\ga_a=\ans{G(b):b\in [a]_{\esn}\cap\aH}$ 
is at most countable.
\ect
\bpf
Otherwise there is a pair of $e\in E$ and $q\in\dqp$ 
and an uncountable set $B\sq [a]_{\esn}\cap\aH$ 
such that $q(b)=q$ and $e(b)=e$ for any $b\in B$ and 
$G(b')\ne G(b)$ for any two different $b,b'\in B.$ 
Note that any $G(b),\msur$ $b\in B,$ is a \dd qgalaxy 
in one and the same set $W(a)=W(b)=W_e,$ therefore, 
if ${b\ne b'}\in B$ then $b'\nin G(b)$ and 
$\sm(b,b')\ge q.$ 
On the other hand, as $B\sq [a]_{\esn},$ we have 
$\sm(a,b)<\piy$ for all $b\in B,$ hence, there is  
$m$ and a still uncountable set $B'\sq B$ such that 
$\smy m (a,b)<q/2$ for all $b\in B'.$ 
Now take a pair of ${b\ne b'}\in B'$ with 
${b\res \ir0m}={b'\res\ir0m}:$ then  
$\sm(b,b')<q,$ contradiction.
\epF{Claim}

It follows that $x\mapsto G(\vt(x))$ maps any 
\dde class into a countable set of galaxies $G(a).$ 
To code the galaxies by single points, let 
$S(a)=\bigcup_m\ens{b\res m}{b\in G(a)}.$ 
Thus $S(a)\sq\bse$ codes the Polish topological
closure of the galaxy $G(a)$.

\bct
\label{cl2}
If\/ $a,\,b\in\aH$ and\/ $\neg\;a\esn b$ then\/ $b$ 
does not belong to the (topological) closure of\/ 
$G(a),$ in particular,\/ $b\res m\nin S(a)$ for some $m$.
\ect
\bpf 
Take $m$ big enough for $\smp0m(a,b)\ge 2.$ 
Then $s={b\res m}$ does not belong to $S(a)$ because 
any $a'\in G(a)$ satisfies $\sm(a,a')<1$.
\epF{Claim}

Elementary computation shows that the sets 
\dm
{\bf G}=\ans{\ang{a,b}:a\in\aH\land b\in G(a)}
\quad\hbox{and}\quad
{\bf S}=\ans{\ang{a,s}:a\in\aH\land s\in S(a)}\,.
\dm
belong to $\is11,$ but this is not enough to claim 
that $a\mapsto S(a)$ is a Borel map. 
Yet we can change it appropriately to get a Borel map 
with similar properties. 
First of all define the following $\is11$ \er\ 
on $\aH$:
\dm
a\rF b\qquad\hbox{iff}\qquad
e(a)=e(b)\land q(a)=q(b)\land G(a)=G(b)\,.
\dm
(To see that $\rF$ is $\is11$ note that here $G(a)=G(b)$ 
is equivalent to $b\in G(a),$ and that $\bf G$ is $\is11.$)
It follows from Claim~\ref{cl2} and \Kres\  
that there is a $\id11$ function $\mu:\aH\ti\aH\to\dN$ 
such that for any pair of $a,\,b\in\aH$ with $a\nsn b$ 
we have $b\res \mu(a,b)\nin S(a).$ 
Then the set
\dm
R(a)=\ans{b\res \mu(a',b):
a',b\in\aH \land {a\rF a'}\land {a'\nsn b})}\,\sq\,\bse
\dm
does not intersect $S(a),$ for any $a\in\aH,$ hence, 
the $\is11$ set 
\dm
{\bf R}=\ans{\ang{a,s}:a\in\aH\land s\in R(a)}
\dm 
does not intersect ${\bf S}.$ 
Note that by definition $\bf R$ is \ddf invariant \vrt\ 
the 1st argument, \ie, if $a,\,a'\in\aH$ satisfy 
$a\rF a'$ then $R(a)=R(a').$ 
It follows from Lemma~\ref{d20} that there is a 
$\id11$ set ${\bf Q}\sq\aH\ti\bse$ with 
${\bf S}\sq{\bf Q}$ but ${\bf R}\cap{\bf Q}=\pu,$ 
\ddf invariant in the same sense. 
Then the map $a\mapsto Q(a)=\ans{s:{\bf Q}(a,s)}$ is $\id11$. 

\bct
\label{cl4}
Suppose that\/ $a,\,b\in\aH.$ 
Then$:$ ${a\rF b}$ implies\/ ${Q(a)=Q(b)}$ and\/ ${a\nsn b}$ 
implies\/ ${Q(a)\ne Q(b)}$.
\ect
\bpf 
The first statement holds just because $Q$ is 
\ddf invariant. 
Now suppose that $a\nsn b.$ 
Then by definition $s=b\res \mu(a,b)\in R(a),$ hence, 
$s\nin Q(a).$ 
On the other hand, $s\in S(b)\sq Q(b)$.
\epF{Claim}

Define $\tau(x)=Q(\vt(x))$ for $x\in\bn,$ so that $\tau$ 
is a $\id11$ map $\bn\to\cP(\bse)$. 

\bct
\label{cl5}
If\/ $x\in\bn$ then\/ 
$T_a=\ans{\tau(y):y\in\eke x}$ is at most countable. 
\ect
\bpf 
Suppose that $y,\,z\in\eke x.$ 
Then $a=\vt(x),\msur$ $b=\vt(y),$ and $c=\vt(z)$ belong to 
$\aH,$ and $b,\,c\in[a]_{\esn}.$ 
It follows from Claim~\ref{cl4} that if $G(b)=G(c),\msur$ 
$e(b)=e(c),$ and $q(b)=q(c),$ then $Q(b)=Q(c).$ 
It remains to note that $G$ takes only countably many 
values on $\aH\cap \ek a{\esn}$ by Claim~\ref{cl3}. 
\epF{Claim}

Finally note that, if ${x\nE y}\in\bn$ then 
$\vt(x),\,\vt(y)$ belong to $\aH$ and satisfy 
$\vt(x)\nsn \vt(y),$ hence, $\tau(x)\ne \tau(y)$ by 
Claim~\ref{cl4}. 
Thus, the Borel map $\tau$ witnesses that the given 
\er\ $\rE$ is 
essentially countable by Lemma~\ref{KT}.

\punk{Case 2}
 
Thus we suppose that the $\is11$ set 
$\aB=\ran\vt\dif C$ is non-empty. 
Note that, by Claim~\ref{cl1}, there is no non-empty 
$\is11$ \grap\ set $A\sq\aB$.\vom 

Let $\cB_s=\ans{a\in\dn:s\su a}$ for $s\in\bse$ and 
$\cN_u=\ans{x\in\dnn:u\su x}$ for $u\in\nse$ 
(basic open nbhds in $\dn$ and $\dnn$). 

If $A,\,B\sq\dn$ and $m,\,k\in\dN,$ then 
$A\xE km B$ will mean that for any $a\in A$ there is $b\in B$ 
with $\smy{k}(a,b)<2^{-m},$ and conversely, for any $b\in B$ 
there is $a\in A$ with $\smy{k}(a,b)<2^{-m}.$ 
This is not a \er, of course, yet the conjunction of 
$A\xE km B$ and $B\xE km C$ implies $A\xE k{m-1} C$.

$0^m$ will denote the sequence of $m$ zeros. 

To prove that $\esn\reb\rE$ in Case 2, we define 
an increasing sequence of natural numbers 
$0=\yk0<\yk1<\yk2<...,$ and also objects 
$A_s,\:g_s,\:v_s$ for any $s\in\bse,$ 
which satisfy the following list of requirements 
\ref{7i} -- \ref x.  

\ben
\tenu{(\roman{enumi})}
\itla{7i}\msur
if $s\in 2^m$ then $g_s\in 2^{\yk m},$ and 
${s\su t}\imp {g_s\su g_t}$; 

\itla2\msur
$\pu\ne A_s\sq\aB\cap\cB_{g_s},$ $A_s$ is $\is11,$ 
and $s\su t\imp A_t\sq A_s$.

\itla3
if $s\in 2^n$ then $A_{0^n}\xE{\yk n}{n+2} A_s$;

\itla5
if $s\in 2^n,\msur$ $m<n,\msur$ $s(m)=0,$ then 
$\smp{\yk{m}}{\yk{m\plo}}(g_s,g_{0^m})<2^{-m-1}$;

\itla6
if $s\in 2^n,\msur$ $m<n,\msur$ $s(m)=1,$ then 
$|\smp{\yk{m}}{\yk{m\plo}}(g_s,g_{0^m})-\frac1{m+1}|<2^{-m-1}$;

\itla7
if $s\yo t\in 2^n\yt m<n\yt s(m)=t(m),$ then 
$|\smp{\yk{m}}{\yk{m\plo}}(g_s,g_{t})|<2^{-m}$;

\vyk{
\itla8 
sets $A_s$ are enough generic;
}

\itla9
if $s\in2^n$ then $v_s\in\dN^n,$ and $s\su t\imp v_s\su v_t$;

\itla x\msur
$A_s\sq 
\ans{a\in\aB:\vt\obr(a)\cap\cN_{v_s}\ne\pu}$.
\een


We can now accomplish Case~2 as follows. 
For any $a\in\dn$ define 
$F(a)=\bigcup_n g_{a\res n}\in\dn$ 
(the only element satisfying $g_{a\res n}\su F(a)$ 
for all $n$) 
\vyk{
(this is where \ref8 is applied to provide the 
non-emptiness of the intersection)
}%
and $\rho(a)=\bigcup_n v_{a\res n}\in\bn.$ 
It follows, by \ref x and the continuity of $\vt,$ 
that $F(a)=\vt(\rho(a))$ for any $a\in\dn.$ 
Thus the next claim proves that $\rho$ is a 
Borel (in fact, here continuous) 
reduction $\esn$ to $\rE$ and ends Case~2. 

\bct
\label{cl8}
The map\/ $F$ reduces\/ $\esn$ to\/ $\esn,$ that is, 
the equivalence\/ 
${a\esn b}\eqv{F(a)\esn F(b)}$ holds for all\/ 
$a,\,b\in\dn$.
\ect
\bpf
By definition 
$\sm(F(a),F(b))=
\tlim_{n\to\iy}\smp0{\yk n}(g_{a\res n},g_{b\res n}).$ 
\vyk{
As $\ga_s=g_{0^n}\sd g_s$ for $s\in 2^n,$ we can 
re-write this as follows:
\dm
\sm(F(a),F(b))=
{\textstyle\tlim_{n\to\iy}}\:\sm(\ga_{a\res n},\ga_{b\res n})\,.
\dm
}
However it follows from \ref5, \ref6, \ref7 that 
\dm
|\smp0{\yk n}(g_{a\res n},g_{b\res n})-
\smp0{n}(a\res n,b\res n)|
\,\le\,
{\textstyle\sum_{m<n}}2^{-m}\,<\,2\,.
\dm
We conclude that 
$|\sm(F(a),F(b)) - \sm(a,b)|\le2,$ as required.
\epF{Claim}

\punk{Construction}

The construction goes on by induction. 
To begin with we set $\yk0=0,\msur$ $g_\La=\La$ and 
$A_\La=\aB.$ 
Suppose that, for some $n,$ we have the objects as 
required for all $n'\le n,$ and extend the construction 
on the level $n+1$.

As $A_{0^n}$ is not \grap\ (see above), there is a 
pair of elements $a^0,\,a^1\in A_{0^n}$ such that 
$|\sm(a^0,a^1)-\frac1{n+1}|<2^{-n-2}.$ 
Note that $a^0\res {\yk n}=a^1\res {\yk n}$ by \ref{7i} and 
\ref2, hence, there is $\yk{n\plo}>\yk n$ such that 
$|\smp{\yk n}{\yk{n\plo}}(a^0,a^1)-\frac1{n+1}|<2^{-n-2}.$ 
According to \ref3, for any $s\in2^n$ there exist 
$b^0_s,\,b^1_s\in A_s$ such that  
and $\smy{\yk n}(a^i,b^i_s)<2^{-n-2}$ for $i=0,1;$ 
we can, of course, assume that $b^i_{0^n}=a^i.$ 
Moreover, the number $\yk{n+1}$ can be chosen big enough 
for the following to hold: 
\dm
\smy{\yk{n\plo}}(b^i_s,a^0)< 2^{-n-3}
\quad\hbox{--- \ for all}\quad s\in 2^n
\quad\hbox{and}\quad i=0,1.
\eqno(1)
\dm

We let $g_{s\we i}=b^i_s\res {\yk{n\plo}}$ for all  
$s\we i\in 2^{n+1}.$ 
This definition preserves \ref{7i}. 
To check \ref5 for $s'=s\we 0\in2^{n+1}$ and $m=n,$ 
note that
\dm
\smp{\yk{n}}{\yk{n\plo}}(g_{s'},g_{0^{n+1}})=
\smp{\yk{n}}{\yk{n\plo}}(b^0_s,a^0)<2^{-n-2}.
\dm 
To check \ref6 for $s'=s\we 1\in2^{n+1}$ and $m=n,$ 
note that
\dm
|\smp{\yk{n}}{\yk{n\plo}}(g_{s'},g_{0^{n+1}})-
{\textstyle\frac1{n+1}}| \,\le \,
\smp{\yk{n}}{\yk{n\plo}}(b^1_s,a^1)
+|\smp{\yk n}{\yk{n\plo}}(a^0,a^1)-{\textstyle\frac1{n+1}}|
\,<\,2^{-n-1}.
\dm 

To fulfill \ref9, choose, for any $s\we i\in 2^{n+1},$ 
a sequence $v_{s\we i}\in\dN^{n+1}$ so that $v_s\su v_{s\we i}$ 
and there is $\cN_{v_{s\we i}}\cap \vt\obr(b^i_{s})\ne\pu$. 

Let us finally define the sets $A_{s'}\sq A_s,$ 
for all $s'=s\we i\in 2^{n+1}$ 
(so that $s\in2^n$ and $i=0,1$).  
To fulfill \ref2 and \ref x, we begin with 
\dm
A'_{s\we i}=\ans{a\in A_s\cap\cB_{g_{s\we i}}:
\vt\obr(a)\cap\cN_{v_{s\we i}}\ne\pu}\,.
\dm 
This is a $\is11$ subset of $A_s,$ containing $b^i_{s}.$ 
To fulfill \ref3, we define $A_{0^{n+1}}$ to be the set of 
all $a\in A'_{0^{n+1}}$ such that 
\dm
\kaz s'=s\we i\in 2^{n+1}\;\sus b\in A'_{s'}\;
\skl\smy{\yk{n\plo}}(a,b)<2^{-n-3}\skp\;;
\dm
this is still a $\is11$ set containing $b^0_{0^n}=a^0$ by 
(1). 
It remains to define, for any $s\we i\ne 0^{n+1},$ 
$A_{s\we i}$ to be the set of 
all $b\in A'_{s\we i}$ such that 
\dm
\sus b\in A_{0^{n+1}}\;
\skl\smy{\yk{n\plo}}(a,b)<2^{-n-3}\skp\;.
\dm
This ends the definition for the level $n+1$.

\vyk{

According to \ref3, for any $s\in2^n$ there exist 
$h^0_s,\,h^1_s\in H_{\yk{n}}(2^{-n-2})$ such that 
\dm
b^0_{s}=\ovp s\sd h^0_s\sd a_0
\quad\hbox{ and }\quad
b^1_{s}=\ovp s\sd h^1_s\sd a_m=
\ovp s\sd h^1_s\sd h\sd a_0
\dm
belong to $A_s.$ 
We can assume that $h^0_{0^n}=h^1_{0^n}=0^\dN$ 
(the constant $0$), so that 
$b^0_{0^n}=a_0,\msur$ $b^1_{0^n}=a_m$. 
Now, as $h^0_s$ and $h^1_s$ belong to 
$H_{\yk{n}}(2^{-n-2}),$ 
there is a number $\yk{n\plo}>\yk{n}$ big enough for 
the following to hold~\footnote
{\ Here Hjorth also adds inequalities 
$\sm_{\ge \yk{n+1}}(h^0_s) < 2^{-n-3}$ and 
$\sm_{\ge \yk{n+1}}(h^1_s) < 2^{-n-3},$ which does 
not seem to be of any use.-- ?}:
\pagebreak[0]%
\ben
\tenu{(\roman{enumi})}
\addtocounter{enumi}8
\itla y
$|\sm(h\res\ir{\yk n}{\yk{n\plo}})-\frac1{n+1}| 
< 2^{-n-2}$.
\een

\vyk{
\dm 
\left. 
\bay{rcll}
\sm_{\ge \yk{n+1}}(h^0_s) &<& 2^{-n-3} &;\\[1ex]

\sm_{\ge \yk{n+1}}(h^1_s) &<& 2^{-n-3} &;\\[1ex]

|\sm(h\res\ir{\yk n}{\yk{n\plo}})-\frac1{n+1}| 
&<& 2^{-n-2} &.
\eay
\right\}
\eqno(\ast)
\dm
}

We let $g_{s\we 0}=b^0_s\res \ir0{\yk{n\plo}}$ and 
$g_{s\we 1}=b^1_s\res \ir0{\yk{n\plo}}.$ 
This definition preserves \ref{7i}. 
To check \ref5 for $s'=s\we 0\in2^{n+1}$ and $m=n,$ 
note that
\dm
\sm(\ga_{s'}\res\ir{\yk n}{\yk{n+1}}) =
\sm(h^0_s\res\ir{\yk n}{\yk{n+1}})<2^{-n-2}
\dm 
since $h^0_s\in H_{\yk{n}}(2^{-n-2}).$ 
To check \ref6 for $s'=s\we 1\in2^{n+1}$ and $m=n,$ 
note that
\dm
|\sm(\ga_{s'}\res\ir{\yk n}{\yk{n+1}})-
{\textstyle\frac1{n+1}}| \,\le \,
\sm(h^1_s\res\ir{\yk n}{\yk{n+1}}) 
+|\sm(h\res\ir{\yk n}{\yk{n+1}})-
{\textstyle\frac1{n+1}}|
\,<\,2^{-n-1}
\dm 
because $h^0_s\in H_{\yk{n}}(2^{-n-2})$ 
and by \ref y.  
The verification of \ref7 for $s',\,t'\in 2^{n+1}$ and 
$m=n$ is similar. 
}

\vtm

\epF{Construction and Theorem~\ref{h}}

\parf{\protect\co equalities}
\label{co}

Suppose that $\stk{X_k}{d_k}$ is a finite metric space 
for each $k\in\dN.$ 
Farah~\cite{f-co} defines an equivalence relation 
$\rD=\rD(\stk{X_k}{d_k})$ on 
$\dX=\prod_{k\in\dN}X_k$ as follows: 
\dm
x\rD y\quad\hbox{ iff }\quad
\tlim_{k\to\iy}\,d_k(x_k,y_k)=0\,.
\dm
\er s of this form are called \co{\it equalities\/}. 
In addition, $\rD(\stk{X_k}{d_k})$ is {\it nontrivial\/} 
if $\tlis_{k\to\iy}\dia(X_k)>0$ 
(otherwise $\rD(\stk{X_k}{d_k})$ makes everything 
equivalent). 
Every \co equality is easily a Borel \er, more exactly, 
of class $\fp03$.

\punk{Some examples and simple results}
\label{co1}

\bex
(1) 
Let $X_k=\ans{0,1}$ with $d_k(0,1)=1$ for all $k.$ 
Then clearly the relation 
$\rD(\stk{X_k}{d_k})$ on $\dn=\prod_k\ans{0,1}$ 
is just $\Eo$.\vom

(2) 
Let $X_{kl}=\ans{0,1}$ with $d_{kl}(0,1)=k\obr$ for all 
$k,\,l\in\dN.$ 
Then the relation 
$\rD(\stk{X_{kl}}{d_{kl}})$ on 
$2^{\dN\ti\dN}=\prod_{k,l}\ans{0,1}$ 
is just $\Et=\rE_{\ofi}$.\vom

(3) 
Generally, if $0=n_0<n_1<n_2<...$ and $\vpi_i$ is a 
submeasure on $\il i{i+1},$ then let $X_i=\pws{\il i{i+1}}$ 
and $d_i(u,v)=\vpi_i(u\sd v)$ for $u,\,v\sq \il i{i+1}.$ 
Then $\rD(\stk{X_{i}}{d_{i}})$ is clearly isomorphic 
to $\rei,$ where 
\dm
\cI=\Exh(\vpi)=
\ans{x\sq\dN:\tlim_{n\to\iy}\vpi(x\cap\ir n\iy)=0}
\dm
and $\vpi(x)=\tsup_i\vpi_i(x\cap\il i{i+1})$.\vom

(4) 
Let $\rdm=\rD(\stk{X_k}{d_k}),$ where 
$X_k=\ans{0,\frac1k,\frac2k,...,1}$ and  
$d_k$ is the distance on $X_k$ inherited from $\dR.$ 
\eex

\bpro[{{\rm Farah~\cite{f-co} with a reference to Oliver}}]
\label{comax}
\ben
\tenu{{\rm(\roman{enumi})}}
\itla{comax1} 
$\rdm\eqb\rzo\,;$  

\itla{comax2} 
if\/ $\rD$ is a\/ \co equality then\/ 
$\rD\reb\rdm,$ moreover, $\rD\rea\rdm$.
\een
\epro
Thus $\rdm$ is a maximal, in a sense, among \co equalities. 
\bpf
\ref{comax1} 
It is clear that $\rdm$ is the same as ${\fco}\res\dX,$ 
where $\dX\sq\rtn$ is defined as in the proof of 
Lemma~\ref{co=d}, where it is also shown that  
${\fco}\eqb{{\fco}\res\dX}$. 

\ref{comax2} 
To prove $\rD\reb\rdm,$ it suffices, by \ref{comax1} and 
Lemma~\ref{co=d}, to show that $\rD\reb{\fco}.$
The proof is based on the following: 

\bct
\label{comax'}
Any finite\/ \dd nelement metric space\/ $\stk Xd$ is 
isometric to an\/ \dd nelement subset of\/ 
$\stk{\dR^n}{\rho_n},$ where\/ $\rho_n$ be the distance 
on\/ $\dR^n$ 
defined by\/ $\rho_n(x,y)=\tmax_{i<n}|x_i-y_i|$.
\ect
\bpc
Let $X=\ans{x_1,\dots,x_n}.$ 
It suffices to prove that for any $k\ne l$ there is a 
set of reals $\ans{r_1,\dots,r_n}$ 
such that $|r_k-r_l|=d(x_k,x_l)$ and 
\bit
\item[$(\ast)$]
$|r_i-r_j| \le d_{ij}=d(x_i,x_j)$ for all $i\zi j$. 
\eit
We can assume that $k=1$ and $l=n.$ 

{\sl Step 1\/}. 
There is a least number $h_1\ge0$ such that 
$(\ast)$ holds for the numbers 
$\ans{\underbrace{0,0,\dots,0}_{n-1\;\,\text{times}},h}$ 
for any $0\le h\le h_1.$ 
Then, for some $k,\,\;1\le k<n,$ we have   
$h_1-0=d_{kn}$ exactly. 
Suppose that $k\ne1;$ 
then it can be assumed that $k=n-1.$\vom 

{\sl Step 2\/}. 
Similarly, there is a least number $h_2\ge0$ such that 
$(\ast)$ holds for the numbers 
$\ans{\underbrace{0,0,\dots,0}_{n-2\;\,\text{times}},h,h_1+h}$ 
for any $0\le h\le h_2.$ 
Then, for some $k,\nu,\,\;1\le k<n-1\le\nu\le n,$ 
we have $h_2-0=d_{k\nu}$ exactly. 
Suppose that $k\ne1;$ 
then it can be assumed that $k=n-2.$\vom 

{\sl Step 3\/}. 
Similarly, there is a least number $h_3\ge0$ such that 
$(\ast)$ holds for the numbers 
$\ans{\underbrace{0,0,\dots,0}_{n-3\;\,\text{times}},
h,h_2+h,h_1+h_2+h}$ 
for any $0\le h\le h_3.$ 
Then again, for some $k,\nu,\,\;1\le k<n-2\le\nu\le n,$ 
we have $h_3-0=d_{k\nu}$ exactly. 
Suppose that $k\ne1;$ 
then it can be assumed that $k=n-3.$\vom 

{\sl Et cetera\/}.\vom

This process ends, after a number $m$ ($m<n$) steps, 
in such  a way that the index $k$ obtained at the final 
step is equal to $1.$
Then $(\ast)$ holds for the numbers 
$\ans{\underbrace{0,0,\dots,0}_{n-m\;\,\text{times}},
r_{n-m+1},r_{n-m+1},\dots,r_n},$ 
where $r_{n-m+j}=h_m+h_{m-1}+\dots+h_{m-j+1}$ for each 
$j=1,\dots m.$ 
Moreover it follows from the construction that there is 
a decreasing sequence $n=k_0>k_1>k_2>\dots>k_\mu=1$ 
($\mu\le m$) 
such that $r_{k_{i}}-r_{k_{i+1}}=d_{k_{i+1},k_{i}}$ 
exactly for any $i.$ 
Then $d_{1n}\le \sum_i r_{k_{i}}-r_{k_{i+1}}$ by the 
triangle inequality. 
But the right-hand side is a part of the sum 
$r_n=h_1+\dots+h_m,$ and hence $r_n\ge d_{1n}.$ 
It follows that, cutting the construction at an 
appropriate step $m'\le m$) 
(and taking an appropriate value of $h\le h_{m'}$), 
we obtain a sequence of numbers 
$r_1=0\le r_2\le\dots\le r_{n-1}\le r_n$ 
still satisfying $(\ast)$ and satisfying 
$r_n=r_n-r_0=d_{1n}.$ 
This ends the proof. 
\epc

Now, to carry out the proof of $\rD\reb{\fco},$   
suppose that $\rD=\rD(\stk{X_k}{d_k})$ is an equivalence 
relation on $\dX=\prod_{k\in\dN}X_k,$ where each 
$\stk{X_k}{d_k}$ is a finite metric space.
Let $n_k$ be the number of elements in $X_k.$ 
Let, by the claim, $\eta_k:X_k\to \dR^{n_k}$ be an 
isometric embedding of $\stk{X_k}{d_k}$ into 
$\stk{\dR^{n_k}}{\rho_{n_k}}.$ 
The map 
$\vt(x)= \eta_0(x_0)\we\eta_1(x_1)\we\eta_2(x_2)\we\dots$  
(from $\dX$ to $\rtn$) 
reduces $\rD$ to $\fco$.
\epf

The structure of \co equalities tend to be  
connected more with the additive reducibility $\rea$ 
(see \prf{reli} on $\rea$ and the associated relations 
$\reas$ and $\eqa$) than with the general 
Borel reducibility. 
In particular, we have 

\ble
\label{<co}
For any\/ \co equality\/ $\rD=\rD(\stk{X_k}{d_k}),$ 
if\/ $\rD'$ is a Borel \er\ on\/ a set\/ 
$\prod_kX'_k$ (with finite nonempty\/ $X'_k$) 
and\/ $\rD'\rea\rD$ then\/ $\rD'$ is a\/ 
\co equality.
\ele
\bpf
Let a sequence $0=n_0<n_1<n_2<...$ and a collection 
of maps $H_i:X'_i\to\prod_{n_i\le k<n_{i+1}}X_k$ 
witness $\rD'\rea\rD.$ 
For $x',\,y'\in X'_i$ put 
\dm
d'_i(x',y')=
\tmax_{n_i\le k<n_{i+1}}d_k(H_i(x')_k,H_i(y')_k)\,.
\dm
Then easily $\rD'=\rD(\stk{X'_k}{d'_k})$. 
\epf

\ble[{{\rm Farah~\cite{f-co} with a reference to Hjorth}}]
\label{coact}
Every \co equality\/ $\rD=\rD(\stk{X_k}{d_k})$ 
is induced by a continuous action of a Polish group.
\ele
(The domain $\dX=\prod_kX_k$ of $\rD$ 
is considered with the product topology.)
\bpf(sketch) \ 
For any $k$ let $S_k$ be the (finite) group of all 
permutations of $X_k,$ with the distance  
$\rho_k(s,t)=\tmax_{x\in X_k}\,d_k(s(x),t(x)).$ 
Then
\dm
\dG=\ans{g\in{\TS\prod_kS_k}:
\tlim_{k\to\iy}\rho_k(g_k,e_k)=0}\,,
\quad\hbox{where $e_k\in S_k$ is the identity}\,,
\dm
is easily a subgroup of $\prod_kS_k,$ moreover, 
the distance  
$d(g,h)=\tsup_k\rho_k(g_k,h_k)$ converts $\dG$ into a 
Polish group, the natural action of which on 
$\dX$ 
(\ie, $(g\app x)_k=g_k(x_k),\msur$ $\kaz k$) 
is continuous and induces $\rD$.
\epf

\punk{Classification}
\label{coclass}

Recall that for a metric space $\stk Ad,$ a rational 
$q>0,$ and $a\in A,$ $\gal qAa$ is the set of all 
$b\in A$ which can be connected with $a$ by a finite 
chain $a=a_0,a_1,...,a_n=b$ with $d(a_i,a_{i+1})<q$ 
for all $i.$ 
Farah defines, for $r>0,$ 
\dm
\da(r,A)=\tinf\,
\ans{q\in\dqp:\sus a\in A\;(\dia(\gal qAa)\ge r)}
\dm 
(with the understanding that here $\tinf\pu=\piy$), 
and
\dm
\Da(A)=\ans{d(a,b):a\ne b\in A}\,,\quad
\text{so that}\quad 
\dia A=\tsup(\Da(A)\cup\ans0)\,.
\dm

Now let $\rD=\rD(\stk{X_k}{d_k})$ be a \co equality 
on $\dX=\prod_{k\in\dN}X_k.$ 
The basic properties of $\rD$ are determined by the 
following two conditions:

\ben
\tenu{{\rm(co\arabic{enumi})}}
\itla{*1}
$\tlii_{k\to\iy}\da(r,X_k)=0$ for some $r>0$.\enuci

\itla{*2}
$\kaz\ve>0\;\sus \ve'\in(0,\ve)\;\susi k\;
\skl\Da(X_k)\cap\ir{\ve'}\ve\ne\pu\skp$.
\een
Easily \ref{*1} implies both the nontriviality 
of $\rD=\rD(\stk{X_k}{d_k})$ and \ref{*2}.

\bte[{{\rm Farah~\cite{f-co}}}]
\label{coT}
Let\/ $\rD=\rD(\stk{X_k}{d_k})$ be a nontrivial 
\co equality. 
Then
\imar
{Comment upon turbulent in \ref{c03}.}%
\ben
\tenu{{\rm(\roman{enumi})}}
\itla{c01}
If\/ \ref{*2}, hence,\/ \ref{*1} fail  then\/ 
$\rD\eqa\Eo,$ hence,\/ $\rD\eqb\Eo\;;$

\itla{c02}
If\/ \ref{*1} fails but\/ \ref{*2} holds then\/ 
$\rD\eqa\Et,$ hence,\/ $\rD\eqb\Et\;;$

\itla{c03}
If\/ \ref{*1}, hence,\/ \ref{*2} hold then\/ 
$\Eo\reas\rD$ and\/ $\rD_1\rea\rD$ for a turbulent 
\co equality\/ $\rD_1$ satisfying\/ $\Et\rea\rD_1$.%
\een
\ete
\bpf
\ref{c01}
To show that $\Eo\rea\rD$ note that, by the 
nontriviality of $\rD,$ there exist: a number 
$p>0,$ an increasing sequence $0=n_0<n_1<n_2<...\;,$ 
and, for any $i,$ a pair of points 
$x_{n_i},\,y_{n_i}\in X_{n_i}$ with 
$d_{n_i}(x_{n_i},y_{n_i})\ge p.$ 
For $n$ not of the form $n_i$ fix an arbitrary 
$x_n\in X_n.$ 
Now, if $a\in\dn,$ then define $\vt(a)\in\prod_kX_k$ 
so that $\vt(a)_n=z_n$ for $n$ not of the form $n_i,$ 
while $\vt(a)_{n_i}=x_{n_i}$ or $=y_{n_i}$ if resp.\ 
$a_i=0$ or $=1.$ 
This map $\vt$ witnesses $\Eo\rea\rD$.

Now prove that $\rD\rea\Eo.$ 
As \ref{*2} fails, there is $\ve>0$ such that for each 
$\ve'$ with $0<\ve'<\ve$ we have only finitely many $k$ 
with the propery that $\ve'\le d_k(\xi,\eta)<\ve$ 
for some $\xi,\,\eta\in X_k.$ 
Let $G_k$ be the (finite) set of all 
\dd{\frac\ve2}galaxies in $X_k,$ and let 
$\vt:\dX=\prod_kX_k\to G=\prod_kG_k$ be defined as 
follows: 
$\vt(x)_k$ is that galaxy in $G_k$ to which $x_k$ 
belongs. 
Let $\rE$ be the \dd Gversion of $\Eo,$ \ie, if 
$g,\,h\in G$ then 
$g\rE h$ iff $g_k=h_k$ for all but finite $k.$ 
As easily $\rE\rea\Eo,$ it suffices to demonstrate that 
$\rD\rea\rE$ via $\vt.$ 
Suppose that $x,\,y\in \dX$ and $\vt(x)\rE\vt(y)$ and 
prove $x\rD y$ (the nontrivial direction). 
Let, on the contrary, $x\nD y,$ so that there is a number 
$p>0$ with $d_k(x_k,y_k)>p$ for infinitely many $k.$ 
We can assume that $p<\frac\ve2.$ 
On the other hand, as $\vt(x)\rE\vt(y),$ there is 
$k_0$ such that $x_k$ and $y_k$ belong to one and the 
same \dd{\frac\ve2}galaxy in $X_k$ for all $k>k_0.$ 
Then, for any $k>k_0$ with $d_k(x_k,y_k)>p$ 
(\ie, for infinitely many values of $k$) 
there exists an element $z_k\in X_k$ in the same 
galaxy such that $p<d_k(x_k,z_k)<\ve,$ but this is a 
contradiction to the choice of $\ve$ 
(indeed, take $\ve'=p$).\vom

\ref{c02} 
Let us show first that if \ref{*2} holds then 
$\Et\rea\rD$ (independently of \ref{*1}).
It follows from \ref{*2} that there exist: 
an infinite sequence $\ve_1>\ve_2>\ve_3>...>0,$ 
for any $i$ an infinite set $J_i,$ 
and for any $j\in J_i$ a pair of elements 
$x_{ij},\,y_{ij}\in X_j$ with 
$d_j(x_{ij},y_{ij})\in\ir{\ve_{i+1}}{\ve_i}.$ 
We may assume that the sets $J_i$ are pairwise 
disjoint. 
Then the \co equality 
$\rD'=\rD(\stk{\ans{x_{ij},y_{ij}}}{d_j}
_{i\in\dN,\:j\in J_i})$ 
satisfies both $\rD'\rea\rD$ and $\rD'\isi\Et$ 
(via a bijection between the underlying sets). 

Now, assuming that, in addition, \ref{*1} fails, 
we show that $\rD\rea\Et.$ 
For all $k,\,n\in\dN$ let $G_{kn}$ be the (finite) 
set of all \dd{\frac1n}galaxies in $X_k.$ 
For any $x\in \dX=\prod_iX_i$ define 
$\vt(x)\in G=\prod_{k,n}G_{kn}$ so that 
$\vt(x)_{kn}$ is that \dd{\frac1n}galaxy in 
$G_{kn}$ to which $x_k$ belongs (for all $k,\,n$). 
The \er\ $\rE$ on $G,$ defined so that 
$g\rE h$ iff $\kaz n\;\kazi k\;(g_{kn}=h_{kn})$ 
($g,\,h\in G$) is easily $\rea\Et,$ so it suffices 
to show that $\rD\rea\rE$ via $\vt.$ 
Suppose that $x,\,y\in \dX$ and $\vt(x)\rE\vt(y)$ 
and prove $x\rD y$ (the nontrivial direction). 
Otherwise there is some $r>0$ with $d_k(x_k,y_k)>r$ 
for infinitely many $k.$ 
As \ref{*1} fails for this $r,$ there is $n$ big 
enough for $\da(r,X_k)>\frac1n$ to hold for almost 
all $k.$ 
Then, by the choice of $r,$ we have 
$\vt(x)_{kn}\ne\vt(y)_{kn}$ for infinitely many $k,$ 
hence, $\vt(x)\nE\vt(y),$ contradiction.\vom

\ref{c03} 
Fix $r>0$ with $\tlii_{k\to\iy}\da(r,X_k)=0.$ 
As for any increasing sequence $n_0<n_1<n_2<...$ 
we have $\rD(\stk{X_{n_i}}{d_{n_i}})\rea\rD,$ 
it can be assumed that $\tlim_{k}\da(r,X_k)=0,$ 
and further that $\da(r,X_k)<\frac1k$ for all $k.$ 
Then every $X_k$ contains a \dd{\frac1k}galaxy 
$Y_k\sq X_k$ of $\dia Y_k\ge r.$ 
As easily $\rD(\stk{Y_k}{d_k})\rea\rD,$ 
the following lemma suffices to prove \ref{c03}.

\blt
\label{coturb}
Suppose that\/ $r>0$ and each\/ $X_k$ is a single\/ 
\dd{\frac1k}galaxy in itself with\/ $\dia(X_k)\ge r.$  
Then\/ $\rD=\rD(\stk{X_k}{d_k})$ is turbulent and\/ 
$\Et\rea\rD$.
\elt
\bpf
We know from the proof of \ref{c03} above 
that $\Et\rea\rD.$
Now prove that the natural action of the Polish group 
$\dG$ defined as in the proof of Lemma~\ref{coact} is 
turbulent under the assumptions of the lemma. 

That every \dd\rD class is dense in $\dX=\prod_kX_k$ 
(with the product topology on $\dX$) is an easy exercise. 
To see that every \dd\rD class $[x]_{\rD}$ also is 
meager in $\dX,$ note that by the assumptions of the lemma 
any $X_k$ contains a pair of elements $x'_k,\,x''_k$ with 
$d_k(x'_k,x''_k)\ge r.$ 
Let $y_k$ be one of $x'_k,\,x''_k$ which is 
\dd{d_k}fahrer than $\frac r2$ from $x_k.$ 
Now the set $Z=\ans{z\in \dX:\susi k\;(z_k=y_k)}$ 
is comeager in $\dX$ and disjoint from $[x]_{\rD}.$ 
It remains to prove that local orbits are somewhere 
dense.

Let $G$ be an open nbhd of the identity in $\dG$ and 
$\pu\ne X\sq\dX$ be open in $\dX.$ 
We can assume that, for some $n,$ $G$ is the 
\dd{\frac1n}ball around the identity in $\dG$ while 
$X=\ans{x\in\dX:\kaz k<n\;(x_k=\xi_k)},$ where 
elements $\xi_k\in X_k,\;k<n,$ are fixed. 
It is enough to prove that all classes of the local 
orbit relation $\ler GX$ are dense in $X.$ 
Consider an open set 
$Y=\ans{y\in\dX:\kaz k<m\;(y_k=\xi_k)}\sq X,$ where 
$m>n$ and elements $\xi_k\in X_k,\;n\le k<m,$ are fixed 
in addition to the above. 

Let $x\in X.$ 
Then $x_k=\xi_k$ for $k<n.$ 
Let $n\le k<m.$ 
The elements $\xi_k$ and $x_k$ belong to $X_k,$ 
which is a \dd{\frac1k}galaxy, therefore, there is 
a chain, of a length $\ell(k),$ 
of elements of $X_k,$ which connects $x_k$ 
and $\xi_k$ so that every step within the chain has 
\dd{d_k}length $<\frac1k.$ 
Then there is a permutation $g_k$ of $X_k$ 
such that $g_k^{\ell(k)}(x_k)=\xi_k,$ 
$g_k(\xi_k)=x_k,$ and 
$d_k(\xi,g_k(\xi))<\frac1k$ for all $\xi\in X_k.$ 
Let $g_k$ be the identity on $X_k$ whenever $k<n$ 
or $k\ge m.$ 
This defines an element $g\in\dG$ which obviously 
belongs to $G,$ moreover, $X$ is \dd ginvariant 
and $g^\ell(x)\in U,$ where 
$\ell=\prod_{n\le k<m}\ell(k),$ 
hence, $x\ler GX g(x),$ as required.
\epF{Lemma}

\epF{Theorem~\ref{coT}}

\bre
\label{coTrem}
Theorem~\ref{coT} shows that any nontrivial 
\co equality $\rD$ \dd\rea contains a turbulent 
\co equality $\rD'$ with $\Et\rea \rD'$   
(and the turbulence of $\rD'$ holds, in particular, 
via the natural action defined in the proof of 
Lemma~\ref{coact}), 
unless $\rD$ is $\eqa$ to $\Eo$ or $\Et,$ 
and that \ref{*1} is necessary for the turbulence 
of $\rD$ itself and sufficient for a turbulent 
\co equality $\rD'\rea\rD$ to exist.
\ere

\punk{LV-equalities}
\label{lv}

By Farah, an {\it\lv equality\/} is a \co equality 
$\rD=\rD(\stk{X_k}{d_k})$ satisfying 

\ben
\tenu{{\rm({\sc lv}\arabic{enumi})}}
\itla{lv1}
$\kaz m\;\kaz\ve>0\;\kazi k\;
\kaz x_0,...,x_m\in X_k\;
\skl
d_k(x_0,x_m)\le\tmax_{j<m}d_k(x_j,x_{j+1})+\ve
\skp$.
\een
In other words, the metrics involved are postulated 
to be ``asymptotically close'' to ultrametrics. 
This sort of \co equalities was first considered by 
Louveau and Velickovic \cite{lv}. 
The following simple fact is analogous to 
Lemma~\ref{<co}.

\ble
\label{<lv}
For any\/ \lv equality\/ $\rD,$ 
if\/ $\rD'$ is a Borel \er\ on\/ a set\/ 
$\prod_kX'_k$ (with finite nonempty\/ $X'_k$) 
and\/ $\rD'\rea\rD$ then\/ $\rD'$ is an\/ 
\lv equality.\qeD
\ele

\bex[{{\rm Louveau and Velickovic \cite{lv}}}]
\label{lvex}
We define $X_k=\ans{1,2,...,2^k}$ and 
$d_k(m,n)=\log(|m-n|+1)/k$ for $1\le m,n\le2^k$.
\eex

\bte[{{\rm Essentially, Louveau and Velickovic \cite{lv}}}]
\label{coAB}
Let\/ $\rD=\rD(\stk{X_k}{d_k})$ be a turbulent\/ 
\lv equality. 
Then we can associate, with each infinite\/ $A\sq\dN,$ 
a\/ \lv equality\/ $\rD_A\rea\rD$ such that for 
all\/ $A,\,B\sq\dN$ the following are equivalent$:$
\ben
\tenu{{\rm(\roman{enumi})}}
\itla{AB1} 
$A\sqa B$ {\rm(\ie, $A\dif B$ is finite)};

\itla{AB2} 
$\rD_A\rea\rD_B\;;$

\itla{AB3}
$\rD_A\reB\rD_B$ 
{\rm(\ie, via a Baire measurable reduction)}.
\een
\ete 

This theorem was the first major application 
of \co equalities.
One of its corollaries is that there exist big 
families of mutually irreducible Borel \er s !

\bpf 
As $\rD$ is turbulent, the necessary turbulence 
condition \ref{*1} of \prf{coclass} holds, moreover, 
as in the proof of Theorem\ref{coT} (case~\ref{c03}), 
we can assume that it takes the following special form 
for some $r>0$: 
%
\ben
\tenu{(\arabic{enumi})}
\itla{lv+}
Each $X_k$ is a single 
\dd{\tmin\ans{\frac r{2},\frac1{k+1}}}galaxy 
of $\dia(X_k)\ge4r$.
\een
The intended transformations (reduction to a certain 
infinite subsequence of spaces $\stk{X_k}{d_k},$ 
and then each $X_k$ to a suitable galaxy $Y_k\sq X_k$) 
preserve \ref{lv1}, of course, moreover, going to 
subsequences once again, we can assume that \ref{lv1} 
holds in the following special form: 
%
\ben
\tenu{(\arabic{enumi})}
\addtocounter{enumi}1
\itla{lv5}\msur
$d_k(x_0,x_{m_k})\le
\tmax_{i<m_k}d_k(x_i,x_{i+1})+\frac1{k+1}$ 
whenever $x_0,...,x_{m_k}\in X_k,$ 
where $m_k=2^{\prod_{j=0}^{k-1}\#(X_j)}$.
\een

We can derive the following important consequence:

\ben
\tenu{(\arabic{enumi})}
\addtocounter{enumi}2
\itla{lv4}\msur
For any 
$k$ there is a set $Y_k\sq X_k$ of  
$\#(Y_k)=m_k$ such that we have 
$d_k(x,y)\ge r$ for all $x\ne y$ in $Y_k$.
\een

To prove this note that by \ref{lv+} there is a 
set $\ans{x_0,...,x_m}\sq X_k$ such that 
$d_k(x_0,x_m)\ge 4r$ but 
$d_k(x_i,x_{i+1})<r$ for all $i.$ 
We may assume that $m$ is the least possible length 
of such a sequence $\ans{x_i}.$ 
Now let us define a subsequence $\ans{y_0,y_1,...,y_n}$ 
of $\ans{x_i},$ the number $n\le m$ will be specified 
in the course of the construction. 
Put $y_0=x_0.$ 
If $y_j=x_{i(j)}$ has been defined, and there is 
$l>i(j),\;l\le m,$ such that $d_k(y_j,x_l)\ge r,$ 
then let $y_{j+1}=x_l$ for the least such $l,$ 
otherwise put $n=j$ and stop the construction. 

By definition $d_k(y_j,y_{j+1})\ge r$ for all $j<n,$ 
moreover, $d_k(y_{j'},y_{j+1})\ge r$ for any 
$j'<j$ by the minimality of $m.$ 
Thus $Y_k=\ans{y_j:j\le n}$ satisfies 
$d_k(x,y)\ge r$ for all $x\ne y$ in $Y_k.$ 
It remains to prove that $n\ge m_k.$ 
Indeed we have $d_k(y_j,y_{j+1})<2r$ by the 
construction, hence, if $n\le m_k$ then we would have 
$d_k(y_0,y_n)\le 3r$ by \ref{lv5}, which implies 
$d_k(y_n,x_m)\ge r,$ a contradiction to the assumption 
that the construction stops with $y_n$,

\vyk{
Requirement \ref{lv1}, together with the 
necessary turbulence condition \ref{*1},  
yield the existence of 
arbitrarily large ``discrete'' sets:

\blt 
\label{lvdis}
In the assumptions of the theorem, 
there is\/ $r>0$ such that 
\ben
\tenu{{\rm({\sc lv}\arabic{enumi})}}
\addtocounter{enumi}1
\itla{lv2}
$\kaz m\;\susi k\;
\sus x_0,...,x_m\in X_k\;\kaz i\ne j\le m\;
\skl d_k(x_i,x_j)\ge r\skp$.
\een
\elt
\bpf

\epF{Lemma~\ref{lvdis}}

We can assume that in fact both \ref{*1} and 
\ref{lv2} hold for 
$r=\frac12,$ if necessary, multiplying 
all distances accordingly.  
Moreover, we can assume that \ref{lv1}, \ref{lv2}, 
and \ref{*1} hold in the following stronger forms, 
where $m_j=\#(X_j),$ 
going to a suitable subsequence of the sequence of 
the spaces $\stk{X_k}{d_k},$ if necessary:
}

\vyk{
\ben
\tenu{(\arabic{enumi})}
\vyk{
\itla{lv5}\msur
$d_k(x_0,x_{m_j})\le
{\DS\tmax_{i<m_j}d_k(x_i,x_{i+1})}+\frac1j$ 
whenever $j<k$ and $x_0,...,x_{m_j}\in X_k$.
}
\itla{lv4}
For any $k$ there is a set $Y_k\sq X_k$ of 
$\#(Y_k)>2^{\prod_{j<k}m_j}$ such that we have 
$d_k(x,y)\ge\frac12$ for all $x\ne y$ in $Y_k$.

\een
}

This said, we proceed to the proof of the theorem. 
First note that 

\blt
\label{ab3ab2}
\ref{AB3} implies that\/ \ref{AB2} 
holds at least for some (infinite) $A'\sq A$.
\imar
{Is it true that for a pair of \co equalities 
$\rD,\,\rD',$ if $\rD\reb\rD'$ then $\rD\rea\rD'$?}
\elt
\bpf 
A Borel reduction can be extracted from a Baire 
measurable one by a version of the ``stabilizers'' 
construction (see proofs of ... .) 
\epF{Lemma~\ref{ab3ab2}}

Thus it remains only to show that \ref{AB2} implies 
\ref{AB1}, even simpler, that, for any disjoint 
infinite sets $A,\,B\sq\dN,$ $\rD_A\rea\rD_B$ 
fails. 
Suppose, towards the contrary, that $\rD_A\rea\rD_B$ 
holds, and let this be witnessed by a reduction 
$\Psi$ defined (as in \prf{reli}) from an increasing 
sequence $\tmin B=n_0<n_1<n_2<...$ of numbers 
$n_i\in B$ and a collection of maps 
$H_k:X_k\to \prod_{j\in\il i{i+1}\cap B}X_j\,,\msur$ 
$k\in A.$ 
Let 
\dm
f_k(\da)=
\tmax_{\xi,\,\eta\in X_k,\:d_k(\xi,\eta)<\da}
\;\;\;
\tmax_{j\in\il i{i+1}\cap B}
\;\;\;
d_j(H_k(\xi)_j,H_k(\eta)_j)\,,
\dm 
for  $k\in\dN$ and $\da>0$  
(with the understanding that $\tmax\pu=0$ if 
applicable).
Then $f(\da)=\tsup_{k\in A}f_k(\da)$ 
is a nondecreasing map $\dR^+\to\ir0\iy$.

\blt
\label{clai3}
$\tlim_{\da\to0}f(\da)=0$.
\elt 
\bpf
Otherwise there is $\ve>0$ such that $f(\da)\ge\ve$ 
for all $\da.$ 
Then the numbers
\dm
\mu_k=
\TS\tmin_{\xi,\,\eta\in X_k\,,\;\xi\ne\eta}\,
d_k(\xi,\eta)
\quad\hbox{(all of them are $>0$)}
\dm
must satisfy $\tinf_{k\in A}\mu_k=0.$ 
This allows us to define a sequence 
$k_0<k_1<k_2<...$ of numbers $k_i\in A,$ and, 
for any $k_i,$ a pair of $\xi_i,\,\eta_i\in X_{k_i}$ 
with $d_{k_i}(\xi_i,\eta_i)\to0,$ and also 
$j_i\in\il{k_i}{k_i+1}\cap B$ such that 
$d_{j_i}(H_{k_i}(\xi_i)_{j_i},H_{k_i}(\eta_i)_{j_i})
\ge\ve.$ 
Let $x,\,y\in\prod_{k\in A}X_k$ satisfy 
$x_{k_i}=\xi_i$ and $y_{k_i}=\eta_i$ for all $i$ and 
$x_k=y_k$ for all $k\in A$ not of the form $k_i.$  
Then easily $x\rD_A y$ holds but 
$\Psi(x)\rD_B\Psi(y)$ fails, which is a contradiction.  
\epF{Lemma~\ref{clai3}}

Let $k\in A,$ and let $Y_k\sq X_k$ be as in \ref{lv4}. 
Then there exist elements $x_k\ne y_k$ 
in $Y_k$ such that $H_k(x_k)\res k=H_k(y_k)\res k.$ 
By \ref{lv+} there is a chain 
$x_k=\xi_0,\xi_1,...,\xi_n=y_k$ of elements 
$\xi_i\in X_k$ with $d_k(z_i,z_{i+1})\le\frac1{k+1}$ 
for all $i<n.$ 
Now  
$H_k(\xi_i)\in\prod_{j\in\il i{i+1}\cap B}X_j$ 
for each $i\le n.$ 
Let $j\in\il i{i+1}\cap B.$ 
If $j>k$ then the elements $y^j_i=H_k(\xi_i)_j,\msur$ 
$i\le n,$ satisfy 
$d_j(y^j_i,y^j_{i+1})\le f_k(\frac1{k+1}).$ 
As clearly $n<m_j,$ we conclude that 
$d_j(H_k(x_k)_j,H_k(y_k)_j)\le 
f_k(\frac1{k+1})+\frac1{j+1}$ 
by \ref{lv5}. 
If $j<k$ then simply $H_k(x_k)_j=H_k(y_k)_j$ 
by the choice of $x_k,\,y_k.$ 
Thus totally 
\ben
\tenu{(\arabic{enumi})}
\addtocounter{enumi}3
\itla{lv=}\msur
$d_j(H_k(x_k)_j,H_k(y_k)_j)\le 
f(\frac1{k+1})+\frac1{k+1}$ 
for all $j\in\il i{i+1}\cap B$.
\een 
(as $k\nin B$).  
Let $x=\sis{x_k}{k\in A}$ and 
$y=\sis{y_k}{k\in A},$ both are elements of 
$\prod_{k\in A}X_k,$ and $x\rD_A y$ fails because 
$d_k(x_k,y_k)\ge r$ for all $k.$ 
On the other hand, we have $\Psi(x)\rD_B\Psi(y)$ by 
\ref{lv=}, because $f(\da)\to 0$ with $\da\to0$ 
by Lemma~\ref{clai3}.
This is a contradiction to the assumption that $\Psi$ 
reduces $\rD_A$ to $\rD_B$.\vtm

\epF{Theorem~\ref{coAB}}

\newpage

\parf{$\rtd$ is not reducible to ...}
\label{rtd}

This section contains a theorem saying that the \er\ $\rtd$ 
of equality of countable sets of the reals is not Borel 
reducible to \er s which belong to a family of 
{\it \PP\/} \er s, including, for 
instance, continuous actions of \cli\ groups and some 
ideals, not only Polishable, and is closed under the 
Fubini product modulo $\ifi.$   
But the {\sl prima facie\/} definition of the family is based 
on a rather metamathematical property which we extracted from 
Hjorth~\cite{h:orb}. 

Recall that $\rtd$ is defined on $\nnn$ as follows: 
$x\rtd y$ iff $\ran x=\ran y$.

Suppose that $X$ is $\fs11$ or $\fp11$ in the universe 
$\dV,$ and an extension $\dvp$ of $\dV$ is considered. 
In this case, let $\di X$ denote what results by the 
definition of $X$ applied in $\dvp.$ 
There is no ambiguity here by Shoenfield, and easily 
$X=\di X\cap\dV$.

\punk{\PPP\ \er s do not reduce $\rtd$}
\label{tagcli}

Fix a Polish space $\dX$ and let $\sis{B_n}{n\in\dN}$ 
be a base of its topology.  
By a {\it Borel code\/} for $\dX$ we shall understand 
\index{Borel code}%
a pair $p=\ang{T,f}$ of a wellfounded tree 
$\pu\ne T=T_p\sq\Ord\lom$ (then $\La\in T$) 
and a map $f:\Max T\to\dN,$ where $\Max T$ 
is the set of all \dd\sq maximal elements of $T.$  
\index{zzMaxT@$\Max T$}%
We define $\bk p(t)\sq\dnn$ for any $t\in T$ 
by induction on the rank of $t$ in $T,$ so that 
\bit
\item
$\bk p(t)=\bk{f(t)}$ for all $t\in\Max T,$ and 

\item
$\bk p(t)=\dop\bigcup_{t\we\xi\in T}\bk p(t\we\xi)$ for 
\index{zzcf@$\bk\tau$}%
$t\in T\dif\Max T\,;$ 

\item
finally, put $\bk p=\bk p(\La)$. 
\eit

For a Borel code $p=\ang{T,F},$ let $\tsup p=\tsup T$ 
be the least ordinal $\ga$ with $T\sq\ga\lom.$ 
A code $p$ is {\it countable\/} 
\index{Borel code!countable}%
if $\tsup p<\omi,$ 
in this case the coded set $\bk p$ is a Borel subset 
of $\dX$. 

\bdf
\label{t2like} 
A $\fs11$ \er\ $\rE$ is {\it\PP\/} if,
\index{equivalence relation, ER!pinned@\PP}%
\index{pointp@\PP}%
for any (perhaps, uncountable) Borel code $p,$ 
\poq{if} $\bk p$ is  
2wise \dd{\drE}equivalent in any generic extension   
of $\dV$ 
and non-empty in some generic extension of $\dV,$  
\poq{then} there is a point $x\in\dom\rE,$ ``pinning'' 
$p$ in the sense that $\bk p\sq \ek{x}{\drE}$ in any  
extension of $\dV$. 
\edf

\bct
\label{t2unlike+}
$\rtd$ is \poq{not} \PP.
\ect
\bpf
Consider a Borel code $p$ for the set  
$\ans{x\in\nnn:\ran x=\dnn\cap\dV},$ 
so that $\bk p\sq(\kon^\dV)\lom.$ 
Then \poq{if} of Definition~\ref{t2like} holds, actually, 
$\bk p$ is a \dd\rtd equivalence class in any universe 
where it is non-empty, but \poq{then} fails.  
\epf

\blt
\label{t2unlike}
\vyk{
{\rm(Hjorth~\cite{h:orb})} \ 
\PPP\/ $\fs11$ \er s\/ $\rE$ do not satisfy\/ 
$\rtd\reb\rE$.
}%
If\/ $\rE,\,\rF$ are\/ $\fs11$ \er s, $\rE\reb\rF,$  and\/ 
$\rF$ is \PP, then so is\/ $\rE$. 
\elt
\bpf
Suppose that, in $\dV,$ $\vt:\dX\to\dY$ is a Borel
reduction of $\rE$ to $\rF,$ where $\dX=\dom\rE$ and 
$\dY=\dom\rF.$    
We can assume that $\dX$ and $\dY$ are just two copies of 
$\dn.$ 
Let $r$ be a (countable) Borel code for $\vt$ as a subset 
of $\dX\ti\dY.$ 
Let $p$ be a Borel code satisfying \poq{if} of 
Definition~\ref{t2like}. 
There is perhaps no Borel code $q$ such that 
$\bk q={\bk r}\ima \bk p$ everywhere, but still 
there is a code $q$ with $\bk q\sq{\bk r}\ima \bk p$ and 
$B_q\ne\pu$ somewhere. 
Indeed, let, in $\dV,$ $\la=\card(\tsup p)$ and 
$\ka=\la^+$ (the next cardinal). 
Consider the formula $A(p,r,y)$ saying:
\bit
\item\msur
$y\in\dY$ and there is a forcing term 
$\tau\in \bL[p,r,y]$ such that the forcing 
$\text{\sc Coll}(\dN,\la)$ forces 
$\tau[\poq G]\in\bk p$ and $y=\bk r(\tau[\poq G])$.
\eit
As it is known, 
\imar{rfrnce~?}%
there is a Borel code $q$ such that 
\dm
\bk q=\ans{y:\bL_{\ka}[p,r,y]\mo A(p,r,y)}
\dm 
in \poq{an}y extension of $\dV.$ 
Then easily $\bk q\sq\bk r\ima \bk p,$ 
hence, $\bk q$ is 2wise \dd{\drF}equivalent in
any universe, in addition, $\bk q$ is nonempty 
somewhere. 
\vyk{
if $x$ is a generic collapse  
$\dN\onto\dnn\cap\dV$ then $y=\bk r(x)$ satisfies 
$A(p,r,y),$ hence, belongs to $\bk q$.  
}

As $\rF$ is \PP, there is, in $\dV,$ a point $y\in\dY$ 
such that $\bk q\sq \ek{y}{\drF}$ holds, in particular, 
in \dd{\text{\sc Coll}(\dN,\la)}generic extension 
$\dvp$ of $\dV,$ where $\bk q\ne\pu,$
\vyk{
Suppose, on the contrary, that $\rE$ is \PP, \ie,  
there is $y\in \dnn\cap\dV$ such that 
$\bk q\sq \ek{y}{\drE}$ holds, in particular, 
in \dd{{\tt Coll}(\dN,\cont^\dV)}generic extension 
$\dvp$ of $\dV,$ where $\bk q\ne\pu,$
}%
hence, there is $x\in\bk p\cap\dvp$ with 
${y\drF\bk r(x)}.$ 
It follows, by Shoenfield, that 
${y\rF \vt(x')}$ for some $x'\in\dX$ in $\dV.$ 
Thus $x\drE x',$ which implies that $x'\in\dV$ 
pins $p,$ as required.
\epf

\punk{Fubini product of \PP\ \er s is \PP}
\label{fbbsbs}

Recall that the Fubini product $\rE=\fps{\ifi}{\rE_k}{k\in\dN}$ 
of \er s $\rE_k$ on $\dnn$ modulo $\ifi$ is a \er\ on $\nnn$ 
defined as follows: $x\rE y$ if $x(k)\rE_k y(k)$ for all 
but~finite~$k$.

\bpro
\label{bsfub}
The family of all \PP\/ $\fs11$ \er s is closed 
under Fubini products modulo\/ $\ifi$.
\epro
\bpf
Suppose that \er s $\rE_k$ on $\dnn$ are \PP; 
prove that the Fubini product 
$\rE=\fps{\ifi}{\rE_k}{k\in\dN}$ 
is \PP. 
Define $x\rF_k y$ iff $x(k)\rE_k y(k):$ $\rF_k$ are 
$\fs11$ \er s on $\nnn$ and $x\rE y$ iff $x\rF_k y$ for 
almost all~$k$.

\bct
\label{bsfub'}
Each\/ $\rF_k$ is \PP.
\ect
\bpf
Consider a Borel code $p$ for a subset of $\nnn,$ 
satisfying \poq{if} of Definition~\ref{t2like} \vrt\ $\rF_k.$  
By the same argument as in the proof of Lemma~\ref{t2unlike}, 
there is a Borel code $q$ for a subset of $\dnn,$ such that 
$\bk q\ne\pu$ in some extension of $\dV$ and 
$\bk q\sq\ans{x(k):x\in\bk p}$ in any extension of $\dV,$ 
hence, $q$ satisfies \poq{if} of Definition~\ref{t2like} 
\vrt\ $\rE_k.$ 
As $\rE_k$ is \PP, there is $a\in\dnn$ such that 
$\bk q\sq\ek{a}{\di\rE_k}$ in any extension, but then easily 
$\bk p\sq\ek{x}{\di\rF_k}$ in any extension, where 
$x\in\nnn\cap\dV$ has only to satisfy $x(k)=a$ for the given 
$k$.
\epF{Claim}

In continuation of the proof of the proposition, 
consider a Borel code $p$ for a subset of $\nnn,$ 
satisfying \poq{if} of Definition~\ref{t2like} \vrt\ $\rE.$ 
Our plan is to find another Borel code $\bap$ with 
$\bk\bap\sq\bk p$ everywhere, which satisfies 
\poq{if} of Definition~\ref{t2like} for almost all $\rE_k.$ 
This involves a forcing by Borel codes.

Let, in $\dV,$ $\la=\tsup p$ and $\ka=\la^+,$ 
thus, $\tsup p<\ka.$  
Let $\dP$ be the set of all Borel codes $q\in\dV$ for  
subsets of $\nnn$ such that $\tsup q<\ka$ 
and $\bk q\ne\pu$ in a generic extension of the 
universe $\dV.$ 
$\dP$ is considered as a forcing, with   
$q\lef p$ ($q$ is stronger) iff $\bk q\sq \bk p$ 
in all generic extensions of $\dV.$  
It is known that $\dP$ forces a point of $\nnn,$ so that  
$\bigcap_{q\in G}\bk q=\ans{x_G}$ for any \dd\dP generic, 
over $\dV,$ set $G\sq\dP.$ 
Let $\dox$ be the name of the generic element of $\nnn$.

By the choice of $p,$ $\ang{p,p}$ \dd{\dP\ti\dP}forces 
$\doxl\drE\doxr,$ hence, there are codes $q,\,r\in\dP$ 
and a number $k_0$ such that $\ang{q,r}$ 
\dd{\dP\ti\dP}forces 
$\doxl\mathbin{\di{\rF_k}}\doxr$ for any $k\ge k_0.$ 
By a standard argument, we have 
$x\mathbin{\di{\rF_k}} y$ for all $k\ge k_0$ 
in any extension of $\dV$ for any two \dd\dP generic, over 
$\dV,$ elements $x,\,y\in\bk q.$ 
We can straightforwardly define in $\dV$ a Borel code 
$\bap$ (perhaps, not a member of $\dP$!) 
such that, in 
any extension of $\dV,$ $\bk\bap$ is the set of all 
\dd\dP generic, over $\dV,$ elements of $\bk q.$ 
Then $\bap$ satisfies \poq{if} of Definition~\ref{t2like} 
\vrt\ any $\rF_k$ with $k\ge k_0.$ 
Hence, by the claim, there is, in $\dV,$ a sequence of points 
$x_k\in\nnn$ such that $\bk\bap \sq\ek{x_k}{\di\rF_k}$ 
in any generic extension of $\dV,$ for any $k\ge k_0.$ 
Define $x\in\nnn\cap\dV$ so that $x(k)=x_k(k)$ for any 
$k\ge k_0,$ then, by the definition of $\rF_k,$ we have  
$\bk\bap\sq\ek{x}{\di\rF_k}$ for all $k\ge k_0$ 
in any extension of $\dV.$ 
Yet 
$\bigcap_{k\ge k_0}\ek{x}{\di\rF_k}\sq\ek{x}{\drE}$.
\epF{Proposition}

\punk{Complete left-invariant actions produce 
\PP\ \er s}
\label{bs:cli}

Recall that a Polish group $\dG$ is 
{\it complete left-invariant\/}, \cli\ for brevity, 
if $\dG$ admits a compatible  
\index{group!cli@\cli}%
\index{group!complete left-invariant}%
\index{cli@\cli}%
left-invariant complete metric. 
Then easily $\dG$ also admits a compatible 
\poq{ri}g\poq{ht}-invariant 
complete metric, which will be practically used.

\bte
\label{clit}
{\rm(Hjorth~\cite{h:orb})} \ 
Suppose that\/ $\dG$ is a Polish \cli\ group continuously 
acting on a Polish space\/ $\dX.$  
Then\/ $\egx$ is \PP, hence, 
$\rtd$ is not Borel reducible to\/ $\egx$. 
\ete
\bpf
Fix a Borel code $\wtau$ satisfying \poq{if} of 
Definition~\ref{t2like} \vrt\ $\egx.$  
Let $\ka$ be a cardinal in $\dV$ satisfying  
$\tsup{\wtau}<\ka.$ 
Define forcing $\dP$ as above, thus,  
$\dP$ forces an element of $\dX.$ 

\vyk{
Say that a set $X\sq\dX$ (of class $\fs11$ or $\fp11$) 
is {\it forcible over $q\in\dP$\/} 
if there is $r\in\dP,\msur$ $r\lef q,$ which forces 
$\dox\in \di X.$ 
For instance, any $\bk q$ is forcible, which is witnessed 
by $q$ itself.

\bct
\label{l2}
If a set\/ $X=\bigcup_nX_n\sq\dX$ is forcible over\/ 
$q\in\dP$ then so is at least one of\/ $X_n$.
\ect
\bpf
By definition some $r\in\dP,\msur$ $r\lef q,$ forces 
$\dox\in \di X.$ 
If none of $X_n$ is forcible then $r$ forces 
$\dox\nin \di{X_n}$ for each $n,$ leading to a contradiction.
\epF{Claim} 
}

Let $\rho$ be a compatible \poq{ri}g\poq{ht}-invariant metric 
on $\dG$. 

For any $\ve>0,$ let $G_\ve=\ans{g\in\dG:\rho(g,1_\dG)<\ve}.$ 
Say  that $q\in\dP$ {\it is of size\/ $\le\ve$\/} 
if $\ang{q,q}$ \dd{(\dP\ti\dP)}forces that 
there is $g\in \di{G_\ve}$ with $\doxl=g\app \doxr.$
In this case, in any generic extension of the universe, 
if $\ang{x,y}\in \bk q\ti \bk q$ is a 
\dd{(\dP\ti\dP)}generic pair then there is 
$g\in \di{G_{\ve}}$ with $y=g\app x.$ 
\vyk{
(Indeed, take $z\in \bk q$ to be generic enough for both 
$\ang{x,z}$ and $\ang{y,z}$ to be \dd{\dP\ti\dP}generic 
pairs. 
Then $x=f\app z$ and $y=h\app z$ for some 
$f,\,h\in \di{G_\ve},$ so that $y=hf\obr \app x.$ 
But $\rho(hf\obr,1_\dG)=\rho(h,f)$ by the right invariance, 
and this is $\le \rho(h,1_\dG)+\rho(f,1_\dG)\le2\ve.$) \ 
By the same argument, if $q\in\dP$ forces 
$\dox\in\di{(G_\ve\app y)}$ for some $y\in\dX$ 
then $q$ is of size $2\ve$. 
}

\blt
\label{l4}
If\/ $q\in\dP,\msur$ $q\lef\wtau,$ 
and\/ $\ve>0,$ then there exists a condition\/ 
$r\in\dP,$ $\msur r\lef q,$ of size\/ $\le \ve$.
\elt
\bpf
Otherwise for any $r\in\dP,\msur$ $r\lef q,$ there 
is a pair of conditions $r',\,r''\in\dP$ stronger than 
$r$ and such that $\ang{r',r''}$ \dd{(\dP\ti\dP)}forces 
that there is no $g\in \di{G_\ve}$ with $\doxl=g\app\doxr.$ 
Applying, in a sufficiently generic extension $\dvp$ of 
$\dV,$ an ordinary splitting construction, we find a perfect 
set $X\sq\bk q$ such that any pair $\ang{x,y}\in X^2$ with 
$x\ne y$ is \dd{(\dP\ti\dP)}generic, 
hence, there is no $g\in \di{G_\ve}$ with $y=g\app x.$ 
Fix $x_0\in X.$
As $X$ is a pairwise \dd\egx equivalent set 
(together with $\bk q$)
we can associate, in $\dvp,$ with each $x\in X,$ an element 
$g_x\in\di G$ such that $x=g_x\app x_0,$ and $g_x\nin\di{G_\ve}$ 
by the above. 
Moreover, we have $g_yg_x\obr\app x=y$ for all $x,\,y\in X,$ 
hence $g_yg_x\obr\nin\di{G_\ve}$ whenever $x\ne y,$ which 
implies $\rho(g_x,g_y)\ge\ve$ by the right invariance. 
But this contradicts the separability of $G$.
\vyk{
Take any $y\in \bk q.$ 
Let $Q$ be a countable dense subset of $\dG.$ 
For any $g\in Q,$ the set $G_\ve g$ is the open 
\dd\ve nbhd of $g$ in $\dG$ by the right invariance, 
hence, $\dG=\bigcup_{g\in Q}G_\ve g.$ 
Then $\bk q\sq\bigcup_{g\in Q}G_\ve g\app y,$ hence, 
by Lemma~\ref{l2}, one of the sets 
$(G_\ve g\app y)\cap \bk q=(G_\ve\app x)\cap \bk q,$ 
where $x=g\app y,$ is forcible over $q,$ so that there 
is a ``condition'' $r\in\dP,$ $\msur r\lef q,$ which 
forces $\dox\in\di{((G_\ve\app x)\cap \bk q)}.$ 
On the other hand, $r$ is of size $\le2\ve$ by the 
remark above.
}%
\epF{Lemma}

It follows that there is, in $\dV,$ a sequence 
of codes $q_n\in\dP$ such that 
$q_0\lef\wtau,$ 
$q_{n+1}\lef q_n,\msur$ $q_n$ has size $\le 2^{-n},$ 
and $\bk{q_n}$ has \dd\dX diameter $\le 2^{-n}$ for any $n.$ 
The only limit point $x$ of the sequence of sets $\bk{q_n}$ 
belongs to $\dV,$ thus, it remains to show that 
$\bk{\wtau}\sq\ek{x}{\di{(\egx)}}$ in any extension 
$\dvp$ of the universe $\dV$.

We can assume that $\dvp$ is rich enough to contain, 
for any $n,$ an element $x_n\in \bk{q_n}$ such that 
each pair $\ang{x_n,x_{n+1}}$ is \dd{(\dP\ti\dP)}generic 
(over $\dV$).   
Then $\tlim_n x_n=x.$ 
Moreover, for any $n,$ both $x_n$ and $x_{n+1}$ 
belong to $\bk{q_n},$ hence, as $q_n$ has size 
$\le 2^{-n-1},$ there is $g_{n+1}\in \di\dG$ with 
$\rho(1,g)\le 2^{-n}$ such that $x_{n+1}=g_{n+1}\app x_n.$ 
Thus, $x_n=h_n\app x_0,$ where $h_n=g_{n}...g_1.$ 
Note that $\rho(h_n,h_{n-1})=\rho(g_n,1_\dG)\le 2^{-n+1}$ by 
the right-invariance of the metric, thus, $\sis{h_n}{n\in\dN}$ 
is a Cauchy sequence in $\di\dG.$ 
Let $h=\tlim_{n\to\iy}h_n\in\di\dG$ be its limit. 
As the action is continuous, we have $x=\tlim_nx_n=h\app x_0.$ 
It follows that $x\egx x_0.$ 
However $x_0\in \bk{q_0}\sq\bk{\wtau},$ therefore, 
$\bk{\wtau}\sq\ek{x}{\di{(\egx)}},$ as required.\vom

\epF{Theorem~\ref{clit}}

\punk{All $\Fs$ ideals are \PP}
\label{s02}

Let us say that a Borel ideal $\cI$ is {\it\PP\/} 
if so is the induced \er\ $\rei.$ 
\index{ideal!pointpin@\PP}%
It immediately follows from Theorem~\ref{clit} that any 
polishable ideal is \PP. 
Yet there are \PP\ ideals among non-polishable ones. 

\bte
\label{fspp}
Any\/ $\Fs$ ideal\/ $\cI\sq\pn$ is\/ \PP.
\ete
\bpf
We have $\cI=\bigcup_nF_n,$ where all sets $F_n\sq\pn$ 
are closed. 
It can be assumed that $F_n\sq F_{n+1},$ moreover, since 
for any closed $F\sq\pn$ the set 
$\Da F=\ans{X\sd Y:x,y\in F}$ 
is also closed (by the compactness of $\pn$), 
it can be assumed that $\Da F_n\sq F_{n+1}$ for all $n$. 

Let $\wtau$ be a Borel code, for a subset of $\pn,$ 
satisfying \poq{if} of Definition~\ref{t2like} \vrt\ 
the induced \er\ $\rei$ on $\pn,$ 
thus, $\wtau\in\dP,$ where $\dP$ is a forcing 
defined as in the proof of Proposition~\ref{bsfub} 
(but now $\dP$ forces a subset of $\pn,$ of course). 
Obviously there exists a pair of conditions $q,\,r\in\dP$ 
with $q,\,r\le\wtau,$ and a number $\nu\in\dN,$ 
such that $\ang{q,r}$ forces that 
$\ang{\doxl,\doxr}\in\di{F_\nu}.$ 
Then $\ang{q,q}$ forces $\doxl\sd\doxr\in\di{F_{\nu+1}}$ 
because $\Da F_\nu\sq F_{\nu+1}.$ 
It follows that, in $\dV,$ there is a sequence of numbers 
$i_0<i_1<i_2<...,$ a sequence 
$q\gef p_0\gef p_1\gef p_2\gef...$ of codes 
in $\dP,$ and, for any $n,$ a set $u_n\sq\ir0n,$ such that 
\ben
\tenu{(\arabic{enumi})}
\itla{'bs1}
each $p_n$ \dd\dP forces $\dox\cap\ir0n=u_n$;

\itla{'bs2}
any \dd\dP generic, over $\dV,$ $x,\,y\in\bk{p_n}$ 
satisfy $x\sd y\in \di{F_{\nu+1}}.$ 
\een
Let, in $\dV,$ $a=\bigcup_nu_n,$ then  
$a\cap\ir0n=u_n$ for all $n.$ 
Prove that $a$ pins $\bk\wtau,$ \ie, 
$\bk\wtau\sq\ek a{\di\rei}$ in any extension of $\dV$.

We can assume that, in the extension, for any $n$ there is 
a \dd\dP generic, over $\dV,$ element $x_n\in\bk{p_n}.$  
Then we have, by \ref{'bs2}, $x_0\sd x_n\in\di{F_{\nu+1}}$ 
for any $n,$ thus, $x_0\sd a\in \di{F_{\nu+1}}$ as well, 
because $\sis{x_n}{}\to a.$ 
We conclude that $x_0\mathbin{\di\rei} a,$ and 
$\bk\wtau\sq \ek a{\di\rei},$ as required.
\epf

\punk{Another family of \PP\ ideals}
\label{bsti}

We here present another family of \PP\ ideals. 
Suppose that $\sis{\vpi_i}{i\in\dN}$ is a sequence of 
lower semicontinuous (\lsc) submeasures on $\dN.$ 
Define
\dm
\Exh_{\sis{\vpi_i}{}}\,=\,
\ans{X\sq\dN:\vpy(X)=0}\,,
\quad\hbox{where}\quad
\vpy(X)=\tlis_{i\to\iy}\vpi_i(X)\,.
\dm
the exhaustive ideal of the sequence of submeasures. 
By Solecki's Theorem~\ref{sol} for any Borel P-ideal $\cI$ 
there is a single \lsc\ submeasure $\vpi$ such that  
$\cI=\Exh_{\sis{\vpi_i}{}}=\Exh_\vpi,$ where 
$\vpi_i(x)=\vpi(x\cap\iry i),$ however, for example, the 
non-polishable ideal $\Ii=\fio$ also is of the form 
$\Exh_{\sis{\vpi_i}{}},$ where for $x\sq\dN^2$ we define 
$\vpi_i(x)=0\,\hbox{ or }\,1$ if resp.\ 
$x\sq\,\hbox{ or }\,\not\sq\ans{0,...,n-1}\ti\dN$. 

\bte
\label{ibs}
Any ideal of the form\/ $\Exh_{\sis{\vpi_i}{}}$ is 
\PP.
\ete
\bpf
Thus let $\cI=\Exh_{\sis{\vpi_i}{}},$ all $\vpi_i$ 
being \lsc\ submeasures on $\dN.$  
We can assume that the submeasures $\vpi_i$ decrease, 
\ie, $\vpi_{i+1}(x)\le\vpi_i(x)$ for any $x,$ for if not 
consider the \lsc\ submeasures 
$\vpi'_i(x)=\tsup_{j\ge i}\vpi_j(x).$ 
Let $\wtau$ be a Borel code, for a subset of $\pn,$ 
satisfying \poq{if} of Definition~\ref{t2like} \vrt\ 
the induced \er\ $\rei$ on $\pn,$ 
thus, $\wtau\in\dP,$ where $\dP$ is a forcing 
defined as in the proof of Proposition~\ref{bsfub} 
($\dP$ forces a subset of $\pn$). 

Using the same arguments as above, we see that for any 
$p\in\dP,\msur$ $p\lef\wtau,$ and $n\in\dN,$ 
there are $i\ge n$ and codes $q,\,r\in\dP$ with 
$q,\,r\lef p,$ 
such that $\ang{q,r}$ \dd{\dP\ti\dP}forces that 
$\vpi_i(\doxl\sd\doxr)\le 2^{-n-1},$ hence, any two 
\dd\dP generic, over $\dV,$ elements $x,\,y\in\bk q$ 
satisfy $\vpi_i(x\sd y)\le 2^{-n}.$ 
It follows that, in $\dV,$ there is a sequence of numbers 
$i_0<i_1<i_2<...,$ a sequence 
$\wtau\gef p_0\gef p_1\gef p_2\gef...$ of codes 
in $\dP,$ and, for any $n,$ a set $u_n\sq\ir0n,$ such that 
\ben
\tenu{(\arabic{enumi})}
\itla{bs1}
each $p_n$ \dd\dP forces $\dox\cap\ir0n=u_n$;

\itla{bs2}
any \dd\dP generic, over $\dV,$ $x,\,y\in\bk{p_n}$ 
satisfy $\vpi_{i_n}(x\sd y)\le 2^{-n}.$ 
\een
Let, in $\dV,$ $a=\bigcup_nu_n,$ then  
$a\cap\ir0n=u_n$ for all $n.$ 
Prove that $a$ pins $\bk\wtau,$ \ie, 
$\bk\wtau\sq\ek a{\di\rE_\cI}$ in any extension of $\dV$.

We can assume that, in the extension, for any $n$ there is 
a \dd\dP generic, over $\dV,$ element $x_n\in\bk{p_n}.$  
Then we have, by \ref{bs2}, $\vpi_{i_n}(x_n\sd x_m)\le 2^{-n}$ 
whenever $n\le m.$ 
It follows that $\vpi_{i_n}(x_n\sd a)\le 2^{-n},$ because 
$a=\tlim_mx_m$ by \ref{bs1}. 
However we assume that the submeasures $\vpi_j$ decrease, 
hence, $\vpy(x_n\sd a)\le 2^{-n}.$ 
On the other hand, $\vpy(x_n\sd x_0)=0$ because 
all elements of $\bk{p_0}$ are pairwise 
\dd{\di\rE_\cI}equivalent. 
We conclude that $\vpy(x_0\sd a)\le 2^{-n}$ 
for any $n,$ 
in other words, $\vpy(x_0\sd a)=0,\msur$  
$x_0\mathbin{\di\rE_\cI} a,$ and 
$\bk\wtau\sq \ek a{\di\rE_\cI},$ as required.
\epf

\bqe
\label{ppq1}
Are all Borel ideals \PP~? 
The expected answer ``yes'' would show that $\rtd$ is not 
Borel reducible to any Borel ideal.  
Moreover, is any orbit \er\ of a Borel action of a Borel 
abelian group pinned~? 
But even this would not fully cover Hjorth's 
Theorem~\ref{clit}.
\eqe

\bqe[{{\rm Kechris}}]
\label{ppq2}
If Question~\ref{ppq1} answers in the positive, is it 
true that $\rtd$ is the \dd\reb least non-\PP\ Borel \er~?
\eqe

\cite{zap}

\parf{Universal analytic \er s and reduction to ideals}
\label{rosen}

\appendix

\parf{Technical introduction}
\label i

\punk{Notation}
\label{not}


\bit
\item\msur
$\dN=\ans{0,1,2,...}:$ natural numbers. 
$\dN^2=\dN\ti\dN.$ 

\item\msur
$\dnn$ is {\it the Baire space\/}. 
If $s\in\dN\lom$ (a finite sequence of natural numbers) 
then $\bon s\dnn=\ens{x\in\dnn}{s\su x},$ a basic clopen 
nbhd in $\dnn$.
\index{zzN@$\dN$}%
\index{zzons@$\bon s\dnn$}%

\item\msur
$X\sqa Y$ means that the difference $X\dif Y$ is finite.
\index{zzzsqa@$\sqa$}

\item
If a basic set $A$ is fixed then $\dop X=\doP X=A\dif X$ 
for any $X\sq A$.
\index{zzcX@$\dop X$}

\item
If $X\sq A\ti B$ and $a\in A$ then 
$\seq Xa=\ens{b}{\ang{a,b}\in X},$ 
a {\it cross-section\/}.
\index{cross-section}%
\index{zzXa@$\seq Xa$}

\item\msur
$\#X=\#(X)$ is the number of elements of a finite set $X$.
\index{zzznx@$\#(X)$}%
\index{number of elements}

\item\msur
$f\ima X=\ens{f(x)}{x\in X\cap\dom f},$ 
the \dd f{\it image\/} of $X$.  
\index{image}%
\index{zzfiX@$f\ima X$}

\item\msur
$\sd$ is the symmetric difference. 
\index{symmetric difference}%
\index{zzda@$\sd$}

\item\msur
$\susi x\;...\;$ means: 
``there exist infinitely many $x$ such that ...'', \\[1mm]
$\kazi x\;...\;$ means: 
``for all but finitely many $x,$ ... holds''.
\index{zzzex@$\susi x$}%
\index{zzzka@$\kazi x$}
\index{quantifiers!zzzex@$\susi x$}%
\index{quantifiers!zzzka@$\kazi x$}%

\item
{\it An ideal on a set\/ $A$\/} is, as usual, any set 
$\pu\ne\cI\sq\pws A,$ closed under $\cup$ and 
satisfying ${x\in\cI}\imp{y\in\cI}$ whenever 
$y\sq x\sq A.$ 
Thus, any ideal contains $\pu.$ 
We'll usually consider only {\it nontrivial\/} ideals, 
\ie, those which contain all 
singletons $\ans{a}\sq A$ and do not contain $A,$ \ie, 
$\pwf A\sq\cI\sneq\pws A$.
\index{ideal}
\index{ideal!nontrivial}

\item
If $\cI$ is an ideal on a set $A$ then let $\rei$ be 
an \eqr\ (\er, for brevity) on $\pws A,$ 
\index{equivalence relation, ER}%
defined as follows: $X\rei Y$ iff ${X\sd Y}\in\cI$. 
\index{zzei@$\rei$}%

\item
If $\rE$ is an \er\ on a set $X$ then   
$\eke y=\ens{x\in X}{y\rE x}$ for any $y\in X$ 
\index{zzye@$\ek y\rE$}%
\index{equivalence class}%
(the \dde{\it class\/} of $x$) and 
$\eke Y=\bigcup_{y\in Y}\eke y$ 
(the \dde{\it saturation\/} of $Y$) for $Y\sq X.$ 
A set $Y\sq X$ is \dde{\it invariant\/} if 
$\eke Y=Y$. 
\index{set!einvariant@\dde invariant}%
\index{saturation}
\index{zzYe@$\eke Y$}%

\item 
If $\rE$ is an \er\ on a set $X$ then a set $Y\sq X$ is 
{\it pairwise \dde equivalent, {\rm resp.}, 
pairwise \dde inequivalent\/}, 
\index{set!pairwiseeeqv@pairwise \dde equivalent}%
\index{set!pairwiseeineqv@pairwise \dde inequivalent}%
if $x\rE y,$ resp., $x\nE y$ holds
for  all $x\ne y$ in $Y$.  

\item
If $X,\,Y$ are sets and $\rE$ any binary relation then 
$X\rE Y$ means that we have both 
\index{zzXEY@$X\rE Y$}%
$\kaz x\in X\:\sus y\in Y\:(x\rE y)$ and 
$\kaz y\in Y\:\sus x\in X\:(x\rE y)$.
\eit

\punk{Descriptive set theory}
\label{DST}

A basic knowledge of Borel and projective hierarchy, 
both classical and effective, in 
{\it the Baire space\/} $\dnn$ and other 
(recursively presented, in the effective case) 
Polish spaces, is assumed. 

A map $f$ 
(between Borel sets in Polish spaces) is 
{\it Borel\/} iff its graph is a Borel set iff all 
\dd fpreimages of open sets are Borel. 
A map $f$ is {\it Baire measurable\/} 
({\it BM\/}, for brevity) iff all 
\dd fpreimages of open sets are Baire measurable sets. 
\index{map!Borel}%
\index{map!Baire measurable, BM}

\punk{Trivia of ``effective'' descriptive set theory}
\label{eff}

Apart of the very common knowledge, the whole 
instrumentarium of ``effective'' descriptive set theory 
employed in the study of reducibility of ideals and 
\er s, can be summarized in a rather short list  
of key ``principles''.
In those below, by a {\it recursively presented\/} Polish
\index{space!recursively presented}%
\index{recursively presented (space)}%
space one can understand any product space of the form
$\dN^m\ti(\dnn)^n$ without any harm for applications
below, yet in fact this notion is much wider.

\bre
\label{relrel}
For the sake of brevity, the results below are formulated 
only for the ``lightface'' parameter-free classes
$\is11\zd \ip11\zd \id11,$ 
but they remain true for $\is11(p)\zd \ip11(p)\zd \id11(p)$
for any fixed real parameter $p$.
\ere
 
\bde
\item[\Redu\ {\rm and} \Sepa:] \ 
\index{separation@\Sepa}%
\index{principle!separation@\Sepa}%
\index{reduction@\Redu}%
\index{principle!reduction@\Redu}%
If $X,\,Y$ are $\ip11$ sets of a recursively 
presented Polish space then there \poq{disjoint} $\ip11$ 
sets $X'\sq X$ and $Y'\sq Y$ with $X'\cup Y'=X\cup Y.$
The sets $X'\zi Y'$ are said to {\it reduce\/} the pair
\index{reduce}%
$X\zi Y$.

If $X,\,Y$ are disjoint $\is11$ sets of a recursively 
presented Polish space then there is a $\id11$ set $Z$ 
with $X\sq Z$ and $Y\cap Z=\pu$.
The set $Z$ is said to {\it separate\/} the $X$ from 
\index{separate}%
$Y$.

\item[\Cpro:] \ 
\index{countablep@\Cpro}%
\index{principle!countablep@\Cpro}%
If $P$ is a $\id11$ subset of the product $\dX\ti\dY$ 
of two recursively presented Polish spaces and for any 
$x\in\dX$ the cross-section $P_x=\ens{y}{P(x,y)}$ is at most 
countable then $\dom P$ is a $\id11$ set in $\dX$.
\ede
It follows that images of $\id11$ sets via 
countable-to-1, in particular, 1-to-1 $\id11$ maps are 
$\id11$ sets, while images via arbitrary $\id11$ maps 
are, generally, $\is11$.

\bde
\item[\Cenu:] \ 
\index{countablee@\Cenu}%
\index{principle!countablee@\Cenu}%
If $P,\,\dX,\,\dY$ are as in  
\Cpro\ then there is a 
$\id11$ map $f:\dom P\ti\dN\to\dY$ such that 
$P_x=\ens{f(x,n)}{n\in\dN}$ for all $x\in\dom P$.

\item[\Cuni:] \ 
\index{countableu@\Cuni}%
\index{principle!countableu@\Cuni}%
If $P,\,\dX,\,\dY$ are as in  
\Cpro\ then $P$ can be uniformized by a $\id11$ set.

\item[\Kres:] \ 
\index{kreisel@\Kres}%
\index{principle!kreisel@\Kres}%
If $\dX$ is a recursively presented Polish space, 
$P\sq\dX\ti\dN$ is a $\ip11$ set, and $X\sq\dom P$ 
is a $\id11$ set then there is a $\id11$ function 
$f:X\to\dN$ such that $\ang{x,f(x)}\in P$ for al 
$x\in X$.
\ede
The proof is surprisingly simple. 
Let $Q\sq P$ be a $\ip11$ set which uniformizes $P.$ 
For any $x\in X$ let $f(x)$ be the only $n$ with 
$\ang{x,n}\in Q.$ 
Immediately, (the graph of) $f$ is $\ip11,$ however, 
as $\ran f\sq\dN,$ we have 
$f(x)=n\eqv \kaz m\ne n\:(f(x)\ne m)$ whenever 
$x\in X,$ which demonstrates that $f$ is $\is11$ as well.

\bde
\item[\Penu:] \ 
\index{denumera@\Penu}%
\index{principle!denumera@\Penu}%
If $\dX$ is a recursively presented Polish space then 
there exist $\ip11$ sets $C\sq\dN$ and 
$W\sq\dN\ti\dX$ and a $\is11$ set $W'\sq\dN\ti\dX$ such that 
$W_e=W'_e$ for
any  $e\in C$ and a set $X\sq\dX$ is $\id11$ iff there is 
$e\in C$ such that $X=W_e=W'_e.$ \ 
(Here $W_e=\ens{x}{W(e,x)}$ and similarly $W'_e.$)
\ede
%
There is a generalization useful for relativised classes 
$\id11(y)$.

\bde
\item[\sf Relativized \Penu:] \ 
\index{denumera@\Penu!relativized}%
If $\dX,\,\dY$ are recursively presented Polish spa\-ces 
then there exist $\ip11$ sets $C\sq\dY\ti\dN$ and 
$W\sq\dY\ti\dN\ti\dX$ and a $\is11$ set 
$W'\sq\dY\ti\dN\ti\dX$ such that
$W_{ye}=W'_{ye}$  for any $\ang{y,e}\in C$ and, 
for any $y\in\dY,$ a set 
$X\sq\dX$ is $\id11(y)$ iff there is 
$e$ such that $\ang{y,e}\in C$ and $X=W_{ye}=W'_{ye}.$ \ 
($W_{ye}=\ens{x}{W(y,e,x)}$ and similarly $W'_{ye}$.)
\ede
%

Suppose that $\dX$ is a recursively presented Polish space.
A set $U\sq\dN\ti\dX,$ 
is a {\it a universal\/ $\ip11$ set\/} if   
\index{set!universal}%
for any $\ip11$ set $X\sq\dX$ there is an 
index $n$ with $X=U_n=\ens{x}{\ang{n,x}\in U},$
and a {\it a ``good'' universal\/ $\ip11$ set\/} if in
\index{set!universal!good}%
addition for any other $\ip11$ set $V\sq\dN\ti\dX$
there is a
recursive function $f:\dN\to\dN$ such that
$V_n=U_{f(n)}$ for all $n$.

The notions of universal and ``good'' universal $\is11$
sets are similar.

\bde
\item[\Uset:] \ 
\index{universalsets@\Uset}%
\index{principle!universalsets@\Uset}%
For any recursively presented Polish space $\dX$
there exist a ``good'' universal\/ $\ip11$ set
$U\sq\dN\ti\dX$ and a ``good'' universal\/ $\is11$ set
$V\sq\dN\ti\dX.$ 
(In fact we can take $V=(\dN\ti\dX)\dif U$.)
\ede

If a ``good'' universal $\ip11$ set $U$ is fixed then
a collection $\cA$ of $\ip11$ sets $X\sq\dX$ is 
{\it $\ip11$ in the codes\/}
\index{in the codes}%
if $\ens{n}{U_n\in\cA}$ is a $\ip11$ set. 
Similarly, if a ``good'' universal $\is11$ set $V$
is fixed then a collection $\cA$ of $\is11$ sets
$X\sq\dX$ is {\it $\ip11$ in the codes\/} if
$\ens{n}{V_n\in\cA}$ is a $\ip11$ set.
These notions quite obviously do not depend on the
choice of ``good'' universal sets.

To show how ``good'' universal sets work, we prove:

\bpro
\label{effred}
Let\/ $\dX$ be a recursively presented Polish space and\/
$U\sq\dN\ti\dX$ a ``good'' universal\/ $\ip11$ set.
Then for any pair of\/ $\ip11$ sets\/ $V,W\sq\dN\ti\dX$
there are recursive functions\/ $f,g:\dN\to\dN$ such
that for any\/ $m,n\in\dN$ the pair of cross-sections\/
$U_{f(m,n)}\zi U_{g(m,n)}$ reduces the pair\/  
$V_m\zi W_n$.
\epro
\bpf
Consider the following $\ip11$ sets in $(\dN\ti\dN)\ti\dX$:
\dm
P=\ens{\ang{m,n,x}}{\ang{m,x}\in V\land n\in\dN},\;\,
Q=\ens{\ang{m,n,x}}{\ang{n,x}\in W\land m\in\dN}.
\dm
By \Redu, there is a pair of $\ip11$ sets $P'\sq P$
and $Q'\sq Q$ which reduce the given pair $P\zi Q.$
Accordingly, the pair $P'_{mn}\zi Q'_{mn}$ reduces
$P_{mn}\zi Q_{mn}$ for any $m,n.$
Finally, by the ``good'' universality there are
recursive functions $f,g$ such that
$P'_{mn}=U_{f(m,n)}$ and $Q'_{mn}=U_{g(m,n)}$ 
for all $m,n$. 
\epf

The following principle is less 
elementary than the results cited above, but it is 
very useful because it allows to ``compress'' some  
sophisticated arguments with multiple applications of 
Separation and Kreisel selection.

\bde
\item[\Refl:] \ 
\index{reflection@\Refl}%
\index{principle!reflection@\Refl}%
Assume that $\dX$ is a recursively presented Polish 
space. 

{$\,$}\hbox{\hspace*{-5ex}$\ip11$ {\sl form\/\fsur}:} \  
Suppose that a collection $\cA$ of $\ip11$ sets $X\sq\dX$  
is $\ip11$ in the codes.
(In the sense of a fixed ``good'' 
universal $\ip11$ set $U\sq\dN\ti\dX.$)  
Then for any $X\in\cA$ there is a $\id11$ set 
$Y\in\cA$ with $Y\sq X$.

{$\,$}\hbox{\hspace*{-5ex}$\is11$ {\sl form\/}\fsur:} \  
Suppose that a collection $\cA$ of $\ip11$ sets $X\sq\dX$   
is $\ip11$ in the codes.
Then for any $X\in\cA$ there is a $\id11$ set 
$Y\in\cA$ with $X\sq Y$.
\ede

One of (generally, irrelevant here) consequences of 
this principle is that the set of all codes of a 
properly $\ip11$ set or properly $\is11$ set 
is never $\ip11$.

\punk{Polish--like families and the 
Gandy -- Harrington topology}
\label{bairef}

The following notion is similar to the Choquet property but 
somewhat more convenient to provide the 
nonemptiness of countable intersections of pointsets. 

\bdf
\label{genb}
A family $\cF$ is {\it Polish--like\/} if there exists a 
\index{Polish--like}%
countable collection $\ens{\cD_n}{n\in\dN}$ of dense subsets 
$\cD_n\sq\cF$ such that we have $\bigcap_nF_n\ne\pu$ 
whenever $F_0\qs F_1\qs F_2\qs...$ is a decreasing sequence 
of sets $F_n\in\cF$ which intersects every $\cD_n.$ \ 
(Here, a set $\cD\sq\cF$ is {\it dense\/} if 
\index{set!dense}%
$\kaz F\in\cF\:\sus D\in\cD\:(D\sq F).$)
\edf

For instance if $\cX$ is a Polish space then the collection 
of all its non-empty closed sets is Polish--like, for take 
$\cD_n$ to be all closed sets of diameter $\le n\obr.$ 

\bte[{{\rm Kanovei~\cite{umn}, Hjorth~\cite{h-ban}}}]
\label{s11pol}
The collection\/ $\cF$ of all non-empty\/ $\is11$ subsets 
of\/ $\dnn$ is Polish--like.\qeD
\ete

\bpf
For any $P\sq\dnn\ti\dnn$ define 
$\pr P=\ens{x}{\sus y\:P(x,y)}$ 
(the projection). 
If $P\sq\dnn\ti\dnn$ and $s,\,t\in\dN\lom$ then let 
$P_{st}=\ens{\ang{x,y}\in P}{s\su x\land t\su y}.$  
Let $\cD(P,s,t)$ be the collection of all $\is11$ sets 
$\pu\ne X\sq\dnn$ such that either $X\cap \pr P_{st}=\pu$ 
or $X\sq\pr P_{s\we i\,,\,t\we j}$ for some $i,\,j.$ 
(Note that in the ``or'' case $i$ is unique but $j$ may be 
not unique.) 
Let $\ens{\cD_n}{n\in\dN}$ be an arbitrary enumeration of all 
sets of the form $\cD(P,s,t),$ where $P\sq\dnn\ti\dnn$ is 
$\ip01.$ 
Note that in this case all sets of the form $\pr P_{st}$ are 
$\is11$ subsets of $\dnn,$ therefore, $\cD(P,s,t)$ is easily a 
dense subset of $\cF,$ so that all $\cD_n\sq \cF$ are dense. 

Now consider a decreasing sequence $X_0\qs X_1\qs ...$ of 
non-empty $\is11$ sets $X_k\sq\dnn,$ which intersects every 
$\cD_n\,;$ prove that $\bigcap_nX_n\ne \pu.$ 
Call a set $X\sq\dnn$ {\it positive\/} 
if there is $n$ such that $X_n\sq X.$ 
For any $n,$ fix a $\ip01$ set $P^n\sq\dnn\ti\dnn$ such that 
$X_n=\pr P^n.$ 
For any $s,\,t\in\dN\lom,$ if $\pr P^n_{st}$ is positive 
then, by the choice of the sequence of $X_n,$ 
there is a unique $i$ and some
$j$ such that  $\pr P^n_{s\we i\,,\,t\we j}$ is also positive.  
It follows that there is a unique $x=x_n\in\dnn$ and some  
$y=y_n\in\dnn$ (perhaps not unique) such that  
$\pr P^n_{x\res k\,,\,y\res k}$ is positive for any $k.$ 
As $P^n$ is closed, we have $P^n(x,y),$ hence, $x_n=x\in X_n$.

It remains to show that $x_m=x_n$ for $m\ne n.$ 
To see this note that if both $P_{st}$ and 
$Q_{s't'}$ are positive then either $s\sq s'$ 
or $s'\sq s$.
\epf

The collection of all non-empty $\is11$ subsets of $\dnn$ 
is a base of the {\it Gandy -- Harrington topology\/}, which has  
many  remarkable applications in descriptive set theory. 
This topology is easily not Polish, even not metrizable at all, 
yet it shares the following important property of Polish 
topologies:

\bcor
\label{s11ba}
The Gandy -- Harrington topology is Baire, \ie, every 
comeager set is dense.
\ecor
\bpf
This can be proved using Choquet property of the topology, 
see \cite{hkl}, however, the Polish--likeness 
(Theorem~\ref{s11pol}) also immediately yields the result.
\epf



\addcontentsline{toc}{section}{\refname}

\mtho
{\small\printindex}

\end{document}